%% file: O-N_quartic_Memoirs_postprint.tex
\documentclass[openany]{memo_amslatex/memo-l}

\input{zonedefs_color}

\usepackage[T1]{fontenc}
\usepackage{mathtools, verbatim,
  amsfonts, amssymb, wasysym, url}
\usepackage[dvipsnames,table]{xcolor}
\usepackage[shortlabels]{enumitem}
\usepackage[matrix,arrow,curve]{xy}
\usepackage[figuresright]{rotating}
\usepackage{tablefootnote, hyperref}

\newtheorem{thm}{Theorem}[chapter]
\newtheorem{thmstar}[thm]{Theorem*}
\newtheorem{lem}[thm]{Lemma}
\newtheorem{prop}[thm]{Proposition}
\newtheorem{cor}[thm]{Corollary}

\newtheorem{conj}[thm]{Conjecture}

\theoremstyle{definition}
\newtheorem{conv}[thm]{Conventions}
\newtheorem{nota}[thm]{Notation}
\newtheorem{defn}[thm]{Definition}
\newtheorem{examp}[thm]{Example}
\theoremstyle{remark}
\newtheorem{rem}[thm]{Remark}

\numberwithin{section}{chapter}
\numberwithin{equation}{chapter}

\DeclareMathOperator{\tr}{tr}
\DeclareMathOperator{\ct}{ct}
\DeclareMathOperator{\ntc}{ntc}
\DeclareMathOperator{\sgn}{sgn}
\DeclareMathOperator{\ch}{char}
\DeclareMathOperator{\coef}{coef}
\DeclareMathOperator{\disc}{disc}

\DeclareMathOperator{\inv}{inv}
\DeclareMathOperator{\ldr}{ldr}
\newcommand{\ur}{\mathrm{ur}}
\DeclareMathOperator{\Res}{res}
\DeclareMathOperator{\Cor}{cor}
\DeclareMathOperator{\Hom}{Hom}
\DeclareMathOperator{\Frac}{Frac}
\DeclareMathOperator{\Stab}{Stab}
\DeclareMathOperator{\Sym}{Sym}
\DeclareMathOperator{\Aut}{Aut}
\DeclareMathOperator{\Gal}{Gal}
\DeclareMathOperator{\End}{End}
\DeclareMathOperator{\Mat}{Mat}
\DeclareMathOperator{\charpoly}{char\,poly}
\newcommand{\GL}{\mathrm{GL}}
\newcommand{\GA}{\mathrm{GA}}
\newcommand{\SL}{\mathrm{SL}}
\newcommand{\PGL}{\mathrm{PGL}}

\newcommand{\PGGamma}{\mathrm{PG\Gamma}}
\newcommand{\SO}{\mathrm{SO}}

\newcommand{\xo}{{}_\times}
\newcommand{\AAA}{\mathbb{A}}
\newcommand{\CC}{\mathbb{C}}
\newcommand{\FF}{\mathbb{F}}

\newcommand{\NN}{\mathbb{N}}
\newcommand{\QQ}{\mathbb{Q}}
\newcommand{\PP}{\mathbb{P}}
\newcommand{\RR}{\mathbb{R}}

\newcommand{\ZZ}{\mathbb{Z}}
\newcommand{\aaa}{\mathfrak{a}}
\newcommand{\bb}{\mathfrak{b}}
\newcommand{\cc}{\mathfrak{c}}
\newcommand{\dd}{\mathfrak{d}}

\newcommand{\mm}{\mathfrak{m}}
\newcommand{\pp}{\mathfrak{p}}

\newcommand{\sss}{\mathfrak{s}}
\newcommand{\ttt}{\mathfrak{t}}
\newcommand{\cA}{\mathcal{A}}
\newcommand{\cB}{\mathcal{B}}
\newcommand{\cC}{\mathcal{C}}
\newcommand{\cD}{\mathcal{D}}
\newcommand{\cF}{\mathcal{F}}
\newcommand{\cG}{\mathcal{G}}
\newcommand{\fI}{\mathfrak{I}}
\newcommand{\cM}{\mathcal{M}}
\newcommand{\cN}{\mathcal{N}}
\newcommand{\OO}{\mathcal{O}}
\newcommand{\fP}{\mathfrak{P}}
\newcommand{\cS}{\mathcal{S}}
\newcommand{\cV}{\mathcal{V}}
\newcommand{\Sm}{\mathfrak{S}}
\newcommand{\fcr}{\mathfrak{r}}
\newcommand{\Ell}{\mathrm{Ell}}
\newcommand{\One}{\mathbf{1}}
\newcommand{\ZERO}{{\mathord{0}}}
\newcommand{\TEE}{{\mathord{\top}}}
\newcommand{\STAR}{{\mathord{\star}}}
\newcommand{\ba}{\overline}
\newcommand{\bs}{\backslash}
\newcommand{\cross}{\times}
\newcommand{\tensor}{\otimes}
\newcommand{\longto}{\mathop{\longrightarrow}\limits}
\newcommand{\textand}{\quad \text{and} \quad}
\newcommand{\textor}{\quad \text{or} \quad}
\newcommand{\tri}[1]{#1^\triangle}
\newcommand{\size}[1]{\lvert #1 \rvert}
\newcommand{\Size}[1]{\left\lvert #1 \right\rvert}
\newcommand{\floor}[1]{\left\lfloor #1 \right\rfloor}
\newcommand{\ceil}[1]{\left\lceil #1 \right\rceil}
\newcommand{\intsec}{\cap}
\newcommand{\union}{\cup}
\newcommand{\nequiv}{\not\equiv}
\newcommand{\isom}{\cong}
\newcommand{\afterhat}{\,\widehat{}}
\newcommand{\ignore}[1]{}
\newcommand{\ds}{\displaystyle}

\newcommand{\quickcol}[1]{%
  \begin{tabular}{@{}c@{}}%
    #1
  \end{tabular}
}

\newcommand{\support}[1]{%
	\quickcol{%
	  \begin{picture}(4,4)%
	    \put(0,0){\framebox(4,4){}}
	    \put(2,0){\line(0,1){4}}
	    \put(0,2){\line(1,0){4}}
	    #1
	  \end{picture}
  }
}
\newcommand{\sh}[2]{%
  \put(#1,#2){\rule{\the\unitlength}{\the\unitlength}}
}

\numberwithin{equation}{section}

\begin{document}

\frontmatter

\title{Reflection theorems of Ohno-Nakagawa type\\
  for quartic rings and pairs of \emph{n}-ary quadratic forms}
\author{Evan M. O'Dorney}

\subjclass[2010]{
  Primary 11R16, 
  11E12, 
  11E08, 
  Secondary 11A15, 
  11R54 
}

\begin{abstract}
We prove a reflection theorem, conjectured by Nakagawa and Ohno, for the
number of quartic rings, or pairs of ternary quadratic forms, with a given
cubic resolvent. Over $\mathbb{Z}$, our results are unconditional; we also
allow the base to be the ring of integers of a general number field,
conditional on some algebraic identities that are Monte Carlo verified. We also
establish a reflection theorem for quartic $11111$-forms and $48441$-forms that
relates them to the number of $3\times 3$ symmetric matrices with given
characteristic polynomial. Along the way, we find elegant new results on Igusa
zeta functions of conics and the average value of a quadratic character over a
box in a local field.

We conjecture that a reflection theorem holds for pairs of $n$-ary quadratic
forms for any odd $n$, and we prove this for odd cubefree discriminant. This
furnishes a more satisfactory answer for a question raised by Cohen, Diaz y
Diaz, and Olivier, namely whether there exists an infinite family of reflection
theorems of Ohno-Nakagawa type.
\end{abstract}

\maketitle

\tableofcontents

\mainmatter

\chapter{Introduction}

\section{Results on pairs of ternary quadratic forms}

Let $\cV(\ZZ)$ be the space of pairs $(A, B)$ of integer-coefficient quadratic forms in three variables. $\cV(\ZZ)$ is a $12$-dimensional lattice whose points parametrize quartic rings by one of Bhargava's celebrated higher composition laws \cite{B3}. It has an action of $\SL_3(\ZZ)$ by simultaneously changing the variables in the two forms. This action has four algebraically independent invariants, which, if we view $A$ and $B$ as $3 \times 3$ matrices with half-integer off-diagonal entries, can be seen as the four coefficients of the resolvent binary cubic
\begin{equation}
  f(x,y) = 4 \det (A x + B y).
\end{equation}
In this paper, we will show the following relation:
\begin{thm}\label{thm:O-N_quartic_Z}
Let $f(x,y)$ be an integer binary cubic form with no repeated root. Denote by $h(f)$ the number of orbits $(A, B) \in \SL_3(\ZZ)\bs \cV(\ZZ)$ with resolvent cubic $f$, weighting each $(A,B)$ by the reciprocal of the order of its stabilizer in $\SL_3(\ZZ)$. Denote by $h_2(f)$, $h^+(f)$, and $h^+_2(f)$ the number of these satisfying certain added conditions, still weighting in the same way:
\begin{itemize}
  \item For $h_2$, we require that the matrices $(A,B)$ have integer entries, that is, the associated symmetric bilinear forms are integral;
  \item For $h^+$, we require that $A$ and $B$, regarded as a pair of conics in $\PP^2(\RR)$, have a common real point (a nontrivial condition only if $f$ has three real roots);
  \item for $h^+_2$, we impose both of the foregoing conditions.
\end{itemize}
Then if $f$ has three real roots, we have the exact identities
\begin{align}
  h^+_2(4f) &= h(f) \\
  h_2(4f) &= 4h^+(f),
\end{align}
while if $f$ has only one real root, we have
\begin{equation}
  h_2(4f) = h^+_2(4f) = 2h(f) = 2h^+(f).
\end{equation}
\end{thm}

This theorem was conjectured by Nakagawa \cite{NakPairs} and proved there in the case that $f$ corresponds under the Delone-Faddeev parametrization to a cubic ring $C$ whose index in its maximal order is squarefree. The purpose of this memoir will be to prove it, and generalizations thereof.

Heuristically, we normally expect class numbers such as $h(f)$ to behave randomly. (We call a count such as $h(f)$ a class number, because it generalizes Gauss's class numbers of binary quadratic forms; $h(f)$ does not naturally arise as the order of any group.) Therefore it was quite astounding when the class numbers that count integer binary cubic forms by discriminant were found to obey an exact identity, which we call the \emph{Ohno--Nakagawa (O-N) reflection theorem:}
\begin{thm}[Nakagawa \cite{Nakagawa}; conjectured by Ohno \cite{Ohno}] 
  For a nonzero integer $D$, let $h(D)$ be the number of $\GL_2(\ZZ)$-orbits of \emph{binary cubic forms}
  \[
  f(x,y) = ax^3 + bx^2y + cxy^2 + dy^3
  \]
  of discriminant $D$, each orbit weighted by the reciprocal of its number of symmetries (i.e.~stabilizer in $\GL_2(\ZZ)$). Let $h_3(D)$ be the number of such orbits $f(x,y)$ such that the middle two coefficients $b, c$ are multiples of $3$, weighted in the same way.

  Then for every nonzero integer $D$, we have the exact identity
  \begin{equation}
    h_3(-27 D) = \begin{cases}
      3 h(D), & D > 0 \\
      h(D), & D < 0.
    \end{cases}
  \end{equation}
\end{thm}
In that case the forms parametrize cubic rings and $3$-torsion in quadratic rings. Hence it was natural to seek an analogue for pairs of ternary quadratic forms, which parametrize quartic rings and $2$-torsion in cubic rings.

In a companion to this paper \cite{OCubic}, the author provided a streamlined proof of O-N that generalizes to the case when the integers $\ZZ$ are replaced by the ring of integers $\OO_K$ of an arbitrary number field $K$. We apply the same techniques to prove Theorem \ref{thm:O-N_quartic_Z}, and they are largely successful over a general number field, where we reduce the theorem to a finite family of identities among rational functions of the form
\begin{equation} \label{eq:fracs_eq_0}
  \sum_{i = 1}^r \frac{1}{f_i(X_1, \ldots, X_k)} \stackrel{?}{=} 0
\end{equation}
where each $f_i$ is an integer polynomial (in fact, a product of binomials). Even in the case where $R/K$ is unramified, we have $r > 9000$ fractions in $k = 7$ variables,
which is out of reach with current computing methods. However, we can probabilistically verify a relation of the form \eqref{eq:fracs_eq_0} by plugging random values from a large finite field $\FF_p$ in place of the variables $X_i$ and summing. If \eqref{eq:fracs_eq_0} does not hold, we expect the result to be nonzero almost all the time, but in view of the high degree of the denominator, we cannot easily prove the identity (even mod $p$) just by plugging in enough different values. Here and in the sequel, we use the label ``Theorem*'' to denote a theorem that is not rigorously proved but is Monte Carlo verified in this manner.
\begin{thmstar}\label{thmstar:O-N_2x3x3}
  Let $K$ be a number field. Let $\aaa$ and $\ttt$ be ideals of $\OO_K$ with $(1) \supseteq \ttt \supseteq (2)$. Denote by $\cV_{\ttt,\aaa}(\OO_K)$ the space of pairs of ternary quadratic forms
  \begin{equation}\label{eq:quartic}
    (A, B) : \OO_K \cross \OO_K \cross \aaa \rightrightarrows \OO_K \cross \aaa,
  \end{equation}
  that are $\ttt$-integral in the sense that the entries belong to the ideals
  \[
  \left(\begin{bmatrix}
    (1) & 2^{-1}\ttt & 2^{-1}\ttt\aaa^{-1} \\
    & (1) & 2^{-1}\ttt\aaa^{-1} \\
    & & \aaa^{-1}
  \end{bmatrix},
  \begin{bmatrix}
    \aaa & 2^{-1}\ttt\aaa & 2^{-1}\ttt \\
    & \aaa & 2^{-1}\ttt \\
    & & (1)
  \end{bmatrix}
  \right).
  \]
  It has a natural action of the group
  \[
  \cG_{\aaa} = \SL(\OO_K \oplus \OO_K \oplus \aaa) \cross \SL(\OO_K \oplus \aaa)
  \]
  with a single invariant, the \emph{discriminant}
  \[
  \disc(A, B) = \disc\bigl(4\det(Ax + By)\bigr) = 256\disc\bigl(\det(Ax + By)\bigr) \in \ttt^8 \aaa^{-2}.
  \]
  Denote by $h_{\ttt, \aaa}(D)$ the number of $\cG_{\aaa}$-orbits in $\cV_{\ttt, \aaa}(\OO_K)$ having discriminant $D$, each orbit weighted by the reciprocal of the order of its stabilizer in $\cG_{\aaa}$. Denote by $h^{\ntc}_{\ttt,\aaa}(D)$ the number of such orbits (weighted in the same way) which are nowhere totally complex, in the sense that at each real place of $K$, the forms $A$ and $B$ have a common zero in $\RR\PP^2$. Then for all discriminants $D$,
  \[
  h_{\ttt,\aaa}(D) = \frac{N(\ttt)^2}{2^{r_\infty}} \cdot h^{\ntc}_{2\ttt^{-1},2^4\ttt^{-8}\aaa}(D),
  \]
  where $r_\infty$ is the number of real places of $K$ at which $D < 0$ plus twice the number of complex places of $K$.

  The label ``Theorem*'' indicates that, in general, the result has not been rigorously shown, but rather has been reduced to a Monte Carlo verification showing that it is probabilistically true. However, the following cases are rigorously known:
  \begin{enumerate}[$($a$)$]
    \item When $K = \QQ$, or more generally, when $K$ is unramified at $2$;
    \item When the reduced discriminant $D \aaa^{2} \ttt^{-8} $ is coprime to $2$.
  \end{enumerate}
\end{thmstar}

The method of Bhargava \cite{B3}, presented by the author \cite[\textsection8]{ORings} (see Theorem \ref{thm:hcl_quartic} below), shows that forms of the type \ref{eq:quartic}, up to a slightly larger group
\[
\GL(\OO_K \oplus \OO_K \oplus \aaa) \cross \GL(\OO_K \oplus \aaa),
\]
parametrize quartic rings over $\OO_K$ of Steinitz class $\aaa$. We may therefore equally well state our reflection theorems in terms of quartic rings.

\begin{defn}
  Call a nondegenerate quartic ring $\OO/\OO_K$  over a number field \emph{nowhere totally complex (ntc)} if there is no real place $\pp$ of $K$ such that $\OO \tensor_{\OO_K} K_\pp \isom \CC \cross \CC$.
\end{defn}

\begin{defn}
  Let $\OO_K$ be a Dedekind domain, $\ch K \neq 2$, and let $\ttt$ be an ideal dividing $(2)$ in $\OO_K$. A quartic ring $\OO$ over $\OO_K$ equipped with a cubic resolvent $C$ is called \emph{$\ttt$-traced} if, for all $x$ and $y$ in $\OO$, the associated bilinear form
  \[
  \Phi(x,y) = \frac{\Phi(x+y) - \Phi(x) - \Phi(y)}{2}
  \]
  whose diagonal restriction is $\Phi(x,x) = \Phi(x)$ takes values in $2^{-1}\ttt C$. We say that $\OO$ is \emph{$\ttt$-traced} if it admits a $\ttt$-traced resolvent.
\end{defn}
If $\OO$ has a $\ttt$-traced resolvent $C$, it is not hard to show (Proposition \ref{prop:traced}\ref{traced:reduced} below) that $C = \OO_K + \ttt^{2}C_\ttt$ for a unique cubic ring $C_\ttt$; we call $C_\ttt$ a \emph{reduced resolvent} of the $\ttt$-traced ring $\OO$.

\begin{thmstar}[\textbf{``Quartic O-N''}] \label{thmstar:O-N_quartic}
  Let $K$ be a number field. Let $C$ be an order in a nondegenerate cubic $K$-algebra $R$, and let $\ttt \subseteq \OO_K$ be an ideal such that $\ttt \mid (2)$. Let $h(C, \ttt)$ count the $\cG(\OO_K)$-orbits on $\cV_{\ttt,C}(\OO_K)$, each orbit weighted by the reciprocal of its point stabilizer. Equivalently, let $h(C, \ttt)$ count the number of $\ttt$-traced quartic rings with reduced resolvent $C$, each weighted by the reciprocal of its number of resolvent-preserving automorphisms. Let $h^{\ntc}(C, \ttt)$ count the subset of the foregoing that are ntc, weighted in the same way. Then
  \[
  h(C, \ttt) = \frac{N(\ttt)^2}{2^{r_\infty}} \cdot h^{\ntc}(C, 2\ttt^{-1}),
  \]
  where $r_\infty$ is the number of real places of $K$ over which $R$ is not totally real plus twice the number of complex places of $K$.

  The following cases are rigorously known:
  \begin{enumerate}[$($a$)$]
    \item When $K = \QQ$, or more generally, when $K$ is unramified at $2$;
    \item When $C$ is maximal at each place lying above $2$.
  \end{enumerate}
\end{thmstar}

\begin{thmstar} \label{thmstar:O-N_quartic_by_disc}
  Let $K$ be a number field, and let $\cD = (\aaa, D)$ be a discriminant.
  Denote by $h^\circ_\ttt(\cD)$ the number of $\ttt$-traced quartic rings $\OO$ over $\OO_K$ having discriminant $\cD$, each $\OO$ weighted by $1/\size{\Aut_K(\OO)}$. Denote by $h_\ttt^{\circ,\ntc}(\cD)$ the number of such that are ntc, weighted in the same way. Then for all discriminants $\cD$,
  \begin{equation}\label{eq:O-N_by_disc}
    h^\circ_\ttt(\ttt^8\cD) = \frac{N(\ttt)^2}{2^{r_\infty}} \cdot h^{\circ,\ntc}_{2\ttt^{-1}}(256\ttt^{-8}\cD),
  \end{equation}
  where $r_\infty$ is the number of real places of $K$ at which $D < 0$ plus twice the number of complex places of $K$.

  The cases in which this theorem is known rigorously are the same as in the previous theorem.
\end{thmstar}

The proofs and verifications for the starred theorems are carried out in the attached file \verb|main.sage|\footnote{Available at \url{https://github.com/emo916math/Quartic-Ohno-Nakagawa.}}. The run times (on a modern Intel i7 running Linux Mint) are shown in Table \ref{tab:runtime}.
\begin{table}[bt]
  \begin{center}
  \begin{tabular}{ccc}
  Spl.t. & Monte Carlo verifications & Provable cases \\ \hline
  ur & 3 min & 2 min \\
  $1^3$ & 4 min & 2 min \\
  $1^2 1$ & 2 hr 52 min\tablefootnote{Since Sage does not have a robust garbage collector, it is advisable to restart Sage at each \texttt{reset()} command in \texttt{main.sage} to avoid running out of memory. The time listed here is the sum of two runs.} & 33 min
\end{tabular}
\end{center}
\label{tab:runtime}
\caption{Run times for the verification of Theorems* \ref{thmstar:O-N_2x3x3}, \ref{thmstar:O-N_quartic}, \ref{thmstar:O-N_quartic_by_disc}, and Theorem \ref{thm:BQ}}
\end{table}

\section{Results on binary quartic forms}

Another natural space to study is the space of integer binary quartic forms, which have been shown by Wood \cite{WoodBQ} to parametrize quartic rings having a monogenized resolvent. If
\[
f(x,y) = a x^4 + b x^3 y + c x^2 y^2 + d x y^3 + e y^4
\]
is a quartic polynomial with integer coefficients, its \emph{resolvent} is
\begin{equation}
  g(z) = z^{3} - c z^{2} + (b d - 4 a e) z + 4 a c e - b^2 e - a d^2;
\end{equation}
equivalently, for a generic
\[
f(x,y) = a(x - x_1y)(x - x_2y)(x - x_3y)(x - x_4y),
\]
we set
\[
g(z) =
\big(z - a (x_1 x_2 + x_3 x_4)\big)
\big(z - a (x_1 x_3 + x_2 x_4)\big)
\big(z - a (x_1 x_4 + x_2 x_3)\big).
\]
It can be verified that if $f$ is changed by a transformation in $\PGL_2(\ZZ)$, the resolvent is changed only by a $\ZZ$-translation of $z$. For this reason, we consider resolvents only up to $\ZZ$-translation. For a monic cubic polynomial $g$, write
\begin{equation}\label{eq:quartic_resolvent_scale}
  g^{[4]}(z)=4^3g(z/4).
\end{equation}
Thus $g^{[4]}$ is the resolvent of the quartic form $4f$ whenever $g$ is the resolvent of $f$.

The reflection theorem we get for binary quartic forms has a somewhat unexpected shape. First, it involves $48441$-forms, that is, forms whose five coefficients are divisible by $4$, $8$, $4$, $4$, and $1$ in that order. Secondly, it relates to another classical problem, namely the enumeration of symmetric matrices with given characteristic polynomial!

\begin{thm}[\textbf{Quartic O-N for binary quartic forms}]\label{thm:BQ}
  Let $g$ be a monic integral binary cubic form. Denote by $h(g)$ the number of integral binary quartic forms of resolvent $g$, up to $\PGL_2(\ZZ)$-equivalence and weighted by inverse of $\PGL_2(\ZZ)$-stabilizer. Denote by $h_4(g)$ the number of binary quartic $48441$-forms of resolvent $g$, up to $\PGGamma^0(2)$-equivalence and weighted by inverse of $\PGGamma^0(2)$-stabilizer, where $\PGGamma^0(2) \subset \PGL_2(\ZZ)$ is the subgroup of matrices $\big[\begin{smallmatrix}
    a & b \\ c & d
  \end{smallmatrix}\big]$ such that $b \equiv 0 \pmod 2$.
  Denote by $s(g)$ the number of integral $3\times 3$ symmetric matrices of characteristic polynomial $g$. Then:
  \begin{itemize}
    \item If $\disc g < 0$, then
    \begin{align}
    \label{eq:BQ_neg_disc}
      4 h(g) &= h_4(g^{[4]}).
    \intertext{\item If $\disc g > 0$, then}
      2 h(g) &= h_4^{\text{indef}}(g^{[4]}) \label{eq:BQ_indef} \\
      4 h^{\text{indef or pos def}}(g) &= h_4^{\text{indef or pos def}}(g^{[4]}) \label{eq:BQ_pos} \\
      4 h^{\text{indef or neg def}}(g) &= h_4^{\text{indef or neg def}}(g^{[4]}) \label{eq:BQ_neg} \\
      24\big(h^{\text{indef}}(g) - h^{\text{def}}(g)\big) &= s(g) \label{eq:BQ_24}
    \end{align}
    where the superscripts instruct one to count only forms satisfying the indicated condition at $\infty$, with the same weighting.
  \end{itemize}
\end{thm}
\begin{cor}\label{cor:BQ}
  Let $g$ be a monic integral binary cubic form with three real roots. Among integral binary quartics with resolvent $g$, at least half are indefinite when we weight by inverse size of $\PGL_2(\ZZ)$-stabilizer, with equality exactly when $g$ is not the characteristic polynomial of an integral $3\times 3$ symmetric matrix.
\end{cor}

\section{Results on \texorpdfstring{$2 \times n \times n$}{2×n×n} boxes}

Finally, we have a result on $2\times n\times n$ symmetric boxes, that is, pairs of $n$-ary quadratic forms, which by a result of Wood \cite{W2xnxn} parametrize $2$-torsion (strictly speaking, self-balanced ideals) in rings parametrized by binary $n$-ic forms. For notational simplicity, we restrict ourselves to the case where the ideals denoted $\aaa$ and $\ttt$ in Theorem* \ref{thmstar:O-N_2x3x3} are principal.
\begin{conj}[\textbf{O-N for $2\times n\times n$ boxes}]\label{conj:2xnxn}
  Let $K$ be a number field, $n \geq 1$ an odd integer, and let $\tau$ be a divisor of $2$ in $\OO_K$. Denote by $h_\tau(f)$ the number of $\Gamma(\OO_K) = \SL_n(\OO_K)$-orbits of pairs $(\cA,\cB)$ of $n\times n$ symmetric matrices whose on- and off-diagonal elements belong to $\OO_K$ and $2^{-1} \tau \OO_K$ respectively, each $(\cA,\cB)$ weighted by the reciprocal of the order of its stabilizer in $\Gamma(\OO_K)$. Denote by $h_\tau^{\ntc}(f)$ the count (by the same weighting) of the subset of orbits whose corresponding self-balanced ideal $(R_f, I, \delta)$ has $\delta > 0$ at every real place of $K$. Then we have a global reflection theorem
  \[
  h_\tau(\tau^{n-1}f) = \frac{\size{N_{K/\QQ}(\tau)}^{n-1}}{2^{r_\infty}} \cdot h^{\ntc}_{2\tau^{-1}}\left((2\tau^{-1})^{n-1} f\right),
  \]
  where
  \[
    r_\infty = \sum_{v \text{ complex}} (n - 1) + \sum_{v \text{ real}} (\text{number of pairs of complex conjugate roots of $f$ at $v$}).
  \]
\end{conj}

\begin{thm}\label{thm:2xnxn}
  Conjecture \ref{conj:2xnxn} holds when $K/\QQ$ is unramified at $2$, and $\disc f$ is cubefree and coprime to $2$ in $\OO_K$.
\end{thm}

\begin{rem}
  The $n = 1$ case is trivial, and the $n = 3$ case reduces to Theorem* \ref{thmstar:O-N_2x3x3}. Note that, for the first time, we have (conjecturally) an infinite family of representation spaces where a nontrivial reflection theorem holds for all nonsingular orbits. By comparison, in their work, Cohen--Rubinstein-Salzedo--Thorne \cite{CohON} also produce an infinite family, but as they look only at field extensions, these reflection theorems are closer in spirit to the original Scholz reflection theorem than to O-N. In the author's opinion, the depth of O-N lies in the mysterious coincidence of counting functions on the local level, which this memoir proves in many more cases but cannot be said to explain.
\end{rem}

\section{Methods}

Several proofs of O-N are now in print \cite{Nakagawa,Marinescu,OOnARemarkable,Gao}, all of which use class field theory together with various approaches to understanding the local behavior of nonmaximal cubic orders. We may thus divide each proof into two steps:
\begin{itemize}
  \item A \emph{global} step that uses global class field theory to understand cubic fields, equivalently $\GL_2(\QQ)$-orbits of cubic forms;
  \item A \emph{local} step to count the rings in each cubic field, equivalently the $\GL_2(\ZZ)$-orbits in each $\GL_2(\QQ)$-orbit, and put the result in a usable form.
\end{itemize}
In generalizing to the quartic case, the global step does not pose especial difficulties, but we will profit from the streamlined approach to the global step introduced in \cite{OCubic} using \emph{composed varieties:}

\begin{defn}[Definition 2.1 of \cite{OCubic}]
  \label{defn:composed_quartic}
  Let $K$ be a field. A \emph{composed variety} over $K$ consists of the data $\bigl(\Gamma, V, X, X_0\bigr)$, where
  \begin{itemize}
    \item $V$ is a finite-dimensional, faithful representation of a connected linear algebraic group $\Gamma$ over $K$;
    \item $X \subseteq V$ is a closed subvariety whose points $X(\bar K)$ form a single $\Gamma(\bar K)$-orbit, such that the stabilizer $M = \Stab_\Gamma(x_0)$ is finite, abelian, and \'etale of order not divisible by $\ch K$;
    \item $X_0 = \Gamma(K)x_0$ is an orbit of a $K$-rational point $x_0 \in X$.
  \end{itemize}
\end{defn}

\begin{prop}[Proposition 2.3 of \cite{OCubic}]
  \label{prop:composed_quartic}
  Let $\bigl(\Gamma, V, X, X_0\bigr)$ be a composed variety. Then:
  \begin{enumerate}[$($a$)$]
    \item\label{compq:psi} There is a natural injection
    \[
      \psi : \Gamma(K)\backslash X(K) \hookrightarrow H^1(K,M)
    \]
    by which the orbits $\Gamma(K)\backslash X(K)$ parametrize some subset of the Galois cohomology group $H^1(K, M)$;
    \item\label{compq:0} We have $\psi(X_0) = 1$;
    \item\label{compq:stab} The $\Gamma(\bar K)$-stabilizer of every $x \in X(K)$ is canonically isomorphic to $M$ as a Galois module;
    \item\label{compq:chg_bspnt} If $X_0$ is replaced by a different base orbit $X_1 = \Gamma(K) x_1$, then the associated parametrization $\psi_{X_1} : \Gamma(K)\backslash X(K) \to H^1(K,M)$ differs only by translation:
    \[
      \psi_{X_1}(Y) = \psi(Y)\psi(X_1)^{-1},
    \]
    under the isomorphism between the stabilizers $M$ established in the previous part;
    \item\label{compq:base_chg} If $K'/K$ is a field extension, then $\bigl(\Gamma', V', X', X_0'\bigr)$ is also a composed variety, where $\Gamma', V', X'$ are formed by extension of scalars and $X_0' = \Gamma(K') x_0$. The corresponding parametrizations $\psi, \psi'$ are related by the commutative diagram
    \[
      \xymatrix{
        \Gamma(K)\bs X(K) \ar@{^{(}->}[r]^\psi \ar[d]
          & H^1(K,M) \ar[d]^\Res \\
        \Gamma'(K')\bs X'(K') \ar@{^{(}->}[r]^{\psi'}
          & H^1(K',M).
      }
    \]
  \end{enumerate}
\end{prop}

\begin{defn}[Definition 2.4 of \cite{OCubic}]~
  \label{defn:fullHasse_quartic}
  \begin{enumerate}[$($a$)$]
    \item A composed variety $X$ is \emph{full} if its cohomology parametrization $\psi$ is surjective, that is, $X$ includes a $\Gamma(K)$-orbit for every coclass $\alpha \in H^1(K,M)$.
    \item If $K$ is a global field, a composed variety is \emph{Hasse} if it satisfies the following version of the Hasse local-global principle: For $\alpha \in H^1(K,M)$, if the localization $\alpha_v \in H^1(K_v,M)$ at each place $v$ lies in the image of the local parametrization
    \[
      \psi_v : \Gamma(K_v)\backslash X(K_v) \to H^1(K_v, M),
    \]
    then $\alpha$ also lies in the image of the global parametrization $\psi$.
  \end{enumerate}
\end{defn}

In proving O-N, we took $V = \Sym^3(K^2)$, and $X$ was the closed variety of binary cubic forms of a fixed nonzero discriminant $D$, and $\Gamma = \SL_2$. The point stabilizer $M$ was cyclic of order $3$.

In this memoir, $V = \Sym^2(K^3) \tensor K^2$, and $X = X_f$ is the closed variety of pairs of ternary quadratic forms over $K$ whose resolvent is a fixed separable binary cubic form $f$, and $\Gamma = \SL_3$ acts naturally. We take $X_0$ to be the orbit corresponding to the reducible algebra $K \cross R_f$. The point stabilizer $M$ is isomorphic to $\cC_2 \cross \cC_2$, with $\Gal(\bar K/K)$ permuting the three non-identity elements as it permutes the three roots of $f$.

In previous work \cite{OCubic}, the author developed weighted lattice counts on composed varieties and proved a local-to-global reflection engine. We recall the statements in the form needed here.

Let $K_v$ be a local field, let $(\Gamma,V,X,X_0)$ be a composed variety over $K_v$, and let $\Lambda_v\subseteq V(K_v)$ be a lattice. Write $\Gamma_{\Lambda_v}$ for its stabilizer in $\Gamma(K_v)$. A \emph{local weighting} is a $\Gamma_{\Lambda_v}$-invariant function
\[
  w_v:\Lambda_v\intsec X(K_v)\longto\CC.
\]
For $\alpha_v\in H^1(K_v,M)$ represented by $x_{\alpha_v}\in X(K_v)$, define the \emph{weighted local lattice count}
\begin{equation}\label{eq:local_lattice_count_quartic}
  \ell_v^{w_v}(\alpha_v)
  =\sum_{\substack{\Gamma_{\Lambda_v}\gamma\in
      \Gamma_{\Lambda_v}\bs\Gamma(K_v)\\
      \gamma x_{\alpha_v}\in\Lambda_v}}
    w_v(\gamma x_{\alpha_v}),
\end{equation}
and set it equal to $0$ if $\alpha_v$ is not represented by an orbit in $X(K_v)$.

For a global composed variety, let $\Lambda=\prod_v\Lambda_v$ be an adelic lattice and let $w=\prod_vw_v$ be a global weighting. Its \emph{total weighted class number} is denoted $h^w(X)$. If $X$ is Hasse, it is given by
\begin{equation}\label{eq:total_weighted_class_quartic}
  h^w(X)=\sum_{\alpha\in H^1(K,M)}\prod_v\ell_v^{w_v}(\alpha_v).
\end{equation}
To relate this to integral orbit counts, put $\Gamma_\Lambda=\Stab_{\Gamma(\AAA_K)}\Lambda$ and, for $\gamma\in\Gamma(\AAA_K)$,
\[
  X_\gamma=X(K)\intsec\gamma^{-1}\Lambda,
  \qquad
  \Gamma_\gamma=\Gamma(K)\intsec\Gamma_{\gamma^{-1}\Lambda}.
\]
Then Lemma 2.23 of \cite{OCubic} gives
\begin{lem}\label{lem:total_cl_num_quartic}
  Let $w$ be a global weighting on $X$. Then
  \begin{equation}\label{eq:total_cl_num_quartic}
    \frac{h^w(X)}{\size{H^0(K,M)}}
    =\sum_{\gamma\in\Gamma_\Lambda\bs\Gamma(\AAA_K)/\Gamma(K)}
      \ \sum_{x\in\Gamma_\gamma\bs X_\gamma}
      \frac{w(\gamma x)}{\size{\Stab_{\Gamma_\gamma}x}}.
  \end{equation}
  In particular, if $\Gamma$ has class number one, the right-hand side is the usual weighted count of integral orbits.
\end{lem}

\begin{defn}[Definition 2.25 of \cite{OCubic}]\label{defn:dual_quartic}
  Let $K_v$ be a local field. For $i=1,2$, let $(\Gamma^{(i)},V^{(i)},X^{(i)},X_0^{(i)})$ be composed varieties over $K_v$ with Tate-dual stabilizers $M^{(1)}$ and $M^{(2)}$, equipped with lattices $\Lambda_v^{(i)}$ and local weightings $w_v^{(i)}$. We say that the two weightings are \emph{dual} with \emph{duality constant} $c_v$ if
  \begin{equation}\label{eq:general_local_refl_quartic}
    \widehat{\ell}_v^{w_v^{(1)}}=c_v\cdot\ell_v^{w_v^{(2)}},
  \end{equation}
  where
  \[
    \widehat f(\beta)=\frac{1}{\size{H^0(K_v,M^{(1)})}}
      \sum_{\alpha\in H^1(K_v,M^{(1)})}f(\alpha)\langle\alpha,\beta\rangle.
  \]
  If both weightings are trivial, we say that the two composed varieties with chosen lattices are \emph{naturally dual}.
\end{defn}

\begin{thm}[\textbf{local-to-global reflection engine}; Theorem 2.26 of \cite{OCubic}]\label{thm:main_compose_quartic}
  For $i=1,2$, let $(\Gamma^{(i)},V^{(i)},X^{(i)},X_0^{(i)})$ be Hasse composed varieties over a global field $K$ whose stabilizers $M^{(1)}$ and $M^{(2)}$ are Tate duals. Equip them with adelic lattices $\Lambda^{(i)}$ and global weightings $w^{(i)}=\prod_vw_v^{(i)}$. Suppose that at every place $v$ the local weightings are dual with constant $c_v$. Then $c_v=1$ for almost all $v$, and
  \begin{equation}\label{eq:main_compose_quartic}
    \frac{h^{w^{(2)}}(X^{(2)})}{\size{H^0(K,M^{(2)})}}
    =\left(\prod_vc_v^{-1}\right)
      \frac{h^{w^{(1)}}(X^{(1)})}{\size{H^0(K,M^{(1)})}}.
  \end{equation}
\end{thm}
This theorem follows from Poisson summation on adelic cohomology, using the middle terms of the Poitou--Tate exact sequence with the product of the local Tate pairings. For more details, see \cite{OCubic}.

Hence the quest for reflection theorems devolves into proving results of the type \eqref{eq:general_local_refl_quartic}. This is the \emph{local step,} which can be carried out at one place of $K$ at a time.

In proving cubic O-N, the local step is already rather involved. In the local step for pairs of ternary quadratic forms, we meet the formidable combinatorics of orders in a quartic extension, which Nakagawa himself studied in a volume of these Memoirs \cite{NakOrders}. Even with the more modest goal of counting by discriminant, a coarser invariant than the resolvent, the resulting generating functions are highly intricate, and there are a large number of cases, above all at the wild prime $2$. This is the main reason why the existing methods have struggled to prove quartic O-N.

The heart of quartic O-N is thus the following local duality:

\begin{thmstar}[\textbf{``Local Quartic O-N''}] \label{thmstar:quartic_local}
  Let $K$ be a nonarchimedean local field, $\ch K \neq 2$, and let $C$ be an order in a cubic \'etale algebra $R$ over $K$. For $\tau$ a divisor of $2$, let $\Lambda_{\tau,C}$ be the lattice of pairs of ternary quadratic forms parametrizing $\tau$-traced orders with reduced resolvent $C$. Then the two chosen lattices $\Lambda_{\tau,C}$ and $\Lambda_{2\tau^{-1},C}$ are naturally dual with duality constant $q^{2v_K(\tau)}$. In other words, their local lattice counts
  \[
    g_{\tau,C}:H^1(K,M_R)\longto\NN
  \]
  satisfy the local reflection theorem
  \[
    \widehat g_{\tau,C}=\Size{\OO_K/\tau\OO_K}^2\,g_{2\tau^{-1},C}.
  \]
  The following cases are rigorously known:
  \begin{enumerate}[$($a$)$]
    \item\label{qloc:tame} When $\ch k_K \neq 2$.
    \item\label{qloc:e_eq_1} When $K$ is unramified over $\QQ_2$.
    \item\label{qloc:max} When $C = \OO_R$ is maximal.
  \end{enumerate}
\end{thmstar}

In this memoir, we use two main approaches to prove Theorem* \ref{thmstar:quartic_local}. In favorable cases, we produce a multijection between the quartic rings counted by each side of the claimed equality. In unfavorable cases, we compute the generating function that records how each side changes as a function of $K$, $\ttt$, and the various parameters determining the resolvent $C$. We then use computer calculations to show that the generating functions of each side are equal.

\section{Side results}
In the quest to prove reflection theorems, our methods turned up a few side results that we deemed sufficiently pretty to mention here. They could also serve as starting points for further generalization.

\subsection{The average value of a quadratic character on a box}
\begin{thm}\label{thm:char_box}
  Let $K \supseteq \QQ_2$ be a finite extension, and let $R \supset K$ be a cubic \'etale extension. Let
  \[
  \chi : R^\cross \to \{\pm 1\}
  \]
  be a character, that is, a group homomorphism, such that $\chi(a) = 1$ for all $a \in K^\cross$. (Any such character is automatically continuous and can be put in the form
  \[
  \chi_\alpha(\xi) = \left\langle\alpha \cdot N_{R/K}(\alpha), \xi\right\rangle,
  \]
  where $\alpha \in R^\cross$ and $\left\langle\bullet, \bullet\right\rangle$ is the Hilbert symbol on $R$.) Let $B$ be an $\OO_K$-sublattice contained in the Jacobson radical of $\OO_R$. (That is, $B$ is a subgroup of $\OO_R$ of finite index closed under multiplication by $\OO_K$, and all elements of $B$ have positive valuation at every field factor of $R$. Note that $B$ need \textsc{not} be closed under multiplication.) Then the average value
  \[
  \frac{1}{\int_B 1} \int_B \chi(1 + \xi) \, \mathrm{d}\xi
  \]
  takes on one of the following values:
  \begin{itemize}
    \item $0$
    \item $1$
    \item $\pm q^{-i}$ for some $i$, $1 \leq i \leq e$. (Here $q = \size{k_K}$ is the size of the residue field, and $e = v_K(2)$ is the absolute ramification index.)
  \end{itemize}
\end{thm}
\begin{rem}
  We assume that $K \supseteq \QQ_2$, since for a $p$-adic field with $p$ odd, $\chi(1 + \xi) = 1$ for all $\xi$ in the Jacobson radical, making the integrand identically $1$.
\end{rem}
\begin{rem}
  The condition that $B$ lie within the Jacobson radical appears to be necessary; integrating $\chi$ over regions containing zero-divisors does not seem to give such simple results.
\end{rem}
\begin{rem}
  By talking about average value, we evade the question of normalizing the Haar measure on $B$.
\end{rem}

\subsection{Conductor rings}

\begin{defn}
  Let $R$ be an \'etale algebra over a local (resp.\ global) field $K$, and let $E/R$ be an abelian extension whose Artin map $\psi_{E/R} : R^\cross \to \Gal(E/R)$ (resp. $\psi_{E/R} : I(R, \mm) \to \Gal(E/R)$) vanishes on the base $K^\cross$ (resp.\ $I(K, \mm)$). Call an $\OO_K$-order $\OO \subseteq R$ an \emph{admissible ring} for $E$ if $\psi(\OO^\cross) = 0$. If there is a unique maximal admissible ring, we call it the \emph{conductor ring} of the extension $E/R$.
\end{defn}

Note that $\OO$ is an admissible ring of $E$ if and only if $E$ is contained in the ring class field of $\OO$. Admissible rings always exist, but a conductor ring need not (see Example \ref{ex:no_cdr_ring} below). However, in the case where $R/K$ is cubic and $E/R$ quadratic, the conductor ring has a neat characterization.

Recall that, in addition to the cubic resolvent $R$, a quartic \'etale algebra $L/K$ has a natural sextic resolvent $S$ coming from the map $\cS_4 \to \cS_6$ that sends a permutation of $\{1,2,3,4\}$ to the corresponding permutation of its $2$-element subsets. $S$ is naturally a quadratic \'etale extension of $R$ with the same Kummer element $\delta \in R^{N=1}/\left(R^{N=1}\right)^2$ that parametrizes $L$ in Proposition \ref{prop:Kummer_quartic}.

\begin{thmstar}[\textbf{Conductor rings}]\label{thmstar:cond_ring}
  Let $L/K$ be a quartic \'etale algebra over a global number field or a $p$-adic field. Assume that Theorem* \ref{thmstar:quartic_local} holds over $K$.
  Let $R$ and $S$ be the cubic and sextic resolvent algebras of $L$, respectively. Then the cubic resolvent ring $C_0 \subseteq R$ of $\OO_L$ is the conductor ring of $S/R$.
\end{thmstar}

\section{Conjectures on more general trace conditions} \label{sec:doubly_traced}
It is an open question to classify, analogously to the work of Osborne on binary cubics \cite[Theorem 2.2]{Osborne}, the lattices in the space $V_K$ of pairs of ternary quadratic forms over $K$ invariant under $\GL_2(\OO_K) \cross \GL_3(\OO_K)$. This problem is not quite the most relevant to us because the factor $\GL_2(\OO_K)$, which changes the coordinates of the resolvent, is not relevant when we count rings with a \emph{fixed} resolvent, as we have done in this memoir. We therefore look at lattices invariant under $\GL_3(\OO_K)$ alone. By inspection we find the following examples, which we conjecture exhaust all of them:
\begin{conj}
  Let $K$ be a number field, and let $\cV(\OO_K)$ be the lattice of pairs of ternary quadratic forms over $\OO_K$.
  A primitive, $\GL_3(\OO_K)$-invariant lattice $L$ in $\cV(\OO_K)$ is $\GL_2(\OO_K)$-equivalent to one whose completions $L_\pp$ are as follows:
  \begin{enumerate}[(a)]
    \item \label{it:intquartic} For $\pp|2$, the lattices $L_{\pp,t,s}$ of $(\pp^t,\pp^s)$-traced pairs of ternary quadratics
    \[
    L_{t,s} = L_{\pp,t,s} = \left\{ \left( \sum_{1 \leq i \leq j \leq 3} a_{ij} x_i x_j, \sum_{1 \leq i \leq j \leq 3} b_{ij} x_i x_j \right) :
    a_{ij} \in \pp^t, b_{ij} \in \pp^s\, (i < j)
    \right\}
    \]
    for $0 \leq t \leq s \leq e = v_p(2)$;
    \item For all other $\pp$, the maximal lattice $V_{\OO_{K,\pp}}$ only.
  \end{enumerate}
  In other words, such $L$ are of the form
  \begin{align*}
  L_{\ttt,\sss} &= \cc \bigcap_{\pp | 2} L_{\pp,\ttt_\pp,\sss_\pp} \\
  &=
  \left\{ \left( \sum_{1 \leq i \leq j \leq 3} a_{ij} x_i x_j, \sum_{1 \leq i \leq j \leq 3} b_{ij} x_i x_j \right) :
  a_{ij} \in \cc\ttt, b_{ij} \in \cc\sss\, (i < j)
  \right\}
  \end{align*}
  for ideals $(1) \supseteq \ttt \supseteq \sss \supseteq (2)$ and arbitrary $\cc$.
\end{conj}
As in the cubic case, the same lattices $L_{\ttt,\sss}$ necessarily appear in the analogues of $\cV(\OO_K)$ with any Steinitz class $\aaa$, appropriately adjusting the ideals that the $a_{ij}$, $b_{ij}$ must lie in.

For now, we look at the ($2$-adic) local case $K/\QQ_2$. An element of $L_{t,s}$ can be visualized as a pair of symmetric matrices
\[
(\cA, \cB) =
\left(
\begin{bmatrix}
  a_{11} & \frac{1}{2} a_{12} & \frac{1}{2} a_{13} \\
  \frac{1}{2} a_{12} & a_{22} & \frac{1}{2} a_{23} \\
  \frac{1}{2} a_{13} & \frac{1}{2} a_{23} & a_{33}
\end{bmatrix},
\begin{bmatrix}
  b_{11} & \frac{1}{2} b_{12} & \frac{1}{2} b_{13} \\
  \frac{1}{2} b_{12} & b_{22} & \frac{1}{2} b_{23} \\
  \frac{1}{2} b_{13} & \frac{1}{2} b_{23} & b_{33}
\end{bmatrix}
\right)
\]
with diagonal entries $a_{ii}, b_{ii} \in \OO_K$ and off-diagonal entries $a_{ij} \in 2^{-1}\pi_K^{t}\OO_K$, $b_{ij} \in 2^{-1}\pi_K^{s}\OO_K$. It is easy to check that the cubic resolvent $f(x,y) = 4 \det (\cA x - \cB y)$ has the form
\[
f(x,y) = \pi^{2 s} a x^3 + \pi^{2 s} b x^2 y + \pi^{2 t} c x y^2 + \pi^{2 t} d y^3
\]
with discriminant $D \in \pi^{4 t + 4 s}$. Now $(A,B)$ parametrizes a quartic ring with a resolvent $\Phi : \OO \to C$, and $f(x,y)$ is the index form of $C$. The divisibility conditions on $f(x,y)$ can be interpreted as a non-maximality condition on $C$: since $t \leq s$,
\[
\frac{f(\pi^{-t}x,\pi^{-t}y)}{\pi^{-2 t}} = \pi^{2 s - 2 t} a x^3 + \pi^{2 s - 2 t} b x^2 y + c x y^2 + d y^3
\]
and
\[
\frac{f(\pi^{-t}x,\pi^{-s}y)}{\pi^{-t-s}} = a x^3 + \pi^{s - t} b x^2 y + c x y^2 + \pi^{s - t} d y^3
\]
are integral and thus are the index forms of certain overrings $C_1$, $C_2$ which we call the \emph{reduced resolvent} and the \emph{reduced coresolvent.} Appropriately lifting the basis $\left\langle\bar\xi, \bar\eta\right\rangle$ of $C/\OO_K$ in which $(A,B)$ and hence $f$ are written, we have $C = \left\langle1, \xi, \eta\right\rangle$, $C_1 = \left\langle1, \pi^{-t}\xi, \pi^{-t}\eta\right\rangle$, and $C_2 = \left\langle1, \pi^{-t}\xi, \pi^{-s}\eta\right\rangle$. Note that $C_1$ is a subring of $C_2$ of index $\pi^{s-t}$ and moreover is a \emph{unidirectional} subring in the sense that $C_1/C_2 \isom \OO_K/\pi^{s-t}$ is generated by one element. The integrality properties of $A$ and $B$ translate readily into relations between the resolvent and the rings $C_1$ and $C_2$.
\begin{defn}
  Let $\OO_K$ be a Dedekind domain, and let $L/K$ be a quartic algebra over its field of fractions with resolvent $\Phi : L \to R$. Let $C_1 \subseteq C_2 \subseteq R$ be a pair of subrings with $C_2/C_1 \isom \OO_K/\dd$ for some divisor $\dd$ of $2$ in $\OO_K$. Let $\ttt$ be a divisor of $2\dd^{-1}$ and let $\sss = \ttt \dd$.

  An order $\OO \subseteq L$ is \emph{$(\ttt, \sss)$-traced with reduced resolvent $C_1$ and reduced coresolvent $C_2$} if
  \begin{enumerate}[(a)]
    \item $C = \OO_K + \ttt^2 C_1$ is a resolvent for $\OO$, that is, $\disc \OO = \disc C$ and $\Phi(\OO) \subseteq C$;
    \item The associated bilinear form
    \[
    \tilde \Phi(x,y) = \Phi(x + y) - \Phi(x) - \Phi(y) = xy' + x'y + x''y''' + x'''y''
    \]
    maps $\OO \cross \OO$ into $\OO_K + \ttt^3 \sss C_2$.
  \end{enumerate}
\end{defn}
If $C$ has a basis $\left\langle1, \xi, \eta\right\rangle$ for which $C_1 = \left\langle1, \ttt^{-1}\xi, \ttt^{-1}\eta\right\rangle$ and $C_2 = \left\langle1, \ttt^{-1}\xi, \sss^{-1}\eta\right\rangle$, as always happens when $\OO_K$ is a PID, then this is easily seen to be equivalent to the condition that $\Phi : \OO \to C$ is a resolvent whose matrix under this basis is in $L_{\ttt,\sss}$.

The functional equation for the Shintani zeta functions on $V_\ZZ$ was stated and proved by Sato and Shintani (\cite{SatoShintani}). It relates pairs of integer-coefficient ternary quadratic forms (over $\ZZ$) with pairs of integer-matrix forms, that is, $L_{(1),(1)}$ to $L_{(2),(2)}$ in our notation. In contrast to the lattice of cubic forms, there is no $\SL_3(K)$-invariant inner product in $V_K$: as a representation of $\SL_3(K)$, $V_K$ is not isomorphic to its dual. However, $V_K$ \emph{is} isomorphic to its dual twisted by the automorphism of $\SL_3(K)$ given by inverse transpose, and under this duality, it is easy to see that $L_{\ttt,\sss}$ is interchanged with $L_{2\sss^{-1},2\ttt^{-1}}$. Therefore, it is a pretty conjecture that the corresponding composed varieties are naturally dual.

\begin{conj} \label{conj:O-N_local_quartic_traced}
  Let $K$ be a $2$-adic local field and $0 \leq t \leq s \leq e = v_K(2)$ be integers. Let $C_1 \subseteq C_2$ be orders in a cubic $K$-algebra $R$ such that $C_2/C_1 \cong \OO_K/\pi^{s-t}$ is an $\OO_K$-module with one generator. For each quartic $K$-algebra $L$ with resolvent $R$, denote by $g(L, C_1, C_2, t, s)$ the number of $(t,s)$-traced quartic rings in $L$ with reduced resolvent and coresolvent $C_1$ and $C_2$, respectively. Then the dual of $g(L, C_1, C_2, t, s)$, considered as a function of $L$, is
  \begin{equation}\label{eq:O-N_local_quartic_wild}
    \hat g(L, C_1, C_2, t, s) = q^{t+s} g(L, C_1, C_2, e-s, e-t)
  \end{equation}
  where $q = \Size{k_K}$ is the order of the residue field.
\end{conj}
Note that the quartic rings counted on either side of \eqref{eq:O-N_local_quartic_wild} actually have resolvents $\OO_K + \pi^{t}C_1$ and $\OO_K + \pi^{e - s}C_1$ and discriminants $\pi^{2 t} \disc C_1$ and $\pi^{2 e - 2 s} \disc C_1$, respectively. If $t + s = e$, this conjecture asserts the self-duality (up to the correct scaling) of the indicated local weighting.

If this conjecture is true, then by our reflection engine, we immediately get the following corresponding global result.
\begin{conj} \label{conj:O-N_quartic_doubly_traced}
  Let $K$ be a number field. Let $C_1 \subseteq C_2$ be orders in a cubic $K$-algebra $R$ such that $C_2/C_1 \cong \OO_K/\dd$ is an $\OO_K$-module with one generator. Let $\ttt,\sss \subseteq \OO_K$ be ideals such that $\sss = \dd \ttt \mid (2)$. Let $h(C_1, C_2, \ttt, \sss)$ count the number of $(\ttt, \sss)$-traced quartic rings with reduced resolvent and coresolvent $C_1$ and $C_2$, respectively, each weighted by the reciprocal of its number of resolvent-preserving automorphisms. Let $h^{\ntc}(C_1, C_2, \ttt, \sss)$ count the subset of the foregoing that are ntc, weighted in the same way. Then
  \[
  h(C_1, C_2, \ttt, \sss) = \frac{N(\ttt\sss)}{2^{r_\infty}} \cdot h^{\ntc}(C_1, C_2, 2\sss^{-1}, 2\ttt^{-1}),
  \]
  where $r_\infty$ is the number of real places of $K$ over which $R$ is not totally real plus twice the number of complex places of $K$.
\end{conj}

A conjecture for general lattices of $2 \times n \times n$ boxes can be similarly stated; we omit the details.

\section{Structure of this work}
This memoir is divided into a ``conceptual'' and a ``computational'' part.

The conceptual part (chapters \ref{cha:glo}--\ref{cha:intro_count}) begins with Chapter \ref{cha:glo}, where we use harmonic analytic tools from \cite{OCubic} to reduce all of the quartic theorems to a single reflection theorem over a local field (Theorem* \ref{thmstar:quartic_local}). In Chapter \ref{cha:tame}, we prove this reflection theorem in the ``tame'' case, that is, away from $2$. In Chapter \ref{cha:2xnxn}, we prove Theorem \ref{thm:2xnxn} on $2\times n\times n$ boxes in the cases claimed.

The computational part (pages \pageref{cha:intro_count}--\pageref{cha:zone_examples}) tackles the problem of enumerating, in sufficient explicitness, the rings in a quartic field having given cubic resolvent. In Chapter \ref{cha:intro_count}, we establish some basic structure of the situation. In Chapter \ref{cha:conic}, we solve conics occurring in the enumeration, which is tantamount to computing the Igusa zeta functions of various quadratic forms in three variables. (Igusa zeta functions of quadratic forms are not easy to compute in general, especially in residue characteristic $2$, and our results go beyond the cases considered in \cite{CKWIgusa}.) In Chapter \ref{cha:boxgroups}, we study \emph{boxgroups,} subgroups of the group of all quartic \'etale algebras with fixed resolvent that share certain families of rings. In Chapters \ref{cha:xi1} and \ref{cha:xi2}, respectively, we count the ways to choose the first and second basis vector of the quartic ring, besides $1$ (as the last basis vector carries no information). We conclude with some miscellaneous considerations linking the attached code to the notational world of the work (Chapter \ref{cha:code}). We include in Appendix \ref{cha:zone_examples} some examples of the zone totals computed in Chapters \ref{cha:xi1}--\ref{cha:xi2}.

\section{Notation}
The notation will in general match that of \cite{OCubic}. Some specific points:

We denote by $\NN$ and $\NN^+$, respectively, the sets of nonnegative and of positive integers.

If $P$ is a statement, then
\[
\One_P = \begin{cases}
  1 & \text{$P$ is true} \\
  0 & \text{$P$ is false}.
\end{cases}
\]
If $S$ is a set, then $\One_S$ denotes the characteristic function $\One_S(x) = \One_{x \in S}$.

An \emph{algebra} will always be commutative and of finite rank over a field, while a \emph{ring} or \emph{order} will be a finite-rank, torsion-free ring over a Dedekind domain, containing $1$. An order need not be a domain.

If $\OO_K$ is a Dedekind domain, $K$ is its field of fractions, $V$ is an $n$-dimensional vector space over $K$ and $A,B \subseteq V$ are two full-rank lattices, we denote by the \emph{index} $[A:B]$ the unique fractional ideal $\cc$ such that
\[
\Lambda^n A = \cc \Lambda^n B
\]
as lattices in the top exterior power $\Lambda^n V$. Note that if $A \supseteq B$ and
\[
A/B \cong \OO_K/\cc_1 \oplus \cdots \oplus \OO_K/\cc_r
\]
as $\OO_K$-modules, then
\[
[A:B] = \cc_1\cc_2 \cdots \cc_r.
\]
If $L$ is a $K$-algebra, $\OO \subseteq L$ is an order, and $\aaa \subseteq L$ is a fractional ideal, the index $[\OO : \aaa]$ is called the \emph{norm} of $\aaa$ and will be denoted by $N_\OO(\aaa)$ or, when the context is clear, by $N(\aaa)$. The norm is multiplicative for $\OO$-ideals that are invertible (i.e.\ locally principal), but not for general $\OO$-ideals.

If $K$ is a field, we denote by $\bar K$ the separable closure of $K$.

If $a,b \in L$ are elements of a local or global field, a separable closure thereof, or a finite product of the preceding, we write $a|b$ to mean that $b = c a$ for some $c$ in the appropriate ring of integers $\OO_L$. If $a|b$ and $b|a$, we say that $a$ and $b$ are \emph{associates} and write $a \sim b$. Note that $a$ and $b$ may be zero-divisors.

We will use the semicolon to separate the coordinates of an element of a product of rings. For instance, in $\ZZ \cross \ZZ$, the nontrivial idempotents are $(1;0)$ and $(0;1)$.

If $n$ is a positive integer, then $\zeta_n$ denotes a primitive $n$th root of unity in $\bar\QQ$, while $\bar\zeta_n$ denotes the $n$th root of unity
\[
\bar\zeta_n = \left(1; \zeta_n; \zeta_n^2; \ldots; \zeta_n^{n-1}\right) \in \bar\QQ^n.
\]

Throughout the proofs of the local reflection theorems, we will fix a local field $K$, its valuation $v = v_K$ (which has value group $\ZZ$), its residue field $k_K$ of order $q$, and a uniformizer $\pi = \pi_K$. Often the choice of uniformizer is immaterial, but not always (for example, a tamely ramified cubic extension will be described without loss of generality as $R = K[\sqrt[3]{\pi}]$.) We have the absolute ramification index $e = e_K = v_K(2)$; note that $e$ is either the absolute ramification index of $K$ or $0$. We fix a separable closure $\bar K$ of $K$ (actually, it is possible to do all proofs within a finite extension, but we do not do so). We let $\mm_K$ denote the maximal ideal of $K$, and likewise let $\mm_{\bar K}$ be the maximal ideal of the ring $\OO_{\bar K}$ of algebraic integers over $K$; note that $\mm_{\bar K}$ is not finitely generated. We also allow $v = v_K$ to be applied to elements of $\bar K$, the valuation being scaled so that its restriction to $K$ has value group $\ZZ$. We use the absolute value bars $\size{\bullet}$ for the corresponding metric (i.e.\ multiplicative valuation), whose normalization will be left undetermined. We use superscript parenthesized numbers or letters to denote coordinates; for instance, a general element of $\bar K^n$ is $x = (x^{(1)}; \ldots; x^{(n)})$. We denote the valuations of such an element as follows:
\begin{align}
  \vec{v}(x) &= \left(v_K(x^{(1)}), \ldots, v_K(x^{(n)})\right) \in \QQ^n \label{eq:vec_v} \\
  v(x) &= \min_i v_K(x^{(i)}) \in \QQ. \label{eq:v_Kn}
\end{align}

If $R$ is a ring, we denote by $R^{N=1}$ the subgroup of units of norm $1$. Implicitly, the group operation is multiplication, so that $R^{N=1}[n]$, for instance, denotes the $n$th roots of unity of norm $1$ (of which there may be more than $n$ if $R$ is not a domain).

If $K$ is a local field, an \emph{$n$-pixel} is a subset of an affine or projective space over $\OO_K$ defined by requiring the coordinates to lie in specified congruence classes modulo $\pi^n$. For instance, in $\PP^2(\OO_K)$, a $0$-pixel is the whole space, which is subdivided into
\[
(q^2+q+1)q^{2n - 2}
\]
$n$-pixels for each $n \geq 1$.

If $f : H^1(K, M) \to \CC$ is a weighting on Galois cohomology over a local field, we normalize its Fourier transform $\hat f : H^1(K, M') \to \CC$, where $M' = \Hom(M, \mu)$ is the Tate dual, by
\[
  \hat f(\beta) = \frac{1}{\size{H^0(K,M)}} \sum_{\alpha \in H^1(K, M)} f(\alpha) \left\langle\alpha, \beta\right\rangle,
\]
where $\left\langle\bullet,\bullet\right\rangle$ is the $\mu$-valued Tate pairing. This ensures that, for $K$ nonarchimedean,
\[
  \hat{\hat{f}}(\alpha) = \frac{\size{H^1(K,M)}}{\size{H^0(K,M)} \cdot \size{H^0(K, M')}} f(-\alpha) = q^{v_K(\size{M})} f(-\alpha)
\]
where the last equality follows from a global Euler-characteristic computation.

\chapter{Global-to-local reduction}
\label{cha:glo}

The methods of this memoir build on those of \cite{OCubic}, in particular in the use of \emph{composed varieties} and a machinery (based on harmonic analysis of adelic cohomology) to reduce global questions to local ones. In this section we briefly review the definitions and then reduce all the various forms of quartic reflection in the introduction to a single lemma.

\section{Higher composition laws}

The space $V$, and hence $V_f$, figure in two higher composition laws: one parametrizing quartic rings with resolvent $R_f$, and one parametrizing self-balanced ideals in $R_f$. Both of these parametrizations will be important for us. We consider each in turn.

We first recall the parametrization of cubic rings.
\begin{thm}[Theorem 5.12 of \cite{OCubic}] \label{thm:hcl_cubic_ring_quartic} Let $A$ be a Dedekind domain with field of fractions $K$, $\ch A \neq 3$.
  \begin{enumerate}[$($a$)$]
    \item Cubic rings $\OO$ over $A$, up to isomorphism, are in bijection with cubic maps
    \[
    \Phi : M \to \Lambda^2 M
    \]
    between a two-dimensional $A$-lattice $M$ and its own Steinitz class, up to isomorphism, in the obvious sense of a commutative square
    \[
    \xymatrix{
      M_1 \ar[r]^\sim_{i} \ar[d]^{\Phi_1} & M_2 \ar[d]^{\Phi_2} \\
      \Lambda^2 M_1 \ar[r]^\sim_{\det i} & \Lambda^2 M_2.
    }
    \]
    The bijection sends a ring $\OO$ to the index form $\Phi : \OO/A \to \Lambda^2(\OO/A)$ given by
    \[
    x \mapsto x \wedge x^2.
    \]
    \item \label{cubic:field_q} If $\OO$ is nondegenerate, that is, the corresponding cubic $K$-algebra $L = K \tensor_A \OO$ is \'etale, then the map $\Phi$ is the restriction, under the Minkowski embedding, of the index form of $\bar{K}^3$, which is
    \begin{equation}\label{eq:k3_phi_q}
      \begin{aligned}
        \Phi : \bar{K}^3 / \bar{K} &\to \bar{K}^2 / \bar{K} \\
        (x;y;z) &\mapsto \big((x - y)(y - z)(z - x), 0\big).
      \end{aligned}
    \end{equation}
    \item \label{cubic:lift_q} Conversely, let $L$ be a cubic \'etale algebra over $K$. If $\bar{\OO} \subseteq L/K$ is a lattice such that
    $\Phi$ sends $\bar{\OO}$ into $\Lambda^2 \bar{\OO}$, then there is a unique cubic ring $\OO \subseteq L$ such that, under the natural identifications, $\OO/A = \bar{\OO}$.
  \end{enumerate}
\end{thm}

The basic method for parametrizing quartic orders is by means of \emph{cubic resolvent rings,} introduced by Bhargava in \cite{B3} and developed by Wood in \cite{WQuartic} and the author in \cite{ORings}.

\begin{defn}[\cite{ORings}, Definition 8.1; also a special case of \cite{WQuartic}, p.~1069]\label{defn:rsv}
  Let $A$ be a Dedekind domain, and let $\OO$ be a quartic algebra over $A$. A \emph{resolvent} for $\OO$ (``numerical resolvent'' in \cite{ORings}) consists of a rank-$2$ $A$-lattice $Y$, an $A$-module isomorphism $\Theta : \Lambda^3 (\OO/A) \to \Lambda^2 Y$, and a quadratic map $\Phi : \OO/A \to Y$ such that there is an identity of biquadratic maps
  \begin{equation}\label{eq:resolvent}
    \Theta(x \wedge y \wedge xy) = \Phi(x) \wedge \Phi(y)
  \end{equation}
  from $\OO \cross \OO$ to $\Lambda^2 Y$.
\end{defn}

We collect some basic facts about these resolvents.

\begin{thm}[\textbf{the parametrization of quartic rings}]\label{thm:hcl_quartic} The notion of resolvent for quartic rings has the following properties.
  \begin{enumerate}[$($a$)$]
    \item \label{quartic:to_ring} If $X$ is a rank-$3$ $A$-lattice and $\Theta : \Lambda^3 X \to \Lambda^2 Y$, $\Phi : X \to Y$ satisfy \eqref{eq:resolvent}, then there is a unique (up to isomorphism) quartic ring $\OO$ equipped with an identification $\OO/A \isom X$ making $(\OO, Y, \Theta, \Phi)$ a resolvent.
    \item \label{quartic:cubic_ring} There is a canonical (in particular, base-change-respecting) way to associate to a resolvent $(\OO, Y, \Theta, \Phi)$ a cubic ring $C$ and an identification $C/A \isom Y$ with the following property: For any element $x \in \OO$ and any lift $y \in C$ of the element $\Phi(x) \in C/A$, we have the equality
    \[
    \Theta(x \wedge x^2 \wedge x^3) = y \wedge y^2.
    \]
    It satisfies
    \[
    \disc C = \disc \OO.
    \]
    (Here the discriminants are to be seen as quadratic resolvent rings, as in \cite{ORings}; this implies the corresponding identity of discriminant ideals.)
    If $\OO$ is nondegenerate, then $C$ is unique.
    \item \label{quartic:exist} Any quartic ring $\OO$ has at least one resolvent.
    \item \label{quartic:max} If $\OO$ is maximal, the resolvent is unique (but need not be maximal).
    \item \label{quartic:count} The number of resolvents of $\OO$ is the sum of the absolute norms of the divisors of the \emph{content} of $\OO$, the smallest ideal $\cc$ such that $\OO = \OO_K + \cc \OO'$ for some order $\OO'$.
    \item \label{quartic:field} Let $(Y, \Theta, \Phi)$ be a resolvent of $\OO$ with associated cubic ring $C$, and let $K = \Frac A$. If the corresponding quartic $K$-algebra $L = K \tensor_A \OO$ is \'etale, then the cubic $K$-algebra $R = C \tensor_A K$ can be characterized as follows: The maps $\Theta$ and $\Phi$ are the restrictions, under the Minkowski embedding, of the unique resolvent of $\bar{K}^4$, which is $\bar{K}^3$ with the maps
    \begin{equation}\label{eq:k4_theta}
      \begin{aligned}
        \Theta : \Lambda^3(\bar{K}^4 / \bar{K}) &\to \Lambda^2(\bar{K}^3 / \bar{K}) \\
        (0;1;0;0) \wedge (0;0;1;0) \wedge (0;0;0;1) &\mapsto (0;1;0) \wedge (0;0;1)
      \end{aligned}
    \end{equation}
    and
    \begin{equation}\label{eq:k4_phi}
      \begin{aligned}
        \Phi : \bar{K}^4 / \bar{K} &\to \bar{K}^3 / \bar{K} \\
        (x;y;z;w) &\mapsto (xy + zw; xz + yw; xw + yz),
      \end{aligned}
    \end{equation}
    and $R$ is the preimage of $\Phi(L)$ in $\bar{K}^3$.
    \item \label{quartic:lift}Conversely, let $L$ be a quartic \'etale algebra over $K$ and $R$ its cubic resolvent. Let
    \[
    \Theta_K : \Lambda^3(L/K) \to \Lambda^2(R/K), \quad \Phi_K : L/K \to R/K
    \]
    be the resolvent data of $L$ as a (maximal) quartic ring over $K$. Suppose $\bar{\OO} \subseteq L/K$, $\bar{C} \subseteq R/K$ are lattices such that
    \begin{itemize}
      \item $\Phi_K$ sends $\bar{\OO}$ into $\bar{C}$,
      \item $\Theta_K$ maps $\Lambda^3 \bar{\OO}$ isomorphically onto $\Lambda^2 \bar{C}$.
    \end{itemize}
    Then there are unique quartic rings $\OO \subseteq L$, $C \subseteq R$ such that, under the natural identifications, $\OO/A = \bar{\OO}$, $C/A = \bar{C}$, and the restrictions of $\Theta_K$ and $\Phi_K$ make $\bar C$ a resolvent of $\bar\OO$.
  \end{enumerate}

\end{thm}
\begin{proof}
  \begin{enumerate}[$($a$)$]
    \item See \cite{ORings}, Theorem 8.3.
    \item See \cite{ORings}, Theorems 8.7 and 8.8.
    \item See \cite{ORings}, Corollary 8.6.
    \item This is a special case of the following part.
    \item See \cite{ORings}, Corollary 8.5.
    \item By base-changing to $K$, we see that $Y \tensor_A K = R/K$ is a resolvent for $L$. Since the resolvent is unique, it suffices to show that the cubic resolvent $R'$ defined in the theorem statement is a resolvent for $L$ also. The maps $\Theta$ and $\Phi$ defined in the theorem statement are seen, by symmetry, to restrict to maps of the appropriate $K$-modules. The verification of \eqref{eq:resolvent} and of the fact that the multiplicative structure on $R'$ is the right one can be checked at the level of $\bar{K}$-algebras.
    \item Letting $X = \bar{\OO}$, $Y = \bar{C}$ in part \ref{quartic:to_ring}, we construct the desired $\OO$ and $C$. By comparison to the situation under base-change to $K$, we see that $\OO$, $C$ naturally inject into $L$ and $R$ respectively. Uniqueness is obvious, as $\OO$ must lie in the integral closure $\OO_L$. \qedhere
  \end{enumerate}
\end{proof}

\begin{rem}
  The datum $\Theta$ of a resolvent carries no information, in the following sense. It is unique up to scaling by $c \in A^\cross$, and the resolvent data $(X, Y, \Theta, \Phi)$ and $(X, Y, c\Theta, \Phi)$ are isomorphic under multiplication by $c$ on $X$ and by $c^{2}$ on $Y$. If $A$ is a PID, indeed, neither $X$ nor $Y$ carries any information, and the entire data of the resolvent is encapsulated in $\Phi$, a pair of $3\cross 3$ symmetric matrices over $A$ (with formal factors of $1/2$ off the diagonal) defined up to the natural action of $\GL_3(A) \cross \GL_2(A)$. This establishes the close kinship with Bhargava's parametrization of quartic rings in \cite{B3}. However, it is useful to keep $\Theta$ around.
\end{rem}

In this memoir we only deal with nondegenerate rings, that is, those of nonzero discriminant, or equivalently, those that lie in an \'etale $K$-algebra. Consequently, all resolvent maps $\Theta$, $\Phi$ that we will see are restrictions of \eqref{eq:k4_theta} and \eqref{eq:k4_phi}. When quartic algebras are parametrized Kummer-theoretically, the resolvent map becomes very explicit and simple.

\begin{prop}[\textbf{explicit Kummer theory for quartic algebras}]\label{prop:Kummer_resolvent_quartic}
  Let $R$ be a cubic \'etale algebra over $K$ ($\ch K\neq2$), and let $\delta\in R^\cross$ have square norm. Let $\sqrt\delta\in\bar K^3$ be an element with $\sqrt{N(\delta)}=N(\sqrt\delta)\in K^\cross$, and define
  \[
    \kappa(\xi)=\left(\tr_{\bar K^3/\bar K}
      (\xi\omega\sqrt\delta)\right)_{
      \omega\in(\bar K^3)^{N=1}[2]}\in\bar K^4.
  \]
  Then $L=K+\kappa(R)$ is the quartic algebra corresponding to the Kummer class of $\delta$, and $\kappa$ maps $R$ bijectively onto the traceless hyperplane in $L$. Its resolvent is given explicitly by
  \begin{equation}\label{eq:Kummer_rsv_theta}
    \begin{aligned}
      \Theta:\Lambda^3(L/K)&\to\Lambda^3(R)\cong\Lambda^2(R/K)\\
      \kappa(\alpha)\wedge\kappa(\beta)\wedge\kappa(\gamma)
      &\mapsto16\sqrt{N(\delta)}\,
        \alpha\wedge\beta\wedge\gamma
    \end{aligned}
  \end{equation}
  and
  \begin{equation}\label{eq:Kummer_rsv_phi}
    \begin{aligned}
      \Phi:L/K&\to R/K\\
      \kappa(\xi)&\mapsto4\delta\xi^2.
    \end{aligned}
  \end{equation}
\end{prop}
\begin{proof}
  Since the resolvent is unique (over a field, any \'etale extension has content $1$), it suffices to prove that \eqref{eq:Kummer_rsv_theta} and \eqref{eq:Kummer_rsv_phi} define a resolvent. This can be done after extension to $\bar{K}$, and then it is enough to prove that \eqref{eq:Kummer_rsv_theta} and \eqref{eq:Kummer_rsv_phi} agree with the standard resolvent on $\bar{K}^4$, given in Theorem \ref{thm:hcl_quartic}\ref{quartic:field}. This is a routine calculation.
\end{proof}

Let $M$ be a $K$-Galois module of underlying group $\cC_2 \cross \cC_2$. Since
\[
\GA(M) \coloneqq M \rtimes \Aut(M) \isom \cS_4,
\]
the cohomology $H^1(K,M)$ parametrizes quartic algebras with fixed cubic resolvent $R$. It can also be described explicitly.
If $\ch K \neq 2$, Proposition \ref{prop:composed_quartic} shows that the orbits of $\Gamma(K)$ on $X(K)$ are naturally parametrized (once a base orbit is fixed) by a subset of $H^1(K, M)$. We take the base orbit to be the one corresponding to the reducible algebra $L = K \cross R$; then the parametrization takes a quite explicit form:
\begin{prop}\label{prop:Kummer_quartic}
  Let $K$ be a field with $\ch K\neq2$, and let $M$ be a $K$-Galois module with underlying group $\cC_2\cross\cC_2$. Let $R$ be the cubic \'etale algebra corresponding to $M$. Then there is a group isomorphism
  \[
    H^1(K,M)\isom R^{N=1}/\bigl(R^{N=1}\bigr)^2.
  \]
  A class represented by $\delta\in R^{N=1}$ corresponds to the quartic algebra $L/K$ generated by the image of the $K$-linear map
  \begin{align*}
    \kappa:R&\to\bar K^4\\
    \xi&\mapsto
      \left(\tr_{\bar K^3/\bar K}
      (\xi\omega\sqrt\delta)\right)_{
      \omega\in(\bar K^3)^{N=1}[2]},
  \end{align*}
  where $\sqrt\delta\in\bar K^3$ is chosen to have norm $1$, and
  \[
    (\bar K^3)^{N=1}[2]
    =\{(1;1;1),(1;-1;-1),(-1;1;-1),(-1;-1;1)\}
  \]
  is the set of square roots of $1$ in $\bar K^3$ having norm $1$. Indeed,
  \[
    L=K+\kappa(R).
  \]
  Moreover, this parametrization coincides with the cohomological parametrization of Proposition \ref{prop:composed_quartic}.
\end{prop}
\begin{proof}
  This follows easily from the classical method for solving the quartic; see Knus and Tignol \cite[Proposition 5.13]{QuarticExercises}.
\end{proof}

We deduce that $X_f$ is \emph{full} in the sense of Definition \ref{defn:fullHasse_quartic} since we get an orbit for every cohomology class.

We next prove some results about $\ttt$-traced rings.
\begin{prop}\label{prop:traced} Let $\OO$ be a quartic ring over a Dedekind domain $\OO_K$.
  \begin{enumerate}[$($a$)$]
    \item\label{traced:conds} $\OO$ is $\ttt$-traced if and only if
    \begin{enumerate}[$($i$)$]
      \item\label{traced:trace} $\ttt^{2}|\tr x$ for all $x\in \OO$;
      \item\label{traced:sq} $x^2 \in \OO_K + \ttt\OO$ for all $x\in \OO$.
    \end{enumerate}
    \item\label{traced:count} If $\OO$ is not an order in the \emph{trivial algebra} $K[\epsilon_1, \epsilon_2, \epsilon_3]/(\epsilon_i \epsilon_j)_{i,j=1}^3$, the number of $\ttt$-traced resolvents of $\OO$ is the sum of the absolute norms of the divisors of its \emph{$\ttt$-traced content,} which is the smallest ideal $\cc$ such that $\OO = \OO_K + \cc \OO'$ and $\OO'$ is also $\ttt$-traced.
    \item\label{traced:reduced} If $(\OO,Y,\Theta,\Phi)$ is a $\ttt$-traced resolvent with associated cubic ring $C$, then $\ttt^{2}|\ct(C)$, that is, $C = \OO_K + \ttt^{2}C_\ttt$ for some cubic ring $C_\ttt$. We call $C_\ttt$ a ``reduced resolvent'' of the $\ttt$-traced ring $\OO$. Also, $\ttt^8 | \disc \OO$.
  \end{enumerate}
\end{prop}
\begin{proof}
  \begin{enumerate}[$($a$)$]
    \item
    Since both statements are local at the primes dividing $2$, we can assume that $\OO_K$ is a DVR, and thus that $\ttt = (t)$ is principal. With respect to bases $(1 = \xi_0, \xi_1,\xi_2,\xi_3)$ for $\OO$ and $(1 = \eta_0, \eta_1, \eta_2)$ for a resolvent $C$, the structure constants $c_{ij}^k$ of the ring $\OO$, defined by
    \[
    \xi_i \xi_j = \sum_k c_{ij}^k \xi_k,
    \]
    are determined by the entries of the resolvent
    \[
    \Phi = \left([a_{ij}], [b_{ij}]\right)
    \]
    via the determinants
    \[
    \lambda^{ij}_{k\ell} = 2^{\One_{i\neq j} + \One_{k\neq \ell}} \begin{vmatrix}
      a_{ij} & a_{k\ell} \\
      b_{ij} & b_{k\ell}
    \end{vmatrix}
    \]
    and a set of formulas appearing in Bhargava \cite[equation (21)]{B3} and over a Dedekind domain by the author \cite[equation (12)]{ORings}:
    \begin{equation}\label{eq:c-lam}
      \begin{aligned}
        c_{ii}^j &= -\epsilon \lambda^{ii}_{ik} \\
        c_{ij}^k &= \epsilon \lambda^{jj}_{ii} \\
        c_{ij}^j - c_{ik}^k &= \epsilon \lambda^{jk}_{ii} \\
        c_{ii}^i - c_{ij}^j - c_{ik}^k &= \epsilon \lambda^{ij}_{ik},
      \end{aligned}
    \end{equation}
    where $(i,j,k)$ denotes any permutation of $(1,2,3)$ and $\epsilon = \pm 1$ its sign. (Here the nonappearance of some of the individual $c_{ij}^k$ on the left-hand side of \eqref{eq:c-lam} stems from the ambiguity of translating each $\xi_i$ by $\OO_K$, which does not change the matrix of $\Phi$.)

    Assume first that $\Phi : \OO/\OO_K \to C/\OO_K$ is $\ttt$-traced. Then
    \begin{equation}\label{eq:traced}
      \lambda^{ij}_{k\ell} \in \ttt^{\One_{i\neq j} + \One_{k\neq \ell}}.
    \end{equation}
    We then prove that the conditions \ref{traced:trace} and \ref{traced:sq} must hold:
    \begin{enumerate}[$($i$)$]
      \item The trace
      \begin{align*}
        \tr(\xi_1) &= c_{11}^1 + c_{12}^2 + c_{13}^3 \\
        &= \lambda^{12}_{13} + 2\lambda^{23}_{11} + 4 c_{13}^3 \\
        &\equiv 0 \mod \ttt^2,
      \end{align*}
      and likewise $\tr(\xi_2), \tr(\xi_3) \in \ttt^2$.
      \item The coefficients $c_{11}^i$ of ${\xi_1}^2$ satisfy:
      \begin{align*}
        c_{11}^2 &= \lambda_{13}^{11} \in \ttt
      \end{align*}
      and likewise for $c_{11}^3$; and then $c_{11}^1 \in \ttt$ also, since the trace $c_{11}^1 + c_{12}^2 + c_{13}^3 = \tr(\xi_1) \in \ttt^2 \subseteq \ttt$. So the desired relation $\xi^2 \in \OO_K + \ttt \OO$ holds when $\xi = \xi_1$, indeed $\xi = a_{1}\xi_{1}$ for any $a_1 \in \OO_K$. The same proof works for $\xi = a_2\xi_2$ or $\xi = a_3\xi_3$. Since the case $\xi = a_0 \in \OO_K$ is trivial and squaring is a $\ZZ$-linear operation modulo $2$, we get the result for all $\xi \in \OO$.
    \end{enumerate}
    Conversely, suppose that \ref{traced:trace} and \ref{traced:sq} hold. We first establish \eqref{eq:traced}. We have
    \begin{itemize}
      \item $\lambda^{11}_{13} = c_{11}^2 \in \ttt$
      \item $\lambda^{23}_{11} = c_{12}^2 - c_{13}^3 = \tr \xi_1 - c_{11}^1 - 2c_{13}^3 \in \ttt$
      \item $\lambda^{12}_{13} = c_{11}^1 - c_{12}^2 - c_{13}^3 = \tr \xi_1 - 2\lambda^{23}_{11} + 4c_{13}^3 \in \ttt^2$.
    \end{itemize}
    Permuting the indices as needed, this accounts for all the $\lambda^{ij}_{k\ell}$ about which \eqref{eq:traced} makes a nontrivial assertion.

    Now we work from the $\lambda^{ij}_{k\ell}$ back to the resolvent $(\cA,\cB)$. We may assume that $C$ is nontrivial (the trivial rings, one for each Steinitz class, are plainly $2$-traced with $(\cA,\cB) = (0,0)$.) Then, in the proof of \cite{ORings}, Theorem 8.4, the author established that there are vectors $\mu_{ij}$ in a two-dimensional vector space $V$ over $K$, unique up to $\GL_2(V)$, such that
    \[
    \mu_{ij} \wedge \mu_{k\ell} = \lambda^{ij}_{k\ell} \cdot \omega
    \]
    for some fixed generator $\omega \in \Lambda^2 V$. (The Pl\"ucker relations needed to recover the vectors from their pairwise wedge products follow from the associative law on $\OO$.) This $V$ is none other than $R/K$, the resolvent module of the quartic algebra $L = \OO \tensor_{\OO_K} K$, which admits the unique resolvent
    \begin{equation}\label{eq:rsv_field}
      \begin{aligned}
        \Phi(a_1\xi_1 + a_2\xi_2 + a_3\xi_3) &= \sum_{i\leq j} a_i a_j \mu_{ij} \\
        \Theta(\xi_1 \wedge \xi_2 \wedge \xi_3) &= \omega.
      \end{aligned}
    \end{equation}
    The resolvents of $\OO$ were found to be exactly the lattices $M$ containing the span $M_0$ of the six $\mu_{ij}$, with the correct index
    \[
    [M : M_0] = \cc = \left(\lambda^{ij}_{k\ell}\right)_{i,j,k,\ell} = \left(c_{ii}^j, c_{ij}^k, c_{ij}^j - c_{ik}^k, c_{ii}^i - c_{ij}^j - c_{ik}^k : i \neq j \neq k \neq i\right),
    \]
    where $\cc$ is the content ideal of $\OO$. By inspection of \eqref{eq:rsv_field}, such a lattice $M$ is $\ttt$-traced if and only if it contains the span $\tilde M_0$ of the six vectors
    \[
    \tilde\mu_{ij} = t^{-\One_{i \neq j}} \mu_{ij}.
    \]
    Condition \eqref{eq:traced} says that the $\tilde\mu_{ij} \wedge \tilde\mu_{k\ell}$ are still integer multiples of $\omega$. The index required of a lattice $M \supseteq \tilde M_0$ is therefore the integral ideal
    \[
    [M : \tilde M_0] = \tilde \cc = \left(\tilde\lambda^{ij}_{k\ell}\right)_{i,j,k,\ell}.
    \]
    Thus such an $M$ exists, finishing the proof of \ref{traced:conds}.

    \item It suffices to prove that $\tilde \cc$ is the $\ttt$-traced content of $\OO$. Locally, write $\OO = \OO_K + a \OO'$. The structure coefficients of $\OO'$ are those of $\OO$ divided by $a$, so the corresponding ideal $\tilde \cc$ is divided by $a$ as well. Thus $\OO'$ is again $\ttt$-traced exactly when $a \mid \tilde \cc$, proving the claim.

    \item We can again reduce to the case that $\OO_K$ is a DVR so $\OO$ has an $\OO_K$-basis. Recall that the index form of the resolvent $C$ is given by
    \[
    f(x,y) = 4 \det (\cA x + \cB y)
    \]
    (\cite{B3}, Proposition 11; \cite{ORings}, Theorem 8.7).
    If $\cA$ and $\cB$ have off-diagonal entries in $2^{-1} \ttt$, it immediately follows that $f$ is divisible by $t^2$, so $t^2 \mid \ct(C)$. Consequently $\disc \OO = \disc C$, being quartic in the coefficients of $f$, is divisible by $t^8$. \qedhere
  \end{enumerate}
\end{proof}

We now turn to self-balanced ideals.

\begin{thm}[\textbf{self-balanced ideals in the quartic setting}] \label{thm:hcl_quartic_sbi}
  Let $\OO_K$ be a Dedekind domain, $\ch K \neq 2$, and let $R$ be a cubic \'etale algebra. A \emph{self-balanced triple} in $R$ is a triple $(C, I, \delta)$ consisting of a cubic order $C \subseteq R$, a fractional ideal $I$ of $C$, and a scalar $\delta \in (K C)^\cross$ satisfying the conditions
  \begin{equation} \label{eq:quartic_sbi}
    \delta I^2 \subseteq C, \quad N(I) = (t) \text{ is principal}, \textand N(\delta) t^2 = 1,
  \end{equation}
  Fix an order $C\subseteq R$ and a scalar $\delta \in R^\cross$ with $N(\delta)$ a square $t^{-2}$. Then the mapping
  \begin{equation}\label{eq:sbi_quartic}
    I \mapsto \OO = \OO_K + \kappa(I)
  \end{equation}
  defines a bijection between
  \begin{itemize}
    \item self-balanced triples of the form $(C, I, \delta)$, and
    \item subrings $\OO \subseteq L$ of the quartic algebra $L = K + \kappa(R)$ corresponding to the Kummer element $\delta$, such that $\OO$ is $2$-traced with reduced resolvent $C$.
  \end{itemize}
\end{thm}
\begin{proof}
  The proof is very similar to that of Proposition 5.13
  of \cite{OCubic}, so we simply summarize the main points. The linear isomorphism $\kappa$ establishes a bijection between lattices $I\subseteq R$ and $\kappa(I) \subseteq L/K$. We wish to prove that the following conditions are equivalent:
  \begin{enumerate}[$($i$)$]
    \item\label{sbi:bal} $(C, I,\delta)$ is balanced;
    \item\label{sbi:tr} $\kappa(I)$ is the projection of a $2$-traced order with reduced resolvent $C$.
  \end{enumerate}

  First note that either of these conditions uniquely specifies
  \[
  [C : I] = (t),
  \]
  the former by the balancing condition $N(I) = (t)$, and the latter by the $\Theta$-condition that $\OO$ have discriminant $256 \disc C$.

  We then prove that both conditions are equivalent to the bilinear map
  \begin{align*}
    \Phi : I \cross I &\to R/K \\
    (\alpha_1,\alpha_2) &\to 4\delta \alpha_1\alpha_2
  \end{align*}
  taking values in $4C/K$.

  For \ref{sbi:bal}, this follows from the parametrization of \emph{balanced} pairs of ideals by $2\times 3\times 3$ boxes performed over $\ZZ$ by Bhargava \cite[Theorem 2]{B2} and over a general base by Wood \cite[Theorem 1.4]{W2xnxn}.

  For \ref{sbi:tr}, the diagonal restriction of $\Phi$ is precisely the resolvent of $\kappa(I)$, by Proposition \ref{prop:Kummer_resolvent_quartic}. That $\Phi(\alpha,\alpha) \in 4C/K$ for each $\alpha \in I$ expresses the one condition remaining for $\kappa(I)$ to lift (by Theorem \ref{thm:hcl_quartic}\ref{quartic:lift}) to a quartic ring $\OO$ with resolvent $\OO_K + 4C$. Then, by definition, this resolvent is $2$-traced exactly when $\Phi$ itself has image in $4C/K$.
\end{proof}
\section{Reducing the theorems to a local theorem}

In this section, we show that Theorem \ref{thm:main_compose_quartic} and Lemma \ref{lem:total_cl_num_quartic} reduce all of Theorem \ref{thm:O-N_quartic_Z}, Theorems* \ref{thmstar:O-N_2x3x3}, \ref{thmstar:O-N_quartic}, \ref{thmstar:O-N_quartic_by_disc}, Theorem \ref{thm:BQ}, and Corollary \ref{cor:BQ} to a single local statement, Theorem* \ref{thmstar:quartic_local}:

\begin{proof}[Proof of Theorem \ref{thmstar:O-N_quartic}]
  Let $\aaa$ be the Steinitz class of $C$, and let $f$ be its corresponding binary cubic form. The higher composition law of Theorem \ref{thm:hcl_quartic} identifies the counts in the theorem with weighted integral orbit counts on the full, hence Hasse, composed variety $X_f$. We apply Theorem \ref{thm:main_compose_quartic} to the two adelic lattices whose finite components are $\Lambda_{\ttt,C}$ and $\Lambda_{2\ttt^{-1},C}$.

  At a finite place $v=\pp$, Theorem* \ref{thmstar:quartic_local} gives
  \[
    c_\pp=\begin{cases}
      N(\pp)^{2v_\pp(\ttt)},&\pp\mid2,\\
      1,&\pp\nmid2.
    \end{cases}
  \]
  At a complex place $c_v=1/4$, at a real place where $R_v\isom\RR\cross\CC$ we have $c_v=1/2$, and at a real place split in $R$, the trivial weighting is dual to the ntc selector with $c_v=1$. Thus
  \[
    \prod_vc_v=\frac{N(\ttt)^2}{2^{r_\infty}}.
  \]
  The group $\cG_\aaa$ has class number one by strong approximation, so Lemma \ref{lem:total_cl_num_quartic} identifies the normalized total class numbers in Theorem \ref{thm:main_compose_quartic} with the desired weighted orbit counts. Since the stabilizer module is the same on both sides, the $H^0$-factors cancel. Rearranging \eqref{eq:main_compose_quartic} gives
  \[
    h(C,\ttt)=\frac{N(\ttt)^2}{2^{r_\infty}}\,h^{\ntc}(C,2\ttt^{-1}),
  \]
  as claimed.
\end{proof}

Taking $K = \QQ$, $\aaa = \ttt = 1$, we get Theorem \ref{thm:O-N_quartic_Z} as a special case.

\begin{proof}[Proof of Theorem* \ref{thmstar:O-N_2x3x3}]
  We sum Theorem* \ref{thmstar:O-N_quartic} over all cubic rings $C$ of discriminant $\cD$, weighting each $C$ by the reciprocal of the number of \emph{orientation-preserving} automorphisms of $C$, which is the stabilizer of the corresponding form in $\SL(\OO_K \oplus \aaa)$. A pair is counted once for each orientation-preserving identification of its cubic resolvent with the chosen representative $C$; by orbit-stabilizer, the two stabilizer weights combine to give the reciprocal of its stabilizer in $\cG_\aaa$.
\end{proof}

\begin{proof}[Proof of Theorem \ref{thmstar:O-N_quartic_by_disc}]
  For this theorem we must do a bit more work owing to the non-uniqueness of the resolvent of a quartic ring.

  In Theorem* \ref{thmstar:O-N_2x3x3}, we studied $h_\ttt(\ttt^8\cD)$, which can be interpreted as the number of quartic rings $\OO$ equipped with a resolvent $C$ and an \emph{orientation,} that is, an identification $\Lambda^4 \OO \isom \aaa$ for which the discriminant is $\ttt^8\cD$. Every quartic ring admits two orientations (there are $\size{\OO_K^\cross}$-many identifications $\Lambda^4 \OO \isom \aaa$, but all but one and its negative yield a $D$ scaled by a different square of a unit). So $\frac{1}{2}h_\ttt(\ttt^8\cD)$ is the number of resolvents $(\OO,C,\Theta,\Phi)$ of discriminant $\ttt^8\cD$, up to isomorphism, each weighted by the reciprocal of its number of automorphisms.

  Let $h^1_\ttt(\ttt^8\cD)$ be the number of quartic rings $\OO$ of discriminant $\ttt^8\cD$ with $\ttt$-traced content $1$, weighted by $1/\lvert\Aut\OO\rvert$. If $\OO$ has $\ttt$-traced content $\cc$, then $\OO=\OO_K+\cc\OO'$, where $\OO'$ has $\ttt$-traced content $1$, discriminant $\cc^{-6}\ttt^8\cD$, and the same automorphism group as $\OO$. Thus
  \begin{equation}\label{eq:ring_count_by_content}
    h_\ttt^\circ(\ttt^8\cD)
    =\sum_{\cc^3\mid D\aaa^2}
      h_\ttt^1(\cc^{-6}\ttt^8\cD).
  \end{equation}
  On the other hand, a ring of $\ttt$-traced content $\cc$ has
  \[
    \sigma_1(\cc)=\sum_{\dd\mid\cc}N_{K/\QQ}(\dd)
  \]
  resolvents as maps out of $\OO$ (as pointed out in \cite{ORings}, end of Section 8). Hence
  \begin{equation}\label{eq:resolvent_count_by_content}
    \frac{1}{2}h_\ttt(\ttt^8\cD)
    =\sum_{\cc^3\mid D\aaa^2}
      \sigma_1(\cc)h_\ttt^1(\cc^{-6}\ttt^8\cD).
  \end{equation}

  Let $\mu_K$ be the M\"obius function on the nonzero ideals of $\OO_K$. Since $\sigma_1$ is the Dirichlet convolution of the norm map with the constant function $1$, M\"obius inversion of \eqref{eq:resolvent_count_by_content} gives
  \begin{equation}\label{eq:ring_from_resolvent_count}
    h_\ttt^\circ(\ttt^8\cD)
    =\frac{1}{2}\sum_{\cc^3\mid D\aaa^2}
      \mu_K(\cc)N(\cc)
      h_\ttt(\cc^{-6}\ttt^8\cD).
  \end{equation}
  The same formula holds with the ntc condition imposed. Plugging this in on both sides of \eqref{eq:O-N_by_disc} and applying Theorem* \ref{thmstar:O-N_2x3x3} term by term proves the result.
\end{proof}

\subsection{Binary quartic forms}

\begin{proof}[Proof of Theorem \ref{thm:BQ}]
We begin with Wood's characterization of quartic rings parametrized by binary quartic forms.
\begin{thm}[\cite{WoodBQ}, Theorem 1.1]\label{thm:WoodBQ}
  There is a natural, discriminant-preserving bijection between the set of
  $\GL_2(\ZZ)$-equivalence classes of binary quartic forms and the set of isomorphism classes of
  pairs $(Q, C)$ where $Q$ is a quartic ring and $C$ is a monogenized cubic resolvent of $Q$ (where
  isomorphisms are required to preserve the generator of $C$ modulo $\ZZ$).
\end{thm}
\begin{proof}
  We send a form $\Phi(x,y) = ax^4 + bx^3y + cx^2y^2 + dxy^3 + ey^4$ to the quartic ring (with associated resolvent) given by the pair of ternary quadratic forms
  \[
  (A_0, B) = \left(
  \begin{bmatrix}
    & & 1/2 \\
    & -1 & \\
    1/2 & &
  \end{bmatrix},
  \begin{bmatrix}
    a & b/2 & \\
    b/2 & c & d/2 \\
    & d/2 & e
  \end{bmatrix}\right).
  \]
  Monogenicity arises because the resolvent form $g(x,y) = 4 \det (A_0 x - B y)$ is monic, since $\det A_0 = 1/4$. Further details will be found in \cite{WoodBQ}.
\end{proof}
To apply this theorem, we need to know the number of automorphisms of the quartic ring corresponding to a given form:
\begin{lem} \label{lem:bq_aut}
  In this bijection, the group of resolvent-preserving automorphisms of a quartic ring is in natural isomorphism with the stabilizer (in $\PGL_2(\ZZ)$) of the corresponding form.
\end{lem}
\begin{proof}
  This follows easily from the method of proof of the preceding theorem. By \cite[Theorem 2.5]{WoodBQ}, we can choose bases for $Q$ and $C$ so that the corresponding pair of ternary quadratic forms has the form $(A_0, B)$ above. A resolvent-preserving automorphism is a change of variables $h \in \SL_3(\ZZ)$ that preserves both $A_0$ and $B$. By \cite{WoodBQ}, Lemma 3.2, $h = \epsilon_{A_0}(\tilde h)$ lies in the image of the embedding
  \begin{align*}
    \epsilon_{A_0} : \PGL_2(\ZZ) &\to \SL_3(\ZZ) \\
    \begin{bmatrix}
      a & b \\
      c & d
    \end{bmatrix}
    &\mapsto
    \frac{1}{ad - bc}
    \begin{bmatrix}
      a^2 & ab      & b^2 \\
      2ac & ad + bc & 2bd \\
      c^2 & cd      & d^2
    \end{bmatrix}
  \end{align*}
  By \cite[Theorem 3.1]{WoodBQ}, $h$ preserves $B$ if and only if $\tilde h$ preserves the binary quartic form $f$. This constructs the desired isomorphism.
\end{proof}

We fix a monic binary cubic form $g(x,y)$ and let $C = C_g = \ZZ[\xi]$ be the corresponding monogenized cubic ring, with generator $\xi$. Then in the notation of Theorem \ref{thm:O-N_quartic_Z},
\[
h(g) = \sum_{\substack{\text{quartic rings } Q \\ \text{with resolvent } C, \\ \text{up to $C$-isom}}} \frac{1}{\Size{\Aut_C Q}}
= \sum_{\substack{\text{binary quartics $f(x,y)$} \\ \text{with resolvent $g$,} \\ \text{up to $\PGL_2(\ZZ)$}}} \frac{1}{\Size{\Stab_{\PGL_2(\ZZ)}} f},
\]
and $h^+(g)$ counts the subset of these forms that have at least one real root. The quantities $h_2(4g)$ and $h_2^+(4g)$ appearing on the opposite side of Theorem \ref{thm:O-N_quartic_Z} are not so straightforward to interpret. Here we are counting pairs $(A,B)$ of integer symmetric matrices with $\det(Ax-By)=g$ (with a certain condition at $\infty$), so we need to classify integer symmetric matrices $A$ with $\det A=1$. There are, up to similarity, two:
\begin{lem}\label{lem:sym_mat}
  Every $3\times 3$ integer symmetric matrix $A$ with $\det A = 1$ is similar to
  \[
  A_1 = \begin{bmatrix}
    & & 1 \\
    & -1 & \\
    1 & &
  \end{bmatrix}
  \textor
  I = \begin{bmatrix}
    1 & & \\
    & 1 & \\
    & & 1
  \end{bmatrix}.
  \]
\end{lem}
\begin{proof}
  Integral unimodular lattices have been studied extensively; see \cite[Table~15.4 (indefinite case), Table~16.7 (definite case)]{ConwaySloane}.
\end{proof}

Hence $h_2(4g)$ decomposes into the pairs $(A,B)$ of ``type $A_1$'' and of ``type $I$'' according to the value of $A$ after an appropriate $\SL_3\ZZ$-transformation.

Pairs of type $A_1$ can be understood using the fact that $A_1 \sim A_0$ over $\QQ$. Namely, the transformation
\[
T = \begin{bmatrix}
  1 & & \\
  & 1 & \\
  & & 2
\end{bmatrix}
\]
satisfies $TA_0T^{\top} = A_1$. Let
\[
(A_1,B) = \left( A_1, \begin{bmatrix}
  a & b & c' \\
  b & c & d \\
  c' & d & e
\end{bmatrix} \right).
\]
Then the pair
\[
\left( T^{-1}A_1\left( T^{-1} \right)^{\top} , T^{-1} B \left( T^{-1} \right)^{\top} \right) = (A_0, B')
\]
is determined up to $\SO(\QQ, A_0)$ by $(A_1, B)$, and the corresponding binary quartic
\[
f(x,y) = a x^4 + 2b x^3 y + (c + c') x^2 y^2 + d x y^3 + \frac{1}{4} e y^4
\]
is determined up to $\PGL_2(\QQ)$. Clearing the fractions, we get a $48441$-form
\[
4f(x,y) = 4ax^4 + 8b x^3 y + 4(c + c') x^2 y^2 + 4d x y^3 + e y^4
\]
of resolvent $g^{[4]}$.

The $48441$-forms do not naturally have an action by $\PGL_2\ZZ$, but rather by a group that we can reveal as $\SO(\ZZ, A_1)$:

\begin{lem}
  We have $\SO(\QQ, A_1) \isom \PGL_2(\QQ)$ via the isomorphism
  \[
  \epsilon_{A_1}\colon
  M =
  \begin{bmatrix}
    a & b \\
    c & d
  \end{bmatrix} \mapsto T\epsilon_{A_0}(M)T^{-1}
  = \frac{1}{ad - bc}\begin{bmatrix}
    a^2 & a b & \frac{1}{2} b^2 \\
    2 a c & a d + b c & b d \\
    2 c^2 & 2 c d & d^2.
  \end{bmatrix}
  \]
  Under this map, the subgroup corresponding to $\SO(\ZZ, A_1)$ is $G = \PGGamma^0(2) \sqcup \tau\PGGamma^0(2)$, where
  \[
  \PGGamma^0(2) = \left\{\begin{bmatrix}
    a & b \\
    c & d
  \end{bmatrix} \in \PGL_2(\ZZ) : b \equiv 0 \pmod 2\right\}
  \]
  is a congruence subgroup, and
  \[
  \tau = \begin{bmatrix}
    & 2 \\
    1 &
  \end{bmatrix}.
  \]
\end{lem}
\begin{proof}
  The first statement follows easily from considering the action of an element of $\SO(\QQ, A_1)$ on the locus of isotropic points for $A_1$, a conic in $\PP^2$ that is rationally isomorphic to $\PP^1$. Note that $\epsilon_{A_1}$ is compatible with the map $\epsilon_{A_0}$ found earlier:
  \begin{equation}
    \epsilon_{A_1} = T \epsilon_{A_0} T^{-1}.
  \end{equation}
  As for the second statement, if
  \[
  M = \begin{bmatrix}
    a & b \\
    c & d
  \end{bmatrix} \in \PGL_2(\QQ)
  \]
  is given such that $\epsilon_{A_1}(M)$ is integral, we first left-multiply by $\tau^{-1}$ if needed to make $v_2(\det M)$ even, and then scale $M$ so that $a,b,c,d$ are coprime integers. If a prime $p$ were to divide $ad-bc$, it must divide each of $a$,$b$,$c$,$d$ by the integrality of $\epsilon_{A_1}(M)$, which is a contradiction. Thus $ad-bc=\pm1$, and the upper-right entry $b^2/2(ad-bc)$ shows that $b$ is even, so either $M$ or $\tau^{-1} M$ lies in $\PGGamma^0(2)$, as claimed.
\end{proof}

We have now mapped each $\SL_3$-orbit of pairs of integral symmetric matrices of type $A_1$ to a $G$-orbit of binary quartic $48441$-forms; indeed, it is not hard to show that $G$ is in fact the subgroup of $\PGL_2(\QQ)$ that preserves the lattice of forms of this shape, and we have, by an argument similar to Lemma \ref{lem:bq_aut},
\[
h^{\text{Type }A_1}\big(C,(2)\big) = \sum_{\substack{\text{quartic rings } Q \\ \text{with resolvent } C \\ \text{of $A_1$ type,} \\ \text{up to $C$-isom}}} \frac{1}{\Size{\Aut_C Q}}
= \sum_{\substack{48441\text{-forms } f(x,y) \\ \text{with resolvent $g^{[4]}$,} \\ \text{up to $G$}}} \frac{1}{\Size{\Stab_G f}}.
\]
The subgroup $\PGGamma^0(2)$ has index $2$ in $G$, so its weighted class number is twice the weighted class number for $G$, as used in the statement of the theorem.

We now turn to type $I$. Following the same method, we can write
\[
h^{\text{Type }I}\big(C,(2)\big) = \sum_{\substack{\text{quartic rings } Q \\ \text{with resolvent } C \\ \text{of $I$ type,} \\ \text{up to $C$-isom}}} \frac{1}{\Size{\Aut_C Q}}
= \sum_{\substack{\text{$(I,B)$ with resolvent $g$} \\ \text{up to $\SO(\ZZ,I)$}}} \frac{1}{\Size{\Stab_{\SO(\ZZ,I)} B}}.
\]
At this point we make two striking observations:
\begin{itemize}
  \item The resolvent condition $\det(xI - yB) = g$ is equivalent to $B$ having \emph{characteristic polynomial} $g(x,1)$, so we have connected counting quartic rings to another classical problem, namely counting symmetric matrices of given characteristic polynomial;
  \item Since $\SO(\ZZ,I)$ is a finite group, isomorphic to $\cS_4$ (in its representation as the group of rotations of a cube), there is no need to count \emph{orbits} of symmetric matrices; the matrices themselves will be finite in number.
\end{itemize}
Thus
\[
h^{I\text{-type}}\big(C,(2)\big)
=\frac{1}{24}\Size{\left\{B\in\Mat_{3\times3}(\ZZ):
B^{\top}=B,\ \charpoly(B)=g(x,1)\right\}}.
\]

We must now take the prime at $\infty$ into account. Over $\RR$ there are only two nondegenerate cubic algebras, $\RR \cross \CC$ and $\RR \cross \RR \cross \RR$. If $C_g \tensor \RR \isom \RR \cross \CC$ (that is, $g$ has only one real root), then there is only one quartic algebra with resolvent $C_g$ up to $C_g$-isomorphism, so it does not make sense to impose local conditions at the infinite place. If, on the other hand, $C_g \tensor \RR \isom \RR \cross \RR \cross \RR$, then the three factors of $\RR$ are non-interchangeable, being labeled by the three real roots of $g$, and there are four non-$C_g$-isomorphic quartic algebras with resolvent $C_g$ (one isomorphic to $\RR \cross \RR \cross \RR \cross \RR$ and three to $\CC \cross \CC$), parametrized by the four Kummer elements $\delta \in C_g^{N=1}/(C_g^{N=1})^2.$ The following is not hard to verify:
\begin{lem}
  Let $g$ be a monic binary cubic form over $\RR$ whose dehomogenization has three real roots $\xi_1 < \xi_2 < \xi_3$. Identify the corresponding $\RR$-algebra $C_g$ with $\RR \cross \RR \cross \RR$ with the coordinates ordered so that $\xi \mapsto (\xi_1,\xi_2,\xi_3)$. Let $L/\RR$ be a quartic algebra with resolvent $C_g$. The corresponding pair $(A,B)$ of real symmetric matrices is related to the sign of the corresponding Kummer element $\delta \in C_g^{N=1}/(C_g^{N=1})^2$ in the following way:
  \begin{itemize}
    \item If $\sgn \delta = (+,+,+)$, then $(A,B)$ is of type $A_1$ and yields an indefinite binary quartic form with four real roots.
    \item If $\sgn \delta = (+,-,-)$, then $(A,B)$ is of type $A_1$ and yields a positive definite binary quartic form.
    \item If $\sgn \delta = (-,-,+)$, then $(A,B)$ is of type $A_1$ and yields a negative definite binary quartic form.
    \item If $\sgn \delta = (-,+,-)$, then $(A,B)$ is of type $I$.
  \end{itemize}
\end{lem}

Now Theorem \ref{thm:O-N_quartic_Z} yields \eqref{eq:BQ_neg_disc} for negative discriminant, \eqref{eq:BQ_indef} and \eqref{eq:BQ_24} for positive discriminant. For the rest, we modify Theorem \ref{thm:O-N_quartic_Z} by changing the conditions imposed at $\infty$. Note that since $H^1(\RR, \cC_2\cross \cC_2) \isom \cC_2 \cross \cC_2$ and the Tate pairing is alternating, any order-$2$ subgroup will be maximal isotropic. In particular, this is true of the subgroup $\{(1,1,1), (1,-1,-1)\}$ giving the forms that are indefinite or positive definite, proving \eqref{eq:BQ_pos}. Similarly for \eqref{eq:BQ_neg}.
\end{proof}

Corollary \ref{cor:BQ} is immediate.

We state this theorem and corollary unconditionally because they apply only to the number field $K = \QQ$. We do not attempt to generalize to other number fields. While the quartic rings of types corresponding to soluble conics ($A_0$ and $A_1$ in our notation) continue to correspond to binary quartic forms, the number of insoluble types grows with the degree of $K$, and they are not expected to have such a nice interpretation as occurred for type $I$.

\subsection{Conductor rings and the proof of Theorem* \texorpdfstring{\ref{thmstar:cond_ring}}{1.19}}
To end this section, we deduce Theorem* \ref{thmstar:cond_ring} from the results preceding it. First, note that an abelian extension of local fields may fail to have a conductor ring, even if we impose the appropriate vanishing condition on the base.

\begin{examp}\label{ex:no_cdr_ring}
  Let $p \geq 5$ be a prime. Take $K = \QQ_p$, $R = \QQ_p^4$ and let $\chi : K^\cross \to \FF_p^\cross$ be any multiplicative homomorphism extending the natural projection from $\ZZ_p^\cross$. Define $\psi : R^\cross \to \FF_p^\cross$ by
  \[
  \psi(a;b;c;d) = \chi\left(\frac{a b}{c d}\right).
  \]
  This is the Artin map of a certain $\FF_p^\cross \isom \ZZ / (p-1)\ZZ$-torsor $E/R$. By construction, $\psi$ vanishes on $\QQ_p^\cross$, and the orders
  \begin{align*}
    \OO_1 &= \{(a;b;c;d) \in \ZZ_p^4 : a \equiv c, b \equiv d \pmod p\} \\
    \OO_2 &= \{(a;b;c;d) \in \ZZ_p^4 : a \equiv d, b \equiv c \pmod p\}
  \end{align*}
  are admissible rings for $\psi$. However, no ring strictly containing either $\OO_1$ or $\OO_2$ (of which there are very few) is an admissible ring for $\psi$; in particular, $\OO_1 \union \OO_2$ generates the whole $\OO_R$, which is certainly not an admissible ring for $\psi$. Thus $\psi$ has no conductor ring.
\end{examp}

However, in two special cases the conductor ring not only exists but has a striking characterization: it is the \emph{resolvent ring} of a certain maximal order. These cases are those of the resolvents corresponding to the maps $\cS_3 \to \cS_2$ and $\cS_4 \to \cS_3$. In the former case, since the orders in a quadratic $p$-adic algebra are totally ordered, a conductor ring must exist, and it is a simple exercise to show that it is the following:

\begin{prop}
  Let $L/K$ be a cubic \'etale algebra over a global or local field. Let $T$ be its quadratic resolvent torsor and $E = L T$ its $\cS_3$-closure. Then the conductor ring of $E/T$ is the quadratic resolvent ring $\OO \subseteq T$ of $L$.
\end{prop}

However, the quartic case is deeper, and our proof will make use of the local reflection of Theorem* \ref{thmstar:quartic_local}.
\begin{proof}
  By general formalism, the global case reduces immediately to the local one, so let $L/K$ be a quartic \'etale algebra over a nonarchimedean local field, and let $R$ and $S$ be its cubic and sextic resolvent algebras. To see the vanishing of $\phi_{S/R}$ on $K^\cross$, write $S=R(\sqrt{\delta})$, where $N_{R/K}(\delta)=1$. Then for $a\in K^\cross$,
  \[
    \phi_{S/R}(a)=\left\langle\delta,a\right\rangle
    =\left\langle a,\delta\right\rangle
    =\phi_{R(\sqrt a)/R}(\delta)
    =\phi_{K(\sqrt a)/K}\left(N_{R/K}(\delta)\right)=1.
  \]

  The claim now has two parts:
  \begin{enumerate}[(a)]
    \item\label{it:cyes} $C_0$ is an admissible ring for $S/R$;
    \item\label{it:cno} any admissible ring for $S/R$ is contained in $C_0$.
  \end{enumerate}
  We prove both using Theorem* \ref{thmstar:quartic_local}.

  Let $g(L',C,\ttt)$ denote the number of $\ttt$-traced orders in a quartic algebra $L'$ with reduced resolvent $C$. Theorem* \ref{thmstar:quartic_local}, with $\ttt=(2)$, states that
  \begin{equation}\label{eq:xz}
    \widehat g(L,C,(2))=c\,g(L,C,(1)),
    \qquad c=\Size{\OO_K/2\OO_K}^{2}.
  \end{equation}
  Now $g(L,C,(1))$ is the number of orders in $L$ with resolvent $C$. In particular, it is $1$ if $C=C_0$ and $0$ if $C\not\subseteq C_0$. By Theorem \ref{thm:hcl_quartic_sbi}, $g(L_{\delta'},C,(2))$ counts the fractional ideals $I$ for which $(C,I,\delta')$ is self-balanced. Hence
  \begin{align}
    \widehat g(L,C,(2))
      &=\frac{1}{\size{H^0(K,M_R)}}
        \sum_{\delta'\in H^1(K,M_R)}
        \left\langle\delta,\delta'\right\rangle
        g(L_{\delta'},C,(2)) \nonumber\\
      &=\frac{1}{\size{H^0(K,M_R)}}
        \sum_{\substack{(I,\delta')\\(C,I,\delta')\text{ self-balanced}}}
        \phi_{S/R}(\delta'). \label{eq:xsum}
  \end{align}

  Suppose first that $C_0$ is not admissible. Then there is a $u\in C_0^\cross$ such that $\phi_{S/R}(u)=-1$. Replacing $u$ by $u^3/N_{R/K}(u)$, we may take $N_{R/K}(u)=1$. The involution
  \[
    (I,\delta')\longmapsto(I,u\delta')
  \]
  permutes the self-balanced pairs for $C_0$ and reverses the sign of every summand in \eqref{eq:xsum}. Thus $\widehat g(L,C_0,(2))=0$, contradicting \eqref{eq:xz}, since $g(L,C_0,(1))=1$. This proves \ref{it:cyes}.

  Now suppose that an admissible ring $C\not\subseteq C_0$ exists, and choose $C$ maximal with this property. Divide the self-balanced pairs in \eqref{eq:xsum} into two cases:
  \begin{itemize}
    \item If $I$ is invertible in $C$, write $I=\alpha C$. Then
    \[
      \epsilon=\delta'\alpha^2\in C^\cross,
      \qquad \delta'=\alpha^{-2}\epsilon,
    \]
    so $\phi_{S/R}(\delta')=\phi_{S/R}(\epsilon)=1$. There is at least one such term, namely $(I,\delta')=(C,1)$.
    \item If $I$ is not invertible in $C$, then $C'=\End(I^2)$ strictly contains $C$. Indeed, if $\End(I^2)=C$, Lemma \ref{lem:sqrs_are_inv} would imply that $I^2$, and hence $I$, is invertible. By maximality, $C'$ is not admissible, so, as above, there is a norm-one $u\in(C')^\cross$ with $\phi_{S/R}(u)=-1$. The involution $(I,\delta')\mapsto(I,u\delta')$ pairs the terms with the same $C'$ and reverses their signs.
  \end{itemize}
  Consequently $\widehat g(L,C,(2))>0$, contradicting \eqref{eq:xz}, since $g(L,C,(1))=0$. This proves \ref{it:cno}.
\end{proof}

\section{The Tate pairing}
\label{sec:Tate}
Let $K$ be a nonarchimedean local field. In this section, we seek to understand the local Tate pairing, given by the cup product followed by the invariant map:
\[
  \langle\ ,\ \rangle_T:H^1(K,M)\cross H^1(K,M)
  \longto H^2(K,\mu_2)\xrightarrow{\inv_K}\mu_2.
\]
\begin{prop}\label{prop:Tate_quartic}
  If $M\cong(\ZZ/2\ZZ)^2$ has Galois structure corresponding to a cubic \'etale algebra $R$ over $K$, then the Tate pairing is the restriction of the Hilbert pairing on $R$.
\end{prop}

\begin{proof}
  Let $R_1,\ldots,R_\ell$ be the field factors of $R$; each field factor $R_i$ corresponds to an orbit $G_K\cdot x_i$ for $x_i\in M\bs\{0\}$. Identifying $M$ with $M^\vee$ via the unique alternating pairing $q:M\tensor M\to\mu_2$, each $x_i$ yields a character $\chi_i:M\to\mu_2$.

  Let $\sigma,\tau\in H^1(K,M)$ correspond to Kummer elements $\alpha,\beta\in R^\cross$. The class $\chi_{i*}\Res^K_{R_i}\sigma\in H^1(R_i,\mu_2)$ corresponds to the $R_i$-coordinate $\alpha^{(i)}$, and similarly for $\beta$. We therefore have
  \begin{align*}
    \left\langle\sigma,\tau\right\rangle_{\mathrm{Hilb}}
      &=\prod_{R_i}\left\langle\alpha^{(i)},\beta^{(i)}\right\rangle_{R_i}\\
      &=\prod_{R_i}\inv_{R_i}\Big(
        \chi_{i*}\Res^K_{R_i}\sigma
        \cup \chi_{i*}\Res^K_{R_i}\tau\Big).
  \end{align*}
  Since $\inv_{R_i}=\inv_K\circ\Cor^K_{R_i}$, this becomes
  \begin{align*}
    \left\langle\sigma,\tau\right\rangle_{\mathrm{Hilb}}
      &=\inv_K\sum_{R_i}\Cor^K_{R_i}\Big(
        (\chi_i\tensor\chi_i)_*\Res^K_{R_i}(\sigma\cup\tau)\Big).
  \end{align*}
  We use the following standard transfer identity.

  \begin{lem}\label{lem:cor_res}
    Let $H\subseteq G$ be a subgroup of finite index. Let $X$ and $Y$ be $G$-modules, and let $f:X\to Y$ be $H$-linear. Define the $G$-linear map
    \[
      \widetilde f(x)=\sum_{gH\in G/H}gfg^{-1}(x).
    \]
    Then for $\rho\in H^n(G,X)$,
    \[
      \Cor_H^G\bigl(f_*\Res_H^G\rho\bigr)=\widetilde f_*\rho.
    \]
  \end{lem}
  \begin{proof}
    By dimension shifting it is enough to check $n=0$, where the identity follows directly from the definitions.
  \end{proof}

  Applying the lemma with $G=G_K$, $H=G_{R_i}$, $X=M\tensor M$, $Y=\mu_2$, and $f=\chi_i\tensor\chi_i$, we obtain
  \begin{align*}
    \left\langle\sigma,\tau\right\rangle_{\mathrm{Hilb}}
      &=\inv_K\sum_{R_i}
        \Bigg(\sum_{gG_{R_i}\in G_K/G_{R_i}}
          g(\chi_i\tensor\chi_i)g^{-1}\Bigg)_*(\sigma\cup\tau)\\
      &=\inv_K\Bigg(
        \sum_{\chi\in M^\vee\bs\{0\}}\chi\tensor\chi
        \Bigg)_*(\sigma\cup\tau).
  \end{align*}
  The Tate pairing is $\inv_Kq_*(\sigma\cup\tau)$, and
  \[
    \sum_{\chi\in M^\vee\bs\{0\}}\chi\tensor\chi=q.
  \]
  Indeed, both sides are $\Aut(M)$-invariant, and a direct calculation distinguishes $q$ from $0$.
\end{proof}

\chapter{Tame quartic rings with not totally split resolvent, by multijection}
\label{cha:tame}

There are many cases of local Quartic O-N (Theorem* \ref{thmstar:quartic_local}). We begin with the case that the residue characteristic of $K$ is not $2$ and the resolvent $K$-algebra $R$ is not isomorphic to $K \cross K \cross K$. Here the local reflection theorem says that
\begin{equation}\label{eq:tame_dual}
  g_{(1), C} = \hat g_{(1), C}.
\end{equation}
 Our methods parallel the corresponding case of cubic O-N in \cite[\textsection5.3]{OCubic}: we prove that the functions on each side of \eqref{eq:tame_dual} are equal at $0 \in H^1(K, M_R)$, and then we deduce the full equality by symmetry.

\section{Invertibility of ideals in orders}
We begin with a technical inquiry that has interest in its own right. It is well known that every $\ZZ$-lattice $\aaa$ in a quadratic field is invertible with respect to \emph{some} order, namely its endomorphism ring $\End \aaa$. In a cubic or higher-degree field this is not so, but still, there is some structure. The following two lemmas are exercises in commutative algebra; the second is similar to a theorem in Neukirch \cite[Theorem I.12.4]{Neukirch}.

\begin{lem} \label{lem:order_is_prod}
Let $K$ be a local field, and let $\OO$ be an order in a finite-rank \'etale algebra $L$ over $K$. Then there is a decomposition
\[
  L = L_1 \cross \cdots \cross L_s,
\]
each $L_i$ being the product of some field factors of $L$, with the following properties:
\[
  \OO = \OO_1 \cross \cdots \cross \OO_s
\]
is the product of orders in the $L_i$, and each $\OO_i$ has only a single prime $\pp_i$ above the valuation ideal $\pp$, so that every element of $\OO_i$ not lying in $\pp_i$ is a unit.
\end{lem}
\begin{lem}\label{lem:inv_eq_pri}
  Let $K$ be a local field, and let $\OO$ be an order in a finite-rank \'etale algebra $L$ over $K$. Then a fractional ideal of $\OO$ is invertible if and only if it is principal.
\end{lem}
Now for an apparently new result:
\begin{lem}\label{lem:sqrs_are_inv}
If $\cc \subseteq R$ is a lattice in a cubic algebra over a local field $K$, then $\cc^2$ is invertible in its endomorphism ring $\End(\cc^2)$.
\end{lem}
\begin{proof}
First, $\cc\OO_R$ is an invertible $\OO_R$-ideal, which, since $\OO_R$ is a product of PID's, we can scale to be $\OO_R$.

We first claim that $\cc$ contains a unit, or else has a special form for which $\cc^2=\OO_R$. Let $\pi$ be a uniformizer for $\OO_K$, let $k=\OO_K/\pi\OO_K$, and let $\ba{\OO_R}=\OO_R/\pi\OO_R$. The non-units of $\ba{\OO_R}$ are the union of at most three proper subspaces, coming from the valuation ideals of the field factors. The image $\bar\cc$ cannot lie in any of these subspaces since $\cc\OO_R=\OO_R$. Since a vector space over $k$ cannot be the union of fewer than $\size{k}+1$ proper subspaces, $\cc$ contains a unit unless $k=\FF_2$ and $R=K\cross K\cross K$. In the exceptional case, if $\bar\cc$ contains no unit, then, after permuting the factors,
\[
  \bar\cc=\{(a;b;c)\in\FF_2^3:a+b+c=0\}.
\]
The products of elements of this plane span $\FF_2^3$, so Nakayama's lemma gives $\cc^2=\OO_R$.

We may therefore assume that $\cc$ contains a unit, which we scale to equal $1$. Smith normal form supplies an $\OO_K$-basis $\{1,\alpha,\beta\}$ of $\OO_R$ such that $\{1,\pi^i\alpha,\pi^j\beta\}$ is a basis of $\cc$ for some $i,j\geq0$. If $\alpha\beta=a+b\alpha+c\beta$, replacing $\alpha$ by $\alpha-c$ and $\beta$ by $\beta-b$ does not change $\cc$ and makes $\alpha\beta\in\OO_K$. Then
\[
  \cc^3=\cc^2+\left\langle
    \pi^{3i}\alpha^3,\
    \pi^{2i+j}\alpha^2\beta,\
    \pi^{i+2j}\alpha\beta^2,\
    \pi^{3j}\beta^3
  \right\rangle.
\]
The two mixed terms lie in $\cc\subseteq\cc^2$ because $\alpha\beta=t$. The characteristic polynomials of the integral elements $\alpha$ and $\beta$ show, after multiplying by the displayed powers of $\pi$, that the other two terms lie in $\cc^2$ as well. Hence $\cc^3=\cc^2$.

It follows that $\cc^4=\cc^2$, so $\cc^2$ is an order. Since an element stabilizing $\cc^2$ must itself belong to $\cc^2$, we have $\End(\cc^2)=\cc^2$, and $\cc^2$ is invertible in its endomorphism ring.
\end{proof}
Although Lemma \ref{lem:sqrs_are_inv} is simple to state, we have not found it anywhere in the literature. In general, we suspect the following:
\begin{conj}
  If $\cc$ is a lattice in an algebra $L$ of rank $n$ over a Dedekind domain $\OO_K$, then $\cc^{n-1}$ and all higher powers of $\cc$ are invertible in their common endomorphism ring.
\end{conj}
That the exponent $n-1$ is sharp is seen from the cute example
\[
  \cc = \ZZ_p + (0;1;\ldots;n-1)\ZZ_p + p\ZZ_p^n \subseteq \ZZ_p^n.
\]
The power $\cc^i$ consists of all sequences $(a_0;\ldots;a_{n-1}) \in \ZZ_p^n$ that are congruent modulo $p$ to the values
\[
  (f(0); f(1); \ldots; f(n-1))
\]
of a polynomial $f$ of degree at most $i$ with coefficients in $k$. If $p > n$, then this power stabilizes to the whole of $\ZZ_p^n$ only for $i \geq n - 1$.

\section{Self-duality of the count of quartic orders}
\begin{thm}[\textbf{Local quartic O-N in the tame, not totally split case}] \label{thm:O-N_quartic_local_tame}
Assume $K$ is a local field of residue characteristic not $2$. Let $C \subseteq R$ be an order in a cubic \'etale algebra $R$ that is \emph{not} totally split. Then the assignment $g_C$ to each $L$ of the number of orders $\OO \subseteq L$ with resolvent $C$ is self-dual.
\end{thm}
\begin{proof}
As in \cite[\textsection5.3]{OCubic}, the proof proceeds by reduction to the zero case (i.e.\ that $f(0) = \hat{f}(0)$).

Corresponding to the cubic \'etale algebra $R$, there is a Galois module $M_R$ whose underlying group is $\cC_2 \cross \cC_2$. We will work extensively with $H^1(K, M_R)$, which we abbreviate to $H^1$. Likewise we will abbreviate $H^0(K, M_R)$ to $H^0$.

We have
\[
  \size{H^0} = \frac{\size{R^\cross[2]}}{2} = \begin{cases}
    1 & \text{if $R$ is a field} \\
    2 & \text{if $R \cong K \cross K_2$ for some quadratic field $K_2$} \\
    4 & \text{if $R \cong K \cross K \cross K$.}
  \end{cases}
\]
Moreover, since the unramified cohomology is self-orthogonal, we have
\[
  \size{H^1} = \size{H^0}^2 = \frac{\size{R^\cross[2]}^2}{4}.
\]
Since we are excluding the case $R \isom K \cross K \cross K$, there are just two possibilities:
\begin{itemize}
  \item If $R$ is a field, then $\size{H^1} = 1$, and there is nothing to prove, as any function on $H^1$ is self-dual.
  \item If $R$ is the product of two fields, then $\size{H^1} = 4$. Pick an $\FF_2$-basis $\left\langle\sigma_1,\sigma_2\right\rangle$. The Tate pairing is given by the unique alternating pairing on $H^1$. The space of functions on $H^1$ is four-dimensional, and a basis is
  \[
    \One_{\left\langle\sigma_1\right\rangle}, \One_{\left\langle\sigma_2\right\rangle}, \One_{\left\langle\sigma_1 + \sigma_2\right\rangle}, \One_{\{0\}}.
  \]
  Note that the first three basis elements are self-dual, while the fourth differs from its dual even at $0$. Just as in case (c) of Section 5.3 of \cite{OCubic}, this proves that if $f$ is a function on $H^1$ with $f(0) = \hat f(0)$, then $f = \hat f$.
\end{itemize}
So we have reduced local O-N to the following lemma:
\end{proof}
\begin{lem}\label{lem:tame0}
Assume $K$ is a local field of residue characteristic not $2$. Let $C \subseteq R$ be an order in a cubic \'etale algebra. Then the assignment $g_C$ to each $L$ of the number of orders $\OO \subseteq L$ with resolvent $C$ satisfies self-duality at $0$:
\[
  \hat g_C(0) = g_C(0).
\]
\end{lem}
\begin{proof}
As in the cubic case, the proof is by explicit multijection.

On the one hand,
\[
  \size{H^0} \cdot \hat g_C(0) = \sum_\sigma g_C(\sigma)
\]
counts all quartic orders with cubic resolvent $C$, and using Theorem \ref{thm:hcl_quartic_sbi}, these can be parametrized by self-balanced ideals $(C, \cc, \delta)$, where $\delta$ ranges over a set of representatives for $R^{N=1}/\big(R^{N=1}\big)^2$. On the other hand, $g_C(0)$ is the number of orders with resolvent $C$ in $K \cross R$. Write such an order as $\OO = \OO_K + 0 \cross \aaa$, where $\aaa \subseteq R$ is a lattice. The conditions that $\OO$ be a ring, and in particular that $\aaa \subseteq \OO_R$, are subsumed (by Theorem \ref{thm:hcl_quartic}) by the resolvent conditions, namely that
\begin{enumerate}[(i)]
  \item $N_C(\aaa) = 1$,
  \item $\Phi(0; \alpha) \in C + K$ for all $\alpha \in \aaa$. If $\alpha$ is a non-zero-divisor (a sufficient case), then $\Phi(0; \alpha)$ can be written explicitly as $N(\alpha)/\alpha$.
\end{enumerate}
The multijection is as follows. First, $\cc$ may not be invertible. Let $C_1 = \End \cc^2$ and $\cc_1 = \cc C_1$, an invertible and thus a principal $C_1$-ideal. Let
\begin{equation} \label{eq:multij_quartic}
  \aaa_1 = \frac{[C_1 : C]}{\delta [\cc_1 : \cc]} \left( \cc_1^2\right) ^{-1}.
\end{equation}
Finally, since $\cc_1$ and $\aaa_1$ are both principal and thus scalar multiples of each other, we can take $\aaa$ to be an ideal that sits inside $\aaa_1$ as $\cc$ sits inside $\cc_1$: that is, if $\cc_1 = \gamma_0 C_1$ and $\aaa_1 = \alpha_0 C_1$, then
\[
  \aaa = \frac{\alpha_0}{\gamma_0} \cc.
\]
Before checking that this $\aaa$ yields a valid ring, we check how many-to-one our multijection is. First note that $\aaa$ determines $C_1 = \End \aaa^2$ and $\aaa_1 = \aaa C_1$, and in particular the index $[\cc_1 : \cc] = [\aaa_1 : \aaa]$. Then, by \eqref{eq:multij_quartic}, the ``shadow'' $\bb = \delta \cc^2 = \delta \cc_1^2$ is determined. Note that $\bb$ is an invertible $C_1$-ideal of norm $[C_1:C]^2/[\aaa_1 : \aaa]^2$, a square. The pairs $(\cc_1, \delta)$ satisfying $\bb = \delta \cc_1^2$, where $\cc_1$ is an invertible $C_1$-ideal and $\delta$ is one of the representatives for $R^{N=1}/\big(R^{N=1}\big)^2$, are found to form $\size{H^0}$ classes under the appropriate equivalence relation, by an argument identical to the cubic case (end of Section 5.3 of \cite{OCubic}). Finally, locating $\cc$ within $\cc_1$ involves the same choice as locating $\aaa$ within $\aaa_1$. So we have a string of many-to-one correspondences
\begin{equation} \label{eq:multij_quartic_dgm}
  (\cc, \delta) \xrightarrow{\text{$n$ to 1}}
  (\cc_1, \delta) \xrightarrow{\text{$\size{H^0}$ to 1}}
  \aaa_1 \xleftarrow{\text{$n$ to 1}}
  \aaa,
\end{equation}
and thus overall there are $\size{H^0}$ times as many $(\cc, \delta)$ as $\aaa$.

It remains to prove that the correspondence \eqref{eq:multij_quartic_dgm} preserves the resolvent and balancing conditions. As for the first condition, regarding the discriminant of the ring, we leave it to the reader to verify that
\begin{equation} \label{eq:xnorm}
  N_C(\aaa) = \frac{1}{N(\delta) N_C(\cc)^2}.
\end{equation}
Now we may assume that both sides of \eqref{eq:xnorm} are $1$. Since $\cc_1$ is principal, we may assume that $\cc_1 = C_1$, adjusting $\delta$ by a square if necessary. Then by the conditions of Theorem \ref{thm:hcl_quartic_sbi}, for a suitable $t\in K^\cross$ we have
\[
  N(\delta)=t^{-2}, \qquad
  N_C(\cc)=t, \qquad
  \frac{[C_1:C]}{[\cc_1:\cc]}=\frac{1}{t};
\]
thus
\[
  \aaa_1 = \frac{1}{t\delta} C_1.
\]

We first prove that if $\cc$ satisfies its resolvent condition, so does $\aaa$. Any $\alpha \in \aaa \subseteq \aaa_1$ has the form $\alpha = \frac{1}{t\delta} \beta$, $\beta \in C_1$, and then
\[
  \alpha' \alpha'' = \frac{1}{t^2\delta' \delta''} \beta' \beta'' = \delta \beta' \beta'';
\]
and we note that if $\beta \in C_1$, then $\beta' \beta'' \in C_1$ as well, by the relation
\begin{equation} \label{eq:resolvent_trick}
  \beta' \beta'' = \underbrace{\beta \beta' + \beta \beta'' + \beta' \beta''}_{{} \in \OO_K}
  - \beta(\underbrace{\beta + \beta' + \beta''}_{{} \in \OO_K} - \beta).
\end{equation}
Thus if $\cc$ satisfies the resolvent condition ($\delta \cc^2 \subseteq C$), then $\aaa$ satisfies the resolvent condition ($\alpha' \alpha'' \in C$ for all $\alpha \in \aaa$). To prove the converse, it suffices to show that
\[
  B = \{\beta' \beta'' : \beta \in \cc\}
\]
spans $C_1$ over $\OO_K$, where $\cc \subseteq C_1$ is a sublattice with $\cc^2 = C_1$.

We first claim that $\cc$ contains a unit. Let $\pi$ be a uniformizer for $\OO_K$ and $k = \OO_K/\pi\OO_K$ the residue field. In the cubic $k$-algebra $\ba{C_1} = C_1/\pi C_1$, the non-units are the union of at most three subspaces. The projection $\bar\cc$ of $\cc$ down to $\ba{C_1}$ cannot lie in any of these subspaces since $\cc^2 = C_1$, so, since $\size{k} \geq 3$, $\bar{\cc}$ must contain a unit, which lifts to a unit in $\cc$. There is no harm in rescaling $\cc$ so that $1 \in \cc$.

Now $1 \in B$. Also, for each $\beta \in \cc$,
\[
  \left\langle B\right\rangle \ni \beta'\beta'' - (\beta' + 1)(\beta'' + 1) + \tr \beta + 1 = \beta,
\]
and thus
\[
  \left\langle B\right\rangle \ni \beta' \beta'' - (\beta \beta' + \beta \beta'' + \beta' \beta'') + \beta(\beta + \beta' + \beta'') = \beta^2.
\]
But since $\ch k \neq 2$, the elements $\beta^2$, for $\beta \in \cc$, generate $\cc^2 = C_1$, completing the proof.
\end{proof}

For the totally split case, the method of proof of Theorem \ref{thm:O-N_quartic_local_tame} fails, because $\dim H^1 = 4$ and $\hat f(0) = f(0)$ is no longer enough to imply $\hat f = f$. But when we count by discriminant instead of resolvent, it can be rescued, due to the following symmetry argument.

\begin{thm}
Fix a cubic algebra $R$ over a local field $K$, and let $g_D : H^1 \to \ZZ$ count the number of orders in a quartic algebra $L \in H^1$ with discriminant $D$, weighted by number of resolvents. That is, $g_D$ is the sum of all the $g_C$'s in Theorem \ref{thm:O-N_quartic_local_tame} over $C \subseteq R$ of discriminant $D$. Then $g_D = \hat g_D$.
\end{thm}
\begin{proof}
All cases are covered by Theorem \ref{thm:O-N_quartic_local_tame} except for the totally split case $R \cong K \cross K \cross K$, where $\dim H^1 = 4$. We can write $H^1 = \left\langle\sigma_1, \tau_1, \sigma_2, \tau_2\right\rangle$, where the $\sigma$'s and $\tau$'s correspond to Kummer elements
  \[
    \sigma_1 : (u,u,1), \quad \tau_1 : (\pi,\pi,1), \quad \sigma_2 : (1,u,u), \quad \tau_2 : (1, \pi, \pi),
  \]
  where $\pi \in \OO_K$ is a uniformizer and $u$ is a non-square unit. The Tate pairing is given, by Proposition \ref{prop:Tate_quartic}, by the $\mu_2$-valued pairing
  \[
    \left\langle\sigma_1,\tau_2\right\rangle = \left\langle\tau_1, \sigma_2\right\rangle = -1, \quad \left\langle\sigma_1,\sigma_2\right\rangle = \left\langle\tau_1,\tau_2\right\rangle = \left\langle\sigma_1,\tau_1\right\rangle = \left\langle\sigma_2,\tau_2\right\rangle = 1.
  \]
  The group $\Aut R \cong \cS_3$ acts on $H^1$, permuting $\sigma_1, \sigma_2, \sigma_1 + \sigma_2$ and $\tau_1, \tau_2, \tau_1 + \tau_2$ in the permutation manner.  There are five orbits, represented by $0$, $\sigma_1$, $\tau_1$, $\sigma_1 + \tau_1$, and $\sigma_1 + \tau_2$. Note that $g_D$ must be constant on each orbit, because its definition is $\cS_3$-invariant. The functions
  \[
    \One_{\left\langle\sigma_1,\sigma_2\right\rangle}, \One_{\left\langle\tau_1,\tau_2\right\rangle}, \One_{\left\langle\sigma_1 + \tau_1, \sigma_2 + \tau_2\right\rangle}, \sum_{\gamma \in \Aut R} \One_{\gamma\left( \left\langle\sigma_2, \sigma_1 + \tau_2\right\rangle\right) }, \One_{\{0\}}
  \]
  form a basis for the $5$-dimensional space of $\cS_3$-invariant functions on $H^1$. The first four are self-dual, while the last differs from its dual even at $0$; so, since $g_D(0) = \hat g_D(0)$ by Lemma \ref{lem:tame0}), the coefficient of the last basis element must be $0$ and $g_D$ is self-dual.
\end{proof}
The weighting by number of resolvents is eliminated by the same ideal M\"obius inversion used in \eqref{eq:ring_from_resolvent_count}.

\chapter{\texorpdfstring{$2 \times n \times n$}{2×n×n} boxes}
\label{cha:2xnxn}

Before proceeding, we adapt the methods of the preceding sections to prove Theorem \ref{thm:2xnxn} on the generalization to symmetric $2\times n\times n$ boxes.

First we must speak a little about rings parametrized by binary $n$-ic forms. Let $\OO_K$ be a Dedekind domain. Let $V \isom \OO_K \oplus \bb$ and $L$ be lattices of ranks $2$ and $1$ respectively, and let $f \in \Sym^n V \tensor L$ be a nonzero tensor, which we can also view as a degree-$n$ map from $V^*$ to $L$. A construction, due to Birch--Merriman \cite{BirchMerriman} and Nakagawa \cite{NakBinForms} over $\ZZ$ and to Wood \cite{WRings} in the general case, associates to $f$ a ring $R_f$ of rank $n$ together with an ideal $I_f$ whose $(n-2)$nd power has the same class as the inverse different of $R_f$. This construction respects base change, and if $S$ is an $\OO_K$-algebra, we let $R_{f,S} = R_f \tensor_{\OO_K} S$ and $I_{f,S} = I_f \tensor_{\OO_K} S$ be the corresponding objects over this base. Symmetric $2 \times n \times n$ boxes have the following useful interpretation:
\begin{prop} Let $\OO_K$ be a Dedekind domain. Let $V \isom \OO_K \oplus \bb$ and $L$ be lattices of ranks $2$ and $1$ respectively, and let $f \in \Sym^n V \tensor L$ be a binary $n$-ic form whose leading coefficient $f_0$ and discriminant are nonzero. Let $U \isom \OO_K^{n}$ be a lattice of rank $n$ with Steinitz class $\Lambda^n U \isom \aaa$ and an orientation $\xi : (\Lambda^n U)^{\tensor 2} \to L$.
\begin{enumerate}[$($a$)$]
  \item $\OO_K$-points of the variety
  \[
  \cV_f(\OO_K) = \left\{(A, B) \in \Sym^2 (U) \tensor V : \xi_* \det(Ax + By) = f \right\}
  \]
  correspond to based self-balanced ideals of $R_f$ with Steinitz class $\aaa$, that is, systems
  \[
    \left(R_f, J = \OO_K\xi_1 \oplus \cdots \oplus \OO_K\xi_{n-1} \oplus \aaa\xi_n, \alpha\right)
  \]
  where $J \subseteq R_{f,K}$ is a fractional ideal and $\alpha \in R_{f,K}$ an element with
  \[
  \alpha J^2 \subseteq I_f^{n-3} \textand N(\alpha) N_{R_f}(J)^2 = f_0^{-n+3},
  \]
  up to the symmetry
  \[
  J \mapsto \lambda J, \quad \xi_i \mapsto \lambda \xi_i, \quad \alpha \mapsto \lambda^{-2} \alpha
  \]
  for $\lambda \in R_{f,K}^\cross$.
  \item The bijection is given by sending $(R_f, J, \alpha)$ to the pair of matrices describing the composition
  \begin{align*}
    J \cross J &\to I_f^{n-3}\xrightarrow{(\pi_{n-2},\pi_{n-1})} V \\
    (\eta_1, \eta_2) &\mapsto \alpha \eta_1\eta_2,
  \end{align*}
  where the $\pi_i$ are maps attached to $f$.
  \item The group $\SL(U) \subseteq \Aut_{\OO_K} U$ of linear transformations with determinant $1$ acts on based self-balanced ideals of $R_f$ by changing the basis.
\end{enumerate}
\end{prop}
\begin{proof}
This was essentially proved by Wood (\cite{W2xnxn}, Theorems 4.1 and 6.3). See Swaminathan \cite[\textsection 2]{SwamAvg2Tors} for a closely related proposition presented in a notation close to the present memoir.
\end{proof}

\begin{lem}\label{lem:2xnxn_composed}
  Over a field $K$ ($\ch K \neq 2$), the variety $\cV_f$ is a composed variety, with point stabilizer $M$ having underlying group $\cC_2^{n-1}$.
\end{lem}
\begin{proof}
We first consider the situation over $\bar{K}$. There, $J = R_{f,\bar K}$ is forced, and all $\alpha$ are equivalent, so there is only one orbit on $\cV_f(\bar K)$. Its stabilizer in $\SL_n(\bar K)$ consists of the linear transformations that act on $R_{f,\bar K} \isom \bar K^n$ by multiplication by a square root of unity of norm $1$; these form a group $M \isom \cC_2^{n-1}$.

We then consider the situation over $K$. Here $J$ carries no information, but $\alpha \in K^\cross / (K^\cross)^2$ can be any element of square norm. So we have proved the lemma.
\end{proof}

Adapting the methods of Section \ref{sec:Tate}, it is not so hard to get the following results.
\begin{lem}
Let $K$ be a field ($\ch K \neq 2$), let $n$ be a positive odd integer, and let $f$ be a binary $n$-ic form of nonzero discriminant over $\OO_K$.
\begin{enumerate}[$($a$)$]
  \item The module $M$ has a basis
  $\{m_1,\ldots,m_{n-1}\}$ such that the basis elements, along with the sum $m_n = m_1+\cdots+m_{n-1}$, are permuted by the Galois group in the same manner as the $n$ roots of $f$ in $\PP^1$. In particular, $M \isom M^\vee$ is canonically isomorphic to its own Tate dual under the pairing
  \begin{equation}
    \left\langle m_i, m_j\right\rangle = \left.\begin{cases}
      1 & i = j \\
      -1 & i \neq j
    \end{cases}\right\} \in \mu_2.
    \label{eq:pairing_2xnxn}
  \end{equation}
  \item There is an isomorphism
  \[
    H^1(K, M) \isom R_{f, K}^{N=1} / \left(R_{f, K}^{N=1}\right)^2
  \]
  that takes the orbit corresponding to $(R_{f, K}, R_{f, K}, \alpha)$ to $\alpha$.
  \item The Tate pairing
  \begin{align*}
    \cup : H^1(K, M) \cross H^1(K,M) \to \mu_2
  \end{align*}
  is given by restriction of the Hilbert pairing on $R_{f,K}^\cross / (R_{f,K}^\cross)^2$. (If $f$ is reducible over $K$, we define this by multiplying the Hilbert pairings over each of the field factors as in \cite[\textsection2.2]{OCubic}.)
\end{enumerate}

\end{lem}

It now remains to study the local duality. At the infinite place, the definition of ntc has been devised to exploit the orthogonality between $H^1$ and $\{0\}$, so we may restrict our attention to the finite places. Recall that if $\tau$ is a divisor of $2$ in $\OO_K$, we let $V_{\tau}$ denote the lattice of symmetric boxes
\[
\cV_{\tau}(\OO_K) = \left\{[a_{ij}, b_{ij}] : a_{ij}, b_{ij} \in (2^{-1}\tau)^{\One_{i \neq j}}\OO_K \right\}.
\]
Note that the resolvent $f$ of such a box is divisible by $\tau^{n-1}$, since
\[
f(x,y) = \tau^{n-1} \cdot 2^{-1}\tau \det\left(2\tau^{-1} (\cA x - \cB y)\right),
\]
and the argument to the determinant is an $\OO_K$-integral matrix that is skew-symmetric, hence singular (being of odd order), modulo $2\tau^{-1}$.

\begin{conj}[\textbf{Local O-N for $2\times n\times n$ boxes}]
  \label{conj:2xnxn_local}
  If $K$ is a local field, $\ch K \neq 2$, then
  \[
  \cV_{\tau, \tau^{n-1}f} \textand
  \cV_{2\tau^{-1}, (2\tau^{-1})^{n-1}f}
  \]
  are naturally dual with duality constant $q^{(n-1)v_K(\tau)}$. In other words, their local orbit counters satisfy
  \[
  \widehat g_{\tau,\tau^{n-1}f}
  = \Size{\OO_K/\tau\OO_K}^{n-1}
    g_{2\tau^{-1},(2\tau^{-1})^{n-1}f}.
  \]
\end{conj}
\begin{rem}
If $\ch k_K \neq 2$ and $C$ is an $n$-ic $\OO_K$-algebra, we can analogously define $g_{C}$ to count the self-balanced ideals corresponding to each $\alpha \in C_K^{N=1}/(C_K^{N=1})^2$. The conjecture then states that $g_C$ is its own Fourier transform if $C$ comes from a binary $n$-ic form. The hypothesis that $C$ come from a binary $n$-ic form is essential; the conjecture fails for
\[
  C = \{(a;b;c;d;e) \in \ZZ_p^5 : a \equiv b \equiv c \mod p\}.
\]
\end{rem}
\begin{lem}
  Conjecture \ref{conj:2xnxn_local} is true in the following cases:
  \begin{enumerate}[$($a$)$]
    \item $\ch k_K \neq 2$, $v_K(\disc f) \leq 2$;
    \item $\ch k_K = 2$, $e_K = 1$, $v_K(\disc f) = 0$.
  \end{enumerate}
(That is, $\disc f$ is ``odd'' and ``cubefree'' to the extent that these terms are meaningful.)
\end{lem}
\begin{proof}
Throughout this proof, we rely on the construction of $R_f$, $I_f$, and the projectors $\pi_{n-2}, \pi_{n-1}$ given explicitly by Wood \cite[\textsection 2]{W2xnxn}.

We first address the case $\ch k_K \neq 2$. If $R_f$ is maximal, then there is at most one solution $J$: we have $J = I_f^{(n-3)/2}$ works when $\alpha$ is a unit, and otherwise there is no solution. So
\[
  g_{\tau,f} = g_{1,f} = \One_{H^1_\ur(K, M)},
\]
which is its own Fourier transform.

Still assuming $\ch k_K \neq 2$, if $R_f$ is nonmaximal, there is only one possibility: we have that $R_f \subset \OO_L \subset L$ has index $\pi$ in the ring of integers of an unramified extension of $\OO_K$. This implies that there is a decomposition
\[
  L = L_2 \cross L'
\]
where $L_2 = K[\sqrt{D}]$ is of degree $2$ (split or not) and
\[
  R_f = \{\xi \in L_2 : \xi \equiv a \mod \pi
  \text{ for some } a \in \OO_K \} \cross L'.
\]
It is not hard to compute that $I_f = R_f \xi_0$ is principal, so letting $J = \xi_0^{(n-3)/2} J'$, we seek ideals $J'$ such that
\begin{equation}\label{eq:x:aJ'2}
  \alpha J'^2 \subseteq R_f \textand N(\alpha) N_{R_f}(J')^2 = 1.
\end{equation}
There are now two kinds of self-balanced ideals. If $J'$ is invertible, then the norm condition forces equality in \eqref{eq:x:aJ'2}. After scaling so that $J' = R_f$, we must have $\alpha \in R_f^\cross$, which restricts its class to a subgroup $A \subset H^1_{\ur}(K,M)$ of index $2$. However, each such class yields two solutions: $J' = R_f$ and $J' = (\sqrt{D}; 1) R_f$. So we get a contribution
\[
  2 \cdot \One_A
\]
to $g_{1,f}$.

If $J'$ is not invertible but rather has endomorphism ring $\OO_L$, then consideration of the norms involved shows that the inclusion in \eqref{eq:x:aJ'2} is of index $\pi$ and thus $\alpha J'^2 = (\pi;1) \OO_L$. In particular, $\alpha$ lies in a coset $(\pi;1) H^1_{\ur}(K,M)$, and within this coset, $J'$ is unique, yielding a contribution
\[
  \One_{(\pi;1) H^1_{\ur}(K,M)}
\]
to the total
\[
  g_{1,f} = 2 \cdot \One_A + \One_{(\pi;1) H^1_{\ur}(K,M)}.
\]
The subgroup
\[
  B = H^1_{\ur}(K,M) \sqcup (\pi;1)H^1_{\ur}(K,M)
\]
is easily verified to be orthogonal to $A$ by computing the relevant Hilbert pairing; hence
\[
  g_{1,f} = 2 \cdot \One_A - \One_{H^1_{\ur}} + \One_{B}
\]
is its own Fourier transform, as desired.

We now turn to the case that $\ch k_K = 2$ and $e_K=1$. Put
\[
  g_1=g_{1,f} \textand g_2=g_{2,2^{n-1}f}.
\]
Assume first that the leading and trailing coefficients $f_0,f_n$ are units; we will lift this assumption at the end. The count $g_2$ is the ordinary integral-symmetric-matrix count attached to the unramified ring $R_f$, which has power basis $\{1,\theta,\ldots,\theta^{n-1}\}$ and satisfies $I_f=R_f$. For each $\alpha$, we are counting ideals $J$ of $R_f$ such that
\begin{equation}\label{eq:x:aJpow2:t_eq_1}
  \alpha J^2 \subseteq R_f, \quad N(\alpha) N_{R_f}(J)^2 = 1.
\end{equation}
Since $J$ is invertible, we have
\[
  1 = N(\alpha) N_{R_f}(J)^2 = N(\alpha J^2),
\]
so equality must occur in the first relation of \eqref{eq:x:aJpow2:t_eq_1}. Now $\alpha$ must be a unit and $J = R_f$, so the local orbit counter is
\[
  g_2=\One_U, \qquad
  U=H^1(K,M)\mathbin{\cap}
  \frac{\OO_L^{\cross}}{(\OO_L^{\cross})^2}.
\]

Now we examine $g_1$. Doubling a half-integral box gives an integral box with even diagonal entries and determinant $2f$. Thus the ring involved is
\[
  R_{2f} = \left\langle1, 2\theta, 2\theta^2, \ldots, 2\theta^{n-1}\right\rangle,
\]
and we are counting ideals $J$ of $R_{2f}$ such that
\begin{equation}\label{eq:x:aJpow2:t_eq_2}
  \alpha J^2 \subseteq I_{2f}^{n-3}, \quad N(\alpha) N_{R_{2f}}(J)^2 = 1
\end{equation}
and that also satisfy the condition that the diagonal entries of the box are divisible by $2$. Let $\{\hat\theta_0, \ldots, \hat\theta_{n-1}\}$ be the basis of $L$ dual to $\{1, \theta,\ldots, \theta^{n-1}\}$; then the condition is that, for each $\xi \in J$, the coefficients
\[
  \hat\theta_{n-1}(\alpha\xi^2) \textand \hat\theta_{n-2}(\alpha\xi^2),
\]
which are divisible by $2$ by \eqref{eq:x:aJpow2:t_eq_2}, are in fact divisible by $4$.

Assume first that $\alpha$ is a unit; we will eventually prove that it must be. Then $J \subseteq R_f$. Let
\[
  W = J / (J \intsec 2R_f),
\]
a subspace of the $n$-dimensional $k$-vector space $R_{f,k}$. The norm condition of \eqref{eq:x:aJpow2:t_eq_2} ensures that
\begin{equation}\label{eq:x:dimW}
  \dim W \geq \frac{n - 1}{2},
\end{equation}
with equality exactly when $J \supseteq 2R_f$. Since
\[
  I_{2f}^{n-3} = \left\langle1, \theta, \theta^2, \ldots, \theta^{n-3}, 2\theta^{n-2}, 2\theta^{n-1}\right\rangle,
\]
we have the integrality conditions that, for all $v,w \in W$,
\begin{equation}\label{eq:x:avw}
  \hat\theta_{n-1}(\alpha v w) = \hat\theta_{n-2}(\alpha v w) = 0 \quad \text{in } k.
\end{equation}
Let $\lambda(\xi) = \hat\theta_{n-1}(\alpha \xi)$, viewed as a nontrivial linear functional on $R_{f,k}$. It is easy to show that there is a constant $c \in \OO_K$ such that, for all $\xi \in L$,
\[
  \hat\theta_{n-2}(\xi) = \hat\theta_{n-1}(\theta\xi + c\xi).
\]
Thus we can write \eqref{eq:x:avw} as
\begin{equation}\label{eq:x:tavw}
  \hat\theta_{n-1}(\alpha v w) = \hat\theta_{n-1}(\alpha \theta v w) = 0 \quad \text{in } k.
\end{equation}
In particular, $W$ is isotropic for the quadratic form $Q : (v, w) \mapsto \hat\theta_{n-1}(\alpha v w)$, which is nondegenerate since $R_{f,k}$ is separable. In particular, equality must hold in \eqref{eq:x:dimW}. Note that $W + \theta W \subseteq W^\perp$ for $Q$, so $W + \theta W$ is of dimension at most $(n + 1)/2$, and in particular, $W \intsec \theta W$ is of codimension at most $1$ in $W$.
Now we find inductively that
\begin{equation*}
  W_k = \{w \in W : \theta^i w \in W \text{ for } 0 \leq i \leq k\}
\end{equation*}
is of dimension at least $(n - 1)/2 - k$. Let $w_0 \in W_{(n-3)/2}$ be nonzero; we get
\begin{equation} \label{eq:x:Ww}
  W \supseteq \left\langle w_0, \theta w_0, \ldots, \theta^{\frac{n-3}{2}}w_0\right\rangle.
\end{equation}
We claim that the right-hand vectors are all linearly independent. If not, there is a nonzero polynomial $p(\theta)$ of degree at most $(n-3)/2$ that annihilates $w_0$. But $p(\theta)$ has at most $(n-3)/2$ coordinates zero, so $w_0$ has at least $(n+3)/2$ coordinates $0$ and
\[
  \left\langle w_0, \theta w_0, \ldots, \theta^{\frac{n-3}{2}}w_0\right\rangle = R_{f,k} w_0.
\]
This is a product of field factors of $R_{f,k}$, which cannot be isotropic for $Q$. So equality holds in \eqref{eq:x:Ww}, and we have
\begin{align*}
  J &= \left\langle w_0, \theta w_0, \ldots, \theta^{\frac{n-3}{2}}w_0\right\rangle + 2\OO_L \\
  &= w_0 I_{2f}^{(n-3)/2} + 2\OO_L \\
  \alpha J^2 &= \alpha w_0^2 I_{2f}^{n-3} + 2\OO_L.
\end{align*}
Since the endomorphism ring of $I_{2f}^{n-3}$ is $R_{2f}$, we must have $\alpha w_0^2 \equiv c \mod 2\OO_L$ for some $c \in \OO_K$. Upon rescaling, we have $\alpha \equiv w_0 \equiv 1 \mod 2\OO_L$. Thus we get a unique solution up to scaling
\[
  J = I_{2f}^{(n-3)/2}
\]
to \eqref{eq:x:aJpow2:t_eq_2}. As to the condition on the diagonal entries, it is nontrivial only for the generators $1, \theta, \ldots, \theta^{(n-3)/2}$ of $J$ that are not multiples of $2$. The conditions
\[
  \hat\theta_{n-1}(\alpha\theta^{2i}) \equiv \hat\theta_{n-2}(\alpha\theta^{2i}) \equiv 0 \mod 4, \quad 0 \leq i \leq \frac{n - 3}{2}
\]
can be consolidated into
\begin{equation}\label{eq:x:a4}
  \hat\theta_{n-1}(\alpha \theta^i) \equiv 0 \mod 4, \quad 0 \leq i \leq n - 2.
\end{equation}
As we have scaled $\alpha \equiv 1 \mod 2$, let $\alpha = 1 + 2\beta$ and observe that \eqref{eq:x:a4} is a family of $n - 1$ linearly independent linear constraints on $(\beta \bmod 2) \in R_{f,k}$. Their common solution space is the one-dimensional space of constants, so $\beta \equiv c \mod 2$ for some $c \in \OO_K$. But $N(\alpha) \equiv 1 + 2c \mod 4$ is a square, forcing $c \equiv 0 \mod 2$. Consequently we have found that
\[
  g_1([\alpha]) = \begin{cases}
   1 & \text{if }\alpha \in 1 + 4\OO_L, \text{ that is, } [\alpha] \in H^1_{\ur}(K, M) \\
   0 & \text{if }\alpha \in \OO_L^\cross \text{ is not a square modulo } 4.
  \end{cases}
\]
We must dispose of the case where $\alpha$ is not a unit. Scale $\alpha$ by squares so that all its valuations are $0$ or $1$; this describes a splitting $f = f_0 f_1, R_f = R_0 \cross R_1$ such that
\begin{equation*}
  \alpha \in R_0^\cross \cross 2 R_1^\cross.
\end{equation*}
Since $N(\alpha)$ is a square, $R_1$ has even dimension $2a$. We still have $J \subseteq R_f$.
Now the quadratic form $Q$ on $R_f/2R_f$ defined by $Q : (v, w) \mapsto \hat\theta_{n-1}(\alpha v w)$ is no longer nondegenerate, but has kernel $R_1/2R_1$. In particular, the space $W = J / (J \intsec 2R_f)$, being isotropic for $Q$, can have at most dimension $(n-1)/2 + a$. But the norm condition in \eqref{eq:x:aJpow2:t_eq_2} implies that $W$ has at least this dimension. So equality holds, and $J$ decomposes as
\[
  J = J_0 \cross R_1
\]
where $(R_0, J_0, \alpha|_{R_0})$ is one of the self-balanced ideals found in the case of $\alpha$ a unit. But now as $\xi \in R_1$ varies, $\alpha \xi^2$ traverses all classes in $2 R_1 / 4 R_1$, whose images under $\hat\theta_{n-1}$ cannot all vanish modulo $4$, so we do not get a valid $J$. Consequently
\[
  g_1 = \One_{H^1_{\ur}(K,M)}.
\]
The subgroup $U$ is the orthogonal complement of $H^1_{\ur}(K,M)$ under the Hilbert pairing, and $\size{H^1_{\ur}(K,M)}=\size{H^0(K,M)}$. Therefore the normalized Fourier transform gives
\[
  \widehat g_1=g_2.
\]
Applying the Fourier transform once more and using the local Euler-characteristic formula gives
\[
  \widehat g_2=q^{n-1}g_1,
\]
which proves the claimed constants.

We assumed here that $f_0$ and $f_n$ are units. If $\size{k_K}>n$, this can be ensured by an appropriate $\SL_2(\OO_K)$-transformation. In general, pass to an unramified extension $K'/K$ of sufficiently large residue field. The foregoing analysis gives, respectively,
\[
  J\tensor_{\OO_K}\OO_{K'}=I_{f,\OO_{K'}}^{(n-3)/2}
  \quad\text{or}\quad
  J\tensor_{\OO_K}\OO_{K'}=I_{2f,\OO_{K'}}^{(n-3)/2},
\]
with the corresponding unit or square-modulo-$4$ condition on $\alpha$. These conditions descend through the faithfully flat extension $\OO_{K'}/\OO_K$, so the formulas for $g_1$ and $g_2$ remain valid.
\end{proof}

Apply Theorem \ref{thm:main_compose_quartic} to the two composed varieties with fibers
\[
  X_{\tau^{n-1}f}
  \textand
  X_{(2\tau^{-1})^{n-1}f}.
\]
Here
\[
  X_g=\{(A,B)\in\Sym^2(K^n)\tensor K^2:
  2^{n-1}\det(Ax+By)=g\}.
\]
They are full by Lemma \ref{lem:2xnxn_composed}, hence Hasse; moreover, $\SL_n$ has class number one. Lemma \ref{lem:total_cl_num_quartic} converts the normalized total class numbers into the desired weighted orbit counts. The local lemma at the finite places and the archimedean calculation give the identity of Conjecture \ref{conj:2xnxn}, proving Theorem \ref{thm:2xnxn}.

\chapter{Introduction to counting quartic rings}
\label{cha:intro_count}

Here end the cases in which a conceptual, bijective argument has been found to suffice for proving local reflection for quartic rings. To approach the remaining cases, we attack a problem that has interest in its own right: counting orders $\OO$ in a quartic algebra $L$ over a local field $K$ whose cubic resolvent ring $C \subseteq R$ is fixed. This is a finer invariant than in the question addressed by Nakagawa \cite{NakOrders}, namely counting $\OO$ of fixed discriminant (equivalently, index in $\OO_L$).

The index $[\OO_L : \OO]$ is fixed by the condition that $\Theta$ be an isomorphism, so we must analyze the well-definedness of the resolvent map $\Phi : L/K \to R/K$, that is, the condition that
\[
  \Phi(\OO/\OO_K) \subseteq C/\OO_K.
\]
With respect to bases of $\OO$ and $C$, $\Phi$ is given by a pair
\[
  (\cM, \cN) = \left( [M_{ij}], [N_{ij}] \right)
\]
of symmetric $3 \times 3$ matrices, so the resolvent condition can be viewed as the $\OO_K$-integrality of the entries, or half-integrality in the case of the off-diagonal elements. By suitably choosing coordinates, we can ensure that only the integrality of the entries
\[
  M_{11}, \quad N_{11}, \quad M_{12}, \textand M_{22}
\]
is in doubt. Of these, the condition on $M_{11}$ is the most challenging. It amounts to a quadratic condition on the first basis vector $\xi_1$ of $C$, that is, a conic on some \emph{pixel} (clopen neighborhood of a point) in $\PP^2(\OO_K)$. The solubility of this conic over $K$ is governed by an appropriate Hilbert symbol, which we analyze. It is very hard in general to tell if any $K$-points of the conic lie in the requisite pixel, but if there is even one point in the pixel, then, using the rational parametrization of a conic with a basepoint, the volume of points in the pixel can be determined explicitly. Accordingly, our approach to solving the $M_{11}$-condition is a three-step one:
\begin{itemize}
  \item Determine the sum of the solution volumes for $\xi_1$ over all quartic algebras $L$.
  \item Find restrictions on what $L$ can yield a nonzero volume and what that volume can be, providing an upper bound (the \emph{bounding step}).
  \item If the sums of these upper bounds agree, deducing that the bound is attained everywhere, i.e.~all conics not already eliminated have a $K$-point in the requisite pixel (the \emph{summing step}).
\end{itemize}
The $M_{12}$ and $M_{22}$ conditions are essentially linear. We use the computer programs Sage and LattE to sum the ring totals over all possible values of the discrete data and verify the local reflection theorem.

Omitted are the adaptations to be made when $\ch k_K > 2$, where, in view of Theorem \ref{thm:O-N_quartic_local_tame}, only splitting type $(111)$ need be considered. It involves only the black, brown, beige, and white zones; the conics are all very easy to solve and yield the same answers as the wild case upon substituting $e = 0$.

\section{The group \texorpdfstring{$H^1$}{H1} of quartic algebras with given resolvent} \label{sec:H1}
As in the proof of Theorem \ref{thm:O-N_quartic_local_tame}, corresponding to the cubic \'etale algebra $R$, there is a Galois module $M_R$ whose underlying group is $\cC_2 \cross \cC_2$. We will work extensively with $H^1(K, M_R)$, which we abbreviate to $H^1$. Likewise we will abbreviate $H^0(K, M_R)$ to $H^0$.

By Proposition \ref{prop:Kummer_quartic}, we can identify $H^1$ naturally with
\[
  R^{N=1}/\left(R^{N=1}\right)^2.
\]
Now there is a natural isomorphism
\begin{equation}\label{eq:H1_is_summand}
  R^\cross/(R^\cross)^2 \isom K^\cross/(K^\cross)^2 \cross  R^{N=1}/\left(R^{N=1}\right)^2.
\end{equation}
Thus for any $\alpha \in R^\cross$, we can talk about the class $[\alpha]$ of $\alpha$ in $H^1$, that is, the class of $\alpha^3/N(\alpha)$.

The structure of $H^1$ is given by a suitable extension of the Shafarevich basis theorem:
\begin{prop}\
\label{prop:Sh_basis_quartic}
\begin{enumerate}[$($a$)$]
  \item\label{Sh:K} If $K$ is a $2$-adic local field, then $K^\cross / (K^\cross)^2$ is an $\FF_2$-vector space of dimension $ef + 2$. It has a basis consisting of
  \begin{itemize}
    \item one nonunit, the uniformizer $\pi$;
    \item one \emph{intimate unit} $1 + 4u$, the discriminant of the unramified quadratic extension, where $u \in \OO_K$ is any element with unit trace down to $\ZZ_2$;
    \item $ef$-many \emph{generic units} $1 + x\pi^{2i + 1}$, where $0 \leq i < e$ and $x$ ranges over $f$ elements whose reductions mod $\pi$ form an $\FF_2$-basis of $k_K$.
  \end{itemize}
  \item\label{Sh:H1}
  If $R/K$ is a cubic \'etale extension, then $H^1 \cong R^{N=1}/\left(R^{N=1}\right)^2$ is an $\FF_2$-vector space of dimension $2 \dim_{\FF_2} H^0 + 2ef$. It has a basis of $\dim_{\FF_2} H^0$ nonunits, $\dim_{\FF_2} H^0$ intimate units,
  and $2ef$ generic units; the generic units can be chosen as follows:
  \begin{enumerate}[$($i$)$]
    \item If $R$ is unramified, we take $2ef$ units of the form $1 + x\pi^{2i + 1}$, where $0 \leq i < e$ and $x$ ranges over $f$ elements whose reductions mod $\pi$ form an $\FF_2$-basis of $k_R^{\tr = 0}$ for each $i$.
    \item If $R$ is totally tamely ramified, we take $2ef$ units of the form $1 + x\pi^{2i+1}$, where $0 \leq i < 3e$ and $i \not\equiv 1 \pmod 3$, and $x$ ranges over an $\FF_2$-basis of $k_K$ for each $i$.
    \item If $R \isom K \cross Q$ is partially wildly ramified, we take $2ef$ units of the form $(1 + y\pi^{i+1}; 1 + x\pi_Q^{2i+1})$ where $0 \leq i < 2e$, where $x$ ranges over an $\FF_2$-basis of $k_K$ for each $i$, and where $y \in \OO_K$ is adjusted for each $(i,x)$.
  \end{enumerate}
\end{enumerate}
\end{prop}
\begin{proof}
Part \ref{Sh:K} can be proved by elementary techniques; see Del Corso and Dvornicich \cite[Corollary 1]{DCD}. Using this, we construct a basis for $R^\cross / (R^\cross)^2$ and can arrange so that each basis element either lies in $K^\cross$ or has norm $1$ down to $K$; we remove all elements of the former type.
\end{proof}

Define a \emph{level space} in $K^\cross /\left(K^\cross\right)^2$ to be a neighborhood of $1$ in the metric inherited from $K$. The distinct level spaces are
\[
\cF_{i} = \begin{cases}
  K^\cross /\left(K^\cross\right)^2 & i = -1 \\
  \{[\alpha] : \alpha \equiv 1 \mod \pi^{2i}\} & 0 \leq i \leq e \\
  \{[1]\} & i = e+1.
\end{cases}
\]
We have $\size{\cF_{i}} = 2q^{e - i}$ for $0 \leq i \leq e$, and a classical result on Hilbert symbols (see Nguyen-Quang-Do \cite[Proposition 3.4.1.1]{Nguyen}) implies that
\begin{equation} \label{eq:Hilb_prod_size}
  \cF_i^\perp = \cF_{e-i}
\end{equation}
for all $i$. We define the \emph{level} $\ell(\alpha)$ of an element $[\alpha] \in K^\cross/\left(K^\cross\right)^2$ as the largest $i \geq 0$ for which $[\alpha] \in \cF_i$. We have $\ell(1) = e+1$. By convention, if $[\alpha] \notin \cF_0$, we set
\[
\ell(\alpha) = -1/2
\]
to shorten some statements.

Similarly, we define levels of elements in $H^1$. If $R$ is unramified over $K$, we define the \emph{level space}
\[
  \cF_{i} = \begin{cases}
    H^1 & i = -1 \\
    \{[\alpha] \in H^1 : \alpha \equiv 1 \mod \pi^{2i}\} & 0 \leq i \leq e \\
    \{[1]\} & i = e+1,
  \end{cases}
\]
noting that $\size{\cF_i} = \size{H^0}q^{2(e-i)}$ for $0 \leq i \leq e$ and that $\cF_i^\perp = \cF_{e-i}$ for all $i$. If $R$ is ramified over $K$, we define the \emph{level space}
\[
  \cF_{i} = \begin{cases}
  H^1 & i = -1 \\
  \{[\alpha] \in H^1 : \alpha \equiv 1 \mod \pi^{i}\} & 0 \leq i \leq 2e \\
  \{[1]\} & i = e'+1,
  \end{cases}
\]
noting that $\size{\cF_i} = \size{H^0}q^{2e-i}$ for $0 \leq i \leq 2e$ and that $\cF_i^\perp = \cF_{2e-i}$ for all $i$ by \eqref{eq:Hilb_prod_size}.

We will occasionally let
\[
  e' = \begin{cases}
    e, & R \text{ unramified} \\
    2e, & R \text{ ramified}
  \end{cases}
\]
to shorten some statements.

We define the \emph{level} $\ell(\alpha)$ of an element $[\alpha] \in H^1$ as the largest $i \geq 0$ for which $[\alpha] \in \cF_i$. We have $\ell(1) = e'+1$. By convention, if $[\alpha] \notin \cF_0$, we set
\[
  \ell(\alpha) = -1/2
\]
to shorten some future statements.

Let $F_i = \One_{\cF_i}$. By the foregoing, we have:
\begin{prop}\label{prop:levperp}
For $-1 \leq i \leq e'+1$, we have $\hat{F}_i = c_i F_{e' - i}$, where
\begin{equation}
  c_i = \frac{\size{\cF_i}}{\size{H^0}} = \begin{cases}
    \size{H^0} q^{2e}, & i = -1 \\
    q^{2e - 2i}, & 0 \leq i \leq e,\quad R \text{ unramified} \\
    q^{2e - i}, & 0 \leq i \leq 2e,\quad R \text{ ramified} \\
    \dfrac{1}{\size{H^0}}, & i = e'+1.
  \end{cases}
\end{equation}
\end{prop}

\section{Reduced bases}\label{sec:red_basis}
Let $R$ be a rank-$n$ \'etale algebra over a local field $K$. We can embed $R$ into $\bar{K}^n$ by the product of its $n$ embeddings into $\bar{K}$ (such a construction is often called the \emph{Minkowski embedding}).
\begin{defn}\label{defn:red_basis}
Let $I$ be an $\OO_K$-lattice in $R$, and let $\omega \in \bar{K}^n$ be a multiplier with the following property:
\begin{enumerate}[$(*)$]
  \item\label{iota:star} If $\iota$, $\iota'$ are two coordinates of the same field factor of $R$, then $\omega^{(\iota)}$ and $\omega^{(\iota')}$ have the same valuation.
\end{enumerate}
A basis $(\rho_1,\ldots,\rho_n)$ for $\omega I$ is called \emph{reduced} if
\begin{enumerate}[$($a$)$]
  \item $v(\rho_1) \leq \cdots \leq v(\rho_n)$;
  \item If $\rho \in \omega R$ is decomposed as
  \[
  \rho = \sum_i c_i \rho_i, \quad c_i \in K,
  \]
  then for each $i$,
  \[
  v(c_i \rho_i) \geq v(\rho).
  \]
\end{enumerate}
\end{defn}

This notion has the following properties:
\begin{prop}\label{prop:rb} Let $\omega I$ be as above.
\begin{enumerate}[$($a$)$]
\item\label{rb:exists} There exists a reduced basis $(\rho_1, \ldots, \rho_n)$ for $\omega I$.
\item\label{rb:max} If $(\rho_1', \ldots, \rho_n')$ is any other basis for $\omega I$, sorted so that $v(\rho_1') \leq \cdots \leq v(\rho_n')$, then for each $k$,
\[
  v(\rho_k') \leq v(\rho_k).
\]
In particular, if both bases are reduced, equality holds.
\item\label{rb:unique} If $(\rho_1', \ldots, \rho_n')$ is another reduced basis for $\omega I$, then
\[
  \rho_i' = \sum_j c_{ij} \rho_j
\]
for some change-of-basis matrix
\begin{equation} \label{eq:change_basis}
  \left[c_{ij}\right] \in \GL_n(\OO_K), \quad v\left(c_{ij}\right) \geq v(\rho_i) - v(\rho_j).
\end{equation}
Conversely, any matrix $[c_{ij}]$ satisfying \eqref{eq:change_basis} yields a new reduced basis $(\rho_i')_i$.
\item\label{rb:span} As $\OO_{\bar K}$-modules,
\[
  \left\langle\pi^{-v(\rho_i)} \rho_i : 1 \leq i \leq n\right\rangle
  = \left\langle\pi^{-v(\rho)} \rho : \rho \in \omega R^\cross\right\rangle.
\]
\end{enumerate}
\end{prop}

\begin{proof}
\begin{enumerate}[$($a$)$]
  \item Choose a basis $(\rho_1, \ldots, \rho_n)$ such that the sum of the valuations $v(\rho_1) + \cdots + v(\rho_n)$ is maximal. This can be done because there are only finitely many possible valuations of primitive vectors in $\omega I$. Sort the $\rho_i$ in increasing order of valuation. We claim $(\rho_i)_i$ is reduced.

  Let $\rho = \sum_i c_i \rho_i$ be given. Let $a$ be the minimal valuation $v(c_i \rho_i)$ of a term, and suppose that $a < v(\rho)$. Then we have a linear dependency
  \[
    \sum_{v(c_i\rho_i) = a} c_i \rho_i \equiv 0 \mod \pi^a \mm_{\bar K}.
  \]
  Choose $j$ such that $v(c_j\rho_j) = a$ and $v(c_j)$ is minimal. Then
  \[
    \rho_j' = \sum_{v(c_i\rho_i) = a} \frac{c_i}{c_j} \rho_i
  \]
  is an element of $\omega I$ whose valuation exceeds $v(\rho_j)$. Since the coefficient of $\rho_j$ in $\rho_j'$ is $1$, replacing $\rho_j$ by $\rho_j'$ does not change the span $\omega I$ but increases the valuation sum $\sum_i v(\rho_i)$, contradicting the choice of basis $(\rho_i)_i$.
  \item Since $(\rho_i)_i$ and $(\rho_i')_i$ are bases for the same module $\omega I$, such a $[c_{ij}] \in \GL_n(\OO_K)$ certainly exists. Applying the reducedness property to each decomposition
  \[
    \rho_i' = \sum_{j} c_{ij} \rho_j
  \]
  yields a bound
  \[
    v(c_{ij}) \geq v(\rho_i') - v(\rho_j).
  \]
  Suppose that $v(\rho_k) < v(\rho_k')$ for some $k$. Then for $j \leq k \leq i$,
  \[
  v(\rho_j) \leq v(\rho_k) < v(\rho_k') \leq v(\rho_i'),
  \]
  so $c_{ij}$ has positive valuation. Thus, when the matrix $[c_{ij}]$ is reduced modulo $\pi$, it has an $(n - k + 1) \cross k$ block of $0$'s, large enough to make the determinant vanish, which is a contradiction.

  \item By the preceding part, $v(\rho'_k) = v(\rho_k)$. So the associated matrix $[c_{ij}]$ must satisfy
  \begin{equation}
  v(c_{ij}) \geq v(\rho_i') - v(\rho_j) = v(\rho_i) - v(\rho_j).
  \end{equation}
  Conversely, if $[c_{ij}]$ is an invertible matrix satisfying this inequality, we get a new basis $\rho_i'$ with $v(\rho_i') \geq v(\rho_i)$. Equality must hold, and now since $\sum_k v(\rho_k')$ achieves the maximal value, $(\rho_k')_k$ is reduced by the proof of part \ref{rb:exists}.
  \item The $\subseteq$ direction is obvious. For the $\supseteq$ direction, let $\rho \in \omega R^\cross$ be given. Since $(\rho_i)_i$ is reduced,
  \[
    \rho = \sum_i c_i \rho_i, \quad v(c_i \rho_i) \geq v(\rho),
  \]
  so
  \[
    \pi^{-v(\rho)} \rho = \sum_i \left(\pi^{v(\rho_i) - v(\rho)}c_i\right) \pi^{-v(\rho_i)} \rho_i,
  \]
  and the parenthesized coefficients belong to $\OO_{\bar K}$, as desired. \qedhere
\end{enumerate}
\end{proof}

We can find reduced bases with added structure.
\begin{defn}
Fix an ordering $R = R^{(1)} \cross \cdots \cross R^{(r)}$ of the field factors of $R$, and choose a uniformizer $\pi_k$ of each $R^{(k)}$. For $\rho \in \omega R$, let $k$ be the minimal index such that $v(\rho^{(k)}) = v(\rho)$. We say that $\rho$ is \emph{$R^{(k)}$-led.} Put
\[
  m_k(\rho)=e_{R^{(k)}/K}\bigl(v(\rho)-v(\omega^{(k)})\bigr) \in \ZZ
\]
and define the \emph{leader} of $\rho$ by
\[
  \ldr(\rho) = \pi_k^{-m_k(\rho)}\left(\frac{\rho}{\omega}\right)^{(k)} \in \OO_{R^{(k)}}.
\]
Thus $\ldr(\rho)$ is primitive in $\OO_{R^{(k)}}$. We say that a reduced basis $(\rho_i)_i$ is \emph{well-led} if the leaders of the $R^{(k)}$-led basis vectors form a reduced basis for $\OO_{R^{(k)}}$ for each $k$.
\end{defn}

\begin{prop} \label{prop:well_led}
Every $\omega I$, as above, admits a well-led basis.
\end{prop}
\begin{proof}
Consider the element
\[
  \eta_\epsilon = (\pi^\epsilon ; \pi^{2 \epsilon} ; \ldots ; \pi^{r \epsilon})
\]
where $\epsilon$ is a positive rational number, small enough that if $a_1 < a_2$ are two valuations of elements in $\omega I$, then $r\epsilon < a_2 - a_1$. By Proposition \ref{prop:rb}\ref{rb:exists}, there is a reduced basis $\rho_1\eta_\epsilon, \ldots, \rho_n\eta_\epsilon$ for $\omega\eta_\epsilon I$. Each basis element $\rho_i \eta_\epsilon$ has some valuation $u + k\epsilon$, $u \in \QQ$, $k \in \{1,\ldots,r\}$, indicating that $\rho_i$ is $R^{(k)}$-led. Since replacing $\epsilon$ by $0$ preserves non-strict inequalities among valuations in $\eta_\epsilon \omega I$, the $\rho_i$ form a reduced basis for $\omega I$, which we claim is well-led.

Given $\alpha \in \OO_{R^{(k)}}$ primitive, decompose $\omega\alpha = \sum_j c_j \rho_j$ as an element of $\omega R$. We have
\[
  v(c_j\rho_j \eta_\epsilon) \geq v(\alpha \eta_\epsilon) = v(\alpha) + k\epsilon,
\]
and for a nonempty subset of $j$, equality must hold and, in particular, $\rho_j$ must be $R^{(k)}$-led. Let $L_k$ be the set of indices $j$ for which $\rho_j$ is $R^{(k)}$-led, and let $B_k = \{\ldr(\rho_j) : j \in L_k\}$. Now we have
\[
  \omega^{(k)}\alpha = \sum_{j \in L_k} c_j \rho_j^{(k)} + \omega^{(k)}\alpha'
\]
with each term of valuation at least $v(\alpha)$, and where $\alpha' \in \OO_{R^{(k)}}$ with $v(\alpha') > v(\alpha)$. Normalizing the $\rho_j^{(k)}$ as in the definition, we can rewrite this as
\[
  \alpha = \sum_{j \in L_k} c_j' \ldr(\rho_j) + \alpha'
\]
with $c_j' \in \OO_K$. We can iteratively decompose $\alpha'$ the same way, and as its valuation goes to infinity, we get a decomposition
\begin{equation} \label{eq:c_j''}
  \alpha = \sum_{j \in L_k} c_j'' \ldr(\rho_j).
\end{equation}
So $B_k$ generates $\OO_{R^{(k)}}$, and in particular, $\size{L_k} \geq [R^{(k)} : K]$. However,
\[
  \sum_{k} [R^{(k)} : K] = n = \sum_k \size{L_k}.
\]
So equality holds and each $B_k$ is a basis for $\OO_{R^{(k)}}$. Since the decomposition \eqref{eq:c_j''} has every term of valuation at least $v(\alpha)$, and $\alpha$ was any primitive vector, $B_k$ is in fact reduced.
\end{proof}
It is evident that reduced indices take a limited number of values modulo $1$. Indeed, we have the following:
\begin{cor}\label{cor:idxs_mod_1}
  If $I \subseteq R$ is a lattice, then the multiset of valuations $\{a_i\}_i = \{v(\rho_i)\}_i$ mod $1$ of reduced basis elements for $\omega I$ depends only on $R$ and $\omega$, not on $I$. It consists of $f_{R_i/K}$ copies of
  \[
    v(\omega^{(R_i)}) + \frac{j}{e_{R_i/K}}, \quad j = 0, 1, \ldots, e_{R_i/K} - 1,
  \]
  where $R_i$ ranges over the field factors of $R = R_1 \cross \cdots \cross R_r$.
\end{cor}
\begin{proof}
Taking a well-led basis and passing to the leaders, we reduce to the case that $R = R_i$ is a field. Then $\omega$ has equal valuations in all coordinates, and we may assume that $\omega = 1$. Let $\rho_i = \pi_R^{e_{R/K}a_i} \xi_i$, where the $\xi_i \in \OO_R^\cross$ differ only by a unit from the normalizations used before. For each congruence class of $a_i$ modulo $1$, note that if the set of corresponding $\xi_i$ is linearly dependent modulo $\pi_R$, then one of the $\rho_i$ could be increased by an $\OO_K$-linear combination of the others to increase its valuation, contradicting the hypothesis that our basis is reduced. So the $\xi_i$ corresponding to each congruence class of $a_i$ modulo $1$ are linearly independent, and in fact must form a basis for $k_R$ over $k_K$ in order for there to be the full number $e_{R/K}f_{R/K}$ of $\xi_i$. This establishes the claimed multiset of valuations.
\end{proof}

A reduced basis for $\omega I$ does not always remain reduced when we extend the ground field $K$. To study this, we make the following definition.
\begin{defn}
An \emph{extender basis} for $\omega I$ is a basis $(\bar\rho_1, \ldots, \bar\rho_n)$ for $\omega I \tensor_{\OO_K} \OO_{\bar{K}}$ such that the vectors
\[
  \bar\xi_i = \pi^{-v(\bar\rho_i)} \bar\rho_i
\]
form an $\OO_{\bar{K}}$-basis for $\OO_{\bar{K}}^n$. The valuations $\bar a_i = v(\bar\rho_i)$ are called the \emph{extender indices} of the basis, and the $\bar\xi_i$ are called the \emph{extender vectors}.
\end{defn}

If it consists of $\bar \rho_i \in \omega I$, an extender basis is easily seen to be reduced. Fortunately, in the cases of tame ramification, this always holds:

\begin{prop}\label{prop:ext_basis_tame}
If $R/K$ is tamely ramified, then any reduced basis of $\omega I$ is an extender basis.
\end{prop}
\begin{proof}
In view of Proposition \ref{prop:rb}\ref{rb:span}, it is enough to show that
\[
  \OO_{\bar K} \left\langle\pi^{-v(\rho)} \rho : \rho \in \omega R^\cross\right\rangle = \OO_{\bar K}^n.
\]
We immediately reduce to the case that $R$ is a field, and then we may assume $\omega = 1$. It suffices to prove that, for some reduced basis $\rho_1, \ldots, \rho_n$ of $R$,
\[
  \det \left[\pi^{-v(\rho_i)}\rho_i^{(j)}\right]_{i,j} \sim 1.
\]
But that follows from the familiar formula for the discriminant of a tamely ramified extension.
\end{proof}

Extending scalars until $R$ is tamely ramified (or even split), we derive that any $\omega I$ admits an extender basis. The extender indices are seen to be unique, and so we make the following definition:
\begin{defn}
  The \emph{discriminant valuation} of $\omega I$ is twice the sum of its extender indices:
  \[
    v_K \disc(\omega I) = 2(\bar a_1 + \cdots + \bar a_n).
  \]
\end{defn}
Observe the following facts:
\begin{lem}\
\label{lem:extender_disc}
  \begin{enumerate}[$($a$)$]
    \item If $\omega I$ has reduced basis
    \[
      \left(\pi^{a_1} \xi_1, \ldots, \pi^{a_n}\xi_n\right)
    \]
    and an extender basis
    \[
      \left(\pi^{\bar a_1} \bar\xi_1, \ldots, \pi^{\bar a_n}\bar\xi_n\right),
    \]
    then
    \[
    \sum_i \bar a_i - \sum_i a_i = d_{\text{wild}}/2,
    \]
  where
  \[
    d_{\text{wild}} = v_K(\disc_K R) - \sum_{R_i\text{ field factor}}f_{R_i/K}\left(e_{R_i/K} - 1\right)
  \]
  is the wild part of the discriminant.
  \item
  If $\left(\pi^{a_1} \xi_1, \ldots, \pi^{a_n}\xi_n\right)$ is a reduced basis of $\omega I$, then
  \[
  \OO_{\bar K} \cdot \xi_1 \wedge \cdots \wedge \xi_n = \pi^{d_{\text{wild}}/2} \Lambda^{n} \OO_{\bar K}^n.
  \]
  \item If $\omega \in R$ and $\omega I = \OO$ is an order in $R$, the discriminant valuation coincides with the valuation of the discriminant of $\OO$ in the usual sense.
\end{enumerate}
\end{lem}

\subsection{Splitting type (\texorpdfstring{$1^21$}{1²1})}
When $R$ is wildly ramified, it is no longer possible for a reduced basis to be an extender basis (again by discriminant considerations). Here the case that concerns us is that $R = K \cross Q$ is a partially ramified cubic extension. Define
\[
  d_0 = v_K(\disc_K R) = v_K(\disc_K Q) \in \{2, 4, \ldots, 2e - 2, 2e, 2e + 1\}.
\]
Let $(\rho_1, \rho_2, \rho_3)$ with valuations $(a_1, a_2, a_3)$ be a reduced basis. Since $R \tensor Q \isom Q \cross Q \cross Q$ is an unramified extension of $Q$, there exists a well-led extender basis $(\bar{\rho}_1, \bar{\rho}_2, \bar{\rho}_3)$ of elements of $\omega I \tensor_K Q$. The extender indices $\bar a_i$ satisfy $\bar a_i \geq a_i$ by Proposition \ref{prop:rb}\ref{rb:max}. Our task in this section is to pin them down precisely.
\begin{lem} \label{lem:indices_wedge}
If $(\bar{\rho}_1, \ldots, \bar{\rho}_n)$ is an extender basis of $\omega I$, then the elementary indices $\bar a_i = v\left(\bar\rho_i\right)$ can be computed from
\[
  \bar a_1 + \cdots + \bar a_k = \min \left\{v(\eta) : \eta \in \Lambda^k(\omega I)\right\},
\]
the wedge power being taken in $\Lambda^k(\bar K^n) \isom \bar K^{\binom{n}{k}}$ with its valuation \eqref{eq:v_Kn}.
\end{lem}
\begin{proof}
Since the elementary vectors form a basis for $\OO_{\bar K}$, their wedge products form bases for $\Lambda^k \OO_{\bar K}$, and the minimum on the right-hand side is seen to equal the value on the left.
\end{proof}
So the first index $\bar a_1 = a_1$ is unchanged by extension (this is rather obvious), while the sum of all the indices is constrained by Lemma \ref{lem:extender_disc}, so
\begin{equation} \label{eq:disc_defect}
  (\bar a_1 + \bar a_2 + \bar a_3) - (a_1 + a_2 + a_3) = \frac{d_0 - 1}{2}.
\end{equation}

It remains to compute $\bar a_2$:

\begin{lem} \label{lem:a2bar}
Let $(\rho_i)_i$ be a well-led reduced basis of $\omega I$.
\begin{enumerate}[$($a$)$]
  \item If $\rho_1$ or $\rho_2$ is $K$-led, then $\bar a_2 = a_2$, and an extender basis for $\omega I$ is
  \[
    \left(\bar\rho_1, \bar\rho_2, \bar\rho_3\right) = \left(\rho_1, \rho_2, \pi^{a_3 + \frac{d_0 - 1}{2}} (0; 0; 1)\right).
  \]
  \item If $\rho_3$ is $K$-led, then
  \begin{equation}
    \label{eq:a2bar}
    \bar a_2 = \min \left\{a_2 + v^{(K)}(\xi_1),\; a_2 + v^{(K)}(\xi_2),\; a_2 + \frac{d_0 - 1}{2},\; a_3\right\},
  \end{equation}
  and there is an extender basis $(\bar\rho_1, \bar\rho_2, \bar\rho_3)$ for $\omega I$ of one of the following shapes:  \begin{enumerate}[$($i$)$]
    \item If $\bar a_2$ is one of the first three arguments to the minimum in \ref{eq:a2bar}, then we can take
    \[
      \bar\rho_1 = \rho_1, \quad \bar\rho_2 = \rho_2 - u\pi^{a_2-a_1} \rho_1
    \]
    for some $u \in \OO_{\bar K}^\cross$.
    \item If $\bar a_2 = a_3$, then we can take
    \[
      \bar\rho_1 = \rho_1, \quad \bar\rho_2 = \rho_3.
    \]
  \end{enumerate}
\end{enumerate}
\end{lem}

\begin{proof}
We denote the three places of $\bar K^3$ by $(K)$, $(Q1)$, and $(Q2)$ according to their relation to $R$. By Lemma \ref{lem:indices_wedge}, $\bar a_1 + \bar a_2$ is the minimum of the valuations
\[
  v(\rho_i^{(k)} \rho_j^{(\ell)} - \rho_i^{(\ell)} \rho_j^{(k)})
\]
of the nine $2 \times 2$ minors of the basis matrix
\[
  \begin{bmatrix}
    \rho_1^{(K)} & \rho_1^{(Q1)} & \rho_1^{(Q2)} \\
    \rho_2^{(K)} & \rho_2^{(Q1)} & \rho_2^{(Q2)} \\
    \rho_3^{(K)} & \rho_3^{(Q1)} & \rho_3^{(Q2)}
  \end{bmatrix}.
\]
If $\rho_1$ is $K$-led, then $\rho_2$ is $Q$-led and, in view of the valuation
\[
  v(\rho_1^{(K)} \rho_2^{(Q1)} - \rho_1^{(Q1)} \rho_2^{(K)})
  = v(\rho_1^{(K)} \rho_2^{(Q1)}) = a_1 + a_2
\]
of the upper left minor, we get $\bar a_2 \leq a_2$. But we know that $\bar a_2 \geq a_2$, so equality holds. In particular, $\rho_1 \wedge \rho_2$ is primitive in $\Lambda^2(\omega I) \tensor_{\OO_K} \OO_{\bar K}$, so we can keep $\rho_1$ and $\rho_2$ as our first two vectors in the extender basis. As for the last basis vector $\pi^{\bar a_3} \xi_3$, we already know that $\bar a_3 = a_3 + (d_0 - 1)/2$, and we can pick any $\xi_3 \in \OO_{\bar K}^3$ coprimitive with $\xi_1$ and $\xi_2$ without changing the span
\[
  \left\langle\rho_1, \rho_2, \bar\rho_3\right\rangle = \left\langle\rho_1, \rho_2, \pi^{\bar a_3}\OO_{\bar K}^3 \right\rangle.
\]
Since $\xi_1$ and $\xi_2$ necessarily lie in the sublattice
\[
  \left\{\left(\xi^{(K)}, \xi^{(Q1)}, \xi^{(Q2)}\right) \in \OO_{\bar K}^3 : \quad \xi^{(Q1)} \equiv \xi^{(Q2)} \mod \pi^{\frac{d_0 - 1}{2}} \right\},
\]
we get that $\xi_3 = (0; 0; 1)$ is a suitable coprimitive vector.

A similar argument works if $\rho_2$ is $K$-led.

If $\rho_3$ is $K$-led, we look at each of the minors in turn:
\begin{itemize}
  \item For $\rho_1^{(K)} \rho_2^{(Qi)} - \rho_1^{(Qi)} \rho_2^{(K)}$, there cannot be any cancellation, because the $K$-valuations differ by an integer (if both are finite) while the $Q$-valuations must differ by an integer plus $1/2$ to make the leaders a reduced basis of $Q$. So
  \begin{align*}
    v\left(\rho_1^{(K)} \rho_2^{(Qi)} - \rho_1^{(Qi)} \rho_2^{(K)}\right)
    &= \min\left\{v\left(\rho_1^{(K)} \rho_2^{(Qi)}\right), v\left(\rho_1^{(Qi)} \rho_2^{(K)}\right)\right\} \\
    &= \min\left\{a_1 + a_2 + v^{(K)}(\xi_1), a_1 + a_2 + v^{(K)}(\xi_2)\right\},
  \end{align*}
  accounting for the first two arguments to the minimum.
  \item The minor $\rho_1^{(Q1)} \rho_2^{(Q2)} - \rho_1^{(Q2)} \rho_2^{(Q1)}$ measures the discriminant of the lattice $\left\langle\rho_1^{(Q)}, \rho_2^{(Q)}\right\rangle$ and hence has valuation
  \[
    a_1 + a_2 + \frac{d_0 - 1}{2},
  \]
  accounting for the third argument to the minimum. The other minors involving the two $Q$-valuations are at least as large and so can be ignored.
  \item Finally, the minor $\rho_1^{(K)} \rho_3^{(Qi)} - \rho_1^{(Qi)} \rho_3^{(K)}$ has valuation
  \[
    v(\rho_1^{(K)} \rho_3^{(Qi)} - \rho_1^{(Qi)} \rho_3^{(K)})
    = v(\rho_1^{(K)} \rho_3^{(Qi)}) = a_1 + a_3,
  \]
  accounting for the last argument to the minimum. The minor $\rho_2^{(K)} \rho_3^{(Qi)} - \rho_2^{(Qi)} \rho_3^{(K)}$ has valuation $a_2 + a_3 > a_1 + a_3$ and can be ignored.
\end{itemize}
This completes the determination of $\bar a_i$. We now go about constructing a suitable extender basis $\left(\bar \rho_i\right)$. We take $\bar \rho_1 = \rho_1$. The heart of the matter lies in finding a suitable $\bar \rho_2$.

If $\rho_1 \wedge \rho_2$ achieves the minimal valuation $\bar a_1 + \bar a_2$, that is, if the minimum corresponds to one of the first three arguments in \eqref{eq:a2bar}, then dividing by $\pi^{a_1 + a_2}$, we get
\[
  \xi_1 \wedge \xi_2 \equiv 0 \mod \pi^{\bar a_2 - a_2}.
\]
Since the $\xi_i$ are primitive, we get that for some $u \in \OO_{\bar K}^\cross$,
\begin{equation} \label{eq:x_a2bar}
  \xi_2 \equiv u \xi_1 \mod \pi^{\bar a_2 - a_2}.
\end{equation}
Then $\bar \rho_2 = \rho_2 - u\pi^{a_2-a_1} \rho_1$ satisfies the needed divisibility
\[
  v\left(\bar \rho_2\right) \geq \bar a_2,
\]
and since $v\left(\bar \rho_1 \wedge \bar \rho_2\right) = \bar a_1 + \bar a_2$, equality must hold, that is, the primitive normalizations
\[
  \bar \xi_1 = \frac{\bar \rho_1}{\pi^{\bar a_1}} \textand
  \bar \xi_2 = \frac{\bar \rho_2}{\pi^{\bar a_2}}
\]
are linearly independent modulo $\mm_{\bar K}$. Then we can take any $\bar \xi_3$ independent from the previous two to complete the basis, as $\bar \rho_3 = \pi^{\bar a_3} \bar\xi_3$ will automatically lie in $\omega \bar I$.

If $a_3$ achieves the minimum in \eqref{eq:a2bar}, we use the same method, with $\rho_3$ in place of $\rho_2$. The exponent $\bar a_2 - a_3$ appearing in \eqref{eq:x_a2bar} vanishes, allowing us to take $u = 0$ and $\bar\rho_2 = \rho_3$.
\end{proof}

\subsection{The reduced basis of a cubic resolvent ring}
Let $C$ be a candidate resolvent for $\ttt$-traced quartic rings, and let $C_{\ttt}$ be the corresponding reduced resolvent, that is, the unique ring such that $C = \OO_K + \ttt^2 C_{\ttt}$. First, look at the reduced basis of $C_{\ttt} \subset R$ as a lattice; and look at its extender basis, a basis of $\bar C_{\ttt} = C_{\ttt} \tensor_{\OO_K} \OO_{\bar{K}}$. Because $1 \in C_{\ttt}$ is an element of minimal valuation, there are not so many cases:
\begin{itemize}
  \item If $R$ is tamely ramified, then the reduced basis
  \[
    C_{\ttt} = \left\langle 1, \pi^{b_1} \theta_1, \pi^{b_2} \theta_2 \right\rangle, \quad 0 \leq b_1 \leq b_2
  \]
  is also an extender basis. If $R$ is unramified, the $b_i$ are of course integers; if $R$ is tamely ramified, then by Corollary \ref{cor:idxs_mod_1}, we have
  \begin{equation} \label{eq:b_i_1pow3}
    \left\{b_1, b_2\right\} \equiv \left\{\frac{1}{3}, \frac{2}{3}\right\} \mod \ZZ.
  \end{equation}
    \item If $R \isom K \cross Q$ is partially wildly ramified, then the application of Lemma \ref{lem:a2bar} is simplified by the fact that $\theta_0 = 1$ is $K$-led. So we have a reduced basis
  \[
    C_{\ttt} = \left\langle 1, \pi^{b_1} \theta_1, \pi^{b_2} \theta_2 \right\rangle
  \]
  and an extender basis
  \[
    \bar{C}_{\ttt} = \left\langle 1, \pi^{b_1} \theta_1, \pi^{\bar b_2} \bar \theta_2 \right\rangle
  \]
  with
  \[
    \bar b_2 = b_2 + \frac{d_0 - 1}{2} \textand \bar \theta_2 = (0;0;1).
  \]
  We have $b_1, b_2, \bar b_2 \in \frac{1}{2} \ZZ$. By Corollary \ref{cor:idxs_mod_1}, the two $b_i$ are of different classes modulo $\ZZ$, so
  \begin{equation} \label{eq:b_i_1pow21}
    \bar b_2 \equiv b_1 + \frac{d_0}{2} \mod \ZZ.
  \end{equation}
\end{itemize}

Because $\theta_1$ is coprimitive to $1$, at most one pair of its three coordinates can be congruent modulo $\mm_{\bar K}$. We define the \emph{idempotency index} of $C$ with basis $\{1, \theta_1, \theta_2\}$ to be
\[
  s = \max_{i \neq j} v\left(\theta^{(i)} - \theta^{(j)}\right)
\]
Note that $s$ is infinite only when two coordinates of $\theta_1$ are exactly equal. Since $\theta_1$ is determined only up to finite precision by $C$, we can, and will, assume that $s$ is finite.

If $s > 0$, then there is a unique coordinate of $\bar K^3 \isom R \tensor \bar K$ at which $\omega_C$ has positive valuation. This defines a splitting $R \isom K \cross Q$ into a linear and a quadratic (possibly split) factor. We denote the three coordinates of $R$ by $(K)$, $(Q1)$, $(Q2)$, where $(K)$ is the distinguished one; thus we can write an element $\xi \in \bar K^3$ as
  \[
  \xi = \left(\xi^{(K)}; \xi^{(Q1)}; \xi^{(Q2)}\right) = \left(\xi^{(K)}; \xi^{(Q)}\right)
  \]
  where $\xi^{(Q)} = \left(\xi^{(Q1)}; \xi^{(Q2)}\right) \in \bar K^2 = Q \tensor_K \bar K$.

\begin{lem}\label{lem:s}
  If finite, the value of $s$ is constrained as follows:
  \begin{itemize}
    \item If $R$ is unramified, then $s$ is an integer.
    \item If $R$ has splitting type $(1^3)$ (residue characteristic $2$), then $s = 0$.
      \item If $R = K \cross Q$ has splitting type $(1^21)$, then either
    \begin{enumerate}[$($a$)$]
      \item\label{s:a} $b_1 \in \ZZ + \frac{1}{2}$ and $s = \frac{d_0 - 1}{2}$ (letter type \ref{type:A} below).
      \item\label{s:bcde} $b_1 \in \ZZ$ and $s = \frac{d_0}{2} + s'$ for some $s' \in \ZZ_{\geq 0}$
    \end{enumerate}
    In particular, $b_1 + s \equiv \frac{d_0}{2} \mod \ZZ$ and $s \equiv \bar b_2 \mod \ZZ$.
  \end{itemize}
\end{lem}
\begin{proof}
  In the tame splitting types this is immediate, knowing that $\theta_0 = 1$ is an extender vector for $C$ of minimal valuation.
    Assume that $R = K \cross Q$ is wildly ramified. Since
  \[
  \theta_1 = \left(\theta_1^{(K)}, \theta_1^{(Q1)}, \theta_1^{(Q2)}\right)
  \]
  is coprimitive to $1$ and the two $Q$-components are already congruent mod $\pi^{\frac{d_0 - 1}{2}}$, we must have
  \[
  v\left(\theta_1^{(K)} - \theta_1^{(Q1)}\right) = 0 \textand s = v\left(\theta_1^{(Q1)} - \theta_1^{(Q2)}\right) \geq \frac{d_0 - 1}{2}.
  \]
  If $b_1 \in \ZZ$, then $\theta_1 \in \OO_R$. Let $\fI$ be the linear functional on $Q$ defined by
  \[
  \fI(\xi) = \frac{\xi^{(Q1)} - \xi^{(Q2)}}{\pi_Q^{(Q1)} - \pi_Q^{(Q2)}},
  \]
  so that
  \[
  \fI(1) = 0, \quad \fI(\pi_Q) = 1.
  \]
  Then
  \[
  s = \frac{d_0}{2} + v\left(\fI(\theta_1)\right) = \frac{d_0}{2} + s',
  \]
  where $s' \geq 0$. If $b_1 \in \ZZ + 1/2$, then $\theta_1 \in \sqrt{\pi} \OO_R$, and in order for $\theta_1$ to be primitive, we must have $\sqrt{\pi}\theta_1^{(Q)} \sim \pi_Q$, so
  \[
  s = \frac{d_0 - 1}{2} + v\left(\fI(\sqrt{\pi}\theta_1^{(Q)})\right) = \frac{d_0 - 1}{2},
  \]
  as desired.
\end{proof}

\begin{lem}\label{lem:impish}
  By choosing a suitable extender basis $(1, \bar\theta_1, \bar\theta_2)$, we can assume that
  \begin{equation}\label{eq:impish}
    s \leq \bar b_2 - b_1.
  \end{equation}
  This $\bar\theta_1$ can be taken to lie in $\pi^{-b_1} R$, and thus be a reduced vector, except in the following situation, which we call the \emph{impish case:} $q = 2$, $R = K\cross K\cross K$, and $C = \OO_K + \pi^{b_1}\OO_R$ is a content ring, that is, $b_1 = b_2$.
\end{lem}
\begin{proof}
Take a reduced basis $(1, \theta_1, \theta_2)$ and convert it to an extender basis $(1, \theta_1, \bar\theta_2)$ with at most the last vector differing, as we explained in the previous paragraphs. Assume that the desired inequality \eqref{eq:impish} does not hold. Then, in particular, $s \neq 0$ and $R = K \cross Q$ admits a splitting.

We first consider the unramified case (splitting type $(111)$ or $(12)$). Here, as $\theta_0 = 1$ and $\theta_1$ both have their two $(Q)$-coordinates congruent modulo $\pi^{s}$, we must have
\[
  \theta_2^{(Q1)} \nequiv \theta_2^{(Q2)} \mod \pi
\]
to keep the three extender vectors linearly independent modulo $\mm_{\bar K}$. We then try the basis vector
\[
  \theta_1' = c \theta_1 + \pi^{b_2 - b_1} \theta_2, \quad c \in \OO_K.
\]
By construction, $\pi^{b_1} \theta_1' \in C_\ttt$. Moreover,
\[
  \theta_1'^{(Q1)} - \theta_1'^{(Q2)} = c\left(\theta_1^{(Q1)} - \theta_1^{(Q2)}\right) + \pi^{b_2 - b_1}\left(\theta_2^{(Q1)} - \theta_2^{(Q2)}\right)
\]
The latter term is exactly divisible by $\pi^{b_2 - b_1}$ and the first term has valuation at least $s > b_2 - b_1$. Thus, for arbitrary $c$, we get a $\theta_1'$ whose $(Q)$-coordinates have a difference of the desired valuation $b_2 - b_1$. If $\pi \nmid c$, $\theta_1'$ is coprimitive to $1$. We can always complete to a reduced basis $\{1, \pi^{b_1}\theta_1', \pi^{b_2}\theta_2'\}$.

We must consider the possibility that some other pair of coordinates $(K)$ and $(Qi)$ of $\theta_1'$ have been made equal modulo $\pi$. This can only occur when $R$ is totally split and (as $\theta_1'$ is coprimitive to $1$) when $b_2 = b_1$. Then the expressions in question are
\[
  \theta_1'^{(K)} - \theta_1'^{(Qi)} = c\left(\theta_1^{(K)} - \theta_1^{(Qi)}\right) + \pi^{b_2 - b_1}\left(\theta_2^{(K)} - \theta_2^{(Qi)}\right).
\]
The coefficient of $c$ is a unit, so this expression can vanish modulo $\pi$ for only one of the $q$ congruence classes modulo $\pi$. This excludes only two congruence classes, so there is a valid choice of $c$ unless $q = 2$. When $q = 2$, the ring
\[
  C_\ttt = \left\langle1, \pi^{b_1}(1;0;0), \pi^{b_1}(0;1;0)\right\rangle = \OO_K + \pi^{b_1} \OO_K^3
\]
is the content ring, and we are in the impish case.

When $R$ has splitting type $(1^2 1)$, we adapt the foregoing proof. We must be in case \ref{s:bcde} of Lemma \ref{lem:s}, as in case \ref{s:a}, the inequality \eqref{eq:impish} is automatic. Therefore $\theta_1 \in \OO_R$, and $\theta_2 \in \sqrt{\pi} R$ with its $(Q)$-coordinates units. Again consider
\[
  \theta_1' = c \theta_1 + \pi^{b_2 - b_1} \theta_2, \quad c \in \OO_K.
\]
Observe that
\[
  \pi^{b_2 - b_1} \left(\theta_2^{(Q1)} - \theta_2^{(Q2)}\right) \sim \pi^{b_2 - b_1} \cdot \pi^{(d_0 - 1)/2} = \pi^{\bar b_2 - b_1}.
\]
So in the expression
\[
\theta_1'^{(Q1)} - \theta_1'^{(Q2)} = c\left(\theta_1^{(Q1)} - \theta_1^{(Q2)}\right) + \pi^{b_2 - b_1}\left(\theta_2^{(Q1)} - \theta_2^{(Q2)}\right),
\]
the latter term is exactly divisible by $\pi^{\bar b_2 - \bar b_1}$ and the first term has valuation at least $s > \bar b_2 - \bar b_1$. So any $c \in \OO_K^\cross$ will work.
\end{proof}

In the proof below, we will ignore the impish case, making note at the end that our answers remain true in it (Section \ref{sec:impish}.)

\subsection{The dual lattice of the resolvent}
A common tool in understanding nonmaximal orders is their duals under the trace pairing. Hence it is fitting that we should understand the element $\omega_C \in \OO_{\bar K}^3$, unique up to scaling by $\OO_{\bar K}^\cross$, that satisfies the relations
\[
  \tr \omega_C = \tr (\theta_1 \omega_C) = 0;
\]
that is, $\bar K \left\langle\omega_C \right\rangle = \bar K \left\langle1, \theta_1\right\rangle^\perp$ under the trace pairing.

One explicit choice of $\omega_C$ is as follows: If $\tilde\theta = \pi^{b_1}\theta_1 = (\tilde\theta^{(1)}, \tilde\theta^{(2)}, \tilde\theta^{(3)}) \in C$ is the second reduced basis vector, then
\[
  \hat{\omega}_C =
  \left(\tilde\theta^{(2)} - \tilde\theta^{(3)}\right) \left(\tilde\theta^{(3)} - \tilde\theta^{(1)}\right) \left(\tilde\theta^{(1)} - \tilde\theta^{(2)}\right) \left(\tilde\theta^{(2)} - \tilde\theta^{(3)}; \tilde\theta^{(3)} - \tilde\theta^{(1)}; \tilde\theta^{(1)} - \tilde\theta^{(2)}\right).
\]
The symmetry ensures that $\hat{\omega}_C \in R$. Note that $N\left(\hat{\omega}_C\right) = (\disc \tilde\theta)^2$ is a square in $K^\cross$. Note also that
\[
  \tr\left(\hat\omega_C\right) = \tr\left(\theta_1\hat\omega_C\right) = 0.
\]
Then the valuations of $\hat \omega_C$ are
\[
\vec{v}(\hat{\omega}_C) = (4b_1+2s, 4b_1+s, 4b_1+s),
\]
so we can take
\[
\omega_C = \pi^{-4b_1-s} \hat\omega_C, \qquad \vec{v}(\omega_C) = (s, 0, 0).
\]
Also let
\begin{equation*}
\widetilde\omega_C = \pi^{-4b_1-2s} \hat\omega_C, \qquad \vec{v}(\widetilde\omega_C) = (0, -s, -s).
\end{equation*}
Observe that $\omega_C$ and $\widetilde\omega_C^{-1}$ are primitive (that is, have valuation $0$) in $\OO_{\bar K}^3$.

\section{Resolvent conditions}
Throughout the local calculation, $C_{\ttt} \subseteq R$ denotes the fixed reduced resolvent and
$C=\OO_K+\ttt^2C_{\ttt}$ the corresponding actual resolvent. Their reduced bases are related by
\[
C = \left\langle1, \pi^{b_1 + 2t} \theta_1, \pi^{b_2 + 2t} \theta_2\right\rangle \textand
C_\ttt = \left\langle 1, \pi^{b_1} \theta_1, \pi^{b_2} \theta_2 \right\rangle, \quad 0 \leq b_1 \leq b_2,
\]
where $t = v_K(\ttt)$. Let $\delta \in R^\cross$ be an element of square norm, and let $L$ be the corresponding quartic algebra with resolvent $R$. As we noted in the proof of Theorem \ref{thm:hcl_quartic}, an order $\OO \subseteq L$ is completely determined by the lattice $I_\OO \subseteq R$ such that $\OO/\OO_K = \kappa(I_\OO)$, where
\begin{align*}
  \kappa : \bar K^3 &\to \bar K^4 \\
  \xi &\mapsto \left(\tr_{\bar K^3/\bar K} \xi \omega \sqrt{\delta}\right)_\omega
\end{align*}
is the map in Proposition \ref{prop:Kummer_quartic}, with the four coordinates indexed by the norm-one $2$-torsion elements $\omega$ of $\bar K^3$. For reasons that will become clear below, we take the reduced and extender bases, not of $I_\OO$ itself, but of
\[
  I = \sqrt{\pi^s\delta\omega_C}\,I_\OO.
\]
Take a reduced basis
\[
  I = \left\langle\pi^{a_1} \xi_1, \pi^{a_2} \xi_2, \pi^{a_3} \xi_3 \right\rangle
\]
and an extender basis
\[
  \bar{I} = \left\langle\pi^{\bar a_1} \bar \xi_1, \pi^{\bar a_2} \bar \xi_2, \pi^{\bar a_3} \bar \xi_3 \right\rangle.
\]
Both $\delta$ and $\omega_C$ satisfy property \ref{iota:star} in Definition \ref{defn:red_basis}, so the work done in Section \ref{sec:red_basis} applies. In the tamely ramified splitting types, the reduced basis is also an extender basis by Proposition \ref{prop:ext_basis_tame}, so we let $\bar a_i = a_i$, $\bar \xi_i = \xi_i$. On the other hand, in splitting type $(1^21)$, only the first basis vector $\pi^{a_1} \xi_1 = \pi^{\bar a_1} \bar \xi_1$ can be carried over in general, and we will need both bases.

We also use the primed lattice
\[
  I' = \sqrt{\delta\hat\omega_C}\,I_\OO = \pi^{2b_1}I,
\]
whose reduced and extender indices are
\[
  a_i' = a_i + 2b_1 \textand \bar a'_i = \bar a_i + 2b_1.
\]


Then we can write the resolvent conditions as follows:
\begin{lem} \label{lem:rsv}
  With respect to the above setup, a lattice $I$ yields a ring $\OO$ with a $\ttt$-traced resolvent to $C$ if and only if the following conditions hold:
  \begin{itemize}
    \item Discriminant $(\Theta)$ condition:
    \begin{equation}
      v_K(\disc I) = v_K(\disc C) + 4s - 8e, \label{eq:disc_quartic}
    \end{equation}
  that is,
    \begin{equation}
      \bar a_{1} + \bar a_{2} + \bar a_{3} = \bar b_{1} + \bar b_{2} + 2s + 4t - 4e \label{eq:disc_cdn_quartic};
    \end{equation}
    \item Resolvent $(\Phi)$ conditions for $\bar\theta_2$-coefficients:
    \[
      \cM_{ij} \colon \quad  \tr(\bar\xi_i \bar\xi_j) \equiv 0 \mod \pi^{\bar m_{ij}}
    \]
    where
    \[
      \bar m_{ij} = \begin{cases}
        \bar b_2 + 2t - 2e - 2\bar a_i + s, & i = j \\
        \bar b_2 + 3t - 3e - \bar a_i - \bar a_j + s, & i \neq j.
      \end{cases}
    \]
    \item Resolvent $(\Phi)$ conditions for $\bar\theta_1$-coefficients:
    \[
      \cN_{ij} \colon \quad \text{All coordinates of} \quad \widetilde\omega_C^{-1} \cdot \bar\xi_i \bar\xi_j \quad \text{are congruent} \mod \pi^{\bar n_{ij}},
    \]
    where
    \[
      \bar n_{ij} = \begin{cases}
      \bar b_1 + 2t - 2e - 2\bar a_i + 2s, & i = j \\
      \bar b_1 + 3t - 3e - \bar a_i - \bar a_j + 2s, & i \neq j.
      \end{cases}
    \]
  \end{itemize}
\end{lem}
\begin{proof}
The conditions that $C$ is a $\ttt$-traced resolvent for $\OO$ are that all coefficients in the coordinate representations of $\Theta$, $\Theta^{-1}$, and $\Phi$ have nonnegative valuation. In particular, it is equivalent to study when $\bar{C} = C \tensor_{\OO_K} \OO_{\bar{K}}$ is a resolvent of $\bar{\OO} = \OO \tensor_{\OO_K} \OO_{\bar{K}}$.

We have
\[
\bar{C} = \left\langle 1, \pi^{2t + \bar b_1} \bar\theta_1, \pi^{2t + \bar b_2} \bar\theta_2 \right\rangle,
\]
so
\[
  \Lambda^2 (\bar{C}/\OO_{\bar{K}}) = \left\langle \pi^{4t + \bar b_1 + \bar b_2} \bar\theta_1 \wedge \bar\theta_2 \right\rangle
  = \pi^{4t + \bar b_1 + \bar b_2} \Lambda^2 (\OO_{\bar{K}}^3 / \OO_{\bar{K}})
\]
and, by the formula for $\Theta^{-1}$ in Proposition \ref{prop:Kummer_resolvent_quartic},
\[
  \Theta^{-1}(\Lambda^2 (\bar{C}/\OO_{\bar{K}})) = \frac{1}{16\sqrt{N(\delta)}} \pi^{4t + \bar b_1 + \bar b_2} \cdot \Lambda^3\kappa(\OO_{\bar{K}}^3)
\]
Meanwhile,
\begin{equation} \label{eq:I_basis_copy}
  \ba{I_\OO} = \frac{1}{\sqrt{\pi^{s} \delta \omega_C}}\bar{I} = \frac{1}{\sqrt{\pi^{s} \delta \omega_C}}\left\langle\pi^{\bar a_1} \bar\xi_1, \pi^{\bar a_2} \bar\xi_2, \pi^{\bar a_3} \bar\xi_3 \right\rangle
\end{equation}
so
\begin{align*}
  \Lambda^3 \ba{I_\OO} &= \frac{1}{\sqrt{\pi^{3 s} N(\delta \omega_C)}} \pi^{\bar a_1 + \bar a_2 + \bar a_3} \left\langle\bar\xi_1 \wedge \bar\xi_2 \wedge \bar\xi_3\right\rangle \\
  &= \frac{1}{\sqrt{\pi^{3s} N(\delta \omega_C)}} \pi^{\bar a_1 + \bar a_2 + \bar a_3} \Lambda^3(\OO_{\bar{K}}^4/\OO_{\bar{K}}).
\end{align*}
So the condition for $\Theta$ to map $\Lambda^3(\bar{\OO}/\OO_{\bar{K}})$ isomorphically onto $\Lambda^2(\bar{C}/\OO_{\bar{K}})$ is that
\begin{align*}
4t + \bar b_1 + \bar b_2 - 4e - \frac{1}{2} v(N(\delta)) = \bar a_1 + \bar a_2 + \bar a_3 - \frac{3}{2}s - \frac{1}{2} v(N(\delta\omega_C)),
\end{align*}
or, since $v(N(\omega_C)) = s$,
\[
  4t - 4e + \bar b_1 + \bar b_2 + 2s = \bar a_1 + \bar a_2 + \bar a_3,
\]
as desired.

Likewise, we use the formula
\[
  \Phi(\kappa(\alpha)) = 4\delta\alpha^2
\]
from Proposition \ref{prop:Kummer_resolvent_quartic} to transform the $\Phi$-condition to the following:
\begin{enumerate}[(i)]
  \item\label{xon} For every $\rho \in \bar{I}$, we have $4 \pi^{-s}\omega_C^{-1} \rho^2 \in \bar{C} + \bar{K}$;
  \item\label{xoff} For every $\rho,\sigma \in \bar{I}$, we have $8 \pi^{-t-s}\omega_C^{-1} \rho\sigma \in \bar{C} + \bar{K}$.
\end{enumerate}
In terms of the basis \eqref{eq:I_basis_copy} for $\bar{I}$, this is to say that the diagonal entries of the matrix representing $\Phi$ belong to $\bar{C} + \bar{K}$ and the off-diagonal entries to $\ttt\bar{C} + \bar{K}$. Hence it suffices to consider \ref{xon} for $\rho$ a basis element and \ref{xoff} for $\rho,\sigma$ distinct basis elements.

To test whether an $\alpha \in \bar{K}^3$ lies in
\[
  \bar{C} = \left\langle 1, \pi^{2t + \bar b_1} \bar\theta_1, \pi^{2t + \bar b_2} \bar\theta_2 \right\rangle,
\]
we can pair it with a basis of the dual lattice $\bar{C}^\vee$ with respect to the trace pairing. Let $(\lambda_0, \lambda_1, \lambda_2)$ be the dual basis to $(\bar\theta_0, \bar\theta_1, \bar\theta_2)$ (that is, $\tr(\bar\theta_i \lambda_j) = \One_{i = j}$). Then $(\lambda_1,\lambda_2)$ is a basis for $(\OO_{\bar K}^3)^{\tr = 0}$, and we have already met $\lambda_2$: it is $\omega_C$, up to a unit. Hence
\begin{align*}
  \bar{C}^\vee &= \left\langle\lambda_0, \pi^{-2t - \bar b_1} \lambda_1, \pi^{-2t - \bar b_2} \lambda_2 \right\rangle \\
  &= \left\langle\lambda_0\right\rangle + \pi^{-2t - \bar b_1}(\OO_{\bar{K}}^3)^{\tr = 0} + \pi^{-2t - \bar b_2} \left\langle{\omega}_C\right\rangle
\end{align*}
We actually wish to test not whether $\alpha \in \bar{C}$, but the weaker condition $\alpha \in \bar{C} + \bar{K}$, so (due to the natural duality between $\bar{C}/\OO_{\bar{K}}$ and $(\bar{C}^\vee)^{\tr = 0}$) we pair only with elements of
\[
  (\bar{C}^\vee)^{\tr = 0} = \pi^{-2t - \bar b_1}(\OO_{\bar{K}}^3)^{\tr = 0} + \pi^{-2t - \bar b_2} \left\langle{\omega}_C\right\rangle.
\]
Hence $\alpha \in \bar{C} + \bar{K}$ if and only if
\begin{itemize}
  \item $\pi^{-2t - \bar b_2}\tr({\omega}_C\alpha)$ is integral, and
  \item $\pi^{-2t - \bar b_1}\tr(\alpha \kappa)$ is integral for $\kappa \in \{(1;-1;0), (0; 1; -1)\}$; that is, all coordinates of $\alpha$ are congruent modulo $\pi^{2t + \bar b_1}$.
\end{itemize}
This is the origin of the $\cM$- and $\cN$-conditions respectively. Applying this to the values
\[
\alpha = 4 \pi^{-s}\omega_C^{-1} \rho^2, \quad \alpha = 8 \pi^{-t-s}\omega_C^{-1} \rho \sigma
\]
derived from the basis above yields the desired form of all the $\Phi$-conditions.
\end{proof}

We say that the condition $\cM_{ij}$ or $\cN_{ij}$ is \emph{active} if its corresponding modulus $\bar m_{ij}$ or $\bar n_{ij}$ is positive. An inactive condition is automatically satisfied (noting that $\widetilde\omega_C^{-1}$ has nonnegative valuations).

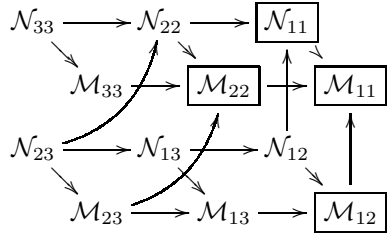
\begin{figure}
\[
\xymatrix@!0{
  \cN_{33} \ar[rr] \ar[rd] & & \cN_{22} \ar[rr] \ar[rd] & & \boxed{\cN_{11}} \ar[rd] \\
  & \cM_{33} \ar[rr] & & \boxed{\cM_{22}} \ar[rr] & & \boxed{\cM_{11}} \\
  \cN_{23} \ar[rr] \ar[rd] \ar@/_3ex/[urur] & & \cN_{13} \ar[rr] \ar[rd] & & \cN_{12} \ar[rd] \ar[uu] \\
  & \cM_{23} \ar[rr] \ar@/_3ex/[urur] & & \cM_{13} \ar[rr] & & \boxed{\cM_{12}} \ar[uu]
}
\]
\caption{Implications among the activity of the resolvent conditions}
\label{fig:inactives}
\end{figure}

Because of the inequalities $\bar a_1 \leq \bar a_2 \leq \bar a_3$, $\bar b_1 + s \leq \bar b_2$, and $0 \leq t \leq e$, the activity of the $\cM_{ij}$ and $\cN_{ij}$ follow a tight array of implications, as shown in Figure \ref{fig:inactives}.
The next lemma limits our concern to the four boxed conditions:
\begin{lem}\
\label{lem:inactives}
    \begin{enumerate}[$($a$)$]
        \item\label{inact:quart} Suppose that the $\bar a_i$ and $\bar\xi_i$ come from the extender decomposition of a quartic ring and a resolvent thereof. Then:
        \begin{itemize}
          \item No conditions are active except $\cM_{11}$, $\cM_{12}$, $\cM_{22}$, $\cN_{11}$, and $\cN_{12}$.
          \item $\cM_{12}$ and $\cM_{22}$ are not both active.
          \item $\cN_{12}$ is \emph{very weakly active,} that is, $\bar n_{12} \leq s/2$.
        \end{itemize}
        \item\label{inact:weak} Suppose that the $\bar a_i$ and $\bar\xi_i$ come from the extender decomposition of \underline{some} lattice $I \subseteq \omega R$. Suppose that the inactivity restrictions from part \ref{inact:quart} hold and that conditions $\cM_{11}$, $\cM_{12}$, $\cM_{22}$, and $\cN_{11}$ are satisfied. Then $\cN_{12}$ is satisfied, and the $\bar a_i$ and $\bar\xi_i$ actually come from a quartic ring. That is, ``$\cN_{12}$ is automatic if it is very weakly active.''
    \end{enumerate}
\end{lem}
\begin{proof}
\begin{enumerate}[$($a$)$]
  \item
Suppose that the $\bar a_i$ and $\bar\xi_i$ come from a quartic ring.

If $\cM_{13}$ is active, so are $\cM_{12}$ and $\cM_{11}$. Since the $\bar\xi_i$ are supposed to form an $\OO_{\bar{K}}$-basis of $\OO_{\bar K}^3$, we obtain for all $\bar\xi \in \OO_{\bar{K}}^3$,
\[
  \tr(\bar\xi_1\bar\xi) \equiv 0 \mod \mm_{\bar K}.
\]
Since $\bar\xi_1$ is primitive, this is a contradiction.

If $\cM_{12}$ and $\cM_{22}$ are active, then so is $\cM_{11}$. We have a $2$-dimensional subspace $V = \left\langle\bar\xi_1, \bar\xi_2\right\rangle$ of the $3$-dimensional space $k_{\bar K}^3$ that is isotropic for the trace pairing. But the trace pairing is nondegenerate, so this is a contradiction.

If $\cM_{33}$ is active, so are $\cM_{22}$ and $\cM_{11}$. If $\ch k_K \neq 2$, then $t = e = 0$ so $\cM_{12}$ is also active, and we have a contradiction as above. If $\ch k_K = 2$, we use that squaring is a linear operation mod $2$ to obtain that for all $\bar\xi \in \OO_{\bar{K}}^3$,
\[
  \tr(\bar\xi^2) \equiv 0 \mod \mm_{\bar K},
\]
which is a contradiction.

If $\cN_{22}$ is active, note that $\ch k_K = 2$ since $\cM_{22}$ is active. There are two cases. If $s = 0$, then $\widetilde\omega_C^{-1}$ is a unit, so the condition
\[
  \cN_{ii} \colon \quad \text{All coordinates of} \quad \widetilde\omega_C^{-1} \cdot \bar\xi_i^2 \quad \text{are congruent} \mod \pi^{\bar n_{ii}}
\]
determines $\bar\xi_i^2$ mod $\mm_{\bar K}$ up to scaling. But since we are in characteristic $2$, square roots are unique, and $\bar\xi_1$ and $\bar\xi_2$ are scalar multiples mod $\mm_{\bar K}$, a contradiction. Now assume $s > 0$, so $\widetilde\omega_C^{-1} \equiv (u;0;0)$ for some unit $u \in \OO_{\bar{K}}^\cross$. Now $\cN_{ii}$ gives $\bar\xi_i^{(K)} \equiv 0 \mod \mm_{\bar K}$. But now $\cM_{ii}$ gives that $\bar\xi_i$ is a unit multiple of $(0;1;1)$ modulo $\mm_{\bar K}$, a contradiction.

Finally, assume that $\cN_{12}$ is active and not very weakly active: that is, $\bar n_{12} > s/2$. Note that $\bar n_{11} > s$ since otherwise $\cN_{22}$ would be active. If $s = 0$, then $\cN_{11}$ implies that $\bar\xi_1^2 \equiv \omega_C \mod \mm_{\bar K}$, up to scaling, and then $\cN_{12}$ implies that $\bar\xi_1\bar\xi_2 \equiv \omega_C \mod \mm_{\bar K}$, up to scaling. Since $\omega_C$ is a unit, this is a contradiction.  So $s > 0$ defines a distinguished coordinate $\bullet^{(K)}$ on $C$. Note that $\omega_C$ is a unit multiple of $(0;1;-1)$ modulo $\pi^{s}$ and that $\widetilde\omega_C^{-1}$ is a unit multiple of $(1;0;0)$ modulo $\pi^{s}$. Now $\cN_{11}$ implies that
\[
  \bar\xi_1^2 \equiv \omega_C \mod \mm_{\bar K}\omega_C = \mm_{\bar K}(\pi^{s};1;1).
\]
So $v_K(\bar\xi_1^{(K)}) = s/2$ exactly. If $v_K(\bar\xi_2^{(K)}) > 0$ also, we get $\bar\xi_1 \equiv \bar\xi_2 \mod \mm_{\bar K}$ by the same argument as when $N_{22}$ is active. So $\bar\xi_2^{(K)}$ is a unit, and the $\cN_{12}$ condition
\[
  \cN_{12} \colon \quad \text{All coordinates of} \quad \widetilde\omega_C^{-1} \cdot \bar\xi_1 \bar\xi_2 \quad \text{are congruent} \mod \pi^{\bar n_{12}}
\]
is unsatisfied, because the $K$-coordinate has valuation $s/2$ and the others have higher valuation.

\item
Note that $\bar n_{12} \leq \frac{1}{2} \bar n_{11}$ because otherwise $\cN_{22}$ would be active. We have $\widetilde\omega_C^{-1}$ a unit multiple of $(1;0;0)$ modulo $\pi^{2\bar n_{12}}$, so $\cN_{11}$ gives $\pi^{2\bar n_{12}} | (\bar\xi_1^2)^{(K)}$, that is, $\pi^{\bar n_{12}} | \bar\xi_1^{(K)}$. Now $\pi^{\bar n_{12}} \mid \widetilde\omega_C^{-1} \bar\xi_1$, so $\cN_{12}$ is satisfied. \qedhere
\end{enumerate}
\end{proof}

Based on this, we will count quartic rings with fixed resolvent. Since a lot will happen with various things being fixed and others varying, it is worthwhile to lay down the following:

\begin{conv}\label{conv:quartic}
We fix variables in the following order:
\begin{itemize}
  \item We always fix the \emph{resolvent data,} which comprise
  \begin{itemize}
    \item a resolvent $C \subseteq R$.
    \item an extender decomposition $C = \left\langle1, \pi^{\bar b_1} \bar\theta_1, \pi^{\bar b_2} \bar\theta_2\right\rangle$, which fixes $\hat\omega_C$ and $s$. We can, and do, assume that $s$ is finite.
    \item a desired trace ideal $\ttt = \mm_K^t$, $0 \leq t \leq e$, which defines a reduced resolvent $C_{\ttt}$.
  \end{itemize}
  \item Then we fix the \emph{discrete data} of a quartic ring, which comprises
  \begin{itemize}
    \item a choice of \emph{coarse coset} $[\delta] \in H^1/\cF_0$, which determines the valuations of $\delta$. There are $\size{H^0}$ cosets $\delta_0 \cF_0$. Then $\delta = \delta_0 \upsilon$, where $\upsilon \in \OO_R^\cross$ can vary.
    \item its extender indices $\bar a_i$,
which are constrained by the conditions modulo $\ZZ$ needed for a sublattice of $\sqrt{\delta \hat\omega_C} R$ and the inactivity inequalities of Lemma \ref{lem:inactives}.
  \end{itemize}
  \item Then we choose $\delta$ and $\xi_1$, which are constrained by the $\cM_{11}$ and $\cN_{11}$ conditions.
  \item Lastly, we choose $\xi_2$ and $\xi_3$, equivalently $\bar\xi_2$ and $\bar\xi_3$. In either case the last basis vector carries no information. The basis vectors are constrained by coprimitivity and by the $\cM_{12}$ and $\cM_{22}$ conditions.
\end{itemize}
Whenever we speak about possibilities for any of the items on this list, it will be implicitly assumed that all the previous items have been fixed in conformity with their respective restrictions.
\end{conv}

Since $\xi_1 \in \pi^{-a_1'} \omega R$, conditions $\cM_{11}$ and $\cN_{11}$ can be viewed in another way, which will be simpler for some purposes:
\begin{lem} \label{lem:beta}
A $\xi_1 \in \pi^{-a_1'}\sqrt{\delta\hat\omega_C}R \intsec \OO_{\bar K}^3$ satisfies the $\cM_{11}$ and $\cN_{11}$ resolvent conditions if and only if
\[
  \beta = \frac{\xi_1^2}{\omega_C} \in \pi^{-2a_1+s}R
\]
is a linear combination of the extender basis vectors of $\bar C_{\ttt}$ of the form
\[
  \beta = x + y\pi^{\bar n_{11}-s}\bar\theta_1
              + z\pi^{\bar m_{11}}\bar\theta_2,
  \quad x\in\bar K,\quad y,z\in\OO_{\bar K},
\]
or equivalently, a linear combination of the reduced basis vectors of $C_{\ttt}$ of the form
\[
  \beta = x + y\pi^{n_{11}-s}\theta_1
              + z\pi^{m_\beta}\theta_2,
  \quad x\in\pi^{-2a_1+s}K,\quad y,z\in\OO_K,
\]
where
\[
  m_\beta = b_2 + 2t - 2e - 2a_1 + s
\]
is computed using $b_2$ instead of $\bar b_2$ as in $\bar m_{11}$. Since the first reduced and extender basis vectors agree, $n_{11}=\bar n_{11}$.
\end{lem}
\begin{proof}
Since $\pi^{a_1'}\xi_1=\pi^{\bar a_1'}\bar\xi_1$ is the first basis element of $\bar I'$, the $\cM_{11}$ and $\cN_{11}$ conditions are equivalent to
\begin{equation*}
\bar K+\bar C\ni
4\delta\left(\frac{\pi^{a_1'}\xi_1}{\sqrt{\delta\hat\omega_C}}\right)^2
=\frac{4\pi^{2a_1'}\xi_1^2}{\hat\omega_C}
=4\pi^{2a_1-s}\beta.
\end{equation*}
Thus
\begin{align*}
  \beta &\in \bar K+\pi^{s-2a_1-2e}\bar C \\
  &=\bar K+\pi^{\bar n_{11}-s}\OO_{\bar K}\bar\theta_1
              +\pi^{\bar m_{11}}\OO_{\bar K}\bar\theta_2 \\
  &=\bar K+\pi^{n_{11}-s}\OO_{\bar K}\theta_1
              +\pi^{m_\beta}\OO_{\bar K}\theta_2.
\end{align*}
The first description is immediate. Since
$(1,\pi^{b_1}\theta_1,\pi^{b_2}\theta_2)$ is an $\OO_K$-basis of $C_{\ttt}$, the second description follows with $x\in\pi^{s-2a_1}K$ and $y,z\in\OO_K$.
\end{proof}

\subsection{Transformation, and ring volumes in the white zone}
\label{sec:white}
We will proceed to compute the volumes of the solution sets in $\omega R$ in which the reduced vectors $\xi_i$ lie. We normalize our measures so that $\OO_R$ has volume $\mu(\OO_R) = 1$, and its projectivization $\PP(\OO_R)$ has volume $\mu(\PP(\OO_R)) = 1 + 1/q + 1/q^2$, so that an $n$-pixel has volume $q^{-2n}$ for $n \geq 1$.

We now transform $\xi_i$ from $\omega R$ to $R$.
\begin{lem}\label{lem:gamma_white}
Fix the discrete data. In particular, $\delta = \upsilon \delta_0$ lies in a fixed coarse coset. For $1 \leq i \leq 3$, there is a $\gamma = \gamma_{i} \in \bar{K}^3$ with the properties that, letting
\[
  \gamma_i = \sqrt{\upsilon}\,\gamma_{i,0},
\]
we have that $\xi_i' = \gamma_{i}^{-1} \xi_i$ is a primitive vector in $\OO_R$.
\end{lem}
\begin{proof}
Note that $\xi_i$ must lie in
\begin{equation} \label{eq:J_for_gamma}
  J = \pi^{-a_i'} \sqrt{\delta \hat\omega_C} R \intsec \OO_{\bar{K}}^3 = \sqrt{\upsilon} J_0,
\end{equation}
where
\[
  J_0 = \pi^{-a_i'} \sqrt{\delta_0 \hat\omega_C} R \intsec \OO_{\bar{K}}^3
\]
is an $\OO_R$-lattice of dimension $1$. As $\OO_R$ is a principal ideal ring (it's a product of DVR's), we obtain that $J_0 = \gamma_{i,0} \OO_R$ for some $\gamma_{i,0}$, clearly not a zero-divisor. Since $\xi_i \in J \setminus \pi J$, we get $\xi_i' = \gamma_i^{-1} \xi_i \in \OO_R \setminus \pi\OO_R$, as desired.
\end{proof}

\begin{lem} \label{lem:vol_to_rings}
  Let $\OO$ be a quartic ring. The set $V$ of triples $(\xi'_1, \xi'_2, \xi'_3)$ in $\PP(\OO_R)^3$ whose associated ring is $\OO$ has a volume determined by the discrete data of $\OO$ alone. It is given by
  \[
  \mu(V) = q^{\ds -\ceil{a_2 - a_1} - \ceil{a_3 - a_2} - \ceil{a_3 - a_1} + \frac{3d_0^\mathrm{tame}}{2} + v_K\big(N(\gamma_1\gamma_2\gamma_3)\big)} \cdot c,
  \]
  where
  \[
  d_0^{\mathrm{tame}} = \sum_{R_i} f_{R_i/K}\left(e_{R_i/K} - 1\right)
  \]
  is the standard lower bound for the discriminant valuation, attained for tame extensions, and
  \[
  c = \begin{cases}
    \left(1 + \frac{1}{q}\right)\left(1 + \frac{1}{q} + \frac{1}{q^2}\right) & \text{if } a_1 = a_2 = a_3 \\
    1 + \dfrac{1}{q} & \text{if } a_1 = a_2 < a_3 \text{ or } a_1 < a_2 = a_3 \\
    1 & \text{if } a_1 < a_2 < a_3. \\
  \end{cases}
  \]
\end{lem}
\begin{proof}
  We have
  \[
  \xi_i, \gamma_i \in \pi^{-a_i'} \sqrt{\delta\hat\omega_C} R, \quad \xi_i' \in \OO_R.
  \]

  To determine whether $(\Xi_1', \Xi_2', \Xi_3') \in V$, there are two conditions: firstly, the generators $\pi^{a_i} \Xi_i$, where $\Xi_i = \Xi_i'\gamma_i$ are the associated reduced basis vectors, belong to the correct lattice
  \begin{equation} \label{eq:Xi_in}
    \pi^{a_i} \Xi_i \in I = \left\langle\pi^{a_1} \xi_1, \pi^{a_2} \xi_2, \pi^{a_3} \xi_3\right\rangle;
  \end{equation}
  and secondly, they generate the whole of $I$.
  Since the $\pi^{a_i} \xi_i$ form a $K$-basis for $I\tensor_{\OO_K}K$, we can write
  \begin{equation}\label{eq:c_to_Xi}
    \Xi_i = \sum_{j = 1}^3 \pi^{a_j - a_i} c_{ij} \xi_j
  \end{equation}
  for some coefficients $c_{ij} \in K$. Condition \eqref{eq:Xi_in} is then equivalent to
  \[
  v\left(c_{ij}\right) \geq \max\{a_i - a_j, 0\},
  \]
  while the condition that the $\pi^{a_i} \Xi_i$ generate the whole of $I$ is equivalent to the change of basis being invertible:
  \[
  v\left(\det\left[c_{ij}\right]\right) = 0.
  \]
  Thus we have parametrized $V$ by the group
  \[
  \Gamma = \left\{[c_{ij}] \in \GL_3(\OO_K) : v\left(c_{ij}\right) \geq a_i - a_j\right\}.
  \]
  More precisely, $V$ is in continuous bijection with the cosets $T \backslash \Gamma$, where $T = \left(\OO_K^\cross\right)^3$ is the subgroup of diagonal matrices, because the $\Xi_i \in \PP(\OO_R)$ are defined only up to scaling.

  Without the invertibility condition, the volume of matrices in $\Mat_3(\OO_K)$ satisfying the valuation restrictions defining $\Gamma$ is
  \[
  q^{\ds -\sum_{1 \leq i,j \leq 3} \max\{0, \ceil{a_{i} - a_{j}}\}} = q^{\ds -\ceil{a_2 - a_1} - \ceil{a_3 - a_2} - \ceil{a_3 - a_1}}.
  \]
  The invertibility depends only on the $c_{ij}$ modulo $\pi$, and the fraction of matrices over $k_K$ of the shapes
  \[
  \begin{bmatrix}
    * & * & * \\
    * & * & * \\
    * & * & *
  \end{bmatrix},
  \begin{bmatrix}
    * & * & * \\
    0 & * & * \\
    0 & * & *
  \end{bmatrix},
  \textand
  \begin{bmatrix}
    * & * & * \\
    0 & * & * \\
    0 & 0 & *
  \end{bmatrix}
  \]
  that are nondegenerate is seen to be
  \[
  \left(1 - \frac{1}{q}\right)^3 \cdot c,
  \]
  accounting for the three cases in the definition of $c$. Projectivizing, $T\backslash \Gamma$ is a subset of $\PP(\OO_K^3)^3$ of volume
  \[
  q^{\ds -\ceil{a_2 - a_1} - \ceil{a_3 - a_2} - \ceil{a_3 - a_1}} \cdot c.
  \]

  It remains to compute how the volume transforms under the bijection $\Psi : T\backslash \Gamma \isom V$ that we have constructed. This map is $K$-linear and is a product of the three maps
  \begin{align*}
    \Psi_i : K \cross K \cross K &\to R \\
    \left(c_{i1}, c_{i2}, c_{i3}\right) &\mapsto \sum_{j=1}^3 c_{ij} \pi^{a_j - a_i} \frac{\gamma_j}{\gamma_i} \xi_j'.
  \end{align*}
  On the domain where it sends primitive vectors to primitive vectors, $\Psi_i$ scales volumes by $q^{-n_i}$, where $n_i$ is the determinant valuation, i.e.
  \begin{equation} \label{eq:x_det_rho}
    \Psi_i\left(\Lambda^3 \OO_K^3\right) = \pi^{n_i} \Lambda^3 \OO_R.
  \end{equation}
  Extending scalars to $\OO_{\bar K}$, the left side of \eqref{eq:x_det_rho} becomes
  \begin{align*}
    \Psi_i\left(\Lambda^3 \OO_{\bar K}^3\right) &= \prod_{j} \pi^{a_j - a_i} \bigwedge_j \frac{\gamma_j \xi_j'}{\gamma_i} \\
    &= \frac{\pi^{\sum_j a_j - 3 a_i}}{N(\gamma_i)} \bigwedge_j \xi_j
  \end{align*}
  When $R$ is tamely ramified, the wedge product of the $\xi_j$ generates the whole of $\Lambda^3 \OO_{\bar K}$, because the $\xi_j$ are an extender basis.
    When $R = K\cross Q$ is partially wildly ramified, the wedge product of the $\xi_j$ generate only $\pi^{(d_0 - 1)/2} \Lambda^3 \OO_{\bar K}$, as we saw in \eqref{eq:disc_defect}.
    Indeed, in all cases, if we let
  \[
  d_0^{\mathrm{wild}} = d_0 - d_0^{\mathrm{tame}},
  \]
  then
  \[
  \left\langle\xi_1 \wedge \xi_2 \wedge \xi_3 \right\rangle = \pi^{d_0^{\mathrm{wild}}/2} \Lambda^3 \OO_{\bar K}^3.
  \]
  Accordingly, we get
  \begin{align*}
    \Psi_i\left(\Lambda^3 \OO_{\bar K}^3\right) &= \frac{\pi^{\sum_j a_j - 3 a_i + \frac{d_0^{\mathrm{wild}}}{2}}}{N(\gamma_i)} \Lambda^3 \OO_{\bar K}.
  \end{align*}
  Meanwhile, the right side of \eqref{eq:x_det_rho} is
  \begin{align*}
    \pi^{n_i} \Lambda^3 \left(\OO_R \tensor_{\OO_K} \OO_{\bar K}\right) &= \pi^{n_i + \frac{d_0}{2}} \Lambda^3 \OO_{\bar K}.
  \end{align*}
  Hence
  \[
  n_i = \sum_{j} a_j - 3a_i - v_K\left(N(\gamma_i)\right) - \frac{d_0^{\mathrm{tame}}}{2}
  \]
  so
  \[
  \sum_i n_i = -v_K\left(N(\gamma_1\gamma_2\gamma_3)\right) - \frac{3d_0^{\mathrm{tame}}}{2}
  \]
  and
  \begin{align*}
    \mu(V) &= \mu\left(T\backslash\Gamma\right) \cdot q^{-\sum_i n_i} \\
    &= q^{\ds -\ceil{a_2 - a_1} - \ceil{a_3 - a_2} - \ceil{a_3 - a_1} + \frac{3d_0^{\mathrm{tame}}}{2} + v_K\left(N(\gamma_1\gamma_2\gamma_3)\right)} \cdot c,
  \end{align*}
  as desired.
\end{proof}

Consequently, we can compute the number of rings with any given discrete data by finding the volume of permissible $(\xi'_1, \xi'_2, \xi'_3)$, which we call the \emph{ring volume,} and dividing by $\mu(V)$. After suitably dividing up the cases, the ring volume will be a product of three \emph{ring volumes} for each $\xi'_i$ given the values of the preceding $\xi_j$. We carry out the computation of this volume in the succeeding sections.

Observe that $\mu(V)$ equals $q^{2a_1 - 2a_3}$ times a correction that depends only on the $a_i$ mod $1$ and whether any $a_i$ are equal. This will simplify the entry of the ring volumes into Sage at the end of the proof.

We may begin with the \emph{white zone} where no $\cM_{ij}$ or $\cN_{ij}$ is active. The ring volumes, which in the case of $\xi'_3$ are applicable in all cases, are computed below.
\subsubsection{Unramified splitting types.}
\begin{itemize}
  \item If $[\delta \hat\omega_C] \in \cF_0$, then  all $a_i$ are integers, and $\gamma_i = \gamma$ is a unit. The three $\xi'$ must form a basis of $\OO_R$. If they are found successively, their ring volumes are respectively $\boxed{1 + 1/q + 1/q^2}$, $\boxed{1 + 1/q}$, and $\boxed{1}$.
  \item If $[\delta \hat\omega_C] \in (1; \pi; \pi) \cF_0$, then  there is one $a_{j_0} \in \ZZ$; there $\gamma_{j_0} \sim (1; \sqrt{\pi}; \sqrt{\pi})$, and
  \[
    \xi'_{j_0} \in \OO_K^\cross \cross \OO_Q,
  \]
  a subset whose projectivization has volume $\boxed{1}$. Meanwhile, two $a_{j_1}, a_{j_2}$ lie in $\ZZ + 1/2$; there $\gamma_{j_k} \sim (\sqrt{\pi}; 1; 1)$, and
  \[
    \xi'_{j_k} \in \OO_K \cross \OO_Q
  \]
  with their $\OO_Q$-coordinates forming a basis of $\OO_Q$; thus $\xi'_{j_1}$ has volume $\boxed{1 + 1/q}$ and $\xi'_{j_2}$ has volume $\boxed{1}$.
\end{itemize}
\subsubsection{Splitting type $(1^3)$.}
\begin{itemize}
  \item Here $[\delta \hat\omega_C] \in \cF_0$. The $a_i$ fill out the classes in $\frac{1}{3}\ZZ / \ZZ$, but all $\gamma_i$ are units and all $\xi'_i$ lie in $\OO_R^\cross$, a subset whose projectivization has volume $\boxed{1}$.
\end{itemize}
  \subsubsection{Splitting type $(1^21)$.}
\begin{itemize}
  \item If $[\delta \hat\omega_C] \in \cF_0$, then by Corollary \ref{cor:idxs_mod_1}, the $a_i'$ form the multiset $\{0, 0, 1/2\}$ modulo $1$. For the two $a_i \in \ZZ$, we have $\gamma_i \sim 1$, and the value of $\xi'_i$ must lie off the $1$-pixel of $(0; \pi_R)$, yielding a ring volume of $\boxed{1 + 1/q}$ for the first such $i$ and $\boxed{1}$ for the second (one $1$-pixel being excluded by linear independence). For the one $a_i \in \ZZ + 1/2$, we have $\gamma_i \sim (\sqrt{\pi} ; 1)$, and $\xi'_i$ must have a unit $Q$-component, yielding a ring volume of $\boxed{1}$.
  \item If $[\delta \hat\omega_C] \in (\pi; \pi_R)\cF_0$, then by Corollary \ref{cor:idxs_mod_1}, the $a_i'$ form the multiset $\{1/4, 1/2, 3/4\}$ modulo $1$. The $\gamma_i$ in the three cases are associate to
  \[
    (\pi^{1/4}; 1; 1), \quad (1; \pi^{1/4}; \pi^{1/4}), \textand (\pi^{3/4}; 1; 1)
  \]
  No matter what $a_i'$ is, it determines whether $\xi'_i$ must have a unit $K$-component or unit $Q$-component, yielding a ring volume of $\boxed{1}$ in all cases.
  In practice, in this splitting type, the ring volumes will have to be diminished, even in the white zone, because of the restrictions on $v^{(K)}(\xi_i)$ that govern the relations between the $a_i$ and $\bar a_i$.
\end{itemize}
For future reference, it is worth noting that, in splitting type $(1^21)$,
\[
  v\big(N(\gamma_i))\big) = \frac{o_i}{4}.
\]
We define the \emph{flavor} of a set of discrete data to be the triple $(o_1, o_2, o_3)$, where $o_i = o \cdot \{a_i\}$ is the fractional part, rescaled by the common denominator appropriate to the splitting type:
\[
  o = \begin{cases}
    1 & \text{in splitting type $(3)$} \\
    2 & \text{in splitting types $(111)$ and $(12)$} \\
    3 & \text{in splitting type $(1^3)$} \\
    4 & \text{in splitting type $(1^21)$}.
  \end{cases}
\]
The allowable flavors are
\begin{equation}\label{eq:flavors}
  \begin{cases}
    000 & \text{in splitting type $(3)$} \\
    000,011,101,110 & \text{in splitting types $(111)$ and $(12)$} \\
    012,021,102,120,201,210 & \text{in splitting type $(1^3)$} \\
    002,020,200;123,132,213,231,312,321 & \text{in splitting type $(1^21)$}.
  \end{cases}
\end{equation}

\chapter{The conic over \texorpdfstring{$ \OO_K $}{OK}}
\label{cha:conic}
As noted above, the $\cM_{11}$ and $\cM_{22}$ conditions have the form of conics. Here we study relevant properties of conics over a local field.

\section{Conics of low discriminant over a local field}

For each nonzero $\alpha \in \OO_R$, the equation
\[
  \tr(\alpha\xi^2) = 0
\]
defines a conic in $\PP(\OO_R)$, viewed as a projective plane over $K$. Its determinant is $D_0 N(\alpha)$, up to squares of units, with respect to any $\OO_K$-basis of $\OO_R$, where $D_0$ is the discriminant of $\OO_R$. As we will find, it is preferable to transform the conic so that its discriminant has as low valuation as possible:

\begin{defn}
  Let $K$ be a local field, $\ch K \neq 2$. By a \emph{conic} over $\OO_K$ we mean a lattice $\Xi$ of dimension $3$ over $\OO_K$ equipped with an integral bilinear form $\cC$, or equivalently an integer-matrix quadratic form $\cC : \Xi \to \OO_K$, up to scaling by $\OO_K^\cross$. We say that $\cC$ is
  \begin{itemize}
    \item \emph{unimodular} if $\det \cC \sim 1$ (note that $\det \cC$ is uniquely defined up to squares of units);
    \item \emph{tiny} if $\det \cC \sim \pi$ and there exists a $\xi \in \Xi$ such that $\cC(\xi)$ is a unit;
    \item \emph{relevant} if it is either unimodular or tiny.
  \end{itemize}
\end{defn}
\begin{rem}
  Between changing basis and rescaling the whole form $\cC$, we can scale $\det \cC$ by any unit: hence we will sometimes assume that $\det \cC$ is exactly $1$ or $\pi$.
\end{rem}

Let $\diamondsuit$ be a generator of the different ideal $\dd_{R/K}$. For instance, we can take
\begin{equation} \label{eq:diamond_choice}
\diamondsuit = \begin{cases}
  1 & R \text{ unramified} \\
  \pi_R^2 & R \text{ totally tamely ramified}
  \\
    (1; \bar{\zeta}_2\sqrt{D_0}), & R \isom K \cross Q \text{ partially wildly ramified.}
\end{cases}
\end{equation}
Then by the definition of the different, the formula
\[
  \lambda^\diamondsuit(\xi) = \tr \frac{\xi}{\diamondsuit}
\]
defines a linear functional $\lambda^\diamondsuit : \OO_R \to \OO_K$ that is \emph{perfect,} that is, $\lambda^\diamondsuit$ generates the dual $\OO_R^\vee$ as an $\OO_R$-module, and hence the pairing
\[
  (x, y) \mapsto \lambda^\diamondsuit(xy)
\]
is a perfect $\OO_K$-linear pairing on $\OO_R$. If $\alpha \in \OO_R$, then the conic
\[
  \lambda^\diamondsuit(\alpha \xi^2)
\]
is $\OO_K$-integral on $\OO_R$ (because the corresponding bilinear form $\lambda^\diamondsuit(\alpha\xi\eta)$ is integral) with determinant $N(\alpha)$. We will put the conic defined by the $\cM_{11}$  and $\cM_{22}$ conditions in this form.

The entities involved in transformation will be marked by the symbol $\odot$ (``odot''). This is the symbol for a circle in Euclidean geometry. After the transformation, the solution volumes depend only on the indicated invariants and not on the choice of basepoint, as follows from the Igusa zeta functions computed in Lemmas \ref{lem:conic_1} and \ref{lem:conic_pi}.

First we isolate the \emph{ivory zone} where $\cM_{11}$ turns into a linear condition.
\begin{lem} \label{lem:conic_linear}\label{lem:ivory}
  Assume that $R$ has splitting type $(1^21)$. If
  \[
  m_{11} \leq \frac{d_0}{2},
  \]
  then $\cM_{11}$ is equivalent to the linear condition
  \[
  \xi_1^{(K)} \equiv 0 \mod \pi^{m_{11}/2}.
  \]
\end{lem}
\begin{proof}
  The $\cM_{11}$-condition is given as
  \[
  \tr(\xi_1^2) \equiv 0 \mod \pi^{m_{11}}.
  \]
  Note that
  \[
  m_{11} = \bar b_2 - 2a_1 - 2e + 2t + s \equiv -2a_1 \equiv \frac{o_1}{2} \mod 1.
  \]
  We have
  \[
  \tr(\xi_1^2) = \lambda^\diamondsuit\left(\diamondsuit \xi^2\right) = \left(\diamondsuit \xi^2\right)^{(K)} + \fI\left(\left(\diamondsuit \xi^2\right)^{(Q)}\right).
  \]
  Our strategy will be to prove that the second term vanishes modulo $\pi^{m_{11}}$, so that $\cM_{11}$ simplifies to
  \begin{align*}
    \diamondsuit^{(K)} {\xi_1^{(K)}}^2 &\equiv 0 \mod \pi^{m_{11}} \\
    {\xi_1^{(K)}}^2 &\equiv 0 \mod \pi^{m_{11}} \\
    \xi_1^{(K)} &\equiv 0 \mod \pi^{m_{11}/2}.
  \end{align*}
  To this end, we analyze
  \[
  \fI\left(\left(\diamondsuit \xi^2\right)^{(Q)}\right) = \fI\left(\bar\zeta_2\sqrt{D_0} {\xi^{(Q)}}^2\right)
  \]
  in the two cases:
  \begin{itemize}
    \item If $o_1$ is even, then ${\xi^{(Q)}}^2 \in \OO_Q$, so the argument to $\fI$ belongs to
    \[
    \bar\zeta_2 \sqrt{D_0} \OO_Q = \pi_Q^{d_0} \OO_Q = \OO_K \left\langle \pi^{\ceil{d_0/2}}, \pi^{\floor{d_0/2}} \pi_Q \right\rangle.
    \]
    The $\fI$-values of the last two generators are $0$ and $\pi^{\floor{d_0/2}}$ respectively, yielding the claimed congruence modulo $\pi^{\floor{d_0/2}}$, which is sufficient since $m_{11} \leq \floor{d_0/2}$.
    \item If $o_1$ is odd, then ${\xi^{(Q)}}^2 \in (\pi_Q/\sqrt{\pi})\OO_Q$, so the argument to $\fI$ belongs to
    \[
    \bar\zeta_2 \sqrt{D_0} \frac{\pi_Q}{\sqrt{\pi}} \OO_Q = \frac{\pi_Q^{d_0 + 1}}{\sqrt{\pi}} \OO_Q = \OO_K \left\langle \pi^{\floor{d_0/2} + 1/2}, \pi^{\ceil{d_0/2} - 1/2} \pi_Q\right\rangle.
    \]
    The $\fI$-values of the last two generators are $0$ and $\pi^{\ceil{d_0/2} - 1/2}$ respectively, yielding the claimed congruence modulo $\pi^{\ceil{d_0/2} - 1/2}$, which is sufficient since $m_{11} \leq \ceil{d_0/2} - 1/2$. \qedhere
  \end{itemize}
\end{proof}

Now, on to the cases where we transform $\cM_{ii}$ to a relevant conic.
\begin{lem}\label{lem:tfm_conic}
  Let $i \in \{1,2\}$. Fix the discrete data such a way that $\cM_{ii}$ is active. In splitting type $1^2 1$, assume that $m_{ii} > d_0/2$. Recall that $\delta = \delta_0\upsilon$ is in a fixed coarse coset. Then there is a multiplier $\gamma^{\odot}_i \in R^\cross$, an offset $p^\odot \in \QQ$, and a value $\delta_0^\odot \in \OO_R$ with the following properties:
  \begin{enumerate}[$($a$)$]
    \item\label{conic:prim} First,
    \[
    \xi_i^\odot = \frac{\xi_i'}{\gamma^{\odot}_i} = \frac{\xi_i}{\gamma_i \gamma^{\odot}_i \sqrt{\upsilon}}
    \]
    is a primitive vector in $\OO_R$, that is, $\gamma_i \gamma_i^\odot$ is an alternative to $\gamma_i$ in Lemma \ref{lem:gamma_white};
    \item\label{conic:tfm} the $\cM_{ii}$ condition $\tr(\xi_i^2) \equiv 0 \mod \pi^{m_{ii}}$ is equivalent to a condition of the form
    \[
    \lambda^\diamondsuit\left(\delta^\odot {\xi^\odot}^2\right) \equiv 0 \mod \pi^{m_{ii}^\odot},
    \]
    where $m_{ii}^\odot = m_{ii} - p^\odot$ is an integer and $\delta^\odot = \delta_0^{\odot}\upsilon \in \OO_R$;
    \item\label{conic:class} The classes of $\delta^\odot_0$ and $\delta^\odot$ in $H^1$ are given by
    \[
    [\delta^\odot_0] = [\delta_0 \hat\omega_C \diamondsuit], \quad [\delta^\odot] = [\delta \hat\omega_C \diamondsuit];
    \]
    \item\label{conic:norm} We have
    \[
      v_K(N(\delta^\odot_0)) \in \{0,1\},
    \]
    so that the conic $\cM^\odot(\xi^\odot) = \lambda^\diamondsuit\left(\delta^\odot {\xi^\odot}^2\right)$ is relevant;
    \item\label{conic:unim} Moreover, $\cM^\odot$ is unimodular exactly when
    \[
    [\delta_0 \hat\omega_C \diamondsuit] \in \cF_0.
    \]
  \end{enumerate}
\end{lem}
The details of the transformation are listed in Table \ref{tab:tfm_conic}, where $h_i$ is a quantity defined in the course of the proof.
\begin{table}[bt]
  \caption{The constants used in transforming the conic $\cM_{11}$}
  \label{tab:tfm_conic}
  \begin{tabular}{cccc|ccc|cc}
    spl.t. & $o_i$ & $h_i$ & $[\delta\hat\omega_C\diamondsuit] \in$ & \small $v^{(K)}(\gamma\gamma^\odot)$ & \small $v^{(Q)}(\gamma\gamma^\odot)$ & $p^\odot$ & $\delta^\odot \sim$ & $\xi^\odot \sim$ \\ \hline
    ur & $0$ & $0$ & $\cF_0$ & $0$ & $0$ & $0$ & $1$ & ? \\
    ur & $1$ & $1$ & $(1;\pi;\pi)\cF_0$ & $1/2$ & $0$ & $0$ & $(\pi ; 1 ; 1)$ & $(? ; 1 ; 1)$ \\
    &&&&&&&& \\
    $(1^3)$ &   & $1$ & $\cF_0$ & $0$ & [n/a] & $2/3$ & $1$ & $1$ \\
    $(1^3)$ &   & $-1$ & $\cF_0$ & $-2/3$ & [n/a] & $-2/3$ & $1$ & $\pi_R^2$ \\
    &&&&&&&& \\
    $(1^21)$ & & $0$ & $\cF_0$ & $d_0/4$ & $0$ & $d_0/2$ & $1$ & $(? ; 1)$ \\
    $(1^21)$ & & $1$ & $(\pi; \pi_Q) \cF_0$ & $\dfrac{d_0 - 1}{4}$ & $0$ & $\dfrac{d_0 - 1}{2}$ & $(1; \pi_Q)$ & $(? ; 1)$ \\
    $(1^21)$ & & $2$ & $\cF_0$ & $\dfrac{d_0 - 2}{4}$ & $-1/2$ & $\dfrac{d_0 - 2}{2}$ & $1$ & $(? ; \pi_Q)$ \\
    $(1^21)$ & & $3$ & $(\pi; \pi_Q) \cF_0$ & $\dfrac{d_0 - 3}{4}$ & $-1/2$ & $\dfrac{d_0 - 3}{2}$ & $(1; \pi_Q)$ & $(? ; \pi_Q)$
      \end{tabular}
\end{table}
\begin{proof}[Proof of Lemma \ref{lem:tfm_conic}]
We have $\xi_i \in \pi^{-a_i - 2b_1} \sqrt{\delta \hat\omega_C} \cdot R$.

Note that whatever $\gamma^{\odot}_i$ we pick, the conic takes the form
\[
  \lambda^\diamondsuit \left(\diamondsuit {\gamma_i^\odot}^2 \gamma_i^2 \upsilon {\xi_i^\odot}^2\right) \equiv 0 \mod \pi^{\bar m_{ii}},
\]
or, for any $p^\odot$,
\[
  \lambda^\diamondsuit\left( \frac{\diamondsuit {\gamma_i^\odot}^2 \gamma_i^2}{\pi^{p^\odot}} \upsilon {\xi_i^\odot}^2 \right) \equiv 0 \mod \pi^{\bar m_{ii} - p^\odot}.
\]
So we seek to pick $\gamma_i^{\odot}$ and $p^\odot$ so that
\[
  \delta^\odot_0 = \frac{\diamondsuit {\gamma_i^\odot}^2 \gamma_i^2 \delta_0}{\pi^{p^\odot}}
\]
lies in $\OO_R$ with norm of valuation $0$ or $1$, and
\[
  m_{ii}^\odot = m_{ii} - p^\odot
\]
is an integer. The second condition is easily seen to follow from the first. Part \ref{conic:class} is automatic. We now divide into cases based on the splitting type of $C$.

\subsection {Unramified.}
If $[\delta_0 \hat\omega_C] \in \cF_0$, then all extender indices are in $\ZZ$ and $\gamma_i$ is a unit. So we can choose $\gamma^{\odot}_i = 1$ and $p^\odot = 0$ and get $\delta^\odot_0$ a unit as well.

If $[\delta_0 \hat\omega_C] \in (1;\pi;\pi)\cF_0$ for some ordering of the coordinates, then we have $[\delta] = [(1;\pi;\pi)\delta']$ for some unit $\delta'$. By Corollary \ref{cor:idxs_mod_1}, either
\begin{itemize}
  \item $a_i' \in \ZZ$ and $\xi_i \in \sqrt{\delta'}(\OO_K^\cross \cross \sqrt{\pi} \OO_Q)$, or
  \item $a_i' \in \ZZ + \frac{1}{2}$ and $\xi_i \in \sqrt{\delta'}(\sqrt{\pi} \OO_K \cross \OO_Q)$.
\end{itemize}
In the first case, $\cM_{ii}$ is unsatisfiable if active, because $\alpha\xi_i^2$ has exactly one coordinate of zero valuation. So we have the second case. Observe that $\vec{v}(\gamma_i) = (1/2,0,0)$, so choosing $\gamma^{\odot}_i = 1$ and $p^\odot = 0$, we get $\vec{v}(\delta^\odot) = (1,0,0)$: the conic has determinant $\sim \pi$. Note that $\cM^\odot \nequiv 0 \mod \pi$ as a quadratic form: after passing to an unramified extension we may assume that $L = K \cross K \cross K$, and then $\cM^\odot$ is diagonal with two of the three coefficients units. So $\cM^\odot$ is tiny.

For convenience, we let
\[
  h_i = \begin{cases}
    0, & [\delta_0 \hat\omega_C] \in \cF_0 \\
    1, & \text{otherwise.}
  \end{cases}
\]

\subsection {Splitting type $(1^3)$.} Here $\OO_R$ is generated by a uniformizer $\pi_R$ with $\pi_R^3 = \pi$ (for a suitably chosen uniformizer $\pi$). Under the Minkowski embedding, $\pi_R = \sqrt[3]{\pi} \cdot \bar{\zeta}_3$, where
\[
  \bar{\zeta}_3 = (1; \zeta_3; \zeta_3^2).
\]
We have $\diamondsuit = \pi_R^2$.

Here $[\delta_0 \hat\omega_C] \in \cF_0$ always. The extender indices are in $\frac{1}{3}\ZZ$, and $\gamma_{i}$ is always a unit. Let $h_i \in \{0,1,-1\}$ be the integer such that
\[
  a_i' \in \ZZ - \frac{h_i}{3}.
\]
Then
\[
  \xi_i \in \pi^{h_i/3} \sqrt{\delta^\odot} R = \bar{\zeta}_{3}^{-h_i} \sqrt{\delta^\odot} R;
\]
indeed, since $\xi_i$ is primitive in $\OO_{\bar{K}}^3$,
\[
  \xi_i \in \bar{\zeta}_{3}^{-h_i} \sqrt{\delta^\odot} \OO_R^\cross.
\]

If $h_i = 0$, then $\cM_{ii}$ is unsatisfiable because the trace of a unit in $\OO_R$ is always a unit.

If $h_i = 1$, the choice $\gamma_i^{\odot} = 1$, $p^\odot = 2/3$ works, making $\delta^\odot$ a unit.

If $h_i = -1$, we can no longer take $\gamma_i^{\odot} = 1$, because the maximal possible value for $p^\odot$ is $1/3$ and the corresponding conic has determinant $\sim \pi^2$. Instead, take $\gamma_i^{\odot} = \pi_R^{-2}$. Then the corresponding values of $\xi^\odot$, instead of being units, have valuation $2/3$ and thus are still primitive in $\OO_R$. Take $p^\odot = -2/3$ and observe that $\delta^\odot$ is again a unit.

\subsection {Splitting type $(1^21)$.} Here things are more delicate because
\[
\diamondsuit = \bar{\zeta}_2\sqrt{D_0} \sim (1; \pi_R^{d_0})
\]
is significantly far from being a unit. Recall that we are assuming that $m_{11} > d_0/2$. Already, Lemma \ref{lem:ivory} ensures that $v_K(\xi_i^{(K)}) \geq (d_0 - 1)/4$. In particular, $\rho_i$ is $Q$-led, so either
\begin{itemize}
  \item $[\delta\hat\omega_C] \in \cF_0$ and $a_i' \in \frac{1}{2}\ZZ$, or
  \item $[\delta\hat\omega_C] \in (\pi ; \pi_Q) \cF_0$ and $a_i' \in \frac{1}{4} + \frac{1}{2}\ZZ$.
\end{itemize}
Let $h_i \in \{0,1,2,3\}$ be the integer such that
\[
  h_i \equiv d_0-o_i \mod 4, \quad \text{that is,} \quad a_i' \in \ZZ + \frac{-h_i + d_0}{4}.
\]
Then we may take
\[
  \delta_0 = (\pi ; \pi_Q)^{h_i - d_0}
\]
and $\delta = \delta_0\upsilon$ with $\upsilon \in \OO_R^\cross$.
We get
\[
  \xi_i \in \pi^{\frac{-d_0 + h_i}{4}} \cdot (\pi ; \pi_Q)^{\frac{-h_i + d_0}{2}} \sqrt{\upsilon} R.
\]
We take $\gamma^{\odot}_i$ so that
\[
  \gamma_i\gamma^{\odot}_i = \pi^{\frac{d_0 - h_i}{4}} \cdot \left(1 ; \pi_Q^{\frac{-d_0 + \One_{2\nmid h_i}}{2}}\right) \sqrt{\upsilon}.
\]
It is evident that $\xi_i^\odot = \gamma_i^{-1}{\gamma^\odot_i}^{-1}\xi_i$ belongs to $R$; we need that it is a primitive vector of $\OO_R$. Since $\xi_i^{(Q)}$ is a unit,
\[
  v_Q({\xi_i^\odot}^{(Q)}) = -v_Q\left((\gamma\gamma^{\odot})^Q\right) = 2\left(\frac{- d_0 + h_i}{4}\right) + \frac{d_0 - \One_{2\nmid h_i}}{2}
  = \begin{cases}
    0, & h_i = 0, 1 \\
    1, & h_i = 2, 3
  \end{cases}
\]
while, since $v_K(\xi_i^{(K)}) \geq \frac{d_0-1}{4}$ as noted above,
\[
  v_K({\xi_i^\odot}^{(K)}) \geq \frac{d_0 - 1}{4} + \frac{-d_0 + h_i}{4} = \frac{h_i - 1}{4} \geq -\frac{1}{4}.
\]
Since ${\xi_i^\odot}^{(K)} \in K$, we can strengthen the bound to $0$. So $\xi^\odot \in \OO_R$ and indeed is primitive on account of the $Q$-component. Pick
\[
  p^\odot = \frac{d_0 - h_i}{2} \equiv -2a_i' \equiv m_{11} \mod \ZZ,
\]
thus getting
\[
  \delta^\odot = \frac{\diamondsuit {\gamma_i^\odot}^2 \gamma_i^2}{\pi^{p^\odot}} \sim \frac{(1; \pi^{d_0/2}) \cdot \pi^{\frac{d_0 - h_i}{2}}\left(1; \pi^{\frac{-d_0 + \One_{2\nmid h_i}}{2}}\right)}{\pi^{\frac{d_0 - h_i}{2}}} \sim (1; \pi_Q^{\One_{2\nmid h_i}}).
\]
So the conic has determinant $\sim 1$ or $\sim \pi$ according as $h_i$ is even or odd, that is, as $[\delta\hat\omega_C\diamondsuit]$ belongs to $\cF_0$ or $(\pi; \pi_Q)\cF_0$. Note that $\cM \nequiv 0 \mod \pi$ as a quadratic form because $\cM(1;0) \sim \lambda^\diamondsuit(1; 0) \sim 1$. So the conic is either unimodular or tiny in the two respective cases, as desired. \qedhere
\end{proof}

The advantage of making the conic's determinant associate to either $1$ or $\pi$ is that we have to solve very few isomorphism types of conics. Although we do not prove the following classification, it animates the choice of what invariants we compute:
\begin{conj}\label{conj:conic_classfn}
  Let $\cC$ be a relevant conic over the ring of integers $\OO_K$ of a local field $K$.
  \begin{enumerate}[$($a$)$]
    \item If $\cC$ is tiny, it is determined up to isomorphism by its \emph{Brauer class} $\epsilon(\cC) \in \{\pm 1\}$, the single bit telling whether $\cC(\vec{x}) = 0$ has a nonzero solution over $K$.
    \item If $\cC$ is unimodular, it is determined up to isomorphism by its Brauer class $\epsilon(\cC)$ and its \emph{squareness level} $\ell(\cC)$, the largest $\ell \in \ZZ, 0 \leq \ell \leq e/2$, such that
    \[
      \cC \equiv c\lambda^2 \mod \pi^{\min\{2\ell + 1, e\}}
    \]
    as a quadratic form, for some constant $c \in \OO_K^\cross$ and linear form $\lambda$. Moreover, all combinations of values $ (\epsilon, \ell) $ occur, except that for $e$ even, $\ell(\cC) = e/2$ implies $\epsilon(\cC) = 1$ by Proposition \ref{prop:conic_1_mod_2} below. Thus there are exactly $e + 1$ isomorphism classes of conics of determinant $1$.
  \end{enumerate}
\end{conj}

Note that if $\cC$ is unimodular and $2\ell + 1 \leq e$, then $\cC$ is congruent to a $c\lambda^2$ modulo $\pi^{2\ell}$ if and only if modulo $\pi^{2\ell + 1}$, as the $x^2$, $y^2$, $z^2$ coefficients have square ratios modulo $\pi^{2\ell + 1}$ and the cross-terms are multiples of $2$ anyway. So the squareness level carries the same amount of information as the \emph{squareness}
\begin{align*}
  \square(\cC) &= \max \left\{i : \cC \equiv c\lambda^2 \mod \pi^i \text{ as a quadratic form}\right\} \\
  &\in \left\{0,2,4, \ldots, 2\floor{\frac{e-1}{2}}\right\} \union \{e\}.
\end{align*}
(We cannot have $\square(\cC) > e$, or the determinant would vanish modulo $\pi$.)

When we use coordinates, we will generally use one of two explicit types of conics: the \emph{diagonal conic}
\[
  aX^2 + bY^2 + cZ^2 = 0
\]
and the \emph{basepoint conic}
\[
  2XZ - Y^2 + aZ^2 = 0,
\]
so called because it passes through the basepoint $[1:0:0]$ and is tangent to the line $Z = 0$ there. We begin with results concerning the diagonal conic.

\section{Diagonal conics}

\begin{lem}\label{lem:conic_diag}
Any relevant conic $\cC$ is diagonalizable, that is, there exists a basis $(\xi_1,\ldots,\xi_3)$ for the given lattice $\Xi$ such that
\[
  \cC(x_1 \xi_1 + x_2 \xi_2 + x_3 \xi_3) = a_1 x_1^2 + a_2 x_2^2 + a_3 x_3^2.
\]
\end{lem}
\begin{proof}
When $\ch k_K \neq 2$, we have that \emph{any} conic is diagonalizable by an easy Gram-Schmidt procedure (in fact more is true: see O'Meara \cite{OMeara}, 92:1). So we assume that $\ch k_K = 2$. Here a quadratic space is not diagonalizable in general, and we must use the restrictions given on $\cC$.

Write the matrix of $\cC$, with respect to any basis $(\xi_1, \xi_2, \xi_3)$, as
\[
  \cC = \begin{bmatrix}
    a & w & v \\
    w & b & u \\
    v & u & c
  \end{bmatrix}.
\]
In the case that $\cC$ is unimodular, we see from
\[
  \det \cC = abc + 2uvw - au^2 - bv^2 - cw^2 \sim 1
\]
that at least one of the diagonal entries---say $a$---is a unit. Then we can use $a$ to eliminate $w$ and $v$ (that, is, add multiples of $\xi_1$ to $\xi_2$ and $\xi_3$). Now if $b$ (or, symmetrically, $c$), is nonzero modulo $\pi$, we use it to eliminate $u$, and we are done, as we have found the requisite diagonal form. However, it is possible that
\[
  \cC \equiv \begin{bmatrix}
    a &   & \\
      &   & u \\
      & u &
  \end{bmatrix} \mod \pi
\]
for units $a$ and $u$. Rescaling $\xi_3$, we can assume that
\[
  \cC \equiv a \cdot \begin{bmatrix}
    1 &   & \\
      &   & 1 \\
      & 1 &
  \end{bmatrix} \mod \pi.
\]
If we try to keep our chosen $\xi_1$, we will fail because the unimodular form
\[
\cC_2 = \begin{bmatrix}
  0 & 1 \\
  1 & 0
\end{bmatrix}
\]
on $\OO_K^2$ is not diagonalizable. However, we can use the identity
\[
  \begin{bmatrix}
    1 & 1 & 1 \\
    1 & 1 &   \\
    1 &   & 1
  \end{bmatrix}
  \begin{bmatrix}
    1 &   &   \\
      &   & 1 \\
      & 1
  \end{bmatrix}
  \begin{bmatrix}
    1 & 1 & 1 \\
    1 & 1 &   \\
    1 &   & 1
  \end{bmatrix}
  \equiv
  \begin{bmatrix}
    1 &   &   \\
      & 1 &   \\
      &   & 1
  \end{bmatrix}
  \mod 2
\]
to change to a basis in which
\[
  \cC \equiv a \cdot
  \begin{bmatrix}
    1 &   &   \\
      & 1 &   \\
      &   & 1
  \end{bmatrix} \mod \pi.
\]
Then the diagonalization proceeds without a hitch.

If $\cC$ is tiny, we proceed similarly. Taking $a$ a unit (since we are given $\cC(\xi_1) \sim 1$ for some $\xi_1$), we can eliminate $w$ and $v$. Then since
\[
  \det \cC = a(bc - u^2) \sim \pi,
\]
we must have at least one of $b$ and $c$ a unit, as otherwise $\det \cC$ would be either a unit (if $u$ is a unit) or a multiple of $\pi^2$ (if $\pi \mid u$). So we can eliminate $u$ and again get the desired diagonalization.
\end{proof}

The following lifting lemma for solutions modulo $4\pi$ will be essential for us.
\begin{lem}\label{lem:conic_lift}
Let $\cC$ be a diagonalized conic on an $\OO_K$-lattice $\Xi$, and assume that $v(\det \cC) \leq 1$. Let $\xi \in \Xi$ be a primitive vector with
\[
  \cC(\xi) \equiv 0 \mod \pi^m, \quad m > 2e.
\]
Then there exists a $\xi' \in \Xi$ such that
\[
  \xi' \equiv \xi \mod \pi^{m-e} \textand \cC(\xi') = 0.
\]
\end{lem}
\begin{proof}
We may write the conic in diagonal form
\[
  \cC(x_1\xi_1 + x_2\xi_2 + x_3\xi_3) = a_1 x_1^2 + a_2 x_2^2 + a_3 x_3^2.
\]
Let $\xi = x_1\xi_1 + x_2\xi_2 + x_3\xi_3$. Since $\xi$ is primitive, not all the $x_i$ are zero modulo $\pi$. We claim that there is an $i$ with
\begin{equation}\label{eq:x_move_coord}
  \pi \nmid a_i \textand \pi \nmid x_i.
\end{equation}
If not, then $\det \cC \sim \pi$, and without loss of generality, $a_1,a_2$ are units while $a_3 \sim \pi$; and $x_3$ is a unit while $x_1,x_2$ are multiples of $\pi$. Summing, we find that
\[
  \cC(\xi) \equiv a_3 x_3^2 \nequiv 0 \mod \pi^2,
\]
a contradiction.

Choose $i$ satisfying \eqref{eq:x_move_coord}. We will construct $\xi'$ by changing only the $x_i$ coordinate of $\xi$ to a different value $x_i'$. The desired condition $\cC(\xi') = 0$ takes the form
\[
  x_i'^2 = y
\]
for some $y \equiv x_i^2$ modulo $\pi^{m}$. Since $m > 2e$, we have that $y$ is also a square and, indeed, has a (unique) square root $x_i'$ satisfying $x_i' \equiv x_i$ mod $\pi^{m-e}$. This constructs the desired $\xi'$.
\end{proof}

Here are two easy corollaries.
\begin{prop}\label{prop:conic_1_mod_2}
If $e$ is even and $\cC$ is a unimodular conic of squareness $e$ (the maximal possible value), then $\epsilon(\cC) = 1$, that is, $\cC$ has a rational point.
\end{prop}
\begin{proof}
We may assume that $\cC$ is diagonal:
\[
  \cC(x_1 \xi_1 + x_2 \xi_2 + x_3 \xi_3) = a_1 x_1^2 + a_2 x_2^2 + a_3 x_3^2.
\]
Then $\epsilon(\cC)$ is a Hilbert symbol,
\[
  \epsilon(\cC) = \left\langle \frac{-a_2}{a_1}, \frac{-a_3}{a_1} \right\rangle.
\]
Both arguments are squares of units modulo $2$. But since $e$ is even, they are actually squares modulo $2\pi$, so the Hilbert symbol is $1$ by \eqref{eq:Hilb_prod_size}.
\end{proof}

\begin{prop}\label{prop:conic_perturb}
Let $\cC$ and $\cC'$ be conics of determinant $1$ and squareness level $\ell$. Suppose that the associated bilinear forms of $\cC$ and $\cC'$ are congruent modulo $\pi^{2e - 2\ell}$. Then $\epsilon(\cC) = \epsilon(\cC')$.
\end{prop}
\begin{proof}
We may assume that $\cC$ is diagonal:
\[
\cC(x_1 \xi_1 + x_2 \xi_2 + x_3 \xi_3) = a_1 x_1^2 + a_2 x_2^2 + a_3 x_3^2.
\]
Although $\cC'$ need not be diagonal with respect to the same basis, the orthogonalization procedure furnished by the proof of Lemma \ref{lem:conic_diag} yields a basis $(\xi_1', \xi_2', \xi_3')$ with $\xi_j' \equiv \xi_j \mod \pi^{2e - 2\ell}$ such that $\cC'$ is diagonal with respect to it,
\[
\cC'(x_1 \xi'_1 + x_2 \xi'_2 + x_3 \xi'_3) = a'_1 x_1^2 + a'_2 x_2^2 + a'_3 x_3^2
\]
with $a_j \equiv a'_j \mod \pi^{2e - 2\ell}$. Now compare
\[
  \epsilon(\cC) = \left\langle \frac{-a_2}{a_1}, \frac{-a_3}{a_1} \right\rangle \textand \epsilon(\cC') = \left\langle \frac{-a_2'}{a_1'}, \frac{-a_3'}{a_1'} \right\rangle.
\]
The arguments to the Hilbert symbol belong to $\cF_{\ell}$, the $\ell$th level space of $K^\cross/(K^\cross)^2$. So perturbing them by multipliers in $\cF_{e - \ell}$ does not change the Hilbert symbol by \eqref{eq:Hilb_prod_size}.
\end{proof}

\section{The solution volume of the conic}

We now use the basepoint form to determine volumes of conics.

\begin{lem}[\textbf{Igusa zeta function of a unimodular conic}]\label{lem:conic_1}
  Let $\cC$ be a conic of determinant $1$ on a $3$-dimensional vector space $\Xi$. Suppose that $\cC$ has Brauer class $\epsilon(\cC) = 1$, that is, it admits a basepoint $\xi_0$ such that $\cC(\xi_0) = 0$.

  Let $U_{m^\odot, n^\odot}$ be the volume of $\xi \in \PP(\Xi)$ (counting the whole $\PP(\Xi)$ to have volume $1 + q^{-1} + q^{-2}$) such that
  \begin{align}
    \xi &\equiv \xi_0 \mod \pi^{n^\odot} \label{eq:conic_1_n} \\
    \cC(\xi) &\equiv 0 \mod \pi^{m^\odot}. \label{eq:conic_1_m}
  \end{align}
  (Note that these $m^\odot$ and $n^\odot$ correspond to the $m_{11}^\odot$ and $n^\odot$ of Lemmas \ref{lem:tfm_conic} and \ref{lem:N11}.)

  Then for $m^\odot$ and $n^\odot$ integers with $m^\odot > 2e$ and $m^\odot \geq 2n^\odot$, the volume $U_{m^\odot,n^\odot}$ depends only on $m^\odot$, $n^\odot$, and the squareness level $\ell = \ell(\cC)$. It is given by
  \[
  U_{m^\odot,n^\odot} = q^{e-m^\odot-n^\odot}, \quad n^\odot \geq e, \quad n^\odot > 0
  \]
  and the recurrence
  \begin{align*}
    \frac{U_{m^\odot,n^\odot}}{U_{m^\odot,n^\odot+1}} &= \begin{cases}
      2 & n^\odot = e - 2\ell - 1 \geq 0 \\
      q & n^\odot \equiv e \mod 2, n^\odot > e - 2\ell - 1, n^\odot > 0 \\
      q+1 & n^\odot = 0, \ell = \frac{e}{2} \\
      1 & otherwise.
    \end{cases}
  \end{align*}
  Explicitly,
  \[
  {U_{m^\odot,n^\odot} = }
  \left\{\begin{tabular}{lll}
    $ q^{e-m^\odot-n^\odot} $ & $ n^\odot \geq e,\ n^\odot>0 $ & (a \emph{black conic}) \\
    $ q^{-m^\odot + \floor{\frac{e - n^\odot}{2}}} $ & $ e \geq n^\odot > e - 2\ell - 1, n^\odot > 0 $ & (a \emph{blue conic}) \\
    $ 2 q^{-m^\odot + \ell} $ & $  0 \leq n^\odot \leq e - 2\ell - 1 $ & (a \emph{green conic}) \\
    $ \left( 1 + \dfrac{1}{q}\right) q^{-m^\odot + e/2}, $ & $ n^\odot = 0, \ell = \dfrac{e}{2}, \text{$ e $ even}  $ & (a \emph{beige conic})
  \end{tabular}
  \right.
  \]
\end{lem}
\begin{rem}
  The colors are chosen to correspond to certain zones of discrete data in which these conics occur. Note that the colors become brighter as $n^\odot$ decreases.
\end{rem}

\begin{proof}
  For the black-conic case, we diagonalize the conic to
  \[
  \cC(X,Y,Z) = aX^2 + bY^2 + cZ^2.
  \]
  Let the basepoint be $\xi_0 = [X_0 : Y_0 : Z_0]$. Note that two coordinates of $\xi_0$, say $X_0$ and $Y_0$, are nonzero modulo $\pi$. We may scale so that $Y_0 = 1$ and so that all solutions we seek have $Y = 1$. This eliminates the issue of scaling ambiguity.

  The $n^\odot$-pixel of $[X:1:Z]$ satisfying \eqref{eq:conic_1_n} has volume $q^{-2n^\odot}$. For fixed $Z$, with $Z \equiv Z_0 \mod \pi^{n^\odot}$ the condition \eqref{eq:conic_1_m} simplifies to $X^2 \equiv u \mod \pi^{m^\odot}$, where $u$ is a unit with $u \equiv X_0^2 \mod \pi^{n^\odot + e}$. Hence its solutions form a congruence class mod $\pi^{m^\odot-e}$, and overall, the solution volume is $q^{n^\odot} \cdot q^{e-m^\odot} = q^{e-m^\odot-n^\odot}.$

  We use this as the base case to prove the recursive formula (and hence also the explicit formula) for $U_{m^\odot,n^\odot}$ by downward induction on $n^\odot$. Our aim is to determine the number $r$ of $(n^\odot+1)$-pixels within the $n^\odot$-pixel of $\xi_0$ that contain a solution to $\cC(\xi) = 0$, or equivalently, to $\cC(\xi) \equiv 0 \mod \pi^{m^\odot}$. Then, by induction, there is a volume $U_{m^\odot,n^\odot+1}$ of solutions in each of those, so $U_{m^\odot,n^\odot} = rU_{m^\odot,n^\odot+1}$ as desired. It remains to compute
  \[
    r = \frac{U_{m^\odot,n^\odot}}{U_{m^\odot,n^\odot+1}}.
  \]

  We now abandon the diagonalized form and choose coordinates such that the basepoint is $\xi_0 = [1:0:0]$ and the tangent line there is $Z = 0$. Then the conic has the form
  \[
  \cC(X,Y,Z) = 2uXY + 2vXZ - cY^2 + 2bYZ + aZ^2.
  \]
  Note that $\pi$ does not divide both $u$ and $v$, for then the conic's determinant
  \[
  \begin{vmatrix}
    0 & u & v \\
    u & -c & b \\
    v & b & a
  \end{vmatrix} = 2buv + cv^2 - au^2
  \]
  would be divisible by $\pi^2$. So, by symmetry, we may assume $v \sim 1$. We scale the conic so that $v = 1$, and then the transformation $Z \mapsto Z - uY$ makes $u = 0$. Also, the transformation $X \mapsto X - bY$ makes $b = 0$. Now $c = 1$ to make the determinant $1$. Thus the conic takes the basepoint form
  \[
  \cC(X,Y,Z) = 2XZ - Y^2 + aZ^2,
  \]
  where $a \in \OO_K$ is the only undetermined coefficient. By definition of squareness level, we know that $a$ is a square modulo $\pi^{\min\{2\ell + 1, e\}}$, but not modulo $\pi^{2\ell + 3}$ if $\ell < \floor{e/2}$. For any $a' \in \OO_K$, the transformation $Y \mapsto Y + a'Z$ can be used to increment $a$ by the square $a'^2$, followed by another $X \mapsto X - b'Y$ to remove the $YZ$ term. Picking $a'$ appropriately, we can arrange so that $v(a)$ \emph{reveals the squareness:} either
  \begin{itemize}
    \item $\ell < \floor{e/2}$ and $v(a) = 2\ell + 1$, or
    \item $\ell = \floor{e/2}$ and $2|a$. In this case, indeed, the transformation $X \mapsto X - (a/2)Z$ makes $a = 0$. (So we have proved one case of Conjecture \ref{conj:conic_classfn}: for $\ell = \floor{e/2}$ and $\epsilon = 1$, the conic takes the fixed form $2XZ - Y^2$.)
  \end{itemize}

  To parametrize $ \cC(\OO_K) $, we use the age-old trick of \emph{stereographic projection,} that is, drawing lines of varying slope through the known basepoint $[1:0:0]$. An easy calculation shows that the second intersection of the line $sZ - tY = 0$ with the conic $\cC$ is $[s^2 - at^2 : 2st : 2t^2]$, yielding an isomorphism
  \begin{align*}
    \PP^1(K) &\isom \cC(K) \\
    [s:t] &\mapsto [s^2 - at^2 : 2st : 2t^2].
  \end{align*}
  If $ [s:t] $ is in lowest terms over $ \OO_K $, the image $ [s^2 - at^2 : 2st : 2t^2] $ need not be in lowest terms over $ \OO_K $, but will have cancellation by $\pi^j$, where $ j = \min\{v_K(s^2 - at^2), e + v_K(st), e + v_K(t^2)\} $. Note that $ j \leq e $ because $ s $ and $ t $ are coprime, so
  \[
  j = \min\{v_K(s^2 - at^2), e\}.
  \] Note also that the resulting point
  \[
  [X : Y : Z] = \left[ \frac{s^2 - at^2}{\pi^j} : \frac{2st}{\pi^j} : \frac{2t^2}{\pi^j} \right]
  \]
  lies in the same $ (e - j) $-pixel as $ [1:0:0] $, but not in the same $ (e - j + 1) $-pixel if $ j > 0 $. Hence the points we are interested in, namely in the $ n^\odot $-pixel of the basepoint but outside the $ (n^\odot + 1) $-pixel, correspond exactly to values of $ [s : t] $ for which $ j = e - n^\odot $. That is, the valuation $ v(s^2 - at^2) $ must be exactly $ e-n^\odot $ (if $ n^\odot > 0 $) or at least $ e $ (if $ n^\odot = 0 $). The cases of $n^\odot$ now become cases for $j$.
  \begin{itemize}
    \item Observe that when $n^\odot < e - 2\ell - 1$, there are no solutions. Also, when $n^\odot > e - 2\ell - 1$ is of the same parity as $e$, there are no solutions, because
    \[
    s^2 - at^2 \equiv s^2 \mod \pi^{\min\{2\ell + 1, e\}}
    \]
    has even valuation if nonzero mod $\pi^{\min\{2\ell + 1, e\}}$. The volumes $ {U_{m^\odot,n^\odot}} $ and $ {U_{m^\odot, n^\odot+1}} $ are thus equal in these cases, as claimed.
    \item Suppose that $ n^\odot = e - 2\ell - 1 $. If $ n^\odot > 0 $, we seek the $[s:t]$ such that $s^2 - a t^2$ attains its maximal valuation $2\ell+1$: this happens when $\pi^{\ell+1} | s$. If $n^\odot = 0$, we seek $\pi^{e}|s^2$, which is still equivalent to $\pi^{\ell+1} | s$. Hence we are looking at the $[s:t]$ with $ t = 1 $, $ s = \pi^{\ell + 1} s' $, where $ s' \in \OO_K$. Then the resulting point on $ \cC(\OO_K) $ is
    \begin{align*}
      [X : Y : Z] &= \left[ \frac{s^2 - at^2}{\pi^{2\ell+1}} : \frac{2st}{\pi^{2\ell+1}} : \frac{2t^2}{\pi^{2\ell+1}} \right] \\
      &= \left[ \frac{\pi^{2\ell+2}s'^2 - a}{\pi^{2\ell+1}} : \frac{2\pi^{\ell+1}s'}{\pi^{2\ell+1}} : \frac{2t^2}{\pi^{2\ell+1}} \right] \\
      &= \left[ \pi s'^2 - \frac{a}{\pi^{2\ell + 1}} : \frac{2s'}{\pi^{\ell}} : \frac{2}{\pi^{2\ell + 1}} \right] \\
      &\equiv \left[ \pi s'^2 - \frac{a}{\pi^{2\ell + 1}} : 0 : \frac{2}{\pi^{2\ell + 1}} \right] \mod \pi^{n^\odot+1} = \pi^{e-2\ell}.
    \end{align*}
    We claim that this is actually the same point modulo $ \pi^{n^\odot+1} $ regardless of $ s' $, that is, the $ \pi s'^2 $ term contributes nothing. If $ 2\ell + 1 = e $, this is clear because $ n^\odot = 0 $. Otherwise, $ a' = a/\pi^{2\ell + 1} $ is a unit, and if we multiply all three coordinates by
    \[
    \frac{a'}{a' - \pi s'^2} \equiv 1 \mod \pi,
    \]
    the last two coordinates do not change mod $ \pi^{n^\odot+1} $ because they are $ 0 $ mod $ \pi^{n^\odot} $. Thus all points $ [X : Y : Z] $ obtained lie in a single $ (n^\odot + 1) $-pixel, and hence the ratio $ {U_{m^\odot,n^\odot}} / {U_{m^\odot, n^\odot+1}} $ is $ 2 $.
    \item Suppose that $ n^\odot \equiv e \mod 2$, $n^\odot > e - 2\ell - 1 $, and $ n^\odot > 0 $. Then $ e - n^\odot = j = 2j' $ is even, with $ j' \leq \ell. $ The pairs $ [s:t] $ yielding $ s^2 - at^2 \sim \pi^j $ are exactly those with $ s \sim \pi^{j'} $. Write $ t = 1 $, $ s = \pi^{j'} s' $, where $ s' \in \OO_K^\cross $. Then the resulting point on $ \cC(\OO_K) $ is
    \begin{align*}
      [X : Y : Z] &= \left[ \frac{s^2 - at^2}{\pi^j} : \frac{2st}{\pi^j} : \frac{2t^2}{\pi^j} \right] \\
      &= \left[ - \frac{a}{\pi^{2j'}} + s'^2 : \frac{2s'}{\pi^{j'}} : \frac{2}{\pi^{2j'}} \right] \\
      &\equiv \left[ s'^2 : 0 : \frac{2}{\pi^{2j'}} \right] \quad \mod \pi^{n^\odot + 1},
    \end{align*}
    where at the last step we multiplied all three coordinates by the unit
    \[
    \frac{s'^2}{s'^2 - \frac{a}{\pi^{2j'}}} \equiv 1 \mod \pi.
    \]
    We get $ q - 1 $ different $ (n^\odot + 1) $-pixels, one for each value of $ s' $ mod $ \pi $. Hence the ratio $ {U_{m^\odot,n^\odot}} / {U_{m^\odot, n^\odot+1}} $ is $ q $.
    \item Finally, suppose that $ n^\odot = 0 $ and $ \ell = \frac{e}{2} $. Write $ t = 1 $, $ s = \pi^{e/2} s' $, where $ s' \in \OO_K $. We get
    \begin{align*}
      [X : Y : Z] &= \left[ \frac{s^2}{\pi^e} : \frac{2st}{\pi^e} : \frac{2t^2}{\pi^e} \right] \\
      &= \left[ s'^2 : \frac{2s'}{\pi^{e/2}} : \frac{2}{\pi^{e}} \right] \\
      &\equiv \left[ s'^2 : 0 : \frac{2}{\pi^{e}} \right] \quad \mod \pi = \pi^{n^\odot+1}.
    \end{align*}
    We get $ q $ different $ 1 $-pixels, one for each value of $ s' $ mod $ \pi $. Hence the ratio $ {U_{m^\odot,n^\odot}} / {U_{m^\odot, n^\odot+1}} $ is $ q + 1 $. \qedhere
  \end{itemize}

\end{proof}

For determinant $\pi$, we use the same method; the number of cases is less:
\begin{lem}[\textbf{Igusa zeta function of a tiny conic}]\label{lem:conic_pi}
  Let $\cC$ be an integer-matrix conic over $\OO_K$ of determinant $\pi$. Suppose that $\cC$ has Brauer class $\epsilon(\cC) = 1$, that is, it admits a basepoint $\xi_0$ such that $\cC(\xi_0) = 0$.

  Let $U_{m^\odot,n^\odot}$ be the volume of $\xi \in \PP(\Xi)$ such that
  \begin{align}
    \xi &\equiv \xi_0 \mod \pi^{n^\odot} \label{eq:conic_pi_n} \\
    \cC(\xi) &\equiv 0 \mod \pi^{m^\odot}. \label{eq:conic_pi_m}
  \end{align}

  Then for $m^\odot > 2e$ and $m^\odot \geq 2n^\odot$, the volume $U_{m^\odot,n^\odot}$ depends only on $m^\odot$ and $n^\odot$:
  \[
  {U_{m^\odot,n^\odot} = }
  \left\{\begin{tabular}{lll}
    $ q^{e-m^\odot-n^\odot} $ & $ n^\odot > e $ & (a \emph{tiny black conic}) \\
    $ 2 q^{-m^\odot} $ & $ 0 \leq n^\odot \leq e $ & (a \emph{tiny green conic})
  \end{tabular}
  \right.
  \]
\end{lem}
\begin{proof}
  Diagonalize the conic to the form
  \[
  \cC(X,Y,Z) = a X^2 + b Y^2 + c \pi Z^2 = 0,
  \]
  where $a,b,c \in \OO_K^\cross$. Observe that $X_0$ and $Y_0$ are units, and scale so that $Y = Y_0 = 1$.

  If $n^\odot > e$, then for each $Z \equiv Z_0 \mod \pi^{n^\odot}$, the condition $\cC(X,Y,Z) \equiv 0 \mod \pi^{m^\odot}$ simplifies to $X^2 \equiv u \mod \pi^{m^\odot}$, where $u \equiv X_0^2$ mod $\pi^{n^\odot+e+1}$, and hence the square roots $X$ with $X \equiv X_0$ mod $\pi^{n^\odot}$ form a single congruence class mod $\pi^{m^\odot-e}$. So the volume is $q^{e-m^\odot-n^\odot}$.

  If $n^\odot = e$, then we use the same method, but now the equation $X^2 \equiv u \mod \pi^{m^\odot}$, where $u \equiv X_0^2$ mod $\pi^{2e+1}$, has as solution set two classes mod $\pi^{m^\odot-e}$, each the negative of the other.

  We claim that these are all the solutions mod $\pi^{m^\odot}$; that is, that the whole conic lies within an $e$-pixel. Suppose there is such an $[X : 1 : Z] \nequiv [X_0 : 1 : Z_0] \mod 2$, and let $v_K(X - X_0) = i$, $v_K(Z - Z_0) = j$. Then
  \begin{align*}
    a (X^2 - X_0^2) = \pi c (Z^2 - Z_0^2).
  \end{align*}
  If $i < e$, then the left side has even valuation $2i$ which cannot be matched by the right side. If $j < e$, then the right side has odd valuation $2j+1$ which cannot be matched by the left side. This completes the proof.
\end{proof}
\begin{rem}
A proof using the basepoint form of the conic is also possible.
\end{rem}

\section{The Brauer class}

\begin{normalsize}
Over $ K $, there are just two isomorphism types of conic, one with points and one without. In this section we understand the Brauer class $\epsilon(\cA)$ of conics of the form
\[
\cA_\alpha : \quad \tr (\alpha \xi^2) = 0.
\]
\begin{lem} \label{lem:H_form}
  If $ \alpha \in R^{N=1} $, define
  \[
  \epsilon(\alpha) = \epsilon(\cA_\alpha) = \begin{cases}
  1 & \text{if $ \tr(\alpha \xi^2) = 0 $ for some nonzero $ \xi \in R $} \\
  -1 & \text{otherwise.}
  \end{cases}
  \]
  Then the map of $ \FF_2 $-vector spaces
  \begin{align*}
  H^1 &\to \mu_2 \\\
  \alpha &\mapsto \epsilon(\alpha) / \epsilon(1)
  \end{align*}
  is a nondegenerate quadratic form whose associated bilinear form is none other than the Tate pairing on $ H^1 $. That is,
  \[
  \epsilon(\alpha\beta) = \epsilon(1) \cdot \epsilon(\alpha) \cdot \epsilon(\beta) \cdot \left\langle\alpha,\beta\right\rangle.
  \]
\end{lem}

\begin{rem}
  $\epsilon$ comes up, in a related context, in the work of Bhargava and Gross (\cite{AIT}, \textsection 7.2), where it is stated to be a quadratic form, at least in the tamely ramified case.
\end{rem}

For the reader's interest, we also include the following fact.
\begin{lem}\label{lem:H_extra}
    $ \epsilon(\omega) = 1 $ if and only if the associated quartic algebra $ L_\omega $ has an element $ \xi \neq 0 $ whose characteristic polynomial has the doubly depressed form
    \[
    \xi^4 - a\xi - b = 0.
    \]
    (Note that $ L_\omega = L_{\omega_0\delta} $ is in general distinct from the algebra $ L_\delta $ in which we have been trying to find orders.)
\end{lem}
\begin{proof}
  Note that $ \xi \in L_\omega $ has a doubly depressed characteristic polynomial if and only if
  \[
  \tr \xi = \tr \xi^2 = 0.
  \]
  Now the traceless elements $ \xi \in L_\omega $ have the form
  \[
  \kappa(x) = \tr_{RL_\omega/L_\omega} (x \sqrt{\omega}).
  \]
  For such $ \xi $, we compute
  \[
  \tr(\xi^2) = \tr_{R/K} (\omega x^2),
  \]
  as desired.
\end{proof}

\begin{rem}
  Over fields of characteristic not $ 2 $, a quadratic form is uniquely determined by its associated bilinear form. However, over $ \FF_2 $, the local Hilbert pairing $ \left\langle\bullet, \bullet\right\rangle_\epsilon $ on $ H^1 $ lifts to $ \size{(H^1)^\vee} = \size{H^1} $ quadratic forms, thanks to the ambiguity by adding a linear functional. It is not hard to show, as a corollary to the lemma, that these quadratic forms are exactly
  \[
  \alpha \mapsto \frac{\epsilon(\omega\alpha)}{\epsilon(\omega)},
  \]
  for each $ \omega \in H^1 $.
\end{rem}

To prove Lemma \ref{lem:H_form}, we must recall some facts about the Grothendieck-Witt ring of a local field. The nondegenerate quadratic forms over a field $ K $ ($ \ch K \neq 2 $), up to isomorphism, form a semiring under the operations of orthogonal direct sum $ \perp $ and tensor product $ \tensor $; it is cancellative (the so-called \emph{Witt cancellation theorem}), and the ring obtained by adjoining formal additive inverses is called the \emph{Grothendieck-Witt ring} $ GW(K) $ of $ K $. In the case that $ K $ is a local field, an element of $GW(K)$ is determined (see O'Meara \cite{OMeara}, Theorem 63:20) by three invariants: its dimension $ n \in \ZZ $, its determinant $ D \in K^\cross /(K^\cross)^2 $, and one other bit of information, the \emph{Hasse symbol} $ S \in \{ \pm 1\} $. For a diagonal form $ f = a_1 x_1^2 + \cdots + a_n x_n^2 $, the determinant is given by $ D(f) = a_1 \cdots a_n $ and the Hasse symbol by
\[
S(f) = \prod_{i<j} \left\langle a_i, a_j\right\rangle_K.
\]
(The determinant is also called the \emph{discriminant}; our choice of terminology is influenced partly by the clash in sign with the discriminant of a binary quadratic form.)

Every combination of dimension, determinant, and Hasse symbol determines a unique element of $ GW(K) $. For dimension at least $3$, all elements are actually realized by a quadratic form; in dimension $ 3 $, the Hasse symbol $S(f)$ carries the same information as the Brauer class $\epsilon(\cC)$ of the corresponding conic, possibly with a sign reversal depending on $D(f)$. The structure of $ GW(K) $ in terms of these invariants is easily computed by reducing to the case of diagonal forms; the formulas are here recorded, as they will be useful to us.
\begin{align}
  \dim (f \perp g) &= \dim f + \dim g \label{eq:dim_perp}\\
  \det (f \perp g) &= \det f \det g \\
  S(f \perp g) &= S(f) S(g) \left\langle\det f, \det g\right\rangle_K \label{eq:Witt_add} \\
  \dim (f \tensor g) &= \dim f \dim g \\
  \det (f \tensor g) &= (\det f)^{\dim g} (\det g)^{\dim f} \\
  S(f \tensor g) &= S(f)^{\dim g} S(g)^{\dim f} \left\langle\det f, -1\right\rangle_K^{\binom{\dim g}{2}} \left\langle\det g, -1\right\rangle_K^{\binom{\dim f}{2}} \cross \label{eq:Hasse_tensor} \\
  &\quad \cross \left\langle\det f, \det g\right\rangle_K^{\dim f \dim g - 1}. \nonumber
\end{align}
We denote by $ \Ell_K $ the unique class in $ GW(K) $ of dimension $ 0 $, determinant $ 1 $, and Hasse symbol $ -1 $ (the ``elliptic class''). Note that if $ f $ and $ g $ are nonisomorphic quadratic forms over $ K $ with the same dimension and determinant, then $ [f] = [g] + \Ell_K $ in $ GW(K) $.

Suppose that $ L/K $ is a field extension. If $ q : W \to L $ is a quadratic form over $ L $, we can view $ W $ as a $ K $-vector space and postcompose with the trace $ \tr_{L/K} $ to get a quadratic form $ \tr_{L/K} q $. Since $ \tr_{L/K} $ respects orthogonal direct sums, it induces a group homomorphism (though not a ring homomorphism) from $ GW(L) $ to $ GW(K) $ (the opposite direction to the more familiar extension-of-scalars morphism). We easily compute that
\[
\dim \tr_{L/K} q = [L:K] \cdot \dim q \textand
\det \tr_{L/K} q = (\disc(L/K))^{\dim q} \cdot N_{L/K}(\det q).
\]
We wish to understand how $ S(\tr_{L/K} q) $ behaves. The following is the most important result needed.
\begin{lem}
  For any extension $ L/K $ of nonarchimedean local fields not of characteristic $ 2 $,
  \[
  \tr_{L/K} \Ell_{L} = \Ell_K.
  \]
  In other words, $ \tr_{L/K} : GW(L) \to GW(K) $ preserves the Hasse symbol on classes of dimension $ 0 $ and determinant $ 1 $.
\end{lem}
\begin{proof}
  We may assume that $ L/K $ is a \emph{primitive} extension, that is, has no nontrivial intermediate extensions, since $ \tr_{E/K} \circ \tr_{L/E} = \tr_{L/K} $ for a tower $ L/E/K $.

  Since $ \tr_{L/K} \Ell_{L} $ is of dimension $ 0 $ and determinant $ 1 $, the only other possibility is that $ \tr_{L/K} \Ell_{L} = 0 $. We prove that this cannot hold.

  Let $ a \in K^\cross $ be an element that does not become a square in $ L $. Such an $ a $ exists because $ L/K $ is primitive; if not, then $ L/K $ would contain both an unramified and a ramified quadratic extension. Then choose $ \theta \in L $ such that $ \left\langle a, \theta\right\rangle_{L} = -1 $. We also have $ \left\langle a, N_{L/K}(\theta)\right\rangle_K = -1 $ by the standard relation $ \left\langle a, \theta\right\rangle_{L} = \left\langle a, N_{L/K}(\theta)\right\rangle_K $ ($ a \in K $, $ \theta \in L $). Consider the following quadratic forms over $ L $:
  \[
  f(x,y) = x^2 - \theta y^2, \quad g(x,y) = a x^2 - a \theta y^2.
  \]
  In other words, $ f = q_1 \perp q_{-\theta} $ and $ g = q_a \perp q_{-a\theta} $ where $ q_a(x) = a x^2 $. Both $ f $ and $ g $ have dimension $ 2 $ and discriminant $ -\theta $, but their Hasse symbols are $ 1 $ and $ \left\langle a, -a \theta\right\rangle = -1 $, respectively. Hence $ [f] = [g] + \Ell_{L} $. But $ g = af $, so $ \tr_{L/K}(g) = a\tr_{L/K}(f) = q_a \tensor \tr_{L/K}(f) $. So by the formula \eqref{eq:Hasse_tensor} for the Hasse symbol of a tensor product,
  \begin{align*}
    S(\tr_{L/K} g) &= S(q_a \tensor \tr_{L/K} f) \\
    &= S(q_a)^{2[L:K]} \cdot S(\tr_{L/K} f) \cdot \left\langle\det q_a, -1\right\rangle^{\binom{2[L:K]}{2}} \cdot \left\langle\det \tr_{L/K} f, -1\right\rangle^{\binom{1}{2}} \cross \\
    & \quad \cross \left\langle\det q_a, \det \tr_{L/K} f\right\rangle^{2 [L:K] - 1} \\
    &= S(\tr_{L/K} f) \cdot \left\langle a, -1\right\rangle^{[L:K]} \cdot \left\langle a, (\disc (L/K))^{2} \cdot N_{L/K}(\det f)\right\rangle^{2 [L:K] - 1} \\
    &= S(\tr_{L/K} f) \cdot \left\langle a, -1\right\rangle^{[L:K]}\cdot \left\langle a, N_{L/K}(-\theta)\right\rangle \\
    &= S(\tr_{L/K} f) \cdot \left\langle a, N_{L/K}(\theta)\right\rangle \\
    &= -S(\tr_{L/K} f),
  \end{align*}
  So $ [\tr_{L/K} f] = [\tr_{L/K} g] + \Ell_{K}, $ yielding the desired conclusion.
\end{proof}

\begin{proof}[Proof of Lemma \ref{lem:H_form}]
  We now relate the $ \epsilon $ of the lemma to the Hasse symbol. Denote by $ q_{R,\alpha} $ the quadratic form
  \[
  q_{R,\alpha}(\xi) = \alpha \xi^2
  \]
  over $ R $. Then $ \epsilon(\alpha) = 1 $ if and only if the form $ \tr_{R/K}(q_{R,\alpha}) $ is isotropic, where $ \tr_{R/K} $ is to be interpreted in the obvious way if $ R $ is not a field. Given $ \alpha,\beta \in R^\cross $ of norm $ 1 $, the forms
  \[
  f = q_{R,1} \perp q_{R,\alpha\beta} \textand
  g = q_{R,\alpha} \perp q_{R,\beta}
  \]
  have dimension $ 2 $ and determinant $ \alpha \beta $ over $ R $, and in $ GW(K), $
  \[
  \epsilon(1)\epsilon(\alpha)\epsilon(\beta)\epsilon(\alpha\beta) = S([\tr_{R/K} f] - [\tr_{R/K} g]).
  \]
  Decompose $ R = \prod_{i = 1}^r R_i $ into its field factors ($ 1 \leq r \leq 3 $), and let $ \alpha = (\alpha_1 ; \ldots ; \alpha_r) $ and $ \beta = (\beta_1; \ldots; \beta_r) $. Then
  \begin{align*}
    [\tr_{R/K} f] - [\tr_{R/K} g] = \sum_{i=1}^r \tr_{R_i/K} ([q_{R_i,\alpha_i\beta_i}] - [q_{R_i,\alpha_i}] - [q_{R_i,\beta_i}] + [q_{R_i, 1}]).
  \end{align*}
  Since the invariants of each $q_{R_i, \gamma}$ are known, the invariants of the class in parentheses can be computed by repeated application of \eqref{eq:dim_perp}--\eqref{eq:Witt_add}. We find that it has dimension $ 0 $, determinant $ 1 $ and Hasse symbol $ \left\langle\alpha_{i}, \beta_{i}\right\rangle_{R_i} $. By the preceding lemma, its trace has the same invariants. Hence the whole sum has Hasse symbol
  \[
  \prod_i \left\langle\alpha_{i}, \beta_{i}\right\rangle_{R_i} = \left\langle\alpha, \beta\right\rangle_R,
  \]
  as desired.
\end{proof}

\end{normalsize}
\section{The squareness}
\label{sec:sqness}

\begin{normalsize}
For the cases in Lemma \ref{lem:tfm_conic} in which the transformed conic $\cM$ is unimodular, we need also to compute its squareness.

\begin{lem}\label{lem:squareness}
  Let $[\heartsuit] \in H^1$ be defined as follows:
\begin{itemize}
  \item If $R = K \cross Q$ is partially ramified, let $\pi_Q$ be any uniformizer that does \emph{not} have trace zero, and let
  \[
    [\heartsuit] = [(1; \tr \pi_Q)].
  \]
  \item In all other cases, we may take
  \[
    [\heartsuit] = 1.
  \]
\end{itemize}
Then the conic
\[
  \cM(\xi^\odot) = \lambda^\diamondsuit \left(\delta^\odot {\xi^\odot}^2\right), \quad \delta^\odot \in \OO_R^\cross
\]
has squareness
\begin{align} \label{eq:sqness}
  \square(\cM) &= \max \left\{ \ell : [\delta^\odot\heartsuit] \equiv 1 \mod \pi^{\ell}; \quad \ell = e, \text{ or } \ell < e \text{ and } \ell \text{ is odd}\right\} \\
  &= \begin{cases}
    \min\left\{2\ell(\delta\hat\omega_C\diamondsuit\heartsuit) + 1, e\right\}, & \text{$R$ unramified} \\
    \min\left\{2\floor{\dfrac{\ell(\delta\hat\omega_C\diamondsuit\heartsuit)}{2} }+ 1, e\right\}, & \text{$R$ ramified}.
  \end{cases}
\end{align}
\end{lem}
\begin{rem}
  As stated, the lemma only requires $\heartsuit$ to be defined modulo $\cF_{\floor{e/2}}$ (unramified types) resp.\ $\cF_{2\floor{e/2}}$ (ramified types). We will mostly use $\heartsuit$ in this way, but when the brown zone occurs, we will need a finer definition and will mention this.
\end{rem}
\begin{proof}[Proof for unramified splitting types]
In unramified splitting type, we first claim that going up to an \emph{unramified} extension $K'/K$ does not change either the left or the right side of \eqref{eq:sqness}. The right-hand side is less than $\floor{e/2}$ only if $[\delta \hat{\omega}_C \heartsuit]$ is represented by a generic unit $\delta^\odot = 1 + \alpha \pi^{2i+1}$, $2i + 1 < e$, and this generic unit remains generic in $R' = K' \tensor_K R$. As to the left side, we can diagonalize the conic $\cM$ to have the form
\[
  \cM(X,Y,Z) = aX^2 + bY^2 + cZ^2, \quad abc = 1.
\]
Then
\[
  \square(\cM) = \min\left\{2\ell(b/a) + 1, 2\ell(c/b) + 1, e\right\},
\]
and this remains invariant over $K'$.

Therefore, we may assume that $R \isom K \cross K \cross K$ is totally split. Let $\delta^\odot = (a; b; c)$. Then $\cM$ is diagonal and
\begin{align*}
  \square(\cM) &= \min\left\{2\ell(b/a) + 1, 2\ell(c/b) + 1, e\right\} \\
  &= \min\left\{2\ell(a) + 1, 2\ell(b) + 1, 2\ell(c) + 1, e\right\} \\
  &= \min\left\{2\ell(\delta^\odot) + 1, e\right\},
\end{align*}
as desired.
\end{proof}
\end{normalsize}
We now assume that $R$ is ramified (splitting type $(1^3)$ or $(1^21)$) and that the conic $\cM^\odot$ is unimodular. By Lemma \ref{lem:conic_1}, the solutions to $\cM^\odot(\xi^\odot) = 0$ lie in a limited number of $1$-pixels, either $1$, $2$, or $q + 1$. Now, since $\cM^\odot$ is a square modulo $\pi$, there are $(q + 1)$-many $1$-pixels of $\xi^\odot$ satisfying
\[
\cM^\odot(\xi^\odot) \equiv 0 \mod \pi.
\]
Of them, we call \emph{special} the $1$-pixel containing $(\pi_R^2)$ (in splitting type $(1^3)$) or $(0; \pi_Q)$ (in splitting type $(1^21)$); the remaining $q$-many $1$-pixels, which consist of elements of valuation $0$, we call \emph{generic}. Note that, for fixed discrete data, we only admit $\xi_i^\odot$ in one type of $1$-pixel depending on the $h_1$ of Lemma \ref{lem:tfm_conic}:
\begin{itemize}
  \item In splitting type $(1^3)$, if $h_i = 1$, then $\xi^\odot$ must be in a generic $1$-pixel.
  \item In splitting type $(1^3)$, if $h_i = -1$, then $\xi^\odot$ must be in the special $1$-pixel.
  \item In splitting type $(1^21)$, if $h_i = 0$, then $\xi^\odot$ must be in a generic $1$-pixel.
  \item In splitting type $(1^21)$, if $h_i = 2$, then $\xi^\odot$ must be in the special $1$-pixel.
\end{itemize}
This renders the following lemma very useful.
\begin{lem}[\textbf{the level parity lemma}]\label{lem:level_parity}
  In ramified splitting types, the solutions to the conic
  \[
  \cM^\odot(\xi^\odot) = 0
  \]
  (if any) are constrained by $\ell_d = \ell(\delta^\odot \heartsuit)$ as follows:
  \begin{itemize}
    \item If $\ell_d < 2\floor{e/2}$ is even, then the conic lies inside a generic $1$-pixel.
    \item If $\ell_d < 2\floor{e/2}$ is odd, then the conic lies inside the special $1$-pixel.
    \item If $e$ is odd and $\ell_d = e-1$, then the conic lies in two generic $1$-pixels.
    \item If $e$ is odd and $\ell_d \geq e$, then the conic lies in two $1$-pixels, one special, one generic.
    \item If $e$ is even and $\ell_d \geq e$, then the conic spans all $q+1$ allowable $1$-pixels.
  \end{itemize}
\end{lem}

Before proving this lemma, we prove some general results about conics of known squareness. First, for the cases when the conic lies in one $1$-pixel:
\begin{lem}\label{lem:sqness_low}
Let $\cA$ be a unimodular conic of squareness $2\ell + 1 < e$. After scaling, we can write
\begin{equation}\label{eq:sqness_low}
  \cA(\xi) \equiv \lambda(\xi)^2 + \pi^{2\ell + 1} \mu(\xi)^2 \mod \pi^{2\ell + 2},
\end{equation}
where $\lambda$ and $\mu$ are linear forms, linearly independent modulo $\pi$. Then the solution set of $\cA(\xi) = 0$ is contained in the unique $1$-pixel where
\[
  \lambda(\xi) \equiv \mu(\xi) \equiv 0 \mod \pi.
\]
Conversely, a conic of the form \eqref{eq:sqness_low} has squareness $2\ell + 1$.
\end{lem}
\begin{proof}
By definition, there is a $\lambda$ such that $\cA(\xi) \equiv \lambda(\xi)^2$ mod $\pi^{2\ell + 1}$. The quadratic form
\[
  Q(\xi) = \frac{\cA(\xi) - \lambda(\xi)^2}{\pi^{2\ell + 1}}
\]
is congruent modulo $\pi$ to a diagonal form $ax^2 + by^2 + cz^2$ (the cross terms being divisible by $\pi^{e - (2\ell + 1)} \equiv 0$ mod $\pi$). Taking square roots for $a$, $b$, and $c$ (which exist modulo $\pi$), we get the coefficients for a linear form $\mu$ satisfying \eqref{eq:sqness_low}. We cannot have $\mu$ a multiple of $\lambda$ modulo $\pi$, as then, after rescaling, $\cA \equiv \lambda(\xi)^2$ mod $\pi^{2\ell+2}$, contradicting the hypothesis about its squareness.

For a $\xi$ to satisfy $\cA(\xi) \equiv 0 \mod \pi^{2\ell+1}$, we must have $\lambda(\xi) \equiv 0 \mod \pi^{\ell+1}$. Then $\lambda(\xi)^2 \equiv 0 \mod \pi^{2\ell+2}$ so $\mu(\xi) \equiv 0 \mod \pi$. Thus $\xi$ lies in the desired $1$-pixel.

The converse result is obvious, as for any $\lambda'$,
\[
  \cA - \lambda'^2 \equiv (\lambda - \lambda')^2 + \pi^{2\ell + 1}\mu^2 \nequiv 0 \mod \pi^{2\ell+2}
\]
as the leading terms of $(\lambda - \lambda')^2$ have even valuation.
\end{proof}

Next, the case when the conic lies in two $1$-pixels:
\begin{lem}\label{lem:e_odd}
Let $\cA$ be a unimodular conic of squareness $e$, with $e$ odd. After scaling, we can write
\begin{equation}\label{eq:e_odd}
  \cA(\xi) \equiv \lambda(\xi)^2 + 2Q(\xi) \mod 2\pi,
\end{equation}
where $Q$ is a quadratic form whose restriction to the line $\lambda(\xi) = 0$ has nonzero discriminant modulo $\pi$. Then the solution set of $\cA(\xi) = 0$, if any, is contained in the two $1$-pixels where
\[
\lambda(\xi) \equiv Q(\xi) \equiv 0 \mod \pi.
\]
\end{lem}
\begin{proof}
By definition, there is a $\lambda$ such that $\cA(\xi) \equiv \lambda(\xi)^2$ mod $\pi^{e}$. Consider the quadratic form
\[
Q(\xi) = \frac{\cA(\xi) - \lambda(\xi)^2}{2}.
\]
Taking suitable coordinates, we can assume that $\lambda((x,y,z)) = x$. Then
\[
  \cA(x,y,z) \equiv x^2 + 2b y^2 + 2c z^2 + 2w x y + 2v x z + 2u y z \mod 2\pi
\]
for certain constants $b,c,u,v,w \in \OO_K$. Because $\cA$ is unimodular, $u$ must be a unit, and then on the line $x = 0$,
\[
  Q(0, y, z) \equiv b y^2 + u y z + c z^2 \mod \pi
\]
has unit discriminant $u^2 - 4 b c \equiv u^2 \mod \pi$.

For a $\xi$ to satisfy $\cA(\xi) \equiv 0 \mod 2\pi$, we must have $\lambda(\xi) \equiv 0 \mod \pi^{(e+1)/2}$. Then $\lambda(\xi)^2 \equiv 0 \mod 2\pi$, so $Q(\xi) \equiv 0 \mod \pi$. As a form on the $\PP^1(\FF_q)$ of $1$-pixels with $\lambda(\xi) = 0$, the form $Q$ has nonzero discriminant so it has roots at either zero or two $1$-pixels, as desired.
\end{proof}

\begin{proof}[Proof of Lemmas \ref{lem:squareness} and \ref{lem:level_parity}]
It now remains to work out, in each case of the conic, where the appropriate forms $\lambda$, $\mu$, $Q$ lie with respect to the special pixel.

In splitting type $(1^3)$, let $i \leq e/2$ be the claimed squareness level. Since $\ell(\delta^\odot) \geq 2i$, we have, after rescaling by a suitable square,
\[
  \delta^\odot \equiv 1 + \pi^{2i} u \pi_R + \pi^{2i + 1} v \pi_R^2 \mod \pi^{2i+2}
\]
for some $u, v \in \OO_K$. If $\ell(\delta^\odot) = 2i$, then $\pi \nmid u$; if $\ell(\delta^\odot) = 2i+1$, then $\pi \mid u$ and $\pi \nmid v$. For $\xi^\odot = x + y\pi_R + z\pi_R^2$,
\begin{alignat*}{2}
  \cM^\odot(\xi^\odot) &= \lambda^{\diamondsuit}\left(\delta^\odot {\xi^\odot}^2\right) \\
  &\equiv \lambda^{\diamondsuit}\left((1 + \pi^{2i} u \pi_R + \pi^{2i + 1} v \pi_R^2)(x + y\pi_R + z\pi_R^2)^2 \right) \\
  &\equiv y^2 + \pi^{2i + 1}(vx^2 + uz^2) + 2xz && \mod \pi^{2i+2}.
\end{alignat*}

If $2i + 1 < e$, the $2xz$ term drops out and we have an expression in the form \eqref{eq:sqness_low} with
\[
  \lambda(\xi^\odot) \equiv y, \quad \mu(\xi^\odot) \equiv v'x + u'z \mod \pi
\]
where $u'$ and $v'$ are square roots for $u$ and $v$, respectively, modulo $\pi$. If at least one of $u'$ and $v'$ is nonzero, then the conic has squareness level $i$, all its rational points lying in the pixel $[u' : 0 : v']$, which coincides with the special pixel $[0 : 0 : 1]$ exactly when $\pi \mid u'$, that is, when $\ell(\delta^\odot) = 2i + 1$.

If $2i + 1 = e$, then we have an expression in the form \eqref{eq:e_odd} with
\[
  \lambda(\xi^\odot) \equiv y, \quad Q(\xi^\odot) \equiv \frac{\pi^e}{2} v x^2 + xz + \frac{\pi^e}{2} u z^2 \mod \pi.
\]
The conic has squareness level $i$, and the special pixel $[0 : 0 : 1]$ is one of the two that contain its points exactly when $\pi \mid u$, that is, when $\ell(\delta^\odot) \geq e$.

Lastly, if $i = e/2$, then $\cM^\odot(\xi^\odot) \equiv y^2$ is a square modulo $2$, and there is no more to prove, as we know by Lemma \ref{lem:conic_1} and \ref{prop:conic_1_mod_2} that the conic is distributed among all $(q + 1)$-many special and generic $1$-pixels.

In splitting type $(1^21)$, we use the same method, but the calculation is considerably more involved. Let $c = \tr \pi_Q$. First, we scale $\heartsuit$, which is restricted by $[\heartsuit] = [(1; c)] = [(1; 1/c)]$, to be a unit. If $d_0' = v(c)$ is even, we take
\[
  \heartsuit = \left(1 ; \frac{\pi^{d_0'}}{c}\right).
\]
If $d_0'$ is odd, this transformation need not be admissible because $\pi^{d_0'}$ may not be a square in $Q$. So we take instead
\[
  \heartsuit = \left(1 ; \frac{\pi^{d_0' + 1}}{\pi_Q^2 c}\right).
\]
Let $i$ be a nonnegative integer, $2i + 1 \leq e$, and assume $\ell(\delta^\odot) \geq 2i$. Then, after scaling $\delta^\odot$ by a suitable square, we can assume that
\[
  \delta^\odot \equiv \heartsuit\left(1; 1 + \pi^{2i}\bar{\pi}_Q u\right) \mod \pi^{2i + 2}
\]
where $\bar{\pi}_Q = c - \pi_Q$ is the algebraic conjugate. Write
\[
  u \equiv \frac{c}{\pi^{d_0'}}(u'^2 + \pi v'^2) \mod \pi^2;
\]
it is easy to see that such $u'$, $v'$ exist. As in the previous splitting type, if $\ell(\delta^\odot) = 2i$, then $\pi \nmid u'$; if $\ell(\delta^\odot) = 2i+1$, then $\pi \mid u'$ and $\pi \nmid v'$.

We have that $N(\pi_Q) = \pi_Q\bar{\pi}_Q$ is a uniformizer for $K$; since $\pi$ is arbitrary, we can assume $N(\pi_Q) = \pi$.
We compute, for $d_0'$ even,
\[
  \cM^\odot(\xi^\odot) \equiv
    (x + \pi^i u' y + \pi^{d_0'/2} z)^2 + \pi^{2i+1}(v'^2 y^2 + u'^2 z^2) + 2\dfrac{\pi^{d_0'}}{c} y z \mod \pi^{2i+2},
\]
and for $d_0'$ odd,

\begin{equation*}
  \cM^\odot(\xi^\odot) \equiv (x + (\pi^{(d_0' - 1)/2} + \pi^{i}u')y)^2 + \pi^{2i + 1}(v'^2 y^2 + u'^2 z^2) + 2\dfrac{\pi^{d_0'}}{c} y z \mod \pi^{2i + 2}.
\end{equation*}
In either case we have
\begin{equation}
  \cM^\odot(\xi^\odot) \equiv \lambda(\xi^\odot)^2 + \pi^{2i + 1} (v'^2 y^2 + u'^2 z^2) + 2\frac{\pi^{d_0'}}{c} y z \mod \pi^{2i+2}
\end{equation}
for a linear form $\lambda$ with $\lambda(\xi^\odot) \equiv x$ mod $\pi$. Thus if $2i + 1 < e$, we have
\[
  \mu(\xi^\odot) \equiv v'y + u'z \mod \pi,
\]
while if $2i + 1 = e$, put $u_0 = \pi^{d_0'}/c \in \OO_K^\cross$. Then we have
\[
  Q(\xi^\odot) \equiv \frac{\pi^e}{2} v'^2 y^2 + u_0yz + \frac{\pi^e}{2} u'^2 z^2 \mod \pi.
\]
and if $i = e/2$, then
\[
  \cM^\odot(\xi^\odot) \equiv \lambda(\xi^\odot)^2 \mod 2.
\]
So in all cases, the conic has the desired squareness level $i$ and its solutions lie in the claimed types of $1$-pixels by the same argument used in splitting type $(1^3)$.
\end{proof}

We close this section with the following corollary of Lemma \ref{lem:squareness}:
\begin{lem}\label{lem:Brauer_const}
  The Brauer class $\epsilon(\delta)$ takes the same value for all $\delta$ in the coset $\hat\omega_C \heartsuit \diamondsuit \cF_{\ceil{e'/2}}$.
\end{lem}
\begin{proof}
  Let $[\delta] = [\kappa \hat\omega_C \heartsuit \diamondsuit]$, where
  \[
  \kappa \equiv 1 \mod \begin{cases}
    \pi^{2\ceil{e/2} + 1}, & \text{$R$ unramified} \\
    \pi^{2\ceil{e/2}}, & \text{$R$ ramified.}
  \end{cases}
  \]
  Then $[\delta^\odot \heartsuit] \in \cF_{\ceil{e'/2}}$, so the conic $\cM^\odot = \cM^\odot_\delta$ has maximal squareness level $e$. If $e$ is even, we know that $\epsilon(\delta) = 1$ by Proposition \ref{prop:conic_1_mod_2}. So we may assume that $e$ is odd.

  Now $\kappa \equiv 1 \mod 2\pi$, so the associated bilinear forms
  \[
  \cM^\odot(\xi,\eta) = \lambda^\diamondsuit(\delta^\odot\xi\eta) \textand \cM^{\odot\prime}(\xi,\eta) = \lambda^\diamondsuit(\kappa\delta^\odot\xi\eta)
  \]
  are congruent modulo $2\pi$. So by Proposition \ref{prop:conic_perturb}, the two conics have the same Brauer class.
\end{proof}

\begin{nota}
  We denote by $\ell_C$ the squareness level of the conic $\cM_1$ corresponding to $\delta = 1$, which plays a fundamental role. We let $\square_C = 2\ell_C + 1$, the squareness of the same conic.
\end{nota}

\section{\texorpdfstring{$\cN_{11}$}{N11}}
In this section, we will transform the $\cN_{11}$-condition, which says that all coordinates of $\widetilde\omega_C^{-1} \cdot \xi_1^2$ are congruent modulo $\pi^{n_{11}}$, into a more manageable form.

We will sometimes need to make some subtle reductions, and thus we make the following definition:

\begin{defn}
  A \emph{first vector problem} $\fP$ consists of a choice of resolvent algebra $R$ and as much of the discrete data as is needed to make $\cM_{11}$ and $\cN_{11}$ meaningful: the coarse coset $\delta_0\cF_0$, the reduced vector $\theta_1 \in \pi^{-\{b_1\}}R$ (which determines $\hat\omega_C$ and $s$), and the moduli $m_{11}$ and $n_{11}$. These are required to satisfy the requisite integrality properties, which essentially say that
  \[
  B_{\theta_1}(m_{11},n_{11}-s) = \pi^{ m_{11}} \OO_K {\theta_1} + \pi^{n_{11}-s} \OO_K \theta_2
  \]
  is a subset of $R$, but are otherwise untethered from a cubic or a quartic ring. The \emph{answer} to a first vector problem is the weighting
  \[
  W_{\fP} = W_{{\theta_1},m_{11},n_{11}} : \delta_0 \cF_0 \to \QQ_{\geq 0}
  \]
  that attaches to each quartic algebra $\delta \in \delta_0 \cF_0$ the volume of $\xi'_1 \in \PP(\OO_R)$ such that the corresponding $\xi_1 = \xi'_1\gamma_{1}$ satisfies the resolvent conditions
  \begin{alignat*}{2}
    \cM_{11} &:& \tr(\xi_1^2) \equiv 0 &\mod \pi^{ m_{11}} \\
    \cN_{11} &:& \quad \text{All coordinates of } \widetilde\omega_C^{-1} \cdot \xi_1^2 \text{ are congruent} &\mod \pi^{n_{11}}.
  \end{alignat*}
  We normalize volumes so that
  \[
  \mu(\OO_R) = 1 \textand \mu(\PP(\OO_R)) = 1 + \frac{1}{q} + \frac{1}{q^2}.
  \]
\end{defn}

We write $W^{\odot}$ instead of $W$ when we wish to speak instead of the volume of possibilities for the vector $\xi_1^\odot = \xi_1'/\gamma^{\odot}$ in Lemma \ref{lem:tfm_conic}. The two scalings are related by
\[
W^{\odot}_{\fP} = q^{v\left(N_{R/K}(\gamma^{\odot})\right)}W_{\fP}.
\]

First vector problems will be sorted into \emph{zones}, given by linear inequalities on $m_{11}$ and $n_{11}$, and having the properties that within each zone, the answer has a uniform description. Zones will be named by colors in such a way that a brightening of the color correlates with a lowering of $m_{11}$ and/or $n_{11}$, which increases the answer. Brightening is governed by the following poset:

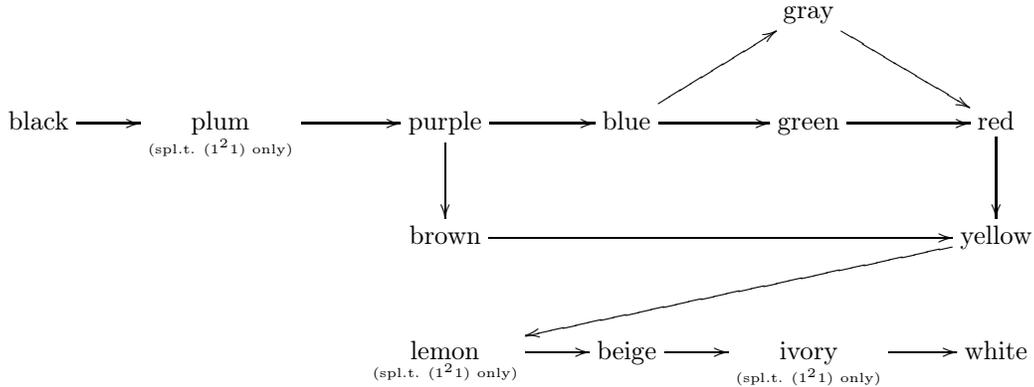
\[
\xymatrix@=0.1in{
  &&&& \text{gray} \ar[rd] \\
  \text{black} \ar[r] &
  \underset{\text{\tiny (spl.t. $(1^21)$ only)}}{\text{plum}}  \ar[r] &
  \text{purple} \ar[r] \ar[d] &
  \text{blue} \ar[r] \ar[ru] &
  \text{green} \ar[r] &
  \text{red} \ar[d] \\
  && \text{brown} \ar[rrr]
  &&& \text{yellow} \ar[llld] \\
  && \underset{\text{\tiny (spl.t. $(1^21)$ only)}}{\text{lemon}}  \ar[r] &
  \text{beige} \ar[r] &
  \underset{\text{\tiny (spl.t. $(1^21)$ only)}}{\text{ivory}}  \ar[r] &
  \text{white}
}
\]

\begin{lem} \label{lem:N11}
  Fix the data of a first vector problem $\fP$ in such a way that $\cN_{11}$ is active with
  \begin{equation}\label{eq:n11_leq_2e}
    0 < n_{11} \leq 2e < m_{11}^{\odot}
  \end{equation}
  and there is a solution $\xi_0$ to $\fP$. Let $\xi^{\odot}_0 = \xi_0 / \gamma_1^\odot$ be its transform. Then there is an $n^\odot \in \ZZ_{\geq 0}$ such that, for any $\xi_1$ satisfying $\cM_{11}$,
  \[
    \xi_1 \text{ satisfies } \cN_{11} \iff \xi_1^\odot \equiv \xi_0^\odot \mod \pi^{n^\odot}.
  \]
  The value of $n^\odot$ is given as follows:
  \begin{itemize}
    \item In unramified splitting types,
    \[
      n^\odot = \ceil{\frac{n_{11}}{2}}.
    \]
    \item In splitting type $(1^3)$,
    \[
      n^\odot = \ceil{\frac{n_{11}}{2} - \frac{h_1}{3}}.
    \]
      \item In splitting type $(1^21)$,
    \[
      n^\odot = \max\left\{\ceil{\frac{n_{11} - d_0}{2}} + \frac{h_1}{2}, 0\right\}.
    \]
  \end{itemize}
\end{lem}
\begin{rem}
  The condition $n_{11} \leq 2e$ (which, as we will see, restricts us to the blue, green, red, yellow, and lemon zones) can be removed, but then our conclusion must be that there is a family $\{\xi_{0(1)}, \ldots, \xi_{0(r)}\}$ of basic solutions, $r \in \{1, 2, 4\}$. The formula for $n^\odot$ becomes more complicated, and we will be able to solve these zones by other means.
\end{rem}
\begin{proof}[Proof of Lemma \ref{lem:N11}]
  In view of Lemma \ref{lem:conic_lift}, we may assume $m_{11} = \infty$, replacing $\xi_1^\odot$ by a value in the same $e$-pixel that satisfies $\cM^\odot(\xi_1^\odot) = 0$ exactly.

  Let $\xi_{0}$ be a fixed solution to $\fP$, and let $\xi_1$ be any solution to $\cM_{11}$. Observe that $\xi_1^2$ and $\xi_{0}^2$ are both traceless, so their wedge product $\xi_1^2 \wedge \xi_{0}^2$ is a scalar multiple of $(1;1;1) \in \OO_{\bar K}^3$. (Here we identify $\Lambda^2 \OO_{\bar K}^3$ with $\OO_{\bar K}^3$ via the trace pairing and standard orientation, so that the wedge product is given by the same formula as the cross product on $\RR^3$.) Let $\{\xi_{0}^2, \alpha\}$ be an $\OO_{\bar K}$-basis for the traceless plane in $\OO_{\bar K}^3$. Write
  \[
  \xi_1^2 = c_0 \xi_{0}^2 + c_1 \alpha.
  \]
  The coefficient $c_1$ controls how far $\xi_1$ deviates from $\xi_{0}$ and thus the satisfaction of $\cN_{11}$:
  \begin{align}
    \cN_{11} &\iff \text{All coordinates of} \quad \widetilde\omega_C^{-1} \cdot \xi_1^2 \quad \text{are congruent} \mod \pi^{n_{11}} \nonumber \\
    &\iff \text{All coordinates of} \quad c_1 \widetilde\omega_C^{-1} \alpha \quad \text{are congruent} \mod \pi^{n_{11}} \nonumber
  \end{align}
  We claim that the element $\widetilde\omega_C^{-1} \alpha \in \OO_{\bar K}^3$ does \emph{not} have all coordinates congruent mod $\mm_{\bar K}$:
  \begin{itemize}
    \item If $s = 0$, then $\widetilde\omega_C$ is a unit so this is equivalent to $\alpha$ and $\xi_0^2$ being linearly independent modulo $\mm_{\bar K}$;
    \item If $s > 0$, then $\widetilde\omega_C^{-1} \sim (1; \pi^{s}; \pi^{s})$, so $\xi_0^{(K)}$ has positive valuation. Hence $\alpha^{(K)}$ and $(\widetilde\omega_C^{-1} \alpha)^{(K)}$ are units, while $(\widetilde\omega_C^{-1} \alpha)^{(Q)}$ is not.
  \end{itemize}

  Consequently
  \begin{align}
    \cN_{11} &\iff c_1 \equiv 0 \mod \pi^{n_{11}} \label{eq:x_c1_tricky} \\
    &\iff \xi_1^2 \wedge \xi_{0}^2 \equiv 0 \mod \pi^{n_{11}}. \label{eq:N11_wedgepow2}
  \end{align}

  Now \eqref{eq:N11_wedgepow2} is advantageous, because the wedge product $\xi_1^2 \wedge \xi_{0}^2$ has all its coordinates equal, so we can test $\cN_{11}$ by looking at any one of them. We have (coordinate indices mod $3$)
  \begin{align*}
    \left(\xi_1^2 \wedge \xi_{0}^2\right)^{(i)} &= {\xi_1^{(i+1)}}^2 {\xi_{0}^{(i-1)}}^2 - {\xi_1^{(i-1)}}^2 {\xi_{0}^{(i+1)}}^2 \\
    &= \left(\xi_1^{(i+1)} \xi_{0}^{(i-1)} - \xi_1^{(i-1)} \xi_{0}^{(i+1)}\right)
    \left(\xi_1^{(i+1)} \xi_{0}^{(i-1)} + \xi_1^{(i-1)} \xi_{0}^{(i+1)}\right)
  \end{align*}
  and the two factors are congruent modulo $2$, so, since $n_{11} \leq 2e$,
  \begin{equation}
    \cN_{11} \iff \left(\xi_1 \wedge \xi_{0}\right)^{(i)} \equiv 0 \mod \pi^{n_{11}/2}. \label{eq:N11_wedge}
  \end{equation}
  We now examine this for each coordinate $i$ in turn, and for each splitting type.

\subsubsection{Unramified.}
We first dispose of the case $h_1 = 1$. Here the conic is tiny, and by Lemma \ref{lem:conic_pi}, all solutions $\xi_1$ lie in a single $e$-pixel:
\begin{align*}
  \xi_1^\odot &\equiv \xi_0^\odot \mod \pi^e \\
  \xi_1 &\equiv \xi_0 \mod \left(\pi^{e + 1/2}; \pi^e; \pi^e\right) \\
  \xi_1^2 &\equiv \xi_0^2 \mod \left(\pi^{2e + 1}; \pi^{2e}; \pi^{2e}\right).
\end{align*}
Since we are assuming $n_{11} \leq 2e$, we find that $\cN_{11}$ is automatic. The stated value $n^\odot=\ceil{n_{11}/2}$ is therefore valid.

Now assume that $h_1 = 0$. Here $\delta^\odot \in \OO_R^\cross$ and $\xi_0, \xi_1 \in \sqrt{\delta^\odot}\OO_R$ are $\OO_{\bar K}$-primitive. We scale $\xi_1$ by $\OO_K^\cross$ to be as close to $\xi_0$ as possible. Then $k = v(\xi_1 - \xi_0)$ is an integer, and
\[
  \xi_0 \textand \frac{\xi_1 - \xi_0}{\pi^k}
\]
are linearly independent elements of $\sqrt{\delta^\odot}\OO_R$. In particular, their wedge product is primitive, so
\[
  v(\xi_0 \wedge \xi_1) = k.
\]
Hence
\begin{align*}
  \cN_{11} &\iff k \geq \frac{n_{11}}{2} \\
  &\iff k \geq \ceil{\frac{n_{11}}{2}} \\
  &\iff \xi_1 \equiv \xi_0 \mod \pi^{\ceil{n_{11}/2}} \\
  &\iff \xi_1^\odot \equiv \xi_0^\odot \mod \pi^{\ceil{n_{11}/2}},
\end{align*}
as desired.
\subsubsection{Splitting type $(1^3)$.} Scale $\xi_1 \in \bar\zeta_3^{-h_1} \sqrt{\delta^\odot} \OO_R^\cross$ to be as close to $\xi_0$ as possible, and consider the valuation $k = v(\xi_1 - \xi_0) \in \frac{1}{3}\ZZ$. If $k \geq e$, then both $\cN_{11}$ and its claimed transformation are easily seen to hold, so assume that $k < e$. Then:
\begin{itemize}
  \item We cannot have $k \in \ZZ$, for then rescaling $\xi_1$ would bring it closer to $\xi_0$.
  \item If $k \in \ZZ + h_1/3$, then $v(\xi_1^2 - \xi_0^2) = 2k \in \ZZ - h_1/3$, and $\xi_1^2 - \xi_0^2$ has its first-order term a multiple of
  \begin{align*}
    &\left(\bar\zeta_3^{-h_1}\sqrt{\delta^\odot}\right)^2 \cdot \pi_R^{2k} \\
    &= \bar\zeta_3^{h_1} \delta^\odot \cdot \pi^{2k} \bar\zeta_3^{-h_1} \\
    &= \delta^\odot \pi^{2k},
  \end{align*}
  which has trace $\sim \pi^{2k}$, contradicting the constraint that both $\xi_1^2$ and $\xi_0^2$ are traceless.
\end{itemize}
Hence $k \in \ZZ - h_1/3$. Note that
\[
  v(\xi_1 \wedge \xi_0) = v\big((\xi_1 - \xi_0) \wedge \xi_0\big) = k,
\]
since the leading terms of $\xi_0$ and $\xi_1 - \xi_0$ are multiples of different powers of $\zeta_3$. So
\begin{alignat*}{2}
  \cN_{11} &\iff &k &\geq n_{11}/2.
\end{alignat*}
Since $k$ belongs to $\ZZ-h_1/3$, this is equivalent to
\[
  k\geq -\frac{h_1}{3}+\ceil{\frac{n_{11}}{2}-\frac{h_1}{3}}.
\]
After the transformation to $\xi_1^\odot$, this gives the claimed
\[
  n^\odot = \ceil{\frac{n_{11}}{2} - \frac{h_1}{3}}.
\]

  \subsubsection{Splitting type $(1^21)$.}
When $h_1 = 1$, the conic is tiny and $n^\odot$ is immaterial for similar reasons to the $h_1 = 1$ case of the unramified splitting types.

When $h_1 = 3$, we claim that $\cM_{11}$ has no solutions. Observe that $n_{11} > 0$, so $m_{11} > d_0/2$ and $m_{11}^\odot = m_{11} - p^\odot = m_{11} - (d_0 - 3)/2 \geq 2$. Then write
\begin{align*}
  \cM^\odot(\xi_1^\odot) &= \left(\delta^\odot {\xi_1^\odot}^2\right)^{(K)} - \fI\left(\left(\delta^\odot {\xi_1^\odot}^2\right)^{(Q)}\right).
\end{align*}
The first term has even valuation since $\delta^{\odot(K)}$ is a unit. But
\[
  \delta^\odot {\xi_1^\odot}^2 \sim \pi_Q^3
\]
has $\fI$-value $\sim \pi$, so $\cM^\odot$ is unsatisfiable mod $\pi^2$. Therefore we may assume $h_1 \in \{0, 2\}$.

We have
\begin{align*}
  \cN_{11} &\iff \xi_1^{(Q1)} \xi_{0}^{(Q2)} - \xi_1^{(Q2)} \xi_{0}^{(Q1)} \equiv 0 \mod \pi^{n_{11}/2} \\
  &\iff \frac{\xi_1^{(Q1)}}{\xi_{0}^{(Q1)}} - \frac{\xi_1^{(Q2)}}{\xi_{0}^{(Q2)}} \equiv 0 \mod \pi^{n_{11}/2} \\
  &\iff \fI\left(\frac{\xi_1^{(Q)}}{\xi_{0}^{(Q)}}\right) \equiv 0 \mod \pi^{\frac{n_{11} - d_0}{2}},
\end{align*}
by the definition of $\fI$. Note that the left-hand side belongs to $\OO_K$, so
\[
  \cN_{11} \iff \fI\left(\frac{\xi_1^{(Q)}}{\xi_{0}^{(Q)}}\right) \equiv 0 \mod \pi^{\ceil{\frac{n_{11} - d_0}{2}} - r}
\]
for any real $r$, $0 \leq r < 1$. We will pick $r$ below. Now, transforming to the coordinate $\xi_1^\odot$,
\begin{align}
  \cN_{11} &\iff \fI\left(\frac{\xi_1^{\odot(Q)}}{\xi_{0}^{\odot(Q)}}\right) \equiv 0 \mod \pi^{\ceil{\frac{n_{11} - d_0}{2}} - r} \nonumber \\
  &\iff \frac{\xi_1^{\odot(Q)}}{\xi_{0}^{\odot(Q)}} \equiv c \mod \pi^{\ceil{\frac{n_{11} - d_0}{2}} + \frac{1}{2} - r} \\
  &\iff \xi_1^{\odot(Q)} \equiv c \cdot \xi_{0}^{\odot(Q)} \mod \pi^{\ceil{\frac{n_{11} - d_0}{2}} + \frac{h_1}{4} + \frac{1}{2} - r}, \quad \text{some } c \in \OO_K^\cross \label{eq:N11_Q}
\end{align}
We may scale so that $c = 1$. We take $r = \frac{2 - h_1}{4} \in \{0, 1/2\}$, ensuring that
\[
  n^\odot = \ceil{\frac{n_{11} - d_0}{2}} + \frac{h_1}{2}
\]
is an integer. We claim it works overall, that is, that we can drop the superscript $(Q)$'s in \eqref{eq:N11_Q}. To do this, we look at the $Q1$-coordinate. We have
\begin{align*}
  \cN_{11} &\iff \xi_1^{(K)} \xi_{0}^{(Q1)} - \xi_1^{(Q1)} \xi_{0}^{(K)} \equiv 0 \mod \pi^{n_{11}/2} \\
  &\iff \xi_1^{(K)} \equiv \xi_{0}^{(K)} \cdot \frac{\xi_1^{(Q1)}}{\xi_{0}^{(Q1)}} \mod \pi^{n_{11}/2} \\
  &\iff \xi_1^{\odot(K)} \equiv \xi_{0}^{\odot(K)} \cdot \frac{\xi_1^{\odot(Q1)}}{\xi_{0}^{\odot(Q1)}} \mod \pi^{\frac{2n_{11} - d_0 + h_1}{4}} \\
  &\implies \xi_1^{\odot(K)} \equiv \xi_{0}^{\odot(K)} \cdot \frac{\xi_1^{\odot(Q1)}}{\xi_{0}^{\odot(Q1)}} \mod \pi^{n^\odot} \\
  &\iff \xi_1^{\odot(K)} \equiv \xi_{0}^{\odot(K)} \mod \pi^{n^\odot},
\end{align*}
as desired, where the lone non-reversible step uses $d_0\geq h_1$ to give
\[
  \frac{2n_{11}-d_0+h_1}{4} \geq \ceil{\frac{n_{11}-d_0}{2}}+\frac{h_1}{2}=n^\odot.
\]
If $n^\odot < 0$, we can obviously set $n^\odot = 0$ without changing the conclusion.

\end{proof}

\chapter{Boxgroups}
\label{cha:boxgroups}
If $\cN_{11}$ is strongly active, then
\[
  \beta = \frac{\xi_1^2}{\omega_C}
\]
is a unit. By Lemma \ref{lem:beta}, the solutions to $\fP$ arise from the $\beta$ that lie in the box
\[
  x + y\pi^{n_{11} - s} \theta_1 + z\pi^{m_{11}} \theta_2, \quad x \in \pi^{-2a_1 + s}K, \quad y,z \in \OO_K.
\]
 Necessarily $-2a_1 + s \in \ZZ$ and $x \in \OO_K^\cross$. Also, $[\delta] = [\beta] \in H^1$. So the support of $\delta$ is bound up with the $H^1$-classes of units in various boxes. Certain boxes have pride of place: those for which the corresponding subset of $H^1$ is a group, which we will call a \emph{boxgroup.}

\section{Signatures of subgroups of \texorpdfstring{$H^1$}{H\textasciicircum 1}}
Recall that in Proposition \ref{prop:Sh_basis_quartic}, we filtered $H^1$ by level spaces $\cF_0 \supset \cF_1 \supset \cdots \supset \cF_{e'}$, where
\[
  e' = \begin{cases}
    e, & \text{$R$ unramified} \\
    2e, & \text{$R$ ramified.}
  \end{cases}
\]
We would like to define some additional subgroups of $H^1$. We use the following notion.

\begin{defn}
If $S \subseteq H^1$ is a subgroup, define the \emph{signature} of $S$ to be the sequence of $e' + 2$ subgroups
\[
  S_i = \frac{S \intsec \cF_i}{S \intsec \cF_{i+1}} \subseteq \cF_i/\cF_{i+1}, \quad -1 \leq i \leq e'.
\]
\end{defn}
The following subgroups $S_i \subseteq \cF_i/\cF_{i+1}$ will occur frequently and will be given names:
\begin{itemize}
  \item $\ZERO$ denotes the zero subgroup $\cF_{i+1}/\cF_{i+1}$;
  \item $\STAR$ denotes the entire group $\cF_{i}/\cF_{i+1}$;
  \item $\TEE$, in unramified resolvent for $0 \leq i < e$,  denotes the order-$q$ subgroup
  \[
    \TEE = \{[1 + a\pi^{2i+1}\theta_1 : a \in \OO_K]\};
  \]
\end{itemize}
Thus, for instance, $\cF_i$ for $0 \leq i \leq e'$ has signature $\ZERO .\ZERO^i \STAR^{e'-i} . \STAR$. We separate the first and last terms of a signature by periods, because there $\cF_i/\cF_{i+1}$ has order $\size{H^0}$ instead of $q$ or $q^2$. In splitting types $(3)$ and $(1^3)$, we can omit these terms since $\size{H^0} = 1$.

$\cF_i$ is the only subgroup with signature $\ZERO .\ZERO^i \STAR^{e'-i} . \STAR$. In general, however, the signature does not uniquely determine the subgroup, though it does determine the \emph{size} of the subgroup, since
\[
  \size{S} = \prod_i \size{S_i}.
\]
Note also that if $S$ has signature $S_{-1}S_0\ldots S_{e'-1}S_{e'}$, then its orthogonal complement $S^\perp$ has signature $S_{e'}^\perp S_{e'-1}^\perp \ldots S_0^\perp S_{-1}^\perp$, since the Tate pairing on $H^1$ induces a perfect pairing between $\cF_{i}/\cF_{i+1}$ and $\cF_{e'-i-1}/\cF_{e'-i}$.

In this section, our aim will be to define a family of \emph{boxgroups} in terms of which the ring totals will be written. These boxgroups will depend on the resolvent data alone. In the unramified splitting types, we will construct, for appropriate $\ell_0,\ell_1,\ell_2$, a boxgroup $T(\ell_0, \ell_1, \ell_2)$ of signature
\[
  \ZERO .\ZERO^{\ell_0} \TEE^{\ell_1} \STAR^{\ell_2} . \STAR.
\]
In splitting type $(1^3)$ they will have signatures
\[
  \ZERO .(\ZERO\ZERO)^{\ell_0} (\ZERO\STAR)^{\ell_1} (\STAR\STAR)^{\ell_2} . \STAR \textand \ZERO .(\ZERO\ZERO)^{\ell_0} (\STAR\ZERO)^{\ell_1} (\STAR\STAR)^{\ell_2} . \STAR, \quad \sum_i \ell_i = e
\]
and will be denoted by $T_{-1}(\ell_0, \ell_1, \ell_2)$ and $T_1(\ell_0, \ell_1, \ell_2)$ respectively, the subscript corresponding to the value of $h$.

Finally, in splitting type $(1^2 1)$, we will use two numbering systems. In the first system, the boxgroups with signatures
\begin{align*}
  &\ZERO .\ZERO^{d_0'}.(\ZERO\ZERO)^{\ell_0} (\ZERO\STAR)^{\ell_1} (\STAR\STAR)^{\ell_2} .\STAR^{d_0'}. \STAR \\
  &\ZERO .\ZERO^{d_0'}.(\ZERO\ZERO)^{\ell_0} (\STAR\ZERO)^{\ell_1} (\STAR\STAR)^{\ell_2} .\STAR^{d_0'}. \STAR
\end{align*}
will be denoted by $T_1(\ell_0, \ell_1, \ell_2)$ and $T_0(\ell_0, \ell_1, \ell_2)$, respectively; the subscript is the value of $h_\eta$, which depends on the letter type of the resolvent. In the second system, the boxgroup with signature
\[
  \ZERO .\ZERO^{g_0 - g_1} (\STAR\ZERO)^{g_1} \STAR^{2e - g_0 - g_1} . \STAR
\]
will be denoted by $T(g_0; g_1)$. Note that the two numbering systems are related by
\begin{equation}\label{eq:l_g_1pow21}
\begin{aligned}
g_0 &= e + \ell_0 - \ell_2 & \ell_0 &= \frac{g_0 - g_1 - d_0' - h_\eta}{2} \\
g_1 &= \ell_1 - h_\eta & \ell_1 &= g_1 + h_\eta \\
&& \ell_2 &= \frac{2e - g_0 - g_1 - d_0' - h_\eta}{2}.
\end{aligned}
\end{equation}

\section{Boxgroups in unramified splitting type}

\begin{lem}\label{lem:boxes_basic_ur} Let $ m, n \in \NN^+ \union \{\infty\} $, $ m \geq n > 0 $. Let $ B_{\theta_1}(m,n) $ be the box
  \[
  B_{\theta_1}(m,n) = \{\pi c_0 + \pi^{n}c_1{\theta_1} + \pi^m c_2 \theta_2 : c_i \in \OO_K \}.
  \]
  \begin{enumerate}[$($a$)$]
    \item\label{boxes_basic_ur:recenter} For every $ \xi \in 1 + B_{\theta_1}(m,n) $,
    \begin{align}
    B_{\theta_1}(m,n) &= \xi B_{\xi^{-1} {\theta_1}}(m,n) \label{eq:recenter_1} \\
    1 + B_{\theta_1}(m,n) &= \xi \left(1 + B_{\xi^{-1} {\theta_1}}(m,n)\right) \label{eq:recenter_2}.
    \end{align}
    \item If
    \begin{equation*}
    m \leq 2n + s,
    \end{equation*}
    then $ B_{\theta_1}(m,n) $ is closed under multiplication and the translated box $ 1 + B_{\theta_1}(m,n) $ is a group under multiplication.
  \end{enumerate}
\end{lem}
\begin{proof}
  \begin{enumerate}[$($a$)$]
    \item Since $ \xi \equiv 1 $ mod $ \pi $, we can take $ \xi^{-1}\theta_2 $ in the role of $ \theta_2 $ for defining $ \xi B_{\xi^{-1} {\theta_1}}(m,n) $, which is a lattice with basis $ [\pi\xi, \pi^n{\theta_1}, \pi^m\theta_2] $, all of which are contained in $ B_{\theta_1}(m,n) $. This proves the reverse inclusion of \eqref{eq:recenter_1}, and equality follows by comparing volumes. To get \eqref{eq:recenter_2}, we add $ \xi $ to both sides and use $ \xi - 1 \in B_{\theta_1}(m,n) $ to simplify the left-hand side.
    \item Observe that $ B_{\theta_1}(m,n) $ is a lattice with basis $ [\pi, \pi^n{\theta_1}, \pi^m\theta_2] $. Since $ m \geq n $, the only product that does not clearly lie in the lattice is $ (\pi^{n} {\theta_1})^2 $, whose $ 1 $- and $ {\theta_1} $-components are divisible by $ \pi^{2n} $, and whose $ \theta_2 $-component is divisible by $ \pi^{2n + s} $. Since $ m \leq 2n + s $, this product lies in the lattice.

    Thus $ B_{\theta_1}(m,n) $ is closed under multiplication and so is $ 1 + B_{\theta_1}(m,n) $. To show the existence of inverses, simply note that the power series
    \[
    \frac{1}{1 + \xi} = 1 - \xi + \xi^2 - \xi^3 + \cdots
    \]
    converges to an element of $ 1 + B_{\theta_1}(m,n) $ for every $ \xi \in B_{\theta_1}(m,n) $.
\end{enumerate}
\end{proof}
\begin{lem} \label{lem:eta_ur}
Let
\[
  \square_C = \min \left\{ 2\ell(\hat\omega_C) + 1, e \right\}
\]
be the squareness of the conic for $\delta = 1$; put $\square_C = 0$ if $[\hat\omega_C] \notin \cF_0$ (i.e.\ $s$ is odd).

If $\theta_1$ is translated by a suitable element of $\OO_K$ and scaled by a suitable element of $\OO_K^\cross$ (neither of which change the associated resolvent $C$), then there is an $\eta \in \OO_R$ such that
\begin{equation} \label{eq:eta_sqrt}
  \eta^2 \equiv \theta_1 \mod \pi^{\square_C + s}
\end{equation}
and
\begin{equation} \label{eq:eta_in}
  \eta \in \left\langle1, \theta_1, \pi^{\ceil{\frac{s}{2}}} \theta_2, 2\theta_2\right\rangle.
\end{equation}
\end{lem}
\begin{proof}
If $s$ is odd, then we can scale and translate $\theta_1$ so that
\[
  \theta_1 \equiv (1; 0; 0) \mod \pi^s.
\]
Then $\eta = \theta_1$ satisfies the desired conditions.

If $s$ is even, then by definition of squareness, there is a linear form $\lambda$ such that
\[
  \tr(1^\odot{\xi^\odot}^2) \equiv c \cdot \lambda(\xi^\odot)^2 \mod \pi^{\square_C}
\]
as functions of $\xi^\odot \in \OO_R$. Here $1^\odot$, the transform of $\delta = 1$ under Lemma \ref{lem:tfm_conic}, is a unit whose class in $H^1$ is $[\hat\omega_C]$; for concreteness, we may take
\[
  1^\odot = \frac{\hat\omega_C}{(\pi^s; 1; 1)}.
\]
Since $\hat\omega_C$ is traceless, the conic has a distinguished basepoint  $\xi'_0 = [\pi^{s/2}; 1; 1]$. Pick a $\xi' = \xi'_1$ in the kernel of $\lambda$ that does \emph{not} lie in the same $1$-pixel as the basepoint. We claim that the choice
\[
  \eta = (1; \pi^{s/2}; \pi^{s/2}) \xi'_1
\]
fulfills the conditions.

The $\theta_2$-coefficient of $\eta^2$ is given by
\[
  \tr(\omega_C \eta^2) = \pi^s \tr(1^\odot \xi_1'^2) \equiv \pi^s c\,\lambda(\xi'_1)^2 = 0 \mod \pi^{s + \square_C}.
\]
Hence there are $a, b \in \OO_K$ such that
\begin{equation}\label{eq:eta_ur_ab}
  \eta^2 \equiv a + b \theta_1 \mod \pi^{s + \square_C}.
\end{equation}
We claim that $\pi \nmid b$. Suppose not. If $s = 0$, we get $\xi'_1 \equiv \xi'_0$ modulo $\pi$, contrary to hypothesis. If $s > 0$, we get $\pi | a$ so $\pi | (\xi'_1)^{(K)}$; we also know that $\pi | (\xi'_0)^{(K)}$. But $1^\odot$ is a unit, so $\lambda$ is a perfect linear functional, and its kernel in $\PP^2(k_K)$ intersects the line $(\xi')^{(K)} \equiv 0$ in only one $1$-pixel, a contradiction.

So $\pi \nmid b$, and hence we can take $a + b \theta_1$ for our new $\theta_1$, so \eqref{eq:eta_ur_ab} proves \eqref{eq:eta_sqrt}.

As for \eqref{eq:eta_in}, since $\theta_1^{(2)} \equiv \theta_1^{(3)} \mod \pi^s$, the span on the right-hand side equals the set of $\xi \in \OO_R$ with
\[
  v_K(\xi^{(2)} - \xi^{(3)}) \geq \min\left\{\frac{s}{2}, e\right\}.
\]
To prove that $\eta$ is such a $\xi$, observe that
\[
  \left(\eta^{(2)} - \eta^{(3)}\right)\left(\eta^{(2)} + \eta^{(3)}\right) = (\eta^{(2)})^2 - (\eta^{(3)})^2 \equiv \theta_1^{(2)} - \theta_1^{(3)} \equiv 0 \mod \pi^s,
\]
with the two factors on the left-hand side congruent modulo $2 \sim \pi^e$.
\end{proof}

\begin{lem} \label{lem:boxgps_ur} Let $0 < n \leq m \leq 2e$ be integers such that
\begin{align}
m &\leq 2n + s && (\text{the gray-red inequality})\label{eq:box_red_ur} \\
m &\leq n + \square_C + s && (\text{the gray-green inequality})\label{eq:box_green_ur} \\
m &\leq e + \frac{n+s+1}{2} && (\text{the gray-blue inequality})\label{eq:box_blue_ur}.
\end{align}
Then the projection $[1 + B_{\theta_1}(m,n)]$ of $ 1 + B_{\theta_1}(m,n) $ onto $H^1$ is a subgroup of signature
\[
  \ZERO .\ZERO^{\ell_0} \TEE^{\ell_1} \STAR^{\ell_2} . \STAR
\]
where
\begin{align*}
\ell_0 &= \floor{\frac{n}{2}} \\
\ell_1 &= \floor{\frac{m}{2}} - \floor{\frac{n}{2}} \\
\ell_2 &= e - \floor{\frac{m}{2}}.
\end{align*}
\end{lem}
\begin{proof}

The gray-red inequality \eqref{eq:box_red_ur} ensures that the projection $T = [1 + B_{\theta_1}(m,n)]$ is a subgroup. It is clear that
  \[
    \cF_{\floor{\frac{m}{2}}} \subseteq T \subseteq \cF_{\floor{\frac{n}{2}}},
  \]
  so the signature of $T$ has the shape $  \ZERO .\ZERO^{\ell_0} ?^{\ell_1} \STAR^{\ell_2} . \STAR $ with $\ell_1$ undetermined terms. These terms contain $\TEE$, because for $i \geq \floor{n/2}$ and for all $a \in \OO_K$, we have
  \[
    [1 + a \pi^{2i + 1} \theta_1] \in T.
  \]
  Thus the signature is at least the one claimed.

  To prove that equality occurs, we fix $ m \leq 2e $ and proceed by downward induction on $ n $. The base case $ n = m $ is clear since $ T = \cF_{\floor{\frac{m}{2}}} $. When moving from $ n + 1 $ to $ n $, note that $ \size{T} $ can grow by at most a factor of
    \[
      [1 + B(m,n) : (1 + B(m, n+1))(1 + \pi^n \OO_K)] = q.
    \]
    If $ n $ is odd, there is nothing to prove, as we claim that $ \size{T} $ actually grows by a factor of $ q $. If $ n $ is even, we claim that $ T $ does not change. It suffices to prove that each of the $ q $ cosets in
    \[
      \frac{(1 + B(m,n))}{(1 + B(m,n+1))(1 + \pi^n \OO_K)}
    \]
    contains a square. For $ c \in \OO_K $, consider
    \begin{align*}
    (1 + \pi^{\frac{n}{2}}c\eta)^2 = 1 + 2\pi^{\frac{n}{2}} c \eta + \pi^n c^2 \eta^2.
    \end{align*}
    The last term is $ \pi^n c^2 {\theta_1} $ up to an error in $ \pi^{n + s + \square_C}\OO_R$, which is in $\pi^m \OO_R$ by the gray-green inequality. We claim that the cross term $2\pi^{n/2} c \eta$ lies in $B(m,n+1)$ also. If $s$ is odd, this is trivial since we took $\eta = \theta_1$. Otherwise, we have
    \[
      2\pi^{n/2} \eta \in 2\pi^{n/2} \left\langle1, {\theta_1}, \pi^{s/2}\theta_2, 2\theta_2\right\rangle
      = B_{\theta_1}\left( e + \frac{n}{2} + \min\left\{\frac{s}{2}, e\right\}, e + \frac{n}{2} \right)
    \]
    We get the needed inequality
    \[
      e + \frac{n}{2} + \frac{s}{2} \geq m
    \]
    from the gray-blue inequality, the difference of whose sides lies in $\ZZ + 1/2$ by parity considerations. So we have found a square in the coset $ 1 + \pi^n c^2 {\theta_1} + B(m,n+1) = (1 + \pi^n c^2 {\theta_1})(1 + B(m,n+1)) $, as desired.
\end{proof}
As a corollary, we have:
\begin{lem}\label{lem:boxgpT_ur}
    For every triple $ (\ell_0, \ell_1, \ell_2) $ of nonnegative integers satisfying
    \begin{align}
    \ell_0 + \ell_1 + \ell_2 &= e \\
    \ell_1 &\leq \ell_0 + \frac{s}{2} + 1 && (\text{the gray-red inequality})\label{eq:boxgp_red} \\
    \ell_1 &\leq \frac{s + \square_C + 1}{2} && (\text{the gray-green inequality}) \label{eq:boxgp_green} \\ 
    \ell_1 &\leq \ell_2 + \frac{s}{2} + 1, && (\text{the gray-blue inequality})\label{eq:boxgp_blue}
    \end{align}
    there is a \emph{boxgroup} $T(\ell_0,\ell_1,\ell_2) \subseteq H^1$ of signature $\ZERO .\ZERO^{\ell_0} \TEE^{\ell_1} \STAR^{\ell_2} . \STAR$, such that, if $m$, $n$ are integers satisfying the conditions of Lemma \ref{lem:boxgps_ur}, then
    \[
    [1 + B_{\theta_1}(m,n)] = T \left( \floor{\frac{n}{2}}, \floor{\frac{m}{2}} - \floor{\frac{n}{2}}, e - \floor{\frac{m}{2}} \right).
    \]
\end{lem}
\begin{proof}
    If $\ell_1 = 0$, take $T(\ell_0, 0, \ell_2) = \cF_{\ell_0}$, the unique subgroup with the correct signature.

    Otherwise, let $ m = 2\ell_0 + 2\ell_1 $, $ n = 2\ell_0 + 1 $ in the preceding lemma. The transformation of the gray-red, gray-green, and gray-blue conditions is routine.

    For the last claim, note that decreasing $m$ or increasing $n$ can only make the conditions of Lemma \ref{lem:boxgps_ur} truer, with the exception of the condition $m \geq n$. If $\floor{m/2} = \floor{n/2}$, then clearly $[1 + B_{\theta_1}(m,n)] = \cF_{\floor{m/2}}$, so we can assume that
    \[
      n \leq 2 \floor{\frac{n}{2}} + 1 < 2 \floor{\frac{m}{2}} \leq m.
    \]
    Clearly
    \[
    [1 + B_{\theta_1}(m,n)] \subseteq \left[1 + B_{\theta_1}\left(2\floor{\frac{m}{2}}, n\right)\right],
    \]
    but both sides have the same signature, so equality holds. Likewise,
    \begin{align*}
    \left[1 + B_{\theta_1}\left(2\floor{\frac{m}{2}}, n\right)\right] &\supseteq \left[1 + B_{\theta_1}\left(2\floor{\frac{m}{2}}, 2\floor{\frac{n}{2}} + 1\right)\right] \\
    &= T\left( \floor{\frac{n}{2}}, \floor{\frac{m}{2}} - \floor{\frac{n}{2}}, e - \floor{\frac{m}{2}} \right),
    \end{align*}
    but both sides have the same signature, so equality holds.
\end{proof}

\subsection{Supplementary boxgroups}\label{sec:suppl_ur}
As thus defined, all boxgroups $T$ satisfy $\cF_e \subseteq T \subseteq \cF_0$. Groups not satisfying these inclusions will be denoted as follows.

If $s > 0$, so that a distinguished splitting $R \isom K \cross Q$ exists, consider the image $\iota(K^\cross)$ of the map
\begin{align*}
  \iota\colon K^\cross/(K^\cross)^2 &\to H^1 \\
  a &\mapsto [(a;1)] = [(1;a)].
\end{align*}
We find that $\iota(K^\cross)$ has signature
\[
  \TEE.\TEE^{e}.\TEE
\]
where the middle $e$-many $\TEE$'s denote the usual subgroups
\[
  \TEE_i = \{1 + \pi^{2i+1} a \theta_1 \} = \{1 + \pi^{2i+1} a (1; 0)\} \subseteq \cF_i/\cF_{i+1}
\]
and where the initial $\TEE_{-1} = \{[1], [(1;\pi)]\} \subseteq H^1/\cF_0$ and the final $\TEE_e = \{1 + 4a : a \in \OO_K\} \subseteq \cF_e$ have size $2$ and $\size{H^0}/2$ respectively. In particular,
\[
  \size{\iota(K^\cross)} = \size{H^0} q^{e} = \sqrt{\size{H^1}}.
\]
From the explicit description in terms of the Hilbert pairing, we find that $\iota(K^\cross)$ is isotropic and hence maximally isotropic for $\left\langle\bullet, \bullet\right\rangle$.

The group $\iota(K^\cross)$ is always important, but it does not behave well with respect to boxgroups unless $s > 2e$, in which case we give it the name $T(\emptyset, e, \emptyset)$.

If $0 \leq \ell_1 \leq e$ and $s > 2\ell_1 $, set
\begin{align*}
T(e - \ell_1, \ell_1, \emptyset) &= \iota(K^\cross) \intsec \cF_{e - \ell_1} = \{[(a ; 1)] \in H^1 : a \equiv 1 \mod \pi^{2e - 2\ell_1 + 1}\} \\
T(\emptyset, \ell_1, e - \ell_1) &= \iota(K^\cross) \cdot \cF_{\ell_1} = \{[(a ; \alpha)] \in H^1 : \alpha \equiv 1 \mod \pi^{2\ell_1+1}\}.
\end{align*}
  (In the first case, the equivalence of the two definitions is established using Lemma \ref{lem:iota}.)
Their signatures are, respectively, $\ZERO.\ZERO^{e-\ell_1}\TEE^{\ell_1}.\TEE$ and $\TEE.\TEE^{\ell_1}\STAR^{e - \ell_1}.\STAR$.

The restrictions on $s$ ensure that these boxgroups satisfy such natural relations as
\begin{align*}
  T(e - \ell_1, \ell_1, \emptyset) \cdot \cF_e &= T(e - \ell_1, \ell_1, 0) \\
  T(\emptyset, \ell_1, e - \ell_1) \intsec \cF_0 &= T(0, \ell_1, e - \ell_1),
\end{align*}
which we will often use without comment.

Finally, in all cases, we let
\begin{align*}
T(e, \emptyset, \emptyset) &= \{1\} \\
T(\emptyset, \emptyset, e) &= H^1.
\end{align*}

It will turn out that $T(\ell_0,\ell_1,\ell_2)$ and $T(\ell_2,\ell_1,\ell_0)$ are orthogonal complements whenever both are defined (Lemma \ref{lem:orth}). Actually, this is simple to prove in the case that one of $\ell_0,\ell_1,\ell_2$ is the symbol $\emptyset$. The sizes of these groups follow immediately from their signatures:
\begin{lem}\label{lem:111_T_size}
  If $T(\ell_0,\ell_1,\ell_2)$ is defined and $\ell_1 \neq \emptyset$, then
  \[
  \size{T(\ell_0,\ell_1,\ell_2)} = \size{H^0}q^{e + \ell_2 - \ell_0},
  \]
  where if $\emptyset$ occurs as either $\ell_0$ or $\ell_2$, it must be replaced by $-1/[k_K : \FF_2] = \log_q(1/2)$.
\end{lem}

\section{Boxgroups in splitting type (\texorpdfstring{$1^3$}{1³})}
Fix the integer $h \in \{1, -1\}$ such that
\[
b_1 \in \ZZ + \frac{h}{3}, \quad \theta_1 \in \bar\zeta_3^h \OO_R^\cross, \quad \theta_2 \in \bar\zeta_3^{-h} \OO_R^\cross.
\]
(Note the tight connection with the $h_i$ of Lemma \ref{lem:tfm_conic}: if $\cN_{11}$ is strongly active, then $\beta$ is a unit in Lemma \ref{lem:beta}, from which we get $m_{11} \in \ZZ$, $a_1 \in \ZZ$, and $h_1 = h$.)

Then fix $\theta_1 \in \bar\zeta_3^{h} \OO_R^\cross$. The last basis vector $\theta_2$, which carries no information, can be fixed at
\[
  \theta_2 = \bar\zeta_3^{-h}.
\]

Let
\[
\square_C = \min \left\{ 2\floor{\frac{\ell(\hat\omega_C)}{2}} + 1, e \right\}
\]
be the squareness of the conic
\[
  \cM^\odot(\xi^\odot) = \lambda^\diamondsuit(1^\odot{\xi^\odot}^2)
\]
that occurs for $\delta = 1$. Note that $\cM^\odot$ is independent of $h_1$, with $[1^\odot] = [\hat\omega_C]$, although we restrict $\xi_1^\odot$ to the generic or special $1$-pixels according as $h_1 = 1$ resp.\ $-1$. Similar to Lemma \ref{lem:eta_ur}, we have:

\begin{lem}\label{lem:eta_1pow3}
Fix a resolvent $C$. If $\theta_1$ is translated and scaled appropriately, there is an $\eta \in \bar\zeta_3^{-h} \OO_R^\cross$ such that
\[
  \eta^2 \equiv \theta_1 \mod \pi^{\square_C - \frac{2 h}{3}}.
\]
\end{lem}
\begin{proof}
By Lemma \ref{lem:squareness}, $\cM^\odot$ has squareness $\square_C$, which means that there is a linear form $\lambda$ and a scalar $c \in \OO_K^\cross$ such that
\[
  \cM^\odot(\xi^\odot) \equiv c \cdot \lambda(\xi^\odot)^2 \mod \pi^{\square_C}
\]
The zero locus of $\cM^\odot$ modulo $\pi$, or equivalently of $\lambda$ modulo $\pi$, consists of $(q+1)$-many $1$-pixels, $q$ generic and one special. There is a basepoint
\[
  \xi^\odot = \pi^{(1 - h)/3}\sqrt{\frac{\diamondsuit\omega_C}{1^\odot}},
\]
due to the tracelessness of $\omega_C$, which lies in a generic pixel if $h = 1$ and the special pixel if $h = -1$. We pick a $\xi^\odot \in \ker \lambda$ from the \emph{other} type of pixel; we have
\[
  \cM^\odot(\xi^\odot) \equiv 0 \mod \pi^{\square_C}.
\]
Then take
\[
  \eta = \pi^{-\frac{h}{3}} \sqrt{\frac{1^\odot}{\hat\omega_C \diamondsuit}} \xi^\odot \in \bar{\zeta}_3^{-h} \OO_R^\cross.
\]
Note that, when expressed in the basis $(1, \theta_1, \theta_2)$,
\begin{align*}
  \coef_{\theta_2}(\eta^2) &= \pi^{-\frac{2h}{3}} \cM^\odot\left(\xi^\odot\right) \\
  &\equiv 0 \mod \pi^{\square_C - \frac{2h}{3}}.
\end{align*}
Hence there exist
\[
  a \in \begin{cases}
    \sqrt[3]{\pi^2}\OO_K, & h = 1 \\
    \sqrt[3]{\pi}\OO_K, & h = -1
  \end{cases}
  \textand b \in \OO_K
\]
such that
\[
  \eta^2 \equiv a + b {\theta_1} \mod \pi^{\square_C - \frac{2h}{3}}.
\]
Since $\eta$ is a unit, we must have $\pi \nmid b$, so $\theta_1$ can be replaced by $a + b\theta_1$.
\end{proof}

The following lemma is proved just like Lemma \ref{lem:boxes_basic_ur}.

\begin{lem}\label{lem:boxes_basic_1pow3} If $m$ and $n$ satisfy
  \[
    0 < n < m < 2e, \quad m \in \ZZ - \frac{h}{3}, \quad n \in \ZZ + \frac{h}{3},
  \]
  let $ B_{\theta_1}(m,n) $ be the box
  \[
  B_{\theta_1}(m,n) = \{\pi c_0 + \pi^{n} c_1 {\theta_1} + \pi^m c_2 \theta_2 : c_0,c_1,c_2 \in \OO_K \} \subseteq \OO_R.
  \]
  \begin{enumerate}[$($a$)$]
    \item For every $ \xi \in 1 + B_{\theta_1}(m,n) $,
    \begin{align}
    B_{\theta_1}(m,n) &= \xi B_{\xi^{-1} {\theta_1}}(m,n) \\
    1 + B_{\theta_1}(m,n) &= \xi \left(1 + B_{\xi^{-1} {\theta_1}}(m,n)\right).
    \end{align}
    \item If
    \begin{equation*}
    m \leq 2n,
    \end{equation*}
    then $ B_{\theta_1}(m,n) $ is closed under multiplication and the translated box $ 1 + B_{\theta_1}(m,n) $ is a group under multiplication.
  \end{enumerate}
\end{lem}
Our goal is to study the projection of $1 + B_{\theta_1}(m,n)$ onto $H^1$. The following yields the conditions under which a useful group is formed thereby:
\begin{lem} \label{lem:boxgps_1pow3} Let $0 < n < m \leq 2e$ be rational numbers with $ m \in \ZZ - \frac{h}{3}, n \in \ZZ + \frac{h}{3}$. Assume that
  \begin{align}
  m &\leq 2n && (\text{the gray-red inequality}) \label{eq:box_red_1pow3} \\
  m &\leq n + \square_C - \frac{2h}{3} && (\text{the gray-green inequality}) \label{eq:box_green_1pow3} \\
  m &\leq \frac{n + 1}{2} + e && (\text{the gray-blue inequality})\label{eq:box_blue_1pow3}.
  \end{align}
  Write
  \[
  \tri m = m - \frac{2h}{3} \in \ZZ, \quad \tri n = n + \frac{2h}{3} \in \ZZ.
  \]
  Then the projection $T = [1 + B_{\theta_1}(m,n)]$ of $ 1 + B_{\theta_1}(m,n) $ onto $H^1$ is a subgroup of signature
  \[
    . (\ZERO\ZERO)^{\ell_0} (\STAR\ZERO)^{\ell_1} (\STAR\STAR)^{\ell_2}.
  \]
  or
  \[
    . (\ZERO\ZERO)^{\ell_0} (\ZERO\STAR)^{\ell_1} (\STAR\STAR)^{\ell_2}.
  \]
  for $h = 1$ and $h = -1$ respectively, where
  \begin{align*}
  \ell_0 &= \floor{\frac{\tri n}{2}} \\
  \ell_1 &= \floor{\frac{\tri m}{2}} - \floor{\frac{\tri n}{2}}  \\
  \ell_2 &= e - \floor{\frac{\tri m}{2}}
  \end{align*}
  except for the case $h = 1$, $\tri m = 2i + 1$, $\tri n = 2i + 2$ ($i \in \ZZ$), where $T = \cF_{2i+1}$.
\end{lem}
\begin{proof}

Suppose $h = 1$. For $a \in \OO_K^\cross$, there are elements in $T$ of the form
\begin{align*}
[1 + a \pi^{2i + \frac{1}{3}}\theta_1], \quad i &\geq \ceil{\frac{3n - 1}{6}} = \floor{\frac{\tri n}{2}} \\
[1 + a \pi^{2i + \frac{5}{3}}\theta_2], \quad i &\geq \ceil{\frac{3m - 5}{6}} = \floor{\frac{\tri m}{2}}.
\end{align*}
Likewise, in the case $h = -1$, there are elements in $T$ of the form
\begin{align*}
[1 + a \pi^{2i + \frac{5}{3}}\theta_1], \quad i &\geq \ceil{\frac{3n - 5}{6}} = \floor{\frac{\tri n}{2}} \\
[1 + a \pi^{2i + \frac{1}{3}}\theta_2], \quad i &\geq \ceil{\frac{3m - 1}{6}} = \floor{\frac{\tri m}{2}}.
\end{align*}
This shows that the signature of $T$ is at least as large as claimed.

To show equality, we fix $m < 2e$ and proceed by downward induction on $n$. The base case $n = m - 1/3$ (for $h = 1$) or $n = m - 2/3$ (for $h = -1$) is clear since $T = \cF_{\tri m - 1}$ or $\cF_{2\floor{\tri m/2} + 1}$ respectively. When moving from $ n + 1 $ to $ n $, note that $ \size{T} $ can grow by at most a factor of
\[
[1 + B(m,n) : (1 + B(m, n+1))(1 + \pi^{\ceil{n}} \OO_K)] = q.
\]
If $ \tri n$ is odd, there is nothing to prove, as we claim that $ \size{T} $ actually grows by a factor of $ q $. If $ \tri n $ is even, we are claiming that $ T $ does not change. It suffices to prove that each of the $ q $ cosets in
\[
\frac{(1 + B(m,n))}{(1 + B(m,n+1))(1 + \pi^{\ceil{n}} \OO_K)}
\]
contains a square. Recall the approximate square root $\eta \in \bar\zeta_3^{-h} \OO_R^\cross$ from Lemma \ref{lem:eta_1pow3}, which satisfies
\[
\eta^2 \equiv \theta_1 \mod \pi^{\square_C - \frac{2 h}{3}}.
\]
For $ c \in \OO_K $, consider
\begin{equation} \label{eq:x_square_1pow3}
(1 + \pi^{\frac{n}{2}}c\eta)^2 = 1 + 2\pi^{\frac{n}{2}} c \eta + \pi^n c^2 \eta^2.
\end{equation}
The last term is $ \pi^n c^2 {\theta_1} $ up to an error in $ \pi^{n + \square_C - \frac{2h}{3}}\OO_R$. To say that this is in $\pi^m \OO_R$, we need the inequalities
\[
  m \leq n + \square_C - \frac{2h}{3} \textand m \leq n + e - \frac{2h}{3}.
\]
The first of these is \eqref{eq:box_green_1pow3}, and the second follows easily from \eqref{eq:box_red_1pow3} and \eqref{eq:box_blue_1pow3}. We claim that the middle term of \eqref{eq:x_square_1pow3} lies in $ \pi^m \OO_R $ also, that is,
\[
  m \leq \frac{n}{2} + e.
\]
This follows from \eqref{eq:box_blue_1pow3} and the fact that $m - n/2 \in \ZZ$. So we have found a square in the coset $ 1 + \pi^n c^2 {\theta_1} + B(m,n+1) = (1 + \pi^n c^2 {\theta_1})(1 + B(m,n+1)) $, as desired.
\end{proof}

As a corollary, just like Lemma \ref{lem:boxgpT_ur}, we get the following:
\begin{lem}\label{lem:boxgpT_1pow3}
For every triple $ (\ell_0, \ell_1, \ell_2) $ of nonnegative integers satisfying
\begin{align}
\ell_0 + \ell_1 + \ell_2 &= e \\
\ell_1 &\leq \ell_0 + 1 - h && (\text{the gray-red inequality})\label{eq:boxgp_red_1pow3} \\
\ell_1 &\leq \frac{\square_C + 1}{2} - h && (\text{the gray-green inequality})\label{eq:boxgp_green_1pow3} \\
\ell_1 &\leq \ell_2 + 1 - h, && (\text{the gray-blue inequality})\label{eq:boxgp_blue_1pow3}
\end{align}
there is a \emph{boxgroup} $T_h(\ell_0,\ell_1,\ell_2) \subseteq H^1$ of signature
\[
. (\ZERO\ZERO)^{\ell_0} (\STAR\ZERO)^{\ell_1} (\STAR\STAR)^{\ell_2}. \quad (h = 1)
\]
or
\[
. (\ZERO\ZERO)^{\ell_0} (\ZERO\STAR)^{\ell_1} (\STAR\STAR)^{\ell_2}. \quad (h = -1)
\]
such that, if $m$, $n$ are rational numbers satisfying the conditions of Lemma \ref{lem:boxgps_1pow3}, then
\[
[1 + B_{\theta_1}(m,n)] = T_h \left( \ell_0', \ell_1', \ell_2' \right)
\]
where $\ell_0', \ell_1', \ell_2'$ are the numbers $\ell_0, \ell_1, \ell_2$ defined in Lemma \ref{lem:boxgps_1pow3}. 
\end{lem}
\begin{proof}
If $\ell_1 = 0$, take $T_h(\ell_0, 0, \ell_2) = \cF_{\ell_0}$. Otherwise, take
\[
  \tri n = 2\ell_0 + 1 \in 2\ZZ + 1, \quad \tri m = 2e - 2\ell_2 \in 2\ZZ,
\]
In other words,
\[
  n = 2\ell_0 + 1 - \frac{2h}{3}, \quad m = 2e - 2\ell_2 + \frac{2h}{3}.
\]
Conditions \eqref{eq:boxgp_red_1pow3}--\eqref{eq:boxgp_blue_1pow3} immediately imply \eqref{eq:box_red_1pow3}--\eqref{eq:box_blue_1pow3}. Just as in Lemma \ref{lem:boxgpT_ur}, we then argue that increasing $m$ by $1$ (resp.\ decreasing $n$ by $1$), if it does not violate \eqref{eq:box_red_1pow3}--\eqref{eq:box_blue_1pow3}, yields a boxgroup of the same signature that is contained in (resp.\ contains) $T_h(\ell_0, \ell_1, \ell_2)$ and thus must equal $T_h(\ell_0, \ell_1, \ell_2)$.
\end{proof}
The subscript ``$h$'' in $T_h(\ell_0, \ell_1, \ell_2)$ is logically superfluous, because ${\theta_1}$ is fixed. But it allows the following manipulation. Define
\begin{equation}\label{eq:neg_boxgps}
  T_{-1}(\ell_0, \ell_1, \ell_2) = T_1 (\ell_0 + \ell_1, -\ell_1, \ell_1 + \ell_2)
\end{equation}
for all $(\ell_0, \ell_1, \ell_2)$ for which either side has been defined. Note that $T_{1}(\ell_0, 0, e - \ell_0) = T_{-1}(\ell_0, 0, e - \ell_0) = \cF_{2\ell_0}$ already fulfills this relation, while allowing
\[
  T_1(\ell_0, -1, e+1 - \ell_0) = T_{-1}(\ell_0 - 1, 1, e - \ell_0) = \cF_{2\ell_0+1}
\]
saves us the trouble of excluding the case $h = 1$, $\tri m = 2i + 1$, $\tri n = 2i + 2$ from Lemma \ref{lem:boxgps_1pow3}. We do not use any other boxgroups with negative $\ell_1$ within this memoir, but in the code we do, using \eqref{eq:neg_boxgps} to convert everything to a $T_{-1}$.

\section{Boxgroups in splitting type (\texorpdfstring{$1^21$}{1²1})}
In splitting type $(1^21)$, the interaction between $d_0$ and the other parameters of the resolvent leads to a complexity that we must address first.

\begin{lem} \label{lem:iota}
  Let $Q = K[\sqrt{D_0}]$ be a quadratic extension, and let $a \in K^\cross$. The level $\ell_Q(a)$ considered as an element of $Q^\cross/(Q^\cross)^2$ is determined by the levels of $a$ and $D_0 a$ in $K^\cross/(K^\cross)^2$ in the following way:
  \begin{enumerate}[$($a$)$]
    \item If $Q/K$ is unramified, then
    \[
    \ell_Q(a) = \begin{cases}
      \ell_K(a), & \ell_K(a) < e \\
      e+1, & \ell_K(a) \geq e.
    \end{cases}
    \]
    \item If $Q/K$ is ramified with $v_K(\disc_K Q) = d_0$, then
    \[
    \ell_Q(a) = \begin{cases}
      e + \ell_K(a), & \ell_K(a) \geq e - d_0/2 \\
      e + \ell_K(D_0 a), & \ell_K(D_0 a) \geq e - d_0/2 \\
      2\ell_K(a) + d_0/2, & -1/2 \leq \ell_K(a) < e - d_0/2.
    \end{cases}
    \]
    (Here we put $\ell_K(a) = -1/2$ if $a \in \pi\OO_K^\cross$.)
    \item\label{iota:custom} In particular,
    \[
    \ell_Q(a) \geq \frac{d_0}{2} + \floor{\ell_K(a)},
    \]
    and for any integer $m$, $0 \leq m \leq e - d_0/2$,
    \[
    \ell_K(a) \geq m \implies \ell_Q(a) \geq 2m + \frac{d_0}{2}.
    \]
  \end{enumerate}
\end{lem}
\begin{proof}
  Although a direct proof is not difficult, it is more illuminating to use the connection between levels and the discriminants of the corresponding quadratic extensions. It is not hard to see that
  \[
  (\disc K[\sqrt{a}]) = \begin{cases}
    (\pi^{2e - 2\ell(a)}), & -1/2 \leq \ell(a) \leq e \\
    (1), & \ell(a) \geq e.
  \end{cases}
  \]
  In particular,
  \[
  \ell_K(D_0) = e - d_0/2.
  \]
  Thus, except for distinguishing levels $e$ and $e+1$, finding $\ell_Q(a)$ is equivalent to determining the relative discriminant $(\disc_Q Q[\sqrt{a}])$. Since
  \begin{equation*}
    N_{Q/K}(\disc_Q Q[\sqrt{a}]) = \frac{\disc_K Q[\sqrt{a}]}{\disc_K Q},
  \end{equation*}
  it is enough (excluding the trivial case $Q = K \cross K$) to compute the absolute discriminant $\disc_K Q[\sqrt{a}] = \disc_K K[\sqrt{D_0}, \sqrt{a}]$.

  Now, since the regular representation of $\cC_2 \cross \cC_2$ is a direct sum of $1$-dimensional representations, we have by an Artin-conductor argument the identity
  \begin{align*}
    \disc_K K[\sqrt{D_0}, \sqrt{a}] = \disc_K K[\sqrt{D_0}] \cdot \disc_K K[\sqrt{a}] \cdot \disc_K K[\sqrt{aD_0}].
  \end{align*}
  This gives the required result in all cases except when $\ell(a) \geq e$ or $\ell(D_0 a) \geq e$, where it gives the incomplete conclusion $\ell_Q(a) \geq e$. Without loss of generality, $\ell(a) \geq e$. If $a$ is actually a square, it remains a square in $Q$, so assume $a = 1 + 4u$ is an intimate unit, $\tr_{k_K/\FF_2} (u \mod \pi) = 1$. If $Q$ is unramified, then $a$ is a square in $Q$: indeed we could have taken $D_0 = a$. But if $Q$ is ramified, then $a$ remains a nontrivial intimate unit in $Q$, so $\ell_Q(a) = 2e$, as claimed.
\end{proof}

\begin{lem}\label{lem:types_1pow21}
The values of $d_0$ and $s'$ constrain the value of $\ell(\hat\omega_C \diamondsuit\heartsuit)$, and hence of $\square_C$, in a different way in each of five cases, which we call \emph{letter types} and label with capital letters:
\begin{enumerate}[(A)]
  \item\label{type:A} If $s' = -1/2$, then $\ell(\hat\omega_C \diamondsuit\heartsuit) = -1/2$, so $\square_C = 0$.
  \item\label{type:B} If $0 \leq s' < d_0/2 - 1$, then $\ell(\hat\omega_C \diamondsuit\heartsuit) = s'$, so $\square_C = 2\floor{s'/2}+1$.
  \item\label{type:C} If $s' = d_0/2 - 1$, then $\hat\omega_C \diamondsuit\heartsuit \in \cF_{d_0/2 - 1}$, so $\square_C \geq 2\floor{(d_0-2)/4} + 1$.
  \item\label{type:D} If $s' \geq d_0/2$ and either $s' \equiv d_0/2$ mod $2$ or $d_0 = 2e + 1$, then $\hat\omega_C \diamondsuit\heartsuit \in \cF_{\floor{d_0/2}}$, so $\square_C \geq 2\floor{d_0/4} + 1$.
  \item\label{type:E} If $s' \geq d_0/2$ and $s' \equiv d_0/2 + 1$ mod $2$, then $\hat\omega_C \diamondsuit\heartsuit \in \spadesuit \cF_{d_0/2}$, so $\square_C = 2\floor{(d_0 - 2)/4} + 1$. Here $\spadesuit$ is the class in $H^1$ of $\iota(\pi) = (1; \pi; \pi)$, well defined up to $\iota(\OO_K^\cross) \subseteq \cF_{d_0/2}$.
\end{enumerate}
\end{lem}
\begin{proof}
Case \ref{type:A} follows immediately from Lemma \ref{lem:tfm_conic}. In the remaining cases, $b_1 \in \ZZ$ and $\theta_1 \in \OO_R$. By translation we may assume that $\theta_1 = (0; c_0 + c_1\pi_Q)$ for some $c_0 \in \OO_K^\cross$, $c_1 \in \OO_K$. We have
\[
  s' = v(c_1).
\]
Let $\iota_Q : Q \to Q$ denote the Galois involution. We compute
\begin{align*}
  \hat\omega_C &= \pi^{4b_1} \big(\theta_1^{(2)} - \theta_1^{(3)}\big) \big(\theta_1^{(3)} - \theta_1^{(1)}\big) \big(\theta_1^{(1)} - \theta_1^{(2)}\big) \big(\theta_1^{(2)} - \theta_1^{(3)}; \theta_1^{(3)} - \theta_1^{(1)}; \theta_1^{(1)} - \theta_1^{(2)}\big) \\
  &= \pi^{4b_1} \left(\theta_1^{(Q)} - \iota_Q\big(\theta_1^{(Q)}\big)\right) N_{Q/K}\big(\theta_1^{(Q)}\big) \left(\theta_1^{(Q)} - \iota_Q\big(\theta_1^{(Q)}\big); \iota_Q\big(\theta_1^{(Q)}\big); -\theta_1^{(Q)}\right) \\
  &= \pi^{4b_1} \cdot N_{Q/K}\left(\theta_1^{(Q)}\right) \cdot c_1\sqrt{D_0} \left(c_1\sqrt{D_0}; \iota_Q(\theta_1^{(Q)}); -\theta_1^{(Q)}\right),
\end{align*}
which is seen to have the same class in $H^1$ as
\[
  \left(1; c_1 \sqrt{D_0} \iota_Q\big(\theta_1^{(Q)}\big); -c_1 \sqrt{D_0} \theta_1^{(Q)}\right).
\]
The last two coordinates are $\iota_Q$-conjugate, as they should be.
Multiplying by
\[
  \diamondsuit = \left(1; \sqrt{D_0}; -\sqrt{D_0}\right) \textand \heartsuit = \big(1; \tr(\pi_Q); \tr(\pi_Q)\big),
\]
we get
\[
  [\hat\omega_C \diamondsuit \heartsuit] = \big(c_1 \tr \pi_Q; \iota_Q(\theta_1); \theta_1\big).
\]
So
\begin{align*}
  \ell(\hat\omega_C \diamondsuit \heartsuit) &= \ell_R(c_1 \tr \pi_Q; \bar\theta_1; \theta_1) \\
  &= \ell_Q\left(c_1 \tr \pi_Q \cdot \theta_1\right) \\
  &= \ell_Q\left(c_1 \cdot \tr \pi_Q \cdot c_0 \cdot \left(1 + \frac{c_1}{c_0} \pi_Q\right)\right).
\end{align*}
The last factor has exact level $s'$ if $s' < 2e$; otherwise it is a square. The remaining factors belong to $K$. By Lemma \ref{lem:iota}, their product is in $\cF_e$ if $d_0 = 2e+1$; otherwise it is in either $\cF_{d_0/2}$ or $\spadesuit \cF_{d_0/2}$ according to the parity of
\[
  v\left(c_1 \cdot \tr \pi_Q \cdot c_0\right) = s' + v(\tr \pi_Q).
\]
It is now easy to work out the level possibilities for each value of $s'$.
\end{proof}

In letter types \ref{type:A} and \ref{type:B}, as well as in the cases when $d_0 = 2e + 1$ is odd, we do not define boxgroups; all ring totals will be expressed in terms of the level spaces alone. Therefore in this subsection we will make the following assumptions:
\begin{itemize}
  \item $d_0 \leq 2e$ is even. We let $d_0 = 2d_0'$.
  \item $s' \geq d_0' - 1$, implying that $\hat\omega_C \diamondsuit\heartsuit \in \cF_{d_0' - 1}$ and $\square_C \geq 2\floor{(d_0'-1)/2} + 1$.
  \item $\theta_1 \in \OO_R$ is fixed. We may fix $\theta_2 = \pi^{-1/2}(0; \pi_Q)$ and $\bar\theta_2 = (0;0;1)$.
\end{itemize}
First we need an analogue of Lemma \ref{lem:eta_ur}:
\begin{lem}\label{lem:eta_1pow21}
Let
\[
\square_C = \min \left\{ 2\floor{\frac{\ell(\hat\omega_C \diamondsuit \heartsuit)}{2}} + 1, e \right\}
\]
be the squareness of the conic for $\delta = 1$. Assume we are in letter type \ref{type:C}, \ref{type:D}, or \ref{type:E} of Lemma \ref{lem:types_1pow21}, and let
\[
  h_\eta = \begin{cases}
    1 & \text{ in type \ref{type:C}} \\
    0 & \text{ in types \ref{type:D} and \ref{type:E}}.
  \end{cases}
\]
Then, if $\theta_1$ is translated by a suitable element of $\OO_K$ and scaled by a suitable element of $\OO_K^\cross$ (neither of which change the associated resolvent $C$), there is an $\eta \in \OO_R$ satisfying the congruence
\begin{equation} \label{eq:eta_main}
  \eta^2 \equiv \pi^{h_\eta}\theta_1 \mod \pi^{m_\eta + h_\eta}
\end{equation}
where
\[
  m_\eta = 2\ceil{\frac{s'}{2}} + \square_C - d_0'.
\]
Also,
\begin{equation} \label{eq:eta_in_1pow21}
  \eta \in \left\langle1, \theta_1, \pi^{\ceil{\frac{s' - d_0'}{2}} + \frac{1}{2}}\theta_2, \pi^{e - d_0' + \frac{1}{2}}\theta_2\right\rangle.
\end{equation}
\end{lem}
\begin{proof}
In letter type \ref{type:E}, we simply scale $\theta_1$ so that
\[
  \theta_1 \equiv (1;0) \mod \pi^{s'}
\]
and take $\eta = \theta_1$. This achieves the desired modulus $m_\eta = s'$, and \eqref{eq:eta_in_1pow21} is also clearly satisfied. We may therefore assume that $\theta_1$ is of letter type \ref{type:C} or \ref{type:D}.

Similar to the other splitting types, we view $\coef_{\theta_2} (\eta^2)$ as being a form of the conic $\cM$ and carry out the transformation of the conic, as detailed in Lemmas \ref{lem:gamma_white} and \ref{lem:tfm_conic}:
\begin{align*}
  \frac{1}{\sqrt{\pi}}\coef_{\theta_2}(\eta^2) &\sim \pi^{- d_0'} \coef_{\bar\theta_2}(\eta^2) \\
  &\sim \pi^{- d_0'} \tr(\omega_C \eta^2) \\
  &= \pi^{- d_0'} \lambda^\diamondsuit(\diamondsuit \omega_C \eta^2) \\
  &= \pi^{s - d_0'} \lambda^\diamondsuit(\diamondsuit \widetilde\omega_C \eta^2) \\
  &= \pi^{s'} \lambda^\diamondsuit(\diamondsuit \widetilde\omega_C \eta^2) \\
  &= \pi^{- \frac{h_1}{2}} \lambda^\diamondsuit(1^\odot {\eta^\odot}^2)\\
  &= \pi^{- \frac{h_1}{2}} \cM(\eta^{\odot}).
\end{align*}
Here $1^\odot \in\OO_R^\cross$ is the coefficient arising in Lemma \ref{lem:tfm_conic}, with
\[
  \frac{\widetilde\omega_C \diamondsuit}{1^\odot} \in (R^\cross)^2,
\]
and
\[
  \eta^{\odot} = \pi^{\ceil{s'/2}} \sqrt{\frac{\widetilde\omega_C \diamondsuit}{1^\odot}} \cdot \eta.
\]
Thus points $\eta^{\odot}$ with $\cM(\eta^{\odot}) \approx 0$ yield values $\eta$ such that $\coef_{\theta_2}(\eta^2) \approx 0$. The conic has a notable basepoint
\[
  \eta^{\odot}_0 = \pi^{\ceil{s'/2}} \sqrt{\frac{\widetilde\omega_C \diamondsuit}{1^\odot}} \sim \left(\pi^{\ceil{s'/2}}; \pi_Q^{h_1/2}\right),
\]
which corresponds to $\eta = 1$. It is evident that this is \emph{not} the $\eta$ we seek. However, it allows us to write the coordinate change succinctly as
\[
  \eta = \frac{\eta^{\odot}}{\eta^{\odot}_0}.
\]

By definition of $\square_C$, there is a linear form $\lambda$ such that
\[
\lambda^\diamondsuit(1^\odot{\eta^\odot}^2) \equiv c \cdot \lambda(\eta^{\odot})^2 \mod \pi^{\square_C}
\]
as functions of $\eta^{\odot} \in \OO_R$. Since $\cM(\eta^{\odot}_0) = 0$, there is no harm in picking $\lambda$ such that $\lambda(\eta^{\odot}_0) = 0$ also.

Our approach will be to complete to a basis $\{ \eta^{\odot}_0, \psi\}$ for $\ker \lambda$ and take $\eta^{\odot} = \pi^i \psi$, where $i$ is chosen such that $\eta = \eta^{\odot}/\eta^{\odot}_0$ is primitive in $\OO_R$.

As in the proof of Squareness Lemma \ref{lem:squareness}, we rescale $1^\odot$ by a scalar in $\OO_K^\cross$ and a square in $\left(\OO_R^\cross\right)^2$ so that the level $\ell(1^\odot \heartsuit)$ is manifest:
\[
1^\odot = \heartsuit (1; (1 + b\pi^{j+1}) + a\pi^j\pi_Q + 2\beta) = \left(1; \pi_Q^{2d_0'} \cdot \frac{1+b\pi^{j+1} + a\pi^j\pi_Q + 2\beta}{\tr \pi_Q}\right),
\]
where $j \in \{0, 1, 2,\ldots, e-1, \infty\}$ controls $\ell(1^\odot\heartsuit)$, and where $a \in \OO_K^\cross$, $b \in \OO_K$, $\beta \in \OO_R$.

Write $N_{Q/K}(\pi_Q)=u\pi$ with $u\in\OO_K^\cross$. We then looked at $\cM(\xi^\odot)$ where $\xi^\odot$ ranges over the basis of $\OO_R$
\[
  (1;0), \quad (0; \alpha_0) = \left(0;\pi_Q^{-d_0'}\pi^{\ceil{\frac{d_0'}{2}}}\right), \quad (0; \alpha_1) = \left(0;\pi_Q^{1-d_0'}\pi^{\floor{\frac{d_0'}{2}}}\right)
\]
and got
\begin{align}
  \cM(1;0) &= \lambda^\diamondsuit(1;0) \nonumber\\
  &= 1 \nonumber\\
  \cM\left(0;\alpha_0\right) &= \pi^{2\ceil{d_0'/2}}\lambda^\diamondsuit\left(0; \frac{1 + b\pi^{j+1} + a\pi^j\pi_Q + 2\beta}{\tr \pi_Q}\right) \nonumber\\
  &\equiv a \cdot \frac{\pi^{2\ceil{\frac{d_0'}{2}}}}{\tr \pi_Q} \cdot \pi^j \mod 2 \label{eq:alpha0} \\
  \cM\left(0;\alpha_1\right) &= \pi^{2\floor{d_0'/2}}\lambda^\diamondsuit\left(0; \pi_Q^2 \cdot \frac{1 + b\pi^{j+1} + a\pi^j\pi_Q + 2\beta}{\tr \pi_Q}\right) \nonumber \\
  &\equiv (1 + \pi^{j+1} b)\pi^{2\floor{\frac{d_0'}{2}}} - au \cdot \frac{\pi^{2\floor{\frac{d_0'}{2}} + 1}}{\tr \pi_Q} \cdot \pi^j + a\pi^{j + 2 \floor{\frac{d_0'}{2}}} \mod 2. \label{eq:bas_alpha1}
\end{align}
We may scale $\lambda$ so that $\lambda(1;0) = 1$. Then
\[
  \ker \lambda = \left\langle\left(-\lambda(0; \alpha_0); \alpha_0\right), \left(-\lambda(0; \alpha_1); \alpha_1\right)\right\rangle;
\]
we choose one of these elements for $\psi$. The details are slightly different in each letter type:
\begin{itemize}
  \item In letter type \ref{type:C}, we have $h_1/2 \equiv s' = d_0' - 1$ mod $2$, so
  \[
    \alpha_0 \sim \pi_Q^{1 - h_1/2} \textand \alpha_1 \sim \pi_Q^{h_1/2}.
  \]
  We can therefore take
  \[
    \psi = \left(-\lambda(0; \alpha_0); \alpha_0\right).
  \]
  Observe that
  \begin{align*}
    v(\cM(0; \alpha_0)) &\geq \min\left\{2\ceil{\frac{d_0'}{2}} - d_0' + j, e\right\} \\
    &= \min\left\{1 - \frac{h_1}{2} + j, e\right\} \\
    &\geq d_0' - \frac{h_1}{2} \\
    &= 2\floor{\frac{d_0' - 1}{2}} + 1.
  \end{align*}
  This lower bound is an odd integer, so we get
  \begin{align*}
    v(\lambda(0; \alpha_0)) &\geq \min\left\{\floor{\frac{d_0' - 1}{2}} + 1, \frac{\square_C + 1}{2}\right\} \\
    &= \floor{\frac{d_0' - 1}{2}} + 1 \\
    &= \floor{\frac{s'}{2}} + 1.
  \end{align*}
  This ensures that
  \[
  \eta^{\odot} = \pi^{h_1/2} \psi
  \]
  yields an $\eta = \eta^{\odot}/\eta^{\odot}_0$ with
  \[
    \pi | \eta_K \textand \eta_Q \sim \pi_Q.
  \]
  We get
  \begin{alignat*}{2}
    \frac{1}{\sqrt{\pi}}\coef_{\theta_2}(\eta^2) &= \pi^{h_1/2} \cM(\psi) \\
    &\equiv 0 &\mod{} &\pi^{h_1/2 + \square_C} \\
    &&={}&\pi^{ 2\ceil{s'/2} - s' + \square_C}\\
    &&={}&\pi^{m_\eta + h_\eta}.
  \end{alignat*}

  Hence $\eta^2$ is congruent modulo $\pi^{m_\eta + h_\eta + 1/2}$ to a value of the form $a + b\theta_1$, where $a,b \in \OO_K$. We have $v^{(K)}(\eta^2) \geq 2$, $v^{(Q)}(\eta^2) = 1$, so $\pi | a, \pi \sim b$. Thus $\eta^2 \equiv \pi\theta_1$ for suitably renormalized $\theta_1$. Note that \eqref{eq:eta_in_1pow21} is trivial, as the right-hand side is all of $\OO_R$. This completes the proof in letter type \ref{type:C}.

  \item In letter type \ref{type:D}, we have $h_1/2 \equiv s' = d_0'$ mod $2$, so
  \[
  \alpha_0 \sim \pi_Q^{h_1/2} \textand \alpha_1 \sim \pi_Q^{1 - h_1/2}.
  \]
  We can therefore take
  \[
  \psi = \left(-\lambda(0; \alpha_1); \alpha_1\right).
  \]
  We have
  \[
  v\left(\cM(0;\alpha_1)\right) = 2 \floor{\frac{d_0'}{2}}:
  \]
  the first term of \eqref{eq:bas_alpha1} dominates due to the bound on $j$ in letter type \ref{type:D}. Since $\square_C \geq 2\floor{d_0'/2} + 1$, we get
  \[
    v\left(\lambda(0;\alpha_1)\right) = \floor{\frac{d_0'}{2}}.
  \]
  Then we get a
  \[
    \psi = (-\lambda(0;\alpha_1); \alpha_1) \sim \left(\pi^{\floor{d_0'/2}}; \pi_Q^{1-h_1/2}\right)
  \]
  and take
  \[
    \eta^{\odot} = \pi^{\ceil{s'/2} - \floor{d_0'/2}} \psi
  \]
  for an $\eta = \eta^{\odot}/\eta^{\odot}_0$ whose $K$-component is a unit. As to the $Q$-component:
  \begin{align*}
    v(\eta^{(Q)}) &= \ceil{\frac{s'}{2}} - \floor{\frac{d_0'}{2}} + v(\psi^{(Q)}) \\
    &= \frac{2s' + h_1}{4} - \frac{2d_0' - h_1}{4} + \frac{2-h_1}{4}\\
    &\geq \frac{2 + h_1}{4}\\
    &\geq \frac{1}{2}.
  \end{align*}
  So $\eta$ is integral. We get
  \begin{alignat*}{2}
    \frac{1}{\sqrt{\pi}}\coef_{\theta_2}(\eta^2) &= \pi^{- \frac{h_1}{2} + 2\left(\ceil{s'/2} - \floor{d_0'/2}\right)} \cM(\psi) \\
    &\equiv 0 &\mod{} &\pi^{- \frac{h_1}{2} + 2\left(\ceil{s'/2} - \floor{d_0'/2}\right) + \square_C} \\
    &&={}&\pi^{ 2\ceil{s'/2} - d_0' + \square_C}\\
    &&={}&\pi^{m_\eta + h_\eta}.
  \end{alignat*}
  Hence $\eta^2$ is congruent modulo $\pi^{m_\eta + h_\eta + 1/2}$ to a value of the form $a + b\theta_1$, where $a,b \in \OO_K$. We have $\pi \nmid b$ because the $K$-component of $\eta^2$ has zero valuation, the $Q$-component positive valuation. So $\eta^2 \equiv \theta_1 \mod \pi^{m_\eta + h_\eta + 1/2}$ for suitably renormalized $\theta_1$, as desired.

  To prove \eqref{eq:eta_in_1pow21}, we observe that, since $\square_C > 0$,
  \[
    \left(\eta^{(Q1)}\right)^2 - \left(\eta^{(Q2)}\right)^2 \equiv \theta_1^{(Q1)} - \theta_1^{(Q2)} \equiv 0 \mod \pi^{s}.
  \]
  Thus
  \[
    \eta^{(Q1)} - \eta^{(Q2)} \equiv 0 \mod \pi^{\min\{s/2, e\}}
  \]
  (recalling that $s = s' + d_0'$), that is,
  \begin{align*}
    v(\fI(\eta^{(Q)})) &\geq \min\left\{\frac{s}{2}, e\right\} - d_0' \\
    &= \min \left\{\frac{s' - d_0'}{2}, e - d_0'\right\}.
  \end{align*}
  Since the left-hand side is an integer, we can take the ceiling and get
  \[
    \eta \in \left\langle1, (1;0), \pi^{\ceil{\frac{s' - d_0'}{2}}}(0; \pi_Q), \pi^{e - d_0'}(0; \pi_Q)\right\rangle.
  \]
  We have $(0; \pi_Q) = \sqrt{\pi}\theta_2$, and since
  \[
    \ceil{\frac{s' - d_0'}{2}} \leq s',
  \]
  we can replace $(1;0)$ by $\theta_1$ to get the desired relation
  \[
    \eta \in \left\langle1, \theta_1, \pi^{\ceil{\frac{s' - d_0'}{2}} + \frac{1}{2}}\theta_2, \pi^{e - d_0' + \frac{1}{2}}\theta_2\right\rangle. \qedhere
  \]
\end{itemize}
\end{proof}

We are now ready to construct boxgroups.
\begin{lem}\label{lem:boxes_basic_1pow21} Assume that $b_1, s' \in \ZZ$, so $\theta_1 \in \OO_R$ is in one of letter types \ref{type:B}--\ref{type:E} of Lemma \ref{lem:types_1pow21}. Let $ m, n \in \NN^+ \union \{\infty\} $, $ m \geq n > 0 $. Let $ B_{\theta_1}(m,n) $ be the box
  \[
  B_{\theta_1}(m,n) = \{\pi c_0 + \pi^{n}c_1\theta_1 + \pi^{m + \frac{1}{2}} c_2 \theta_2 : c_i \in \OO_K \}.
  \]
  \begin{enumerate}[$($a$)$]
    \item\label{boxes_basic_1pow21:recenter} For every $ \xi \in 1 + B_{\theta_1}(m,n) $,
    \begin{align}
      B_{\theta_1}(m,n) &= \xi B_{\xi^{-1} {\theta_1}}(m,n) \label{eq:1pow21_recenter_1} \\
      1 + B_{\theta_1}(m,n) &= \xi \left(1 + B_{\xi^{-1} {\theta_1}}(m,n)\right) \label{eq:1pow21_recenter_2}.
    \end{align}
    \item\label{boxes_basic_1pow21:group} If
    \begin{equation*}
      m \leq 2n + s',
    \end{equation*}
    then $ B_{\theta_1}(m,n) $ is closed under multiplication and the translated box $ 1 + B_{\theta_1}(m,n) $ is a group under multiplication.
  \end{enumerate}
\end{lem}
\begin{rem}
  We will eventually apply this result with
\[
  m = m_\beta - \frac{1}{2} = m_{11} - \frac{d_0}{2}, \quad n = n_{11} - s
\]
of a first vector problem.
\end{rem}
\begin{proof}
The proof proceeds just as in the unramified splitting types (Lemma \ref{lem:boxes_basic_ur}). For \ref{boxes_basic_1pow21:group}, the key is to note that, in the basis $\{1, \theta_1, \sqrt{\pi} \theta_2\}$ for $\OO_R$,
\[
  v\left(\coef_{\sqrt{\pi}\theta_2}(\theta_1^2)\right) = s',
\]
a useful reinterpretation of the (renormalized) idempotency index $s'$.
\end{proof}
\begin{lem}\label{lem:boxgps_1pow21}
Assume that $d_0$ is even and $\theta_1$ is of one of the letter types \ref{type:C}--\ref{type:E} of Lemma \ref{lem:types_1pow21}. Let $0 \leq n \leq m \leq 2e - d_0'$ be integers such that
\begin{align}
m &\geq n + d_0' - 1 \label{eq:gray_domino} \\
m &\leq 2n + s' &&(\text{the gray-red inequality}) \\
m &\leq n + 2\ceil{\frac{s'}{2}} - d_0' + \square_C &&(\text{the gray-green inequality}) \\
m &\leq e + \frac{n + s' - d_0' + 1}{2} &&(\text{the gray-blue inequality})
\end{align}
Then the projection $[1 + B_{\theta_1}(m,n)]$ of $1 + B_{\theta_1}(m,n)$ onto $H^1$ is a subgroup of signature
\begin{align*}
  &\ZERO .\ZERO^{d_0'}.(\ZERO\ZERO)^{\ell_0} (\ZERO\STAR)^{\ell_1} (\STAR\STAR)^{\ell_2} .\STAR^{d_0'}. \STAR \quad (\text{type \ref{type:C}}) \\
  &\ZERO .\ZERO^{d_0'}.(\ZERO\ZERO)^{\ell_0} (\STAR\ZERO)^{\ell_1} (\STAR\STAR)^{\ell_2} .\STAR^{d_0'}. \STAR\quad (\text{types \ref{type:D}--\ref{type:E}})
\end{align*}
where
\begin{align*}
  \ell_0 &= \floor{\frac{n - h_\eta}{2}} \\
  \ell_1 &= \floor{\frac{m - d_0' + h_\eta}{2}} - \floor{\frac{n - h_\eta}{2}} \\
  \ell_2 &= e - \floor{\frac{m + d_0' + h_\eta}{2}}.
\end{align*}
\end{lem}
\begin{proof}
The gray-red inequality ensures that
\[
  T = [1 + B_{\theta_1}(m, n)] = \{[1 + c_1 \pi^{n} \theta_1 + c_2 \pi^{m} (0; \pi_Q)]\}
\]
is a subgroup. As in the other splitting types, we will establish its signature by downward induction on $n$, fixing $m$.

The base case is when \eqref{eq:gray_domino} is an equality. Here $\ell_1 \in \{0,1,-1\}$ is such that the signature simplifies, and we wish to prove that
\[
  T = \cF_{m}.
\]
The $\supseteq$ direction is straightforward:
\[
  T \supseteq \{[\alpha] : \alpha \equiv 1 \mod \pi^{m + 1/2}\} = \cF_{m}.
\]
For the $\subseteq$ direction, we simply need to show that the other generators $1 + c_1\pi^{n} \theta_1$ do not add any new elements to the box. Since $s' \geq d_0' - 1$ (we are in letter types \ref{type:C}--\ref{type:E}), we can scale and translate $\theta_1$ so that
\[
  \theta_1 \equiv (1; 0; 0) \mod \pi^{d_0' - 1/2}.
\]
Then
\[
  [1 + c_1\pi^{n} \theta_1] \equiv [(1 + c_1\pi^{n} ; 1 ; 1)] = \iota(1 + c_1\pi^{n}) \mod \cF_{m}.
\]
The element $1 + c_1\pi^{n}$ has $K$-level at least $\floor{n/2}$ and, by Lemma \ref{lem:iota}\ref{iota:custom}, $Q$-level at least $n + d_0' - 1$, getting the needed bound $\floor{n/2} \leq e - d_0'$ from the bound $m \leq 2e - d_0'$. This completes the proof of the base case.

For the induction step, assume that $[1 + B_{\theta_1}(m, n + 1)]$ has the correct signature. Note that decreasing $n$ by $1$ enlarges the boxgroup at most by a factor of $q$. There are two cases.

\subsubsection{Moving case.} The first case is when
\begin{itemize}
  \item $n$ is odd in letter types \ref{type:D} and \ref{type:E}
  \item $n$ is even in letter type \ref{type:C}.
\end{itemize}
Here we claim that the boxgroup actually does grow by a factor of $q$, and in particular that the new generators $1 + c_1\pi^{n}\theta_1$, $\pi \nmid c_1$, are all of level $n + d_0' - 1$. The method depends slightly on the letter type.

In letter types \ref{type:D} and \ref{type:E}, we have $s' \geq d_0'$ and thus can normalize $\theta_1$ so that
\[
  \theta_1 \equiv (1;0;0) \mod \pi^{d_0' + 1/2}.
\]
Then, for $n$ odd,
\[
[1 + c_1\pi^{n}\theta_1] \equiv [(1 + c_1\pi^{n}; 1; 1)] = \iota(1 + c_1\pi^{n}) \mod \cF_{n + d_0'}.
\]
The element $1 + c_1\pi^{n}$ has $K$-level $(n - 1)/2$. Plugging into Lemma \ref{lem:iota}, we use $n \leq 2e - 2d_0'$ (from $m \leq 2e - d_0'$ and the fact that we are beyond the base case) to nail down the $Q$-level
\[
  \ell\left(\iota(1 + c_1\pi^{n})\right) = 2\left(\frac{n - 1}{2}\right) + d_0' = n + d_0' - 1.
\]

In letter type \ref{type:C}, we have $s' = d_0' - 1$. We may normalize so that
\[
  \theta_1 = \left(1; \pi^{d_0' - 1}u\pi_Q\right), \quad u \in \OO_K^\cross.
\]
Then, for $n$ even,
\begin{align*}
\left[1 + c_1\pi^{n}\theta_1\right] &= \left[\left(1 + c_1\pi^{n}; 1 + c_1\pi^{n + d_0' - 1}u\pi_Q)\right)\right] \\
&= \iota(1 + c_1\pi^{n}) \cdot \left[\left(1; 1 + c_1\pi^{n + d_0' - 1}u\pi_Q)\right)\right].
\end{align*}
The second factor has the desired level $n + d_0' - 1$. As to the first factor, an element of $K$-level at least $n/2$ is fed to $\iota$, yielding $Q$-level at least $n + d_0'$, using the inequality $n \leq 2e - 2d_0'$ again.

As a consequence, the boxgroup grows by a factor of $q$, entirely in level $n + d_0' - 1$: so we get its signature by replacing the $\ZERO$ in that slot by a $\STAR$.

\subsubsection{Stationary case.} We now come to the second and subtler case of the induction step, when
\begin{itemize}
  \item $n$ is even in letter types \ref{type:D} and \ref{type:E}
  \item $n$ is odd in letter type \ref{type:C}.
\end{itemize}
Here our aim is to prove that the boxgroup is unchanged, that is, that the new generators $1 + c_1\pi^{n}\theta_1$ actually lie in the existing group $[1 + B_{\theta_1}(m, n + 1)]$. It suffices to show that there is a square in each of the $q$ cosets
\[
  (1 + c_1^2\pi^{n}\theta_1)\left(1 + B_{\theta_1}(m, n + 1)\right) = 1 + c_1^2\pi^{n}\theta_1 + B_{\theta_1}(m, n + 1), \quad c_1 \in \OO_K
\]
(only $c_1$ mod $\pi$ matters; the two forms of the coset are equal by the gray-red inequality). We produce these squares using the $\eta$-lemma (Lemma \ref{lem:eta_1pow21}).

In letter types \ref{type:D} and \ref{type:E}, we have an $\eta \in \OO_R$ such that (renormalizing $\theta_1$ appropriately)
\[
  \eta^2 \equiv \theta_1 \mod \pi^{m_\eta + 1/2}.
\]
Given $n$ even, we claim that the square
\[
  \left(1 + \pi^{n/2} c_1 \eta\right)^2 = 1 + 2\pi^{n/2} c_1 \eta + \pi^{n} c_1^2 \eta^2
\]
lies in the proper coset, that is,
\[
  1 + 2\pi^{n/2} c_1 \eta + \pi^{n} c_1^2 \eta^2 \equiv 1 + \pi^{n} c_1^2\theta_1 \mod B_{\theta_1}(m, n + 1).
\]
This requires two considerations:
\begin{itemize}
  \item To replace $\eta^2$ by $\theta_1$, we need
  \[
    m_\eta + n \geq m,
  \]
  which is equivalent to the gray-green inequality.
  \item For the cross term, Lemma \ref{lem:eta_1pow21} gives directly
  \[
    2\pi^{n/2} c_1 \eta \in B_{\theta_1}\left(e + \frac{n}{2} + \ceil{\frac{s' - d_0'}{2}}, e + \frac{n}{2}\right) \subseteq B_{\theta_1}(m, n + 1)
  \]
  by the gray-blue inequality, the difference of whose sides is in $\ZZ + 1/2$ in letter type \ref{type:D}. (In letter type \ref{type:E}, since $\eta = \theta_1$, this step is even easier.)
\end{itemize}

Similarly, in letter type \ref{type:C}, we have an $\eta \in \OO_R$ such that (renormalizing $\theta_1$ appropriately)
\[
\eta^2 \equiv \pi\theta_1 \mod \pi^{m_\eta + 1 + 1/2}.
\]
Given $n$ odd, we claim that the square
\[
  \left(1 + \pi^{\frac{n - 1}{2}} c_1 \eta\right)^2 = 1 + 2\pi^{\frac{n - 1}{2}} c_1 \eta + \pi^{n - 1} c_1^2 \eta^2
\]
lies in the proper coset, that is,
\[
  1 + 2\pi^{\frac{n - 1}{2}} c_1 \eta + \pi^{n - 1} c_1^2 \eta^2 \equiv 1 + \pi^{n} c_1^2 \theta_1 \mod B_{\theta_1}(m, n + 1).
\]
This requires two considerations:
\begin{itemize}
  \item To replace $\eta^2$ by $\pi \theta_1$, we need
  \[
    m_\eta + n \geq m
  \]
  (the term $1 = h_\eta$ cancels), which is equivalent to the gray-green inequality.
  \item For the cross term, since $\eta \equiv 0 \mod \sqrt{\pi}$, we have
  \[
    2\pi^{(n-1)/2} c_1 \eta \equiv 0 \mod \pi^{e + \frac{n - 1}{2} + \frac{1}{2}},
  \]
  the exponent being at least $m + 1/2$ by the gray-blue inequality (the difference of whose sides is in $\ZZ + 1/2$). \qedhere
\end{itemize}
\end{proof}

\begin{lem}\label{lem:boxgpT_1pow21}
  For every triple $ (\ell_0, \ell_1, \ell_2) $ of integers satisfying
  \begin{align}
    \ell_0 + \ell_1 + \ell_2 &= e - d_0' \\
    \ell_1 &\leq \ell_0 + \ceil{\frac{s' - d_0'}{2}} + 1 + h_\eta && (\text{the gray-red inequality})\label{eq:boxgp_red_1pow21} \\
    \ell_1 &\leq \frac{\square_C + 1}{2} + \ceil{\frac{s'}{2}} - d_0' + h_\eta && (\text{the gray-green inequality}) \label{eq:boxgp_green_1pow21} \\ 
    \ell_1 &\leq \ell_2 + \ceil{\frac{s' - d_0'}{2}} + 1 + h_\eta, && (\text{the gray-blue inequality})\label{eq:boxgp_blue_1pow21} \\
    \ell_0 &\geq 0, \quad \ell_1 \geq -1 + h_\eta, \quad \ell_2 \geq 0, \nonumber
  \end{align}
  there is a \emph{boxgroup} $T_{h_\eta}(\ell_0,\ell_1,\ell_2) \subseteq H^1$ of signature
  \begin{align*}
    &\ZERO .\ZERO^{d_0'}.(\ZERO\ZERO)^{\ell_0} (\ZERO\STAR)^{\ell_1} (\STAR\STAR)^{\ell_2} .\STAR^{d_0'}. \STAR \quad (\text{type \ref{type:C}}) \\
    &\ZERO .\ZERO^{d_0'}.(\ZERO\ZERO)^{\ell_0} (\STAR\ZERO)^{\ell_1} (\STAR\STAR)^{\ell_2} .\STAR^{d_0'}. \STAR\quad (\text{types \ref{type:D}--\ref{type:E}})
  \end{align*}
  such that, if $m$, $n$ are integers satisfying the conditions of Lemma \ref{lem:boxgps_1pow21}, then
  \[
  [1 + B_{\theta_1}(m,n)] = T_{h_\eta}\bigg( \Big\lfloor{\frac{n - h_\eta}{2}}\Big\rfloor, \Big\lfloor{\frac{m - d_0' + h_\eta}{2}}\Big\rfloor - \Big\lfloor{\frac{n - h_\eta}{2}}\Big\rfloor, e - \Big\lfloor{\frac{m + d_0' + h_\eta}{2}}\Big\rfloor
   \bigg).
  \]
\end{lem}
\begin{proof}
  If $\ell_1 = -1 + h_\eta$, take the appropriate level space, the unique subgroup with the correct signature.

  Otherwise, let
  \[
    m = 2e - 2\ell_2 - d_0' - h_\eta, \quad n = 2\ell_0 + h_\eta + 1
  \]
  in the preceding lemma. When we transform the conditions on $m$ and $n$ to conditions on the $\ell_i$, they become a priori
  \begin{align*}
    \ell_1 &\leq \ell_0 + \left(\frac{s' - d_0' + h_\eta}{2}\right) + 1 + h_\eta && (\text{the gray-red inequality}) \\
    \ell_1 &\leq \frac{\square_C + 1}{2} + \ceil{\frac{s'}{2}} - d_0' + h_\eta && (\text{the gray-green inequality}) \\
    \ell_1 &\leq \ell_2 + \left(\frac{s' - d_0' + h_\eta}{2}\right) + 1 + h_\eta, && (\text{the gray-blue inequality}) \\
  \end{align*}

  In letter types \ref{type:C} and \ref{type:D}, the gray-red and gray-blue inequalities are then put into the desired form using the identity
  \[
    \floor{\frac{s' - d_0' + h_\eta}{2}} = \ceil{\frac{s' - d_0'}{2}},
  \]
  which is easily proved using the interaction between $h_\eta$ and the parity of $s' - d_0'$ in letter types \ref{type:C}, \ref{type:D}. In letter type \ref{type:E} (which corresponds to $s$ odd in the unramified splitting types), the gray-red and gray-blue inequalities are subsumed by the gray-green.

  If the $\ell_i$ are fixed, let $(m_0, n_0)$ be the pair
  \[
    m_0 = 2e - 2\ell_2 - d_0' - h_\eta, \quad n_0 = 2\ell_0 + h_\eta + 1
  \]
  used in the proof of the lemma.

  The other pairs $(m, n)$ yielding the same triple $(\ell_0, \ell_1, \ell_2)$ are
  \[
    (m, n) = (m_0 + 1, n_0), \quad (m_0, n_0 - 1), \textor (m_0 + 1, n_0 - 1).
  \]

  Except for the base-case inequality, which is satisfied as long as $\ell_1 > -1 + h_\eta$, the conditions on $m$ and $n$ only become truer upon decreasing $m$ or increasing $n$. Now
  \[
  [1 + B_{\theta_1}(m,n)] \subseteq \left[1 + B_{\theta_1}\left(m_0, n\right)\right],
  \]
  but both sides have the same signature, so equality holds. Likewise,
  \[
  [1 + B_{\theta_1}(m_0,n)] \supseteq \left[1 + B_{\theta_1}\left(m_0, n_0\right)\right],
  \]
  but both sides have the same signature, so equality holds.
\end{proof}

\begin{nota}
The notation $T(g_0; g_1)$ will simplify the computer step. Note that $T(g_0; 0)$ (and also $T(g_0; -1)$, which inconveniently arises from transforming $T_{1}(\ell_0, 0, \ell_2)$) is the level space $\cF_{g_0}$; we extend the notation
\[
  T(g_0; 0) = \cF_{g_0}
\]
even to values $g_0 < d_0'$ or $g_0 > 2e - d_0'$, for which Lemma \ref{lem:boxgpT_1pow21} and \eqref{eq:l_g_1pow21} do not apply. Accordingly, the $(g_0, g_1)$ for which $T(g_0; g_1)$ is defined forms a nonconvex ``hat'' pictured in Figure \ref{fig:boxgp_hat}.
\begin{figure}
  \begin{center}
  \setlength{\unitlength}{0.1in}
  \begin{picture}(48,7)(-3,-2)%
    \multiput(0,0)(1,0){43}{\circle*{1}}
    \multiput(6,1)(2,0){16}{\circle*{1}}
    \multiput(7,2)(2,0){15}{\circle*{1}}
    \multiput(8,3)(2,0){14}{\circle*{1}}
    \multiput(11,4)(2,0){11}{\circle*{1}}
    \multiput(14,5)(2,0){8}{\circle*{1}}
    \put(5,0){\line(1,1){3}}
    \put(8,3){\line(3,1){6}}
    \put(14,5){\line(1,0){14}}
    \put(28,5){\line(3,-1){6}}
    \put(34,3){\line(1,-1){3}}
    \put(-3,0){\vector(1,0){48}}
    \put(43,0){\makebox(2,2)[r]{$g_0$}}
    \put(0,-2){\vector(0,1){7}}
    \put(0,3){\makebox(2,2)[t]{$g_1$}}
  \end{picture}
  \end{center}
\caption{A ``boxgroup hat'' showing the boxgroups $T(g_0; g_1)$, in the $g$-coordinates, that are defined when $e = 16$, $d_0 = 10$, $s' = 9$, $\ell_C = 4$ (letter type D).}
\label{fig:boxgp_hat}
\end{figure}
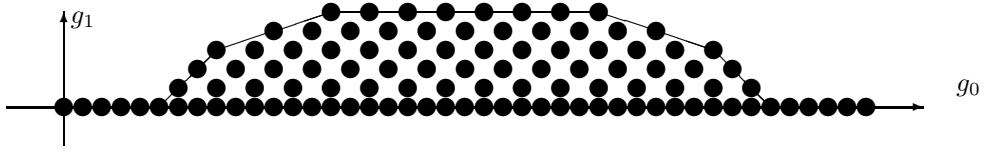
\end{nota}

\subsection{Supplementary boxgroups}\label{sec:suppl_1pow21}
The boxgroups constructed thus far fit in the region
\[
  \cF_{d_0/2} \supseteq T \supseteq \cF_{2e - d_0'}.
\]
Outside this region, there is a less rich variety of groups appearing in the ring totals. They can be constructed using the following:
\begin{lem} \label{lem:iota_image}
Consider the map
\begin{align*}
  \iota\colon K^\cross/(K^\cross)^2 &\to H^1 \\
  a &\mapsto [(a;1)] = [(1;a)].
\end{align*}
\begin{enumerate}
  \item If $d_0 = 2e+1$, then $\iota(K^\cross) = \cF_e$.
  \item If $d_0 = 2d_0'$ is even, then the image $\iota(K^\cross)$ is a group of signature
  \[
  \ZERO.\ZERO^{d_0'-1}\spadesuit.(\STAR\ZERO)^{e - d_0'}.\spadesuit^{\perp}\STAR^{d_0'-1}.\STAR
  \]
  where $\spadesuit \subseteq \cF_{d_0' - 1}/\cF_{d_0'}$ is the $2$-element subgroup generated by $\iota(\pi)$, and $\spadesuit^{\perp} \subseteq \cF_{2e - d_0'}/\cF_{2e - d_0' + 1}$ is its orthogonal complement.
\end{enumerate}
\end{lem}
\begin{proof}
Use Lemma \ref{lem:iota} to determine the size of the signature component at each level. Since $\iota(K^\cross)$ is self-orthogonal, we get the term $\spadesuit^{\perp}$ at level $2e - d_0'$.
\end{proof}

Now assume that $d_0 = 2d_0'$ is even. The group $\iota(K^\cross)$ is always important, but it does not behave well with respect to boxgroups unless $s' > 2e - d_0$, in which case we give it the name $T(\emptyset, e - d_0', \emptyset)$.

If $0 \leq \ell_1 \leq e - d_0'$ and $s' \geq d_0' + 2\ell_1$ (so we are in letter type \ref{type:D} or \ref{type:E}), we let
\begin{align*}
  T(e - d_0' - \ell_1, \ell_1, \emptyset) &= \iota(K^\cross) \intsec \cF_{2e - d_0' - 2\ell_1} = \iota(1 + \pi^{2(e - d_0' - \ell_1) + 1}\OO_K) \\
  T(\emptyset, \ell_1, e - d_0' - \ell_1) &= \iota(K^\cross) \cdot \cF_{d_0' + 2\ell_1}.
\end{align*}
(In the first case, the equivalence of the two definitions is established using Lemma \ref{lem:iota}.) Their signatures are, respectively,
\begin{gather*}
  \ZERO.\ZERO^{d_0'}.(\ZERO\ZERO)^{e - d_0' - \ell_1}(\STAR\ZERO)^{\ell_1}.\spadesuit^\perp\STAR^{d_0' - 1}.\STAR
  \intertext{and}
  \ZERO.\ZERO^{d_0' - 1}\spadesuit . (\STAR\ZERO)^{\ell_1} (\STAR\STAR)^{e - d_0' - \ell_1} . \STAR^{d_0'}.\STAR.
\end{gather*}

The restrictions on $s'$ ensure that these boxgroups satisfy such natural relations as
\begin{align*}
  T(e - d_0' - \ell_1, \ell_1, \emptyset) \cdot \cF_{2e - d_0'} &= T(e - d_0' - \ell_1, \ell_1, 0) \\
  T(\emptyset, \ell_1, e - d_0' - \ell_1) \intsec \cF_{d_0'} &= T(0, \ell_1, e - d_0' - \ell_1),
\end{align*}
which we will often use without comment.

It will turn out that $T(\ell_0,\ell_1,\ell_2)$ and $T(\ell_2,\ell_1,\ell_0)$ are orthogonal complements whenever both are defined. Actually, this is simple to prove in the case that one of $\ell_0,\ell_1,\ell_2$ is the symbol $\emptyset$. The sizes of these groups follow immediately from their signatures:
\begin{lem}\label{lem:1pow21_T_size}
  If $T(\ell_0,\ell_1,\ell_2)$ is defined and $\ell_1 \neq \emptyset$, then
  \[
  \size{T(\ell_0,\ell_1,\ell_2)} = \size{H^0}q^{e + \ell_2 - \ell_0},
  \]
  where if $\emptyset$ occurs as either $\ell_0$ or $\ell_2$, it must be replaced by $-1/[k_K : \FF_2] = \log_q(1/2)$.
\end{lem}

\section{The recentering lemma}
When ${\theta_1}$ is nearly a square, we will sometimes be able to assume that it \emph{is} a square, thanks to the following lemma, which we state separately for each splitting type, although the proof is nearly identical.
\begin{lem}\label{lem:recentering_ur}
  If $C$ is unramified and $s$ is even, then there is an element $ \psi_C \in 1 + B_{\theta_1}(\infty, \square_C) $ and $a \in \OO_K$, $b \in \OO_K^\cross$ such that $ \theta_1' = a + b\psi_C^{-1}\theta_1 = \eta^2 $ is a square in $\OO_R$.
\end{lem}
\begin{proof}
  It is easy to see that tweaking $\theta_1$ by addends in $\OO_K$ or multipliers in $\OO_K^\cross$, which do not change the underlying resolvent ring $C$, do not affect the truth of the lemma either. We may therefore assume, by Lemma \ref{lem:eta_ur}, that there is an $\eta$ such that
  \[
  \eta^2 \equiv \theta_1 \mod \pi^{s + \square_C}.
  \]
  Let $k = \square_C$. In particular, $ \eta^2 \equiv \theta_1 $ mod $ \pi^{k} \OO_K[\theta_1] $. Notice that $ \theta_1' $ can take all values in the orbit
  \[
  \theta_1' = \frac{a_{11} \theta_1 + a_{12}}{a_{21} \theta_1 + a_{22}} = \frac{a_{12}}{a_{22}} + \frac{(a_{11}a_{22} - a_{12}a_{21})\theta_1}{a_{22}(a_{21}\theta_1 + a_{22})}
  \]
  of $ \theta_1 $ under the congruence subgroup
  \[
  \Gamma(k) = \left\{\gamma =
  \begin{bmatrix}
  a_{11} & a_{12} \\ a_{21} & a_{22}
  \end{bmatrix} \in \GL_2(\OO_K): \gamma \equiv
  \begin{bmatrix}
  1 & 0 \\ 0 & 1
  \end{bmatrix}
  \mod \pi^k \right\}.
  \]
  We claim that this orbit is precisely the pixel $ \theta_1 + \pi^k \OO_K[\theta_1] $. The orbit is clearly contained in this pixel and contains all elements of the form
  \begin{align*}
  \gamma(\theta_1) &= \frac{(1 + a_{11}'\pi^k) \theta_1 + a_{12}'\pi^k}{a_{21}' \pi^k \theta_1 + (1 + a_{22}'\pi^k)} \\
  &\equiv \theta_1 + \pi^k(a_{12}' + (a_{11}' - a_{22}')\theta_1 + a_{21}'\theta_1^2) \mod \pi^{2 k} \OO_K[\theta_1].
  \end{align*}
  So at least the orbit contains a point in each congruence class mod $ \pi^{2 k} \OO_K[\theta_1] $ in the claimed pixel. But applying general elements of $ \Gamma(2 k) $ to each of those, we get a point in each congruence class mod $ \pi^{4 k} \OO_K[\theta_1] $, and so on. Hence, the orbit is dense in the pixel, and being compact, it coincides with the pixel, establishing the desired result.
\end{proof}
\begin{lem}\label{lem:recentering_1pow3}
  If $R$ is of splitting type $(1^3)$, there is an element $ \psi_C \in 1 + B_{\theta_1}(\infty, \square_C - 2h/3) $ and $a \in \pi^{-h/3} K \intsec \OO_{\bar K}$, $b \in \OO_K^\cross$ such that $ \theta_1' = a + b\psi_C^{-1}\theta_1 = \eta^2 $ is a square, $\eta \in \bar{\zeta}_3^{-h} \OO_R$.
\end{lem}
\begin{proof}
  It is easy to see that tweaking $\theta_1$ by addends in $\pi^{-h/3}\OO_K$ or multipliers in $\OO_K^\cross$, which do not change the underlying resolvent ring $C$, do not affect the truth of the lemma either. We may therefore assume, by Lemma \ref{lem:eta_1pow3}, that there is an $\eta$ such that
  \[
  \eta^2 \equiv \theta_1 \mod \pi^{k}
  \]
  where $k = \square_C - 2h/3$. Notice that $ \theta_1' $ can take all values in the orbit
  \[
  \theta_1' = \frac{a_{11} \theta_1 + a_{12}}{a_{21} \theta_1 + a_{22}} = \frac{a_{12}}{a_{22}} + \frac{(a_{11} a_{22} - a_{12} a_{21})\theta_1}{a_{22}(a_{21}\theta_1 + a_{22})}
  \]
  of $ \theta_1 $ under the congruence subgroup
  \[
  \Gamma(k) = \left\{\gamma =
  \begin{bmatrix}
  a_{11} & a_{12} \\ a_{21} & a_{22}
  \end{bmatrix} \in \GL_2(\OO_{\bar K}): a_{ij} \in \pi^{\frac{h}{3}(i - j)} K, \gamma \equiv
  \begin{bmatrix}
  1 & 0 \\ 0 & 1
  \end{bmatrix}
  \mod \pi^k \right\}.
  \]
  We claim that this orbit is precisely the pixel
  \[
  \{ \theta_1' \in \bar{\zeta}_3^h \OO_R^\cross : \theta_1' \equiv \theta_1 \mod \pi^k \}.
  \]
  The orbit is clearly contained in this pixel and contains all elements of the form
  \begin{align*}
    \gamma(\theta_1) &= \frac{(1 + a_{11}'\pi^k) \theta_1 + a_{12}'\pi^k}{a_{21}' \pi^k \theta_1 + (1 + a_{22}'\pi^k)} \\
    &\equiv \theta_1 + \pi^k(a_{12}' + (a_{11}' - a_{22}')\theta_1 + a_{21}'\theta_1^2) \mod \pi^{2 k} \OO_R,
  \end{align*}
  where the $a'_{ij}$ are integral in the appropriate groups $\pi^{\frac{h}{3}(i-j) - k} K$.

  So at least the orbit contains a point in each congruence class mod $ \pi^{2 k} \OO_K[\theta_1] $ in the claimed pixel. But applying general elements of $ \Gamma(2 k) $ to each of those, we get a point in each congruence class mod $ \pi^{4 k} \OO_K[\theta_1] $, and so on. Hence, the orbit is dense in the pixel, and being compact, it coincides with the pixel, establishing the desired result.
\end{proof}
  \begin{lem}\label{lem:recentering_1pow21}
  If $C$ has splitting type $(1^21)$ and letter type \ref{type:C} or \ref{type:D}, then there is an element $ \psi_C \in 1 + B_{\theta_1}(\infty, \square_C - d_0' + h_1/2) $ and $a \in \OO_K$, $b \in \OO_K^\cross$ such that $ \theta_1' = a + b\psi_C^{-1}\theta_1 $ has the property that $\pi^{h_\eta} \theta_1' = \eta^2 $ is a square in $\OO_R$.
\end{lem}
\begin{proof}
  It is easy to see that tweaking $\theta_1$ by addends in $\OO_K$ or multipliers in $\OO_K^\cross$, which do not change the underlying resolvent ring $C$, do not affect the truth of the lemma either. We may therefore assume, by Lemma \ref{lem:eta_1pow21}, that there is an $\eta$ such that
  \[
  \frac{\eta^2}{\pi^{h_\eta}} \equiv \theta_1 \mod \pi^{s' + \frac{1}{2} + k}
  \]
  where $k = m_\eta - s' = \square_C - d_0' + h_1/2$.
  In particular, $ \eta^2/\pi^{h_\eta} \equiv \theta_1 $ mod $ \pi^{k} \OO_K[\theta_1] $. Notice that $ \theta_1' $ can take all values in the orbit
  \[
  \theta_1' = \frac{a_{11} \theta_1 + a_{12}}{a_{21} \theta_1 + a_{22}} = \frac{a_{12}}{a_{22}} + \frac{(a_{11}a_{22} - a_{12}a_{21})\theta_1}{a_{22}(a_{21}\theta_1 + a_{22})}
  \]
  of $ \theta_1 $ under the congruence subgroup
  \[
  \Gamma(k) = \left\{\gamma =
  \begin{bmatrix}
    a_{11} & a_{12} \\ a_{21} & a_{22}
  \end{bmatrix} \in \GL_2(\OO_K): \gamma \equiv
  \begin{bmatrix}
    1 & 0 \\ 0 & 1
  \end{bmatrix}
  \mod \pi^k \right\}.
  \]
  We claim that this orbit is precisely the pixel $ \theta_1 + \pi^k \OO_K[\theta_1] $. The orbit is clearly contained in this pixel and contains all elements of the form
 \begin{align*}
   \gamma(\theta_1) &= \frac{(1 + a_{11}'\pi^k) \theta_1 + a_{12}'\pi^k}{a_{21}' \pi^k \theta_1 + (1 + a_{22}'\pi^k)} \\
   &\equiv \theta_1 + \pi^k(a_{12}' + (a_{11}' - a_{22}')\theta_1 + a_{21}'\theta_1^2) \mod \pi^{2 k} \OO_K[\theta_1].
 \end{align*}
 So at least the orbit contains a point in each congruence class mod $ \pi^{2 k} \OO_K[\theta_1] $ in the claimed pixel. But applying general elements of $ \Gamma(2 k) $ to each of those, we get a point in each congruence class mod $ \pi^{4 k} \OO_K[\theta_1] $, and so on. Hence, the orbit is dense in the pixel, and being compact, it coincides with the pixel, establishing the desired result.
\end{proof}

\section{Charmed cosets}
\begin{lem} \label{lem:charm}
Let $H$ be a finite $2$-torsion group, and let $\epsilon : H \to \{\pm 1\}$ be a nondegenerate quadratic form over $\FF_2$. (By ``nondegenerate'' we mean that the associated bilinear form is a nondegenerate symplectic form on $H$; in particular, $\dim_{\FF_2} H$ is even.) Let $W \subseteq H$ be a subspace that is \emph{coisotropic;} that is, $W^\perp$ is isotropic, or equivalently, $W$ contains a maximal isotropic subspace. Then:
\begin{enumerate}[$($a$)$]
  \item \label{charm:exist} There is exactly one coset $\alpha W \subseteq H$ such that
\[
  \sum_{\delta \in \alpha W} \epsilon(\delta) \neq 0.
\]
Indeed, the sum is
\[
  \pm_{\epsilon} \sqrt{\size{H}},
\]
where $\pm_{\epsilon} \in \{\pm 1\}$ is an invariant of the quadratic space $(H, \epsilon)$. We call $\alpha W$ the \emph{charmed coset} of $W$.
\item \label{charm:constant} On any coset $\beta W^\perp$ inside the charmed coset $\alpha W$, $\epsilon$ is constant. By contrast, on any coset $\beta W^\perp$ outside the charmed coset, $\epsilon$ is equidistributed.
\end{enumerate}
\end{lem}
\begin{proof}
\begin{enumerate}[$($a$)$]
\item Assume first that $W$ is maximal isotropic. Then $\size{W} = \sqrt{\size{H}}$. On each coset $\alpha W$, $\epsilon$ looks like a linear form, that is, there is a $\lambda_\alpha \in W^*$ such that
\[
  \epsilon(\alpha \beta) = \epsilon(\alpha) \cdot \lambda_\alpha(\beta).
\]
Note that the linear form $\lambda_\alpha$ is independent of coset representative, so we have a mapping
\[
  \lambda : H/W \to W^*.
\]
If $\lambda$ takes the same value on two different cosets $\alpha_1 W, \alpha_2 W$, then we see that $\epsilon$ is linear on the union $\alpha_1 W \sqcup \alpha_2 W$. Then the associated bilinear form $\left\langle,\right\rangle_\epsilon$ is isotropic on the space $W \oplus \{0, \alpha_2/\alpha_1\}$, which is too big to be isotropic. So $\lambda$ is injective. Comparing sizes, we see that $\lambda$ is surjective also. So there is one coset $\alpha W$ on which $\epsilon$ is identically $1$ or $-1$. This is the charmed coset. On the remaining cosets, the values of $\epsilon$ are those of a nontrivial linear functional on $W$ and hence are equidistributed between $1$ and $-1$.

For a general $W$, take $N \subseteq W$ maximal isotropic. Every $W$-coset decomposes into $N$-cosets, and only the one containing the charmed coset of $N$ will yield a nonzero sum for $\epsilon$, namely $\pm \size{N} = \pm \sqrt{\size{H}}$.

A priori the sign $\pm$ of the sum on the charmed coset depends on both $\epsilon$ and $W$. But if $N \subseteq W$ are coisotropic, then it is easy to see that $N$ and $W$ yield the same sign. Then, taking $W = H$, we obtain that one sign holds for all coisotropic subspaces.

\item Take $N \subseteq W$ maximal isotropic. Then the charmed coset $\alpha W$ of $W$ is the one containing a charmed coset $\alpha N$ of $N$ on which $\epsilon = \pm_{\epsilon} 1$ is constant.

If $\beta/\alpha \in W$, then for all $\gamma \in W^\perp$, we have $\left\langle\alpha\beta, \gamma\right\rangle = 1$ so
\begin{align*}
  \epsilon(\beta \gamma) &= \epsilon(\beta) \epsilon(\alpha) \epsilon(\alpha \gamma) \left\langle\alpha\beta, \gamma\right\rangle \\
  &= \epsilon(\beta),
\end{align*}
since $\alpha, \alpha \gamma \in \alpha N$. By contrast, if $\beta/\alpha \notin W$, there exists $\gamma \in W^\perp$ such that $\left\langle\alpha\beta, \gamma\right\rangle = -1$. Then $\epsilon(\beta \gamma) = -\epsilon(\beta)$, so $\epsilon$ is nonconstant on $\beta W^\perp$. But $W^\perp$ is isotropic, so $\epsilon$ is a linear form (plus a constant) on $\beta W^\perp$ and is therefore equidistributed.
\end{enumerate}
\end{proof}

As you might expect, we apply this lemma to the space $H = H^1$ with its quadratic form $\epsilon$. We will eventually find that $\pm_{\epsilon} = 1$ (it is ``positively charmed,'' one might say) though this is not obvious.

Let $\epsilon_C$ be the following translation of $\epsilon$: for $[\alpha] \in H^1$,
\[
  \epsilon_C(\alpha) = \epsilon\left(\hat\omega_C \alpha\right) = \begin{cases}
    1, & \coef_{\bar\theta_2}(\alpha \xi^2) = 0 \text{ for some } \xi \in R^\cross \\
    -1, & \text{otherwise.}
  \end{cases}
\]
Note that $\epsilon_C$ is still a quadratic form on $H^1$ whose associated bilinear form is the Hilbert symbol $\left\langle\bullet,\bullet\right\rangle_R$. Note that $\epsilon_C(1) = 1$. Lemma \ref{lem:Brauer_const} can be interpreted as saying that
\[
  \heartsuit \diamondsuit \cF_{\ceil{e'/2}}
\]
is charmed for $\epsilon_C$.

\begin{nota}
  If $W \subseteq H^1$ is a subspace, we let
\[
  F_W = \One_{W}
\]
be its characteristic function. If $W$ is coisotropic, we denote by $G_{\epsilon_a,W}$, resp.\ $G_{\epsilon_C,W}$ the characteristic function of its charmed coset with respect to one of the quadratic forms $\epsilon_a, \epsilon_C$ whose associated bilinear form is the Hilbert pairing. The $\epsilon_a$ will be omitted if clear. 
\end{nota}
 The following results will power the computation of Fourier transforms of ring totals, a necessary step in our desired reflection theorems.
\begin{lem} \label{lem:FT_charm}
  Let $W$ be a subspace of $H^1$.
  \begin{enumerate}[(a)]
    \item\label{FT:F}
    \[
      \widehat{F_W} = \frac{\size{W}}{\size{H^0}} F_{W^\perp}.
    \]
    \item\label{FT:G} If $W$ is coisotropic, then
    \[
      \widehat{G_{\epsilon_a,W}} = \frac{\size{W}}{\size{H^0}} \epsilon_a(1) \epsilon_a F_{W^\perp}.
    \]
    \item\label{FT:eG} If $W$ is coisotropic, then
    \[
      \widehat{\epsilon_a G_W}(\delta) = \pm_{\epsilon_a} q^e \epsilon_a(1) \epsilon_a(\delta) F_W(\delta).
    \]
  \end{enumerate}
\end{lem}
\begin{proof}
  Part \ref{FT:F} is a standard property of the Fourier transform. For parts \ref{FT:G} and \ref{FT:eG}, let $N \subseteq W$ be a maximal isotropic subspace, so
  \[
    W \supseteq N = N^\perp \supseteq W^\perp.
  \]
  Let $\beta N$ be the charmed coset of $N$, so $\beta W$ is the charmed coset of $W$.

  For \ref{FT:G}, we compute
  \begin{align*}
    \widehat{G_W}(\delta) &= \frac{1}{\size{H^0}} \sum_{\alpha \in W} \left\langle\alpha \beta, \delta\right\rangle \\
    &= \frac{1}{\size{H^0}} \left\langle\beta, \delta\right\rangle \sum_{\alpha \in W} \left\langle\alpha, \delta\right\rangle \\
    &= \frac{\size{W}}{\size{H^0}} \left\langle\beta, \delta\right\rangle F_{W^\perp} (\delta).
  \end{align*}
  Hence it remains to prove that, for $\delta \in W^\perp$,
  \[
    \left\langle\beta, \delta\right\rangle = \epsilon_a(1)\epsilon_a(\delta).
  \]
  By the definition of the associated bilinear form,
  \[
    \left\langle\beta, \delta\right\rangle = \epsilon_a(1) \epsilon_a(\delta) \epsilon_a(\beta) \epsilon_a(\delta \beta).
  \]
  But $\epsilon_a(\beta) = \epsilon_a(\delta \beta)$ since both arguments lie in the charmed isotropic coset $\beta N$. This establishes \ref{FT:G}.

  For \ref{FT:eG}, we compute
  \begin{align*}
    \widehat{\epsilon_a G_W}(\delta) &= \frac{1}{\size{H^0}} \sum_{\alpha \in W} \epsilon_a(\alpha\beta) \left\langle\alpha\beta, \delta\right\rangle \\
    &= \frac{1}{\size{H^0}} \sum_{\alpha \in W} \epsilon_a(1) \epsilon_a(\delta) \epsilon_a(\alpha\beta\delta) \\
    &= \frac{\epsilon_a(1)\epsilon_a(\delta)}{\size{H^0}} \sum_{\alpha \in W} \epsilon_a(\alpha\beta\delta).
  \end{align*}
  The last sum equals $\pm_{\epsilon_a} \sqrt{\size{H^1}}$ if $\beta\delta W$ is charmed, $0$ otherwise. But $\beta W$ is charmed, so the relevant condition is that $\delta \in W$. So
  \begin{align*}
    \widehat{\epsilon_a G_W}(\delta) &= \frac{\epsilon_a(1)\epsilon_a(\delta)}{\size{H^0}} \cdot \pm_{\epsilon_a} \sqrt{\size{H^1}} F_W(\delta) \\
    &= \pm_{\epsilon_a} q^e \epsilon_a(1)\epsilon_a(\delta) F_W(\delta),
  \end{align*}
  as desired.
\end{proof}

\section{The projectors}

On the space of complex- (or even rational-) valued functions on $H^1$, we can define certain projectors that divide up the work to be done. First look at the cosets of $\cF_0$. Let
\[
  I = \begin{cases}
    \cF_0 \union (1;\pi;\pi)\cF_0, & s > 0 \\
    \cF_0, & s = 0
  \end{cases}
\]
be the union $\cF_0 \union \alpha\cF_0$ of up to two cosets, using the distinguished splitting $R \isom K \cross Q$ if $s > 0$. Let $I_0$, $I_1$, $I_2$ be the restriction operators that restrict the support of a function to $\cF_0$, $I \setminus \cF_0$, and $H^1 \setminus I$, respectively. They are orthogonal idempotents ($I_1$ and/or $I_2$ may vanish). Let $J_i$ be the conjugate of $I_i$ under the Fourier transform. Each $J_i$ is convolution by a certain function supported on $\cF_{e'}$; $J_0$ is none other than the \emph{smear operator} $\Sm_{e'}$ which will occur below. Since $\cF_{e'} \subseteq \cF_0$, each $I_i$ commutes with each $J_j$, so the $I_i J_j$ form a system of nine orthogonal idempotents. For orderliness of presentation, we transform all ring answers to a sum of terms each in the image of one idempotent. The Fourier transform interchanges the images of $I_i J_j$ and $I_j J_i$. In splitting type $(111)$, all nine idempotents are nonzero, although $I_1 J_2$ and $I_2 J_1$ will be found to annihilate every ring total. In the remaining splitting types, some of the idempotents vanish, and correspondingly some terms of our answers can be ignored.

\begin{defn}\label{defn:x's_ur} In unramified splitting types, we define the use of a symbol $\xo$ as follows, where the $\ell_i$ are such that the relevant boxgroups are well defined:
  \begin{itemize}
    \item $\xo F(0, \ell_1, \ell_2) = I_1(F(\emptyset, \ell_1,\ell_2))$, so that
    \[
      F(\emptyset, \ell_1,\ell_2) = F(0, \ell_1, \ell_2) + \xo F(0, \ell_1, \ell_2).
    \]
    \item $\xo\xo F(0, 0, e) = I_2(F(\emptyset, \emptyset, e))$, so that
    \[
      F(\emptyset, \emptyset, e) = F(0, 0, e) + \xo F(0, 0, e) + \xo\xo F(0, 0, e).
    \]
    \item $F\xo(\ell_0, \ell_1, 0) = 2 \cdot J_1(F(\ell_0, \ell_1, \emptyset))$, so that, if there is a distinguished coarse coset,
    \[
      2 F(\ell_0, \ell_1, \emptyset) = F(\ell_0, \ell_1, 0) + F\xo(\ell_0, \ell_1, 0).
    \]
    \item $F\xo\xo(e, 0, 0) = \size{H^0} \cdot J_2(F(e, \emptyset, \emptyset))$, so that
    \[
      \size{H^0} F(e, \emptyset, \emptyset) = F(e, 0, 0) + F\xo(e, 0, 0) + F\xo\xo(e, 0, 0).
    \]
    \item $\xo F\xo(0, e, 0) = 2 \cdot I_1 J_1(F(\emptyset, e, \emptyset))$ (note there must be a distinguished coarse coset for this to be meaningful), so that
    \begin{align*}
      2 F(\emptyset, e, \emptyset) &= 2 F(0, e, \emptyset) + 2 \xo F(0, e, \emptyset) \\
      &= F(\emptyset, e, 0) + F\xo(\emptyset, e, 0) \\
      &= F(0, e, 0) + \xo F(0,e,0) + F\xo(0,e,0) + \xo F\xo(0,e,0).
    \end{align*}
    \item $\xo G(0,\ell_1, \ell_2)$, $\xo\xo G(0,0,e)$, and $\xo G\xo(0,e,0)$, analogously.
  \end{itemize}
\end{defn}
Observe that these definitions are crafted so that the following cute rule applies:
\begin{lem}\label{lem:FT_x}
Let $N$ be one of the symbols $F\xo, \xo F$, $F\xo\xo$, $\xo\xo F$, $\xo F\xo$, and let $N'$ be the symbol made by spelling $N$ backward. If the $\ell_i$ are integers such that $N(\ell_0, \ell_1, \ell_2)$ is meaningful, then
\[
  \widehat{N(\ell_0, \ell_1, \ell_2)} = q^{\ell_2 - \ell_0} N'(\ell_2, \ell_1, \ell_0).
\]
(This lemma also holds for $N = F$, though this will be proved later: see Lemma \ref{lem:orth}.)
\end{lem}
This will allow us to write the ring totals for all three unramified splitting types in a uniform way and verify reflection for them simultaneously.

For the unshifted quadratic form $\epsilon$, $\cF_i$ is charmed for $i \geq e'/2$, an easy consequence of Lemma \ref{lem:Brauer_const}, from which we derive:

\begin{lem} \label{lem:e_projectors}
  The function $\epsilon_C = \epsilon_C F(\emptyset, \emptyset, e)$ has at most the following projections nonzero:
  \begin{enumerate}[$($a$)$]
    \item If $[\hat{\omega}_C] \in \cF_0$, then $I_0 J_0$, $I_1 J_1$, and $I_2 J_2$.
    \item If $[\hat{\omega}_C] \notin \cF_0$, then $I_0 J_1$, $I_1 J_0$, and $I_2 J_2$.
  \end{enumerate}
\end{lem}
\begin{proof}
It suffices to prove that
\[
  I_i \left(\epsilon_C F(\emptyset, \emptyset, e)\right)
  = \epsilon_C({}_{\times^i} F)(0, 0, e)
\]
is in the image of $J_j$, for each pair $(i,j)$ mentioned in the lemma. Now $\cF_0$ is charmed for the unshifted quadratic form $\epsilon$, an easy consequence of Proposition \ref{prop:conic_perturb}. Hence $\hat{\omega}_C \cF_0$ is charmed for $\epsilon_C$, from which we get
\[
  \epsilon_C({}_{\times^i} F)(0, 0, e)
  = \epsilon_C({}_{\times^j} G)(0, 0, e).
\]
Taking the Fourier transform by Lemma \ref{lem:FT_charm}\ref{FT:eG}, we obtain a nonzero scalar multiple of
\[
  \epsilon_C({}_{\times^j} F)(0, 0, e),
\]
which lies in the image of $I_j$. Hence the original $\epsilon_C({}_{\times^i} F)(0, 0, e)$ lies in the image of $J_j$, as desired.
\end{proof}

\section{Notation}

If $T(\ell_0, \ell_1, \ell_2)$ is a boxgroup, we let $F(T) = F(\ell_0, \ell_1, \ell_2)$ be its characteristic function, and $G(T) = G(\ell_0, \ell_1, \ell_2)$ be the characteristic function of its charmed coset, if applicable. Let $T^\cross(\ell_0, \ell_1, \ell_2)$ denote the subset of elements of $T(\ell_0, \ell_1, \ell_2)$ having minimal level, assuming this level is less than $e'$:
\[
  T^\cross(\ell_0, \ell_1, \ell_2) = \begin{cases}
    T(\ell_0, \ell_1, \ell_2) \setminus T(\ell_0 + 1, \ell_1 - 1, \ell_2), & \ell_0 \geq 0, \ell_1 \geq 1 \\
    T(\ell_0, \ell_1, \ell_2) \setminus T(\ell_0 + 1, 0, \ell_2 - 1), & \ell_0 \geq 0, \ell_1 = 0 \\
    T(\emptyset, \ell_1, \ell_2) \setminus T(0, \ell_1, \ell_2), & \ell_0 = \emptyset
  \end{cases}
\]
Let $F^\cross(\ell_0, \ell_1, \ell_2)$ and $G^\cross(\ell_0, \ell_1, \ell_2)$ be the characteristic functions of the subsets $T^\cross(\ell_0, \ell_1, \ell_2)$ and $\psi_C T^\cross(\ell_0, \ell_1, \ell_2)$, respectively. (We will eventually prove that the latter coset is charmed.) This will provide enough notation to write the ring totals in the succeeding sections.

\chapter{Ring volumes for \texorpdfstring{$\xi'_1$}{xi'1}}
\label{cha:xi1}
In this section, we will compute the volume of vectors $\xi'_1 \in \OO_R$ satisfying the $\cM_{11}$ and $\cN_{11}$ conditions. Examples of our answers are tabulated in Appendix \ref{cha:zone_examples}.

Because $\bar\xi_1 = \xi_1$, we freely omit the bar on $m_{11} = \bar m_{11}$ and $n_{11} = \bar n_{11}$.

\section{The smearing lemma}
\label{sec:smear}
We begin with a simple lemma that allows us to reduce to the case $m_{11}$ large.
\begin{lem} \label{lem:smear}
  For a first vector problem $\fP$, define $u_1$ to be the unique value such that, when the conic is transformed to minimal discriminant in accordance with Lemma \ref{lem:tfm_conic}, $\beta^\odot = \delta^\odot {\xi_1^\odot}^2$ must lie in $\pi^{u_1} \OO_R \setminus \pi^{u_1 + 1} \OO_R$, to wit:
  \[
  u_1 = \begin{cases}
    1 & \text{if } (\sigma, h_1) \in \{((1^3), -1)
      , ((1^21), 2), ((1^21), 3)
    \} \\
    0 & \text{otherwise.}
  \end{cases}
  \]
  Assume that the $m_{11}$ and $n_{11}$ of $\fP$ satisfy
  \[
  m_{11}^\odot = m_{11} - p^\odot > u_1
  \]
  (so $\cM_{11}$ is active even after transformation) and $m_{11}^\odot > n_{11}$. Let $m_{11}^\sharp > m_{11}$ lie in the same class mod $\ZZ$, and let $m_{11}^{\odot\sharp} = m_{11}^\sharp - p^{\odot}$. Let $\fP^\sharp$ be the first vector problem with $m_{11} = m_{11}^\sharp$ and the rest of the data the same.

  Then the answer to $\fP$ can be computed from that of $\fP^\sharp$ by the following formula:
  \[
  W_{\fP} = q^{m_{11}^\sharp - m_{11}} \Sm_{r} W_{\fP^\sharp}
  \]
  Here the \emph{smear operator} $\Sm_{r}$ is defined by the following convolution:
  \[
  \Sm_{r} W(\delta) = \frac{1}{\size{\cF_r}} \sum_{\alpha \in \cF_r} W(\alpha\delta),
  \]
  and $r$ is the level for which
  \[
  \cF_r = \left\{ [\eta] : \eta \equiv 1 \mod \pi^{m_{11}^\odot - u_1} \right\}.
  \]
\end{lem}
We call this the \emph{smearing lemma} because it states that the function $W_{\fP}$ can be obtained from $W_{\fP^\sharp}$ by averaging over the cosets of $\cF_{r}$, like blurring a picture by averaging over larger pixels. Here the symbol $\sharp$ is used to mark the ``sharper'' image given by the solutions of $\fP^\sharp$.

The level $r$ is given explicitly as follows:
\begin{itemize}
  \item In unramified type,
  \[
  r = \floor{\frac{m_{11}}{2}}, \quad 0 < m_{11} \leq 2e
  \]
  \item In splitting type $(1^3)$,
  \[
  r = \floor{m_{11}}, \quad 0 < m_{11} \leq 2e.
  \]
    \item In splitting type $(1^21)$,
    \[
    r = m_{11} - \frac{d_0}{2} + \left\{\frac{h_1}{2}\right\},
    \quad \frac{d_0 + 1}{2} \leq m_{11} \leq 2e + \frac{d_0}{2}.
    \]
\end{itemize}

\begin{proof}
  We will prove the identity by computing in two ways the volume of the set
  \begin{equation}
    \cS = \left\{
    (\eta, \xi'):
    \begin{aligned}
      \eta &\equiv 1 \mod \pi^{m_{11}^\odot - u_1}, \\
      \lambda^\diamondsuit(\eta \delta^\odot \xi'^2) &\equiv 0 \mod \pi^{m_{11}^{\odot\sharp}}, \\
      \widetilde\omega_C^{-1} \gamma^2 \xi'^2 &\equiv a \mod \pi^{n_{11}} \text{ for some } a \in \OO_K
    \end{aligned}
    \right\}
    \subseteq \OO_R^\cross \cross \PP(\OO_R).
  \end{equation}
  First, fix $\eta \in 1 + \pi^{m_{11}^\odot - u_1}\OO_R$. The conditions on $\xi'$ are seen to be the $\cM_{11}$ and $\cN_{11}$ conditions for $\delta$ replaced by $\eta \delta$; the omission of the $\eta$ factor in $\cN_{11}$ makes no difference, since the left-hand side is a multiple of $\pi^{u_1}$ and the addition would have valuation at least
  \[
  (m_{11}^\odot - u_1) + u_1 = m_{11}^\odot \geq n_{11}.
  \]
  So the volume of $\xi'$ for fixed $\eta$ is $W_{\fP^\sharp} (\delta\eta)$, and since $[\eta]$ takes all classes in $\cF_{r}$ equally often while ranging in a pixel of volume $q^{-3(m_{11}^\odot - u_1)}$, we get
  \[
  \mu(\cS) = \frac{q^{-3(m_{11}^\odot - u_1)}}{\size{\cF_r}} \sum_{\eta \in \cF_{r}} W_{\fP^\sharp}(\delta\eta).
  \]
  On the other hand, a fixed $\xi' \in \OO_R$ has a chance of being the second coordinate of a pair in $\cS$ only if
  \begin{itemize}
    \item it satisfies the $\cN_{11}$ condition $\widetilde\omega_C^{-1} \gamma^2 \xi'^2 \equiv a \mod \pi^{n_{11}}$ for some $a \in \OO_K$, and
    \item it satisfies the $\cM_{11}$ condition for some $\eta = 1 + \pi^{m_{11}^\odot - u_1} \eta'$, $\eta' \in \OO_R$; in particular,
    \begin{align*}
      0 &\equiv \lambda^\diamondsuit(\eta \delta^\odot \xi'^2) \\
      &= \lambda^\diamondsuit(\delta^\odot \xi'^2) + \pi^{m_{11}^\odot - u_1} \lambda^\diamondsuit(\eta' \delta^\odot \xi'^2) \\
      &= \lambda^\diamondsuit(\delta^\odot \xi'^2) + \pi^{m_{11}^\odot} \lambda^\diamondsuit\left(\eta' \pi^{-u_1} \delta^\odot \xi'^2 \right) \\
      &\equiv \lambda^\diamondsuit(\delta^\odot \xi'^2) \mod \pi^{m_{11}^\odot}.
    \end{align*}
  \end{itemize}
  The volume of $\xi'$ satisfying these conditions is, by definition, none other than $W_{\fP}$. For fixed $\xi'$, the value of $\eta'$ is constrained by $\cM_{11}$ alone:
  \begin{align*}
    \pi^{m_{11}^\odot - u_1} \lambda^\diamondsuit(\eta' \delta^\odot \xi'^2) &\equiv 0 \mod \pi^{m_{11}^{\odot\sharp}} \\
    \lambda^\diamondsuit\left(\eta' \cdot \pi^{-u_1}\delta^\odot\xi'^2\right) &\equiv 0 \mod \pi^{m_{11}^{\odot\sharp} - m_{11}^\odot} = \pi^{m_{11}^\sharp - m_{11}}.
  \end{align*}
  Since $\lambda^\diamondsuit$ is a perfect linear functional and $\pi^{-u_1}\delta^\odot\xi'^2$ is a primitive vector in $\OO_R$, the volume of $\eta'$ satisfying this congruence is $q^{m_{11} - m_{11}^\sharp}$, which makes a volume of $q^{-3(m_{11}^\odot - u_1) + m_{11} - m_{11}^\sharp}$ for $\eta$. So
  \[
  \mu(\cS) = q^{-3(m_{11}^\odot - u_1) + m_{11} - m_{11}^\sharp} W_{\fP}(\delta).
  \]
  Comparing the two expressions for $\mu(\cS)$, the result follows.
\end{proof}

\section[The zones when \texorpdfstring{$\cN_{11}$}{N11} is strongly active]{The zones when \texorpdfstring{$\cN_{11}$}{N11} is strongly active (black, plum, purple, blue, green, and red)}\label{sec:strong}
In this section, we solve first vector problems in which $\cN_{11}$ is strongly active, that is, $n_{11} > s$. In view of the smearing lemma, we assume that $m_{11} > 2e$. Let $n_c$ (``n for the colorful zones'') be $n_{11} - s$.

The following little symmetry will be occasionally useful:
\begin{lem}\label{lem:recenter}
Let $\psi = a + b \theta_1 \in \OO_R$ with $a \in \OO_K^\cross$, $b \in \pi^{n_c} \OO_K^\cross$, and let $\fP'$ be the first vector problem derived from $\fP$ by replacing the pertinent extender vector $\theta_1$ by $\theta_1' = \psi^{-1} \theta_1$. Then:
\begin{enumerate}[$($a$)$]
  \item \label{recenter:w_C} $\hat\omega_{\fP'} = \psi^{-1} \hat\omega_\fP$.
  \item \label{recenter:e_C} $\epsilon_{\fP'}(\delta) = \epsilon_\fP(\psi^{-1} \delta)$.
  \item \label{recenter:W} $W_{\theta_1',m_{11},n_{11}}(\delta) = W_{\theta_1,m_{11},n_{11}}(\psi\delta)$.
\end{enumerate}
\end{lem}
\begin{proof}
The left-hand side is the volume of $\xi'$ for which
\begin{align*}
  \delta \left(\frac{\xi'}{\xi'_0}\right)^2 &\in \OO_K^\cross + B_{\psi^{-1}\theta_1}(m_{11},n_c) \\
  &= \psi^{-1} \cdot \left(\OO_K^\cross + B_{\theta_1}(m_{11},n_c)\right)
\end{align*}
by Lemma \ref{lem:boxes_basic_ur}\ref{boxes_basic_ur:recenter}, since $\psi \in \OO_K^\cross + B_{\theta_1}(m_{11},n_c)$. So the condition on $\xi'$ can be written as
\[
  \delta\psi \left(\frac{\xi'}{\xi'_0}\right)^2 \in \OO_K^\cross + B_{\theta_1}(m_{11},n_{11}),
\]
of which the solution volume is seen to be $W_{\theta_1,m_{11},n_{11}}(\psi\delta)$. Since $H^1$ is an elementary $2$-group, $[\psi^{-1}] = [\psi]$, so this agrees with the inverse appearing in parts \ref{recenter:w_C} and \ref{recenter:e_C}.
\end{proof}

The following formula for the sum of the values of $W_{m_{11},n_{11}}$ will be essential. It can be thought of as an extreme version of the smearing lemma where we sum over the whole of $\cF_0$.
\begin{lem} \label{lem:sum_strong}
If the values of $m_{11}$ and $n_{11}$ make $\cN_{11}$ strongly active, then $W_{m_{11},n_{11}} = 0$ unless the chosen coarse coset is $\cF_0$. In this case
\[
  \sum_{\delta \in \cF_0} W_{m_{11},n_{11}}(\delta) = \size{H^0} q^{2e - m_{11} - n_{11} + \frac{s}{2} + \frac{d_0}{2} + v(N(\gamma))}
\]
where $d_0 = v_K(\disc_K R)$; equivalently,
\[
  \sum_{\delta \in \cF_0} W^{\odot}_{m_{11},n_{11}}(\delta) = \size{H^0} q^{2e - m_{11} - n_{11} + \frac{s}{2} + \frac{d_0}{2} + v(N(\gamma\gamma^{\odot}))}.
\]
\end{lem}
\begin{proof}
By Lemma \ref{lem:beta}, the support of $W_{m_{11},n_{11}}$ consists of the classes $[\delta] = [\beta]$ in $H^1$ of elements $\beta$ of the box
\[
  \OO_K^\cross + B_{m_{11},n_c} = \{x + y \pi^{n_c} \theta_1 + z \pi^{m_{11}} \theta_2, \quad x \in \OO_K^\cross, \quad y,z \in \OO_K\}.
\]
Our method is to show that $W_{m_{11},n_{11}}(\delta)$ can be interpreted as the volume of $\beta$ of class $[\delta]$ in the box, up to a scalar. Then since each $\beta$ belongs to just one square-class, the sum is known.

Since $m_{11},n_c > 0$, $\beta$ must be a unit, explaining why $[\delta] \in \cF_0$. For fixed $\delta \in \OO_R^\cross$, as $\xi'$ ranges over the solution set of its transformed conditions, $\beta$ ranges over the elements of its box of class $[\delta]$, up to scaling. The correspondence is given by a relation of the form
\begin{equation}
  \beta = \delta \left(\frac{\xi'}{\xi'_0}\right)^2,
\end{equation}
where
\[
  \xi'_0 = \frac{1}{\gamma}\sqrt{\delta\hat\omega_C}
\]
is an element of $\OO_R$ constructed from the $\gamma$ of Lemma \ref{lem:gamma_white} whose valuations represent a lower bound on the valuations of $\xi'$.

Now we compare the projective volumes of $\beta$ and $\xi'$ satisfying the $\cM_{11}$ and $\cN_{11}$ conditions. In the sequence
\[
  \xi' \longmapsto \frac{\xi'}{\xi'_0} \longmapsto \left(\frac{\xi'}{\xi'_0}\right)^2 \longmapsto \delta \left(\frac{\xi'}{\xi'_0}\right)^2 = \beta,
\]
each member is a primitive vector in $\OO_R$, so we can speak of projective volumes.

Dividing by $\xi'_0$ is a one-to-one operation that scales both affine and projective volumes by
\[
  q^{v_K(N_{R/K}(\xi'_0))} = q^{s/2 - v(N(\gamma))}.
\]
So there is a volume $q^{s/2 - v(N(\gamma))} W_{m_{11},n_{11}}(\delta)$ of $\xi/\xi'_0$.

Squaring, on units, multiplies small projective volumes by $q^{-2e}$ (since it takes the $i$-pixel about $1$ to the $(i+e)$-pixel for $i > e$). But it is $\size{H^0}$-to-one since there are $\size{H^0}$-many square roots of $1$ in $R$, up to scaling by $\pm 1$. Since the resolvent conditions are invariant under multiplying $\xi'$ by a square root of $1$, the volume of the squared image $(\xi/\xi'_0)^2$ is
\[
  \frac{q^{-2e + s/2 - v(N(\gamma))}}{\size{H^0}} W_{m_{11},n_{11}}(\delta).
\]

Lastly, $\delta$ is a unit, so multiplying by it does not change volumes. Hence
\[
  \frac{q^{-2e + s/2 - v(N(\gamma))}}{\size{H^0}} \sum_{\delta \in \cF_0} W_{m_{11},n_{11}}(\delta)
\]
is the volume of the box $\PP(\OO_K^\cross + B_{m_{11},n_c})$. Thus it suffices to prove that
\[
  \mu(\PP(\OO_K^\cross + B_{m_{11},n_c})) = q^{-m_{11}-n_{11}+s + \frac{d_0}{2}}.
\]
Converting to affine volumes, with $\mu(R) = 1$,
\begin{align*}
  \mu(\PP(\OO_K^\cross + B_{m_{11},n_c})) &= \frac{1}{1 - \frac{1}{q}} \mu(\OO_K^\cross + B_{m_{11},n_c}) \\
  &= \mu(\OO_K + B_{m_{11},n_c}) \\
  &= q^{-i}
\end{align*}
where $i$ is the integer such that
\begin{align*}
  \Lambda^3(\OO_K + B_{m_{11},n_c}) &= \pi^i \Lambda^3 \OO_R \\
  \Lambda^3(\OO_K\left\langle1, \pi^{n_c} \theta_1, \pi^{m_{11}} \theta_2\right\rangle) &= \pi^i \Lambda^3 \OO_R \\
  n_c + m_{11} &= i + \frac{d_0}{2}.
\end{align*}
So
\begin{align*}
  \sum_{\delta \in \cF_0} W_{m_{11},n_{11}}(\delta) &= \size{H^0} \cdot q^{2e - s/2 + v(N(\gamma))} \cdot q^{-m_{11} - n_c + \frac{d_0}{2}}\\
  &= \size{H^0} \cdot q^{2e - m_{11} - n_{11} + \frac{s}{2} + \frac{d_0}{2} + v(N(\gamma))}
\end{align*}
as desired. The $\odot$ version is proved identically.
\end{proof}

Drawing the relevant information from Table \ref{tab:tfm_conic}, we can explicitly compute some of the quantities involved in the foregoing theorem: see Table \ref{tab:sum_strong}.
\begin{table}[ht]
  \begin{tabular}{rc|cl}
    spl.t. & $h_1$ & $\xi^\odot_0 \sim$ & $v(N(\gamma\gamma^{\odot}))$ \\ \hline
    ur & $0$ & $(\pi^{s/2};1;1)$, $s$ even & $0$ \\
    ur & $1$ & $(\pi^{(s-1)/2} ; 1 ; 1)$, $s$ odd & $1/2$ \\
    $(1^3)$ & $1$ & $1$ & $0$ \\
    $(1^3)$ & $-1$ & $\pi_R^2$ & $-2$ \\
    $(1^21)$ & $0$ & $(\pi^{s'/2} ; 1)$, $s'$ even & $d_0/4$ \\
    $(1^21)$ & $1$ & $1$, \quad $s' = -1/2$ & $(d_0 - 1)/4$ \\
    $(1^21)$ & $2$ & $(\pi^{(s'+1)/2} ; \pi_Q)$, $s'$ odd & $(d_0-6)/4$ \\
    $(1^21)$ & $3$ & --- & ---
  \end{tabular}
\caption{Some constants appearing in Lemma \ref{lem:sum_strong}}
\label{tab:sum_strong}
\end{table}
As shown, we obtain as a by-product that $s$ determines $h_1$.
  The last row is left blank because the integrality of the $K$-valuation yields $s' \in 1/2 + 2\ZZ$, contradicting the known fact $s' \in \{-1/2\} \union \ZZ_{\geq 0}$ from Lemma \ref{lem:s}.

We now come to our main lemma, which computes $ W_{m_{11},n_{11}} $ for large $m_{11}$.

\subsection{Unramified}
\begin{lem}\label{lem:111_strong_zones}
  Suppose $R$ is unramified over $K \supseteq \QQ_2$.
  Let $m_{11}$ and $n_c$ be integers with $m_{11}>2e$ and $n_c>0$. Let $\square_C =  \min\left\{2\ell(\hat\omega_C) + 1, e\right\}$ be the squareness of the $\delta = 1$ conic if $s$ is even; let $\square_C = 0$ if $s$ is odd. If $s$ is even, let
  \[
  \tilde{n} = \floor{\frac{2e - s - n_c + 2}{4}} = \ceil{\frac{e - \frac{s}{2} - \ceil{\frac{n_c}{2}}}{2}}.
  \]
  Then $ W_{m_{11},n_{11}} $ is given as follows:
  \begin{enumerate}[$($a$)$]
    \item If $n_c > 2e$ (black zone), then
    \begin{align*}
    W_{m_{11},n_{11}} &= \size{H^0} q^{2e - m_{11} - n_c - \floor{\frac{s}{2}}} F\left( e, \emptyset, \emptyset \right) \\
    &= q^{2e - m_{11} - n_c - \floor{\frac{s}{2}}} [F(e,0,0) + F\xo(e,0,0) + F\xo\xo(e,0,0)]
    \end{align*}
    \item If $ 2e - s < n_c \leq 2e $ (purple zone), then
    \begin{align*}
    W_{m_{11},n_{11}} &= 2q^{e - m_{11} - \ceil{\frac{n_c}{2}} - \floor{\frac{s}{2}}} F\left( \floor{\frac{n_c}{2}}, e - \floor{\frac{n_c}{2}}, \emptyset \right) \\
     &= q^{e - m_{11} - \ceil{\frac{n_c}{2}} - \floor{\frac{s}{2}}} \Big[F\left(\floor{\frac{n_c}{2}}, e - \floor{\frac{n_c}{2}}, 0\right) + F\xo\Big(\floor{\frac{n_c}{2}}, e - \floor{\frac{n_c}{2}}, 0\Big)\Big]
    \end{align*}
    \item If $ 2e - s - 2\square_C < n_c \leq 2e - s $ and $ n_c \geq \dfrac{2e - s}{3} $ (blue zone), then $s$ is even and
    \[
    W_{m_{11},n_{11}} = q^{-m_{11} + \floor{\frac{2e - s - n_c}{4}}} F\left( \floor{\frac{n_c}{2}}, e - \tilde{n} - \floor{\frac{n_c}{2}}, \tilde{n} \right)
    \]
    \item If $ \square_C < n_c \leq 2e - s - 2\square_C $ (green zone), then
    \[
    W_{m_{11},n_{11}} = q^{-m_{11} + \ceil{\ell_C}}(1 + \epsilon_C)F\left( \floor{\frac{n_c}{2}}, \ell_C + \frac{s}{2} + \One_{2\nmid n_c}, e - \ell_C - \frac{s}{2} - \ceil{\frac{n_c}{2}}\right)
    \]
    where $\ell_C = \frac{\square_C - 1}{2} \in \{-1/2\} \union \ZZ_{\geq 0}$.
    \item If $n_c < \dfrac{2e - s}{3}$ and $n_c \leq \square_C$ (red zone), then $s$ is even and
    \begin{align*}
    W_{m_{11},n_{11}} &= \sum_{\floor{\frac{n_c}{2}} \leq \ell < \tilde{n}} q^{-m_{11} + \ell}(1 + \epsilon_C) G^\cross\left( \ell, \ceil{\frac{n_c}{2}} + \frac{s}{2}, e - \ceil{\frac{n_c}{2}} - \frac{s}{2} - \ell \right)  + {}\\
    & \quad {} + q^{-m_{11}+ \floor{\frac{2e - s - n_c}{4}}} G(\tilde{n}, e - 2\tilde{n}, \tilde{n}).
    \end{align*}
  \end{enumerate}
\end{lem}
\begin{proof}
First note that if $s \geq 2e$, then the blue, green, and red zones are empty, and if $s < 2e$ is odd, then since $\ell_C = -1/2$, the blue and red zones are empty. This ensures that the answers are at least well defined.

By Lemma \ref{lem:beta}, the support of $W_{m_{11},n_{11}}$ consists of the classes $[\delta] = [\beta]$ in $H^1$ of elements $\beta$ of the box
\[
1 + B_{m_{11},n_c} = \{1 + y \pi^{n_c} \theta_1 + z \pi^{m_{11}} \theta_2 : y,z \in \OO_K\}.
\]
Since $m_{11} > 2e$, the $z$ term has no effect on $[\beta]$, and we ignore it. In particular, the conic
\[
  \coef_{\theta_2} (\beta \xi^2) = 0, \quad \text{that is,} \quad \tr(\hat\omega_C \beta \xi^2) = 0
\]
has a solution $\xi = 1$, so $\epsilon_C(\beta) = \epsilon(\hat\omega_C \beta) = 1$ for all $\beta$ in the box.

  \subsubsection{Black zone.} In the black zone, we have $\beta \equiv 1 \mod 4\pi$, so $[\beta] = 1$. Hence only $\delta = 1$ yields a nonzero volume, which is, by Lemma \ref{lem:sum_strong},
  \[
    \size{H^0} q^{2e - m_{11} - n_c + \frac{d_0}{2} - (s/2 - v(N(\gamma)))} = \size{H^0} q^{2e - m_{11} - n_c - \floor{\frac{s}{2}}}.
  \]

  For the remaining zones, let $ W_{m_{11},n_{11}}' $ denote the claimed value of $W_{m_{11},n_{11}}$ in each case. Our proof method will consist of two steps:
  \begin{itemize}
    \item We prove that $W_{m_{11},n_{11}}(\delta) \leq W_{m_{11},n_{11}}'(\delta)$ for every $\delta$ (the \emph{bounding step}).
    \item We check that
    \[
    \sum_{\delta \in H^1} W'_{m_{11},n_{11}}(\delta) = \size{H^0} q^{2e - m_{11} - n_c} = \sum_{\delta \in H^1} W_{m_{11},n_{11}}(\delta)
    \]
    (the \emph{summing step}), implying that equality must hold for every $\delta$.
  \end{itemize}

  \subsubsection{Purple zone.} In the purple zone, $s > 0$ defines a splitting $R = K \cross Q$. If we translate $\theta_1$ so that $\theta_1^Q \equiv 0 \mod \pi^s$, then we get $\beta^Q \equiv 1 \mod 4\pi$. Also, $\beta \equiv 1$ mod $\pi^{n_c}$, and indeed, $\beta^{(K)}$ can achieve any value $\equiv 1 \mod \pi^{n_c}$, each congruence class modulo $4\pi$ achieved equally often. So $[\delta] = [\beta]$ ranges uniformly over $T\left(\floor{\frac{n_c}{2}}, e - \floor{\frac{n_c}{2}}, \emptyset\right)$, and for each class $[\delta]$ that is attained,
  \[
    W_{m_{11},n_{11}}(\delta) = \frac{\size{H^0} q^{2e - m_{11} - n_c - \floor{\frac{s}{2}} }}{\size{T\left(\floor{\frac{n_c}{2}}, e - \floor{\frac{n_c}{2}}, \emptyset\right)}}
    = 2q^{e - m_{11} - \ceil{\frac{n_c}{2}} - \floor{\frac{s}{2}}},
  \]
  as claimed.

  \subsubsection{Blue zone.} Note that $s < 2e$ and that $s$ is even (as the bounds imply $\square_C > 0$). Our strategy is to note that
  \[
    \beta \in 1 + B_{m_{11},n_c} \subseteq 1 + B_{m', n_c}
  \]
  for some $m' \leq m_{11}$ for which $[1 + B_{m',n_c}]$ is a boxgroup. Here, we find that the gray-blue condition \eqref{eq:box_blue_ur} in Lemma \ref{lem:boxgps_ur} is the most stringent one, so we take
  \[
    m' = \ceil{\frac{n_c}{2}} + e + \frac{s}{2}
  \]
  and get
  \begin{align*}
  \beta &\in \left[1 + B_{\theta_1}\left( \ceil{\frac{n_c}{2}} + e + \frac{s}{2}, n_c\right)\right] \\
   &= T\left( \floor{\frac{n_c}{2}}, \floor{\frac{\ceil{\frac{n_c}{2}} + e + \frac{s}{2}}{2}} - \floor{\frac{n_c}{2}}, e - \floor{\frac{\ceil{\frac{n_c}{2}} + e + \frac{s}{2}}{2}} \right) \\
  &= T\left( \floor{\frac{n_c}{2}}, e - \tilde{n} - \floor{\frac{n_c}{2}}, \tilde{n} \right).
  \end{align*}
  For each such $\delta$, the value of
  \[
    W_{m_{11},n_{11}}(\delta)
  \]
  is controlled by the conic via Lemma \ref{lem:conic_1}, once we know the level
  \[
    \ell = \min\left\{\ell(\delta\hat\omega_C\diamondsuit\heartsuit), \floor{\dfrac{e}{2}}\right\}.
  \]
  We claim that all these conics are blue in the sense of Lemma \ref{lem:conic_1}; this requires
  \[
    2 \ell + 1 \stackrel{?}{\geq} e - n' = e - \ceil{\frac{n_{11}}{2}} = e - \ceil{\frac{n_c + s}{2}},
  \]
  that is,
  \[
    \delta \hat\omega_C \diamondsuit \heartsuit \stackrel{?}{\in} \cF_{\ceil{\frac{e - \ceil{\frac{n_c + s}{2}}}{2}}} = \cF_{\ceil{\frac{e - \floor{\frac{n_c+s+1}{2}}}{2}}} = \cF_{\ceil{\frac{2e - n_c - s - 1}{4}}}.
  \]
  When $\delta = 1$, the required relation
  \[
    \hat\omega_C \diamondsuit \heartsuit \in \cF_{\ceil{\frac{2e - n_c - s - 1}{4}}}
  \]
  follows from the given inequality $n_c > 2e - s - 2\square_C$. So it suffices to show that
  \[
    \beta \in \cF_{\ceil{\frac{2e - n_c - s - 1}{4}}}.
  \]
  Since $\beta \equiv 1 \mod \pi^{n_c}$, it suffices to show that
  \[
    \floor{\frac{n_c}{2}} \stackrel{?}{\geq} \ceil{\frac{2e - n_c - s - 1}{4}}.
  \]
  But the given red-blue inequality $ n_c \geq \frac{2e - s}{3} $ gives
  \[
    \frac{n_c - 1}{2} \geq \frac{2e - n_c - s - 1}{4},
  \]
  from which the desired inequality follows by taking ceilings. So all conics are blue, and for every $\delta$, $W_{\theta_1,m_{11},n_{11}}(\delta) = q^{-m_{11} + \floor{\frac{2e - s - n_c}{4}}}$ if nonzero. The summing step is now straightforward:
  \begin{align*}
  \sum_{\delta \in H^1} W'_{m_{11},n_{11}}(\delta)
  &= q^{-m_{11} + \floor{\frac{2e - s - n_c}{4}}} \Size{T\left( \floor{\frac{n_c}{2}}, e - \tilde{n} - \floor{\frac{n_c}{2}}, \tilde{n} \right)} \\
  &= \size{H^0} \cdot q^{-m_{11} + \floor{\frac{2e - s - n_c}{4}} + e - \tilde{n} - \floor{\frac{n_c}{2}} + 2\tilde{n}} \\
  &= \size{H^0} \cdot q^{-m_{11} + e - \floor{\frac{n_c}{2}} + \floor{\frac{2e - s - n_c}{4}} + \floor{\frac{2e - s - n_c + 2}{4}}} \\
  &= \size{H^0} \cdot q^{-m_{11} + e - \floor{\frac{n_c}{2}} + \floor{\frac{2e - s - n_c}{2}}} \\
  &= \size{H^0} \cdot q^{-m_{11} + e - \floor{\frac{n_c}{2}} + e - \frac{s}{2} - n_c + \floor{\frac{n_c}{2}}} \\
  &= \size{H^0} \cdot q^{2e - m_{11} - n_c - \frac{s}{2}}.
  \end{align*}
  This completes the proof, and in particular shows that $\epsilon = 1$ identically on $\hat\omega_C \cdot T\left( \floor{\frac{n_c}{2}}, e - \tilde{n} - \floor{\frac{n_c}{2}}, \tilde{n} \right)$. This result will be important in proving the remaining zones.

  \subsubsection{Green zone.} Again, we write
\[
\beta \in 1 + B_{m_{11},n_c} \subseteq 1 + B_{m', n_c}
\]
where $m'$ is as large as possible to make a boxgroup. This time, we find that the gray-green inequality \eqref{eq:box_green_ur} is the most stringent of the conditions in Lemma \ref{lem:boxgps_ur}, so we take $m' = 2 \ceil{\frac{n_c}{2}} + \square_C + s$ (noting that $m'$ is an odd integer) and find that the support of $W_{m_{11},n_{11}}$ is contained in
\[
\left[1 + B_{\theta_1}\left( 2 \ceil{\frac{n_c}{2}} + \square_C + s, n_c\right)\right] = T\bigg( \Big\lfloor{\frac{n_c}{2}}\Big\rfloor, \ell_C + \frac{s}{2} + \One_{2\nmid n_c}, e - \ell_C - \frac{s}{2} - \Big\lceil\frac{n_c}{2}\Big\rceil\bigg)\text{.}
\]
We claim all conics are green of the same squareness $\square_C$. The zone boundaries easily imply $\ell_C < \floor{e/2}$, so
\[
[\hat\omega_C \diamondsuit \heartsuit] \in \cF_{\floor{\ell_C}} \setminus \cF_{\floor{\ell_C} + 1}.
\]
We then note that $[\beta] \in \cF_{\floor{\ell_C} + 1}$, because $\beta \equiv 1 \mod \pi^{n_c}$ and we have the inequality $n_c > \square_C$. So $[\beta\hat\omega_C \diamondsuit \heartsuit]$ is also of exact level $\ell_C$. Thus all conics are green of the same squareness, and by Lemmas \ref{lem:conic_1} and \ref{lem:conic_pi}, we have the bound
\begin{align*}
W_{m_{11},n_{11}} &\leq 2 q^{-m_{11} + \ceil{\ell_C}} F\left( \floor{\frac{n_c}{2}}, \ell_C + \frac{s}{2} + \One_{2\nmid n_c}, e - \ell_C - \frac{s}{2} - \ceil{\frac{n_c}{2}}\right),
\end{align*}
indeed
\begin{align*}
W_{m_{11},n_{11}} &\leq q^{-m_{11} + \ceil{\ell_C}} (1 + \epsilon_C) F\left( \floor{\frac{n_c}{2}}, \ell_C + \frac{s}{2} + \One_{2\nmid n_c}, e - \ell_C - \frac{s}{2} - \ceil{\frac{n_c}{2}}\right).
\end{align*}
Here $\epsilon_C(\delta) = \epsilon(\hat\omega_C \delta)$. This completes the bounding step.

To perform the summing step, we need to compute the sum
\[
q^{-m_{11} + \ceil{\ell_C}} \sum_{\delta \in T\left( \floor{\frac{n_c}{2}}, \ell_C + \frac{s}{2} + \One_{2\nmid n_c}, e - \ell_C - \frac{s}{2} - \ceil{\frac{n_c}{2}}\right)} (1 + \epsilon_C(\delta)).
\]
The term $1$ is found to sum to the desired total $\size{H^0} \cdot q^{2e - m_{11} - n_c - \floor{\frac{s}{2}}}$.
We claim that
\[
\sum_{\delta \in T\left( \floor{\frac{n_c}{2}}, \ell_C + \frac{s}{2} + \One_{2\nmid n_c}, e - \ell_C - \frac{s}{2} - \ceil{\frac{n_c}{2}}\right)} \epsilon_C(\delta) = 0,
\]
in other words that $\epsilon$ is \emph{equidistributed} between $1$ and $-1$ in this boxgroup. This follows from Lemma \ref{lem:charm}\ref{charm:constant}: because $\cF_{\floor{\ell_C} + 1} \neq \hat\omega_C \cF_{\floor{\ell_C} + 1}$ is an uncharmed coset, we have $\epsilon$ equidistributed on cosets of $\cF_{e - \floor{\ell_C} - 1}$.

\subsubsection{Red zone.}

We first recenter. Changing ${\theta_1}$ to the element $a + b\psi_C^{-1}{\theta_1}$ from Lemma \ref{lem:recentering_ur}, keeping the rest of the resolvent data fixed, gives us a new first vector problem $\fP'$ whose associated cubic ring $C'$ has first extender vector $\theta_1'=\eta^2$, a square. By Lemma \ref{lem:recenter}, the ring-count function $W_{\fP}$ simply shifts by $\psi_C$. Note that all boxgroups $T(\ell_0, \ell_1, \ell_2)$ in claimed totals satisfy the gray-green inequality
\[
  \ell_1 \leq \frac{s + \square_C + 1}{2}
\]
Since $\theta_1' \equiv {\theta_1} \mod \pi^{s + \square_C}$, the replacement does not change any of the boxgroups, and their charmed cosets merely translate by $[\psi]$ along with the quadratic form $\epsilon_{C'}(\delta) = \epsilon_C(\psi\delta)$.

So it suffices to prove the result in the case that ${\theta_1}$ is a square. Note that $[\hat\omega_C] \in \cF_{\floor{e/2}}$ since $\omega_C \equiv \theta_1 \mod 2$.

We will actually prove something stronger:
\begin{equation}\label{eq:red_F}
  \begin{aligned}
  W_{m_{11},n_{11}} &= \sum_{\floor{\frac{n_c}{2}} \leq \ell < \tilde{n}} q^{-m_{11} + \ell}(1 + \epsilon_C) F^\cross\bigg( \ell, \Big\lceil\frac{n_c}{2}\Big\rceil + \frac{s}{2}, e - \Big\lceil\frac{n_c}{2}\Big\rceil - \frac{s}{2} - \ell \bigg) \\
  &\quad {} + q^{-m_{11} + \floor{\frac{2e - s - n_c}{4}}} F(\tilde{n}, e - 2\tilde{n}, \tilde{n}),
\end{aligned}
\end{equation}
in which we have replaced all $G$'s by $F$'s. We take a moment to realize why this is actually stronger. The claim that $T = T(\tilde{n}, e - 2\tilde{n}, \tilde{n})$ is in the support of $W_{m_{11},n_{11}}$ implies the following:
\begin{itemize}
  \item $\epsilon_C$ is identically $1$ on $T$, and hence
  \item $T$ is maximal isotropic for the Hilbert pairing, and also
  \item the identity coset $T$ is charmed, so $F_T = G_T$, and
  \item the quadratic form $\epsilon_C$ is positively charmed.
\end{itemize}
All the remaining terms use $F_V$ where $V \supseteq T$, so $F$ is interchangeable with $G$ there too. We now prove \eqref{eq:red_F}.

Begin with an arbitrary element in the box,
\[
\beta = 1 + \pi^{n_c} b {\theta_1}.
\]
Assume first that $\pi^{n_c} b$ is \emph{not} a square modulo $\pi^{e - \frac{s}{2} - \ceil{\frac{n_c}{2}} + 1}$, and let $k$ be the largest integer such that $\pi^{n_c} b$ is a square modulo $\pi^{2k+1}$. Note that $\floor{\frac{n_c}{2}} \leq k < \tilde{n}$. We may write
\[
\pi^{n_c} b = \left(\pi^{\ceil{\frac{n_c}{2}}}a\right)^2 + \pi^{2k + 1} c, \quad \pi \nmid c.
\]
Let $\zeta = 1 + \pi^{\ceil{\frac{n_c}{2}}}a \eta$. We claim that
\[
\frac{\beta}{\zeta^2} \in 1 + B_{\theta_1}\left( 2\ceil{\frac{n_c}{2}} + s + 2k + 1, 2k + 1 \right),
\]
implying that $[\beta] \in T\left( k, \ceil{\frac{n_c}{2}} + \frac{s}{2}, e - \ceil{\frac{n_c}{2}} - \frac{s}{2} - k \right)$ (compare the $k$th term of the sum). Write
\begin{align*}
\frac{\beta}{\zeta^2}
&= \frac{1 + \left( \pi^{\ceil{\frac{n_c}{2}}}a \right)^2{\theta_1} + \pi^{2k + 1} c {\theta_1}}{\zeta^2} \\
&= \frac{\left[ 1 + \left( \pi^{\ceil{\frac{n_c}{2}}}a \right)^2 {\theta_1}\right] \left( 1 + \pi^{2k+1} c {\theta_1} \right) - \pi^{2\ceil{\frac{n_c}{2}} + 2k + 1} a^2 c \theta_1^2}{\zeta^2} \\
&= \frac{\left[ \left( 1 + \pi^{\ceil{\frac{n_c}{2}}}a \eta \right)^2 - 2\pi^{\ceil{\frac{n_c}{2}}} a \eta\right] \left( 1 + \pi^{2k+1} c {\theta_1} \right) - \pi^{2\ceil{\frac{n_c}{2}} + 2k + 1} a^2 c {\theta_1}^2}{\zeta^2} \\
&= \frac{\left( \zeta^2 - 2\pi^{\ceil{\frac{n_c}{2}}} a \eta \right) \left( 1 + \pi^{2k+1} c \theta_1 \right) - \pi^{2\ceil{\frac{n_c}{2}} + 2k + 1} a^2 c \theta_1^2}{\zeta^2}.
\end{align*}

We first claim that the denominator $\zeta^2$ belongs to $\OO_K[\theta_1] = \OO_K + B_{\theta_1}(s,0)$. Since
\[
\zeta^2 = \left( 1 + \pi^{\ceil{\frac{n_c}{2}}}a \eta\right)^2 = 1 + \pi^{2\ceil{\frac{n_c}{2}}}a^2 \theta_1 + 2 \pi^{\ceil{\frac{n_c}{2}}}a \eta,
\]
only the last term is in question, and since $v_K\left(\coef_{\theta_2} \eta\right) \geq s/2$ (as we saw in Lemma \ref{lem:eta_ur}), the inequality needed is
\[
e + \ceil{\frac{n_c}{2}} + \frac{s}{2} \geq s,
\]
a consequence of $s < 2e$.

Therefore the last term of $\beta/\zeta^2$ is
\begin{align*}
\frac{-\pi^{2\ceil{\frac{n_c}{2}} + 2k + 1} a^2 c \theta_1^2}{\zeta^2} &\in \pi^{2\ceil{\frac{n_c}{2}} + 2k + 1} \OO_K[\theta_1] \\
&= B_{\theta_1}\left(2\ceil{\frac{n_c}{2}} + s + 2k + 1, 2\ceil{\frac{n_c}{2}} + 2k + 1\right) \\
&\subseteq B_{\theta_1}\left(2\ceil{\frac{n_c}{2}} + s + 2k + 1, 2k + 1\right).
\end{align*}
Thus it is enough to prove that
\begin{align*}
\frac{\left( \zeta^2 - 2\pi^{\ceil{\frac{n_c}{2}}} a \eta \right) \left( 1 + \pi^{2k+1} c \theta_1 \right) }{\zeta^2} &\in 1 + B_{\theta_1}\left(2\ceil{\frac{n_c}{2}} + s + 2k + 1, 2k + 1\right).
\end{align*}
Since the right-hand side is a group and contains $1 + \pi^{2k+1} c \theta_1$, it is enough to show that
\[
\frac{\zeta^2 - 2\pi^{\ceil{\frac{n_c}{2}}} a \eta}{\zeta^2} \in 1 + B_{\theta_1}\left(2\ceil{\frac{n_c}{2}} + s + 2k + 1, 2k + 1\right),
\]
that is,
\[
\frac{2\pi^{\ceil{\frac{n_c}{2}}} a \eta}{\zeta^2} \in B_{\theta_1}\left(2\ceil{\frac{n_c}{2}} + s + 2k + 1, 2k + 1\right).
\]
But, since $\eta$ belongs to the ring $B_{\theta_1}(s/2, 0)$ and $\zeta^2$ is a unit in that ring,
\begin{align*}
\frac{2\pi^{\ceil{\frac{n_c}{2}}} a \eta}{\zeta^2} &\in \pi^{e + \ceil{\frac{n_c}{2}}} B_{\theta_1}\left( \frac{s}{2}, 0 \right) \\
&= B_{\theta_1} \left( e + \ceil{\frac{n_c}{2}} + \frac{s}{2}, e + \ceil{\frac{n_c}{2}} \right) \\
&\subseteq B_{\theta_1}\left(2\ceil{\frac{n_c}{2}} + s + 2k + 1, 2k + 1\right),
\end{align*}
where the last step uses $k < \tilde{n}$. This establishes the claim that
\[
[\beta] \in T\left( k, \ceil{\frac{n_c}{2}} + \frac{s}{2}, e - \ceil{\frac{n_c}{2}} - \frac{s}{2} - k \right).
\]
To replace the $T$ by $T^\cross$, note that
\[
\frac{\beta}{\zeta^2} \equiv 1 + \pi^{2k+1} c \theta_1 \mod \pi^{2k + 2},
\]
a generic unit of exact level $k$. Hence
\[
\beta \in T^\cross\left( k, \ceil{\frac{n_c}{2}} + \frac{s}{2}, e - \ceil{\frac{n_c}{2}} - \frac{s}{2} - k \right).
\]
Also $\ell(\beta\hat\omega_C) = k$, so at $\delta = \beta$ there is a green conic of squareness $k$ and
\[
W_{m_{11},n_{11}}(\beta) = 2q^{-m_{11} + k}.
\]

If it so happens that $b\pi^{n_c}$ \emph{is} a square modulo $\pi^{e - \frac{s}{2} - \ceil{\frac{n_c}{2}} + 1}$, then writing $\pi^{n_c} b = \left(\pi^{\ceil{\frac{n_c}{2}}}a\right)^2 + \pi^{2\tilde{n}} c$ and carrying out the above computations, \emph{mutatis mutandis,} shows that
\[
[\beta] \in T(\tilde{n}, e - 2\tilde{n}, \tilde{n})
\]
and
\[
\ell(\beta\hat\omega_C) \geq \tilde n,
\]
so there is a blue conic at $\beta$ and
\[
W_{m_{11},n_{11}} = q^{-m_{11} + \floor{\frac{2e - s - n_c}{4}}}.
\]
Overall,
\begin{align*}
W_{m_{11},n_{11}} &\leq \sum_{\floor{\frac{n_c}{2}} \leq \ell < \tilde{n}} q^{-m_{11} + \ell}(1 + \epsilon_C)F^\cross\left( \ell, \ceil{\frac{n_c}{2}} + \frac{s}{2}, e - \ceil{\frac{n_c}{2}} - \frac{s}{2} - \ell \right) \\
&\quad {} + q^{-m_{11} + \floor{\frac{2e - s - n_c}{4}}} F(\tilde{n}, e - 2\tilde{n}, \tilde{n})
\end{align*}
This completes the bounding step. For the summing step, we note that the support $T^\cross\left( \ell, \ceil{\frac{n_c}{2}} + \frac{s}{2}, e - \ceil{\frac{n_c}{2}} - \frac{s}{2} - \ell \right)$ is a union of cosets of $\cF_{e - \ell - 1}$ that do not lie in the charmed coset $\cF_{\ell + 1}$, so $\epsilon_C$ is equidistributed. The summation then proceeds routinely.
\end{proof}

\subsubsection{Further remarks on the red zone.} For general $C$, we end up proving that the coset $\psi_C T$ is charmed for each $T = T(\ell_0, \ell_1, \ell_2)$ appearing (either positively or negatively) in the sum. Now $\psi_C \in 1 + B_{\theta_1}(\infty, \square_C)$, which, if
\begin{equation} \label{eq:G_is_F}
  \ell_0 \leq \ell_C,
\end{equation}
is contained in the box defining $T$. Thus the identity coset is charmed and we can replace $G$ by $F$.

\subsection{Splitting type (\texorpdfstring{$1^3$}{1³})}

\begin{lem}\label{lem:1pow3_strong_zones}
  Suppose $R$ has splitting type $(1^3)$.
  Let $ m_{11} $ and $ n_{11} $ be rational numbers with
  \[
    m_{11} \in \ZZ - \frac{h}{3}, \quad n_{11} \in \ZZ + \frac{h}{3}, \quad m_{11} > 2e, \quad n_{11} > 0.
  \]
  Write
  \[
    \tri m = m_{11} - \frac{2h}{3} \in \ZZ, \quad \tri n = n_{11} + \frac{2h}{3} \in \ZZ.
  \]
  Let $\square_C = \min\left\{2\floor{\ell(\hat\omega_C) / 2} + 1, e\right\}$ be the squareness of the $\delta = 1$ conic. Also let
  \[
  \tilde{n} = \floor{\frac{2e - \tri n + 2 + 2h}{4}}.
  \]
  Then $ W_{m_{11},n_{11}} $ is given as follows:
  \begin{enumerate}[$($a$)$]
    \item If $n_{11} > 2e$ (black zone), then
    \[
      W_{m_{11},n_{11}} = q^{2e - m_{11} - n_{11} + 1} F\left( e, 0, 0 \right).
    \]
    \item If $ 2e - 2\square_C + 1 + \dfrac{4h}{3} \leq n_{11} \leq 2e $ and $ n_{11} \geq \dfrac{2e}{3} $ (blue zone), then
    \[
    W_{m_{11},n_{11}} = q^{-\tri m + 1 - h + \floor{\frac{2e - \tri n + 2h}{4}}} F_h\left( \floor{\frac{\tri n}{2}} , e - \floor{\frac{\tri n}{2}} - \tilde{n}, \tilde{n} \right)
    \]
    \item If $ \square_C + 1 - \dfrac{2h}{3} \leq n_{11} \leq 2e - 2\square_C + \dfrac{4h}{3} $ (green zone), then
    \[
    W_{m_{11},n_{11}} = q^{-\tri m + \ell_C + 1 - h}(1 + \epsilon_C)F_h\bigg( \Big\lfloor{\frac{\tri n}{2}}\Big\rfloor, \ell_C - h + \One_{2\nmid \tri n}, e - \Big\lceil{\frac{\tri n}{2}}\Big\rceil - \ell_C + h\bigg)
    \]
    where $\ell_C = \frac{\square_C - 1}{2} \in \ZZ_{\geq 0}$.
    \item If $n_{11} < \dfrac{2e}{3}$ and $n_{11} \leq \square_C - \dfrac{2h}{3}$ (red zone), then
    \begin{align*}
    W_{m_{11},n_{11}} &= \sum_{\floor{\frac{\tri n}{2}} \leq \ell < \tilde{n}} q^{-\tri m + 1 - h + \ell}(1 + \epsilon_C)G_h^\cross\left( \ell, \ceil{\frac{\tri n}{2}} - h, e - \ceil{\frac{\tri n}{2}} - \ell + h \right) \\
    &\quad{} + q^{-\tri m + 1 + \floor{\frac{2e - \tri n - 2h}{4}}} G_h(\tilde{n}, e - 2\tilde{n}, \tilde{n}) \\
    \end{align*}
  \end{enumerate}
\end{lem}
\begin{proof}

  By Lemma \ref{lem:beta}, the support of $W_{m_{11},n_{11}}$ consists of the classes $[\delta] = [\beta]$ in $H^1$ of elements $\beta$ of the box
  \[
  1 + B_{m_{11},n_{11}} = \{1 + y \pi^{n_{11}} \theta_1 + z \pi^{m_{11}} \theta_2 : y,z \in \OO_K\}.
  \]
  Since $m_{11} > 2e$, the $z$ term has no effect on $[\beta]$, and we ignore it. In particular, the conic
  \[
  \coef_{\theta_2} (\beta \xi^2) = 0, \quad \text{that is,} \quad \tr(\hat\omega_C \beta \xi^2) = 0
  \]
  has a solution $\xi = 1$, so $\epsilon(\hat\omega_C \beta) = 1$ for all $\beta$ in the box.

  In this case, the summation lemma (Lemma \ref{lem:sum_strong}) yields
  \[
    \sum_{\delta \in \cF_0} W_{m_{11},n_{11}}(\delta) = q^{2e - m_{11} - n_{11} + 1},
  \]
  in particular verifying the total in the black zone.

  Again, let $ W_{m_{11},n_{11}}' $ denote the claimed value of $W_{m_{11},n_{11}}$ in each case. Our proof method will consist of bounding and summing, as in the preceding splitting types.

  \subsubsection{Blue zone.} Our strategy is to note that
  \[
  \beta \in 1 + B_{m_{11},n_{11}} \subseteq 1 + B_{m', n_{11}}
  \]
  for some $m' \leq m_{11}$ for which $[1 + B_{m',n_{11}}]$ is a boxgroup. Here, we find that the gray-blue condition \eqref{eq:box_blue_1pow3} in Lemma \ref{lem:boxgps_1pow3} is the most stringent one, so we take
  \[
    \tri{m'}= \ceil{\frac{\tri n + 1}{2}} + e - h \textand m' = \tri{m'} + \frac{2h}{3}
  \]
  We easily verify that
  \[
    e - \floor{\frac{\tri{m'}}{2}} = \tilde{n},
  \]
  so $[1 + B_{\theta_1}(m', n_{11})]$ is exactly the support of $W'_{m_{11},n_{11}}$. We claim all conics are blue. This requires that the squareness
  \[
    \square = \min \left\{2\floor{\frac{\ell(\beta\hat\omega_C)}{2}}+1, e\right\}
  \]
  satisfy
  \[
    \square \stackrel{?}{\geq} e - \ceil{\frac{\tri n}{2}} + h,
  \]
  which simplifies to
  \[
    \beta \hat\omega_C \in \cF_{\floor{\frac{2e - \tri n + 2h}{4}}}.
  \]
  When $\beta = 1$, we have $\ell = \ell_C$, and the required relation follows from the given blue-green inequality $n_{11} \geq 2e - 2\square_C + 1 + \dfrac{4h}{3}$. So it suffices to show that
  \[
  \beta \in \cF_{\floor{\frac{2e - \tri n + 2h}{4}}}.
  \]
  Using the known relation $\ell(\beta) \geq 2\floor{\tri n / 2}$
  and the blue-red inequality $n_{11} \geq 2e/3$, this is not hard to prove. So all conics are blue, and the only possible nonzero value of $W_{m_{11},n_{11}}(\delta)$ is
  \[
    q^{-\tri m + \floor{\frac{2e - \tri n + 2h}{4}}}.
  \]
  This completes the bounding step. The summing step is routine. In particular, this shows that $\epsilon = 1$ identically on boxgroups of the shape in the lemma. This result will be important in proving the remaining zones.

  \subsubsection{Green zone.} Again, we write
  \[
  \beta \in 1 + B_{m_{11},n_{11}} \subseteq 1 + B_{m', n_{11}}
  \]
  where $m'$ is as large as possible to make a boxgroup. This time, we find that the gray-green inequality \eqref{eq:box_green_1pow3} is the most stringent of the conditions in Lemma \ref{lem:boxgps_1pow3}, so we take $m' = n_{11} + \square_C - \frac{2h}{3}$, that is,
  \[
    \tri{m'} = m' - \frac{2h}{3} = \tri n + \square_C - 2h \in \ZZ.
  \]
  and find that the support of $W_{m_{11},n_{11}}$ is contained in
  \[
  \left[1 + B_{\theta_1}\left(m', n_{11}\right)\right] = T_h\left( \floor{\frac{\tri n}{2}}, \ell_C - h + \One_{2\nmid \tri n}, e - \ceil{\frac{\tri n}{2}} - \ell_C + h\right).
  \]
  We claim all conics are green of the same squareness $\square_C$. It is easy to prove that $\square_C < e - 1$ in this zone, so
  \[
  [\hat\omega_C \diamondsuit \heartsuit] \in \cF_{2 \ell_C} \setminus \cF_{2 \ell_C + 2}.
  \]
  We then note that $[\beta] \in \cF_{2 \ell_C + 2}$, because $[\beta] \in \cF_{2\floor{\tri n/2}}$ and we have the inequality $\tri n \geq 2\ell_C + 2$. So $[\beta\hat\omega_C \diamondsuit \heartsuit]$ is also of exact level $2\ell_C$ or $2\ell_C + 1$. Thus all conics are green of the same squareness, and by Lemmas \ref{lem:conic_1} and \ref{lem:conic_pi}, we have the bound
  \begin{align*}
  W_{m_{11},n_{11}} &\leq 2 q^{-\tri m + 1 - h + \ell_C} F_h\left( \floor{\frac{\tri n}{2}}, \ell_C - h + \One_{2\nmid \tri n}, e - \ceil{\frac{\tri n}{2}} - \ell_C + h\right),
  \end{align*}
  indeed
\[
  W_{m_{11},n_{11}} \leq q^{-\tri m + \ell_C + 1 - h}(1 + \epsilon_C)F_h\bigg( \Big\lfloor{\frac{\tri n}{2}}\Big\rfloor, \ell_C - h + \One_{2\nmid \tri n}, e - \Big\lceil{\frac{\tri n}{2}}\Big\rceil - \ell_C + h\bigg).
\]
This completes the bounding step.

  To perform the summing step, we need to compute the sum of $q^{-\tri m + \ell_C} (1 + \epsilon_C)$ over the stated boxgroup. The term $1$ is found to sum to the desired total $\size{H^0} \cdot q^{2e - m_{11} - n_{11} + h}$. We claim that
  \[
  \sum_{\delta \in T\left( 2\floor{\frac{\tri n}{2}} + \frac{1 - h}{2}, 2\ell_C - 1 - h + \One_{2\nmid \tri n}, 2e - 2\ell_C - 2\ceil{\frac{\tri n}{2}} + \frac{3h + 1}{2}\right)} \epsilon_C(\delta) = 0,
  \]
  in other words that $\epsilon$ is equidistributed between $1$ and $-1$ in this boxgroup. This follows from Lemma \ref{lem:charm}\ref{charm:constant}: because $\cF_{2\ell_C + 2} \neq \hat\omega_C \cF_{2\ell_C + 2}$ is an uncharmed coset, we have $\epsilon$ equidistributed on cosets of $\cF_{2e - 2\floor{\ell_C} - 2}$, of which the boxgroup in question is a union by the green-blue inequality.

  \subsubsection{Red zone.} Considerations of space prevent us from writing out the proof, which is like that in the unramified splitting types with the following changes:
  \begin{itemize}
    \item We reduce to the case that $\theta_1 = \eta^2$ is a square using Lemma \ref{lem:recentering_1pow3}, and there we will prove the result with the $G$'s replaced by $F$'s. We begin with an arbitrary
    \[
    \beta = 1 + \pi^{n_{11}} b \theta_1 = 1 + \pi^{\tri{n}}b \cdot \pi^{-2h/3}\theta_1, \quad b \in \OO_K.
    \]
    \item We assume first that $\pi^{\tri n} b$ is \emph{not} a square modulo $\pi^{e - \ceil{\frac{\tri n}{2}} + 1 + h}$, and let $k$ be the largest integer such that $\pi^{\tri n} b$ is a square modulo $\pi^{2k+1}$. We find that
    \[
      [\beta] \in T^\cross\left( k, \ceil{\frac{\tri n}{2}} - h, e - \ceil{\frac{\tri n}{2}} - k + h \right)
    \]
    is in the support of the $k$th term of the claimed answer and that there is a green conic of level $k$ there.
    \item If $\pi^{\tri n} b$ in fact \emph{is} a square modulo $\pi^{e - \ceil{\frac{\tri n}{2}} + 1 + h}$, then we find that $[\beta]$ is in the support of the last term of the claimed answer and that there is a blue conic there.
  \end{itemize}
The rest of the proof, including the summing step, is completely like the unramified splitting types. This completes the proof.
\end{proof}

The $G$'s can be replaced by $F$'s when the index $k$ is at most $\ell_C$, for reasons just like those named above.

\subsection{Splitting type (\texorpdfstring{$1^21$}{1²1})}
In this splitting type, we first note that, if $\cN_{11}$ is strongly active, it forces all the valuations of $\xi_1 \in \pi^{-a_1}\sqrt{\pi^s\delta\omega_C}$ to match those of $\omega_C \in \pi^{-4b_1 - s} R$. Since a $\beta \in \pi^{i}R$ cannot have $\beta \sim 1$ unless $i \in \ZZ$, we obtain the class of $a_{1}$ mod $\ZZ$:
\[
  a_1 \equiv -\frac{s}{2} \mod \ZZ.
\]
In particular, this means that $m_{11}$ and $n_c = n_{11} - s$ are also known mod $1$:
\begin{equation}
  m_{11} \equiv \frac{d_0}{2} + s', \quad n_c \equiv s' \mod \ZZ.
\end{equation}
In particular, $n_c$ belongs to $\ZZ + 1/2$ in letter type \ref{type:A}, and to $\ZZ$ in letter types \ref{type:B}--\ref{type:E}.
\begin{lem}\label{lem:1pow21_strong_zones}
Suppose $R$ has splitting type $(1^21)$.
Let $ m_{11} $ and $ n_{11} $ be rational numbers whose translates
\[
  m_c = m_{11} - \frac{d_0}{2} \textand n_c = n_{11} - s
\]
satisfy $m_c \geq 2e$ and $n_c > 0$ and lie in the correct congruence classes modulo $1$ for the letter type of $\theta_1$. Let
\[
\square_C = \min\left\{2\floor{\ell(\hat\omega_C\diamondsuit\heartsuit) / 2} + 1, e\right\}
\]
be the squareness of the $\delta = 1$ conic. Also, for $d_0$ even, in letter types \ref{type:C} and \ref{type:D}, let
\begin{align*}
  \tilde{n} &= \floor{\frac{2e - s' - d_0' - n_c - 2h_\eta + 2}{4}}
  = \ceil{\frac{e - \ceil{\frac{s' + d_0' + n_c}{2}} - h_\eta}{2}} \\
  \tilde{n}^- &= \floor{\frac{2e - s' - d_0' - n_c - 2h_\eta}{4}}
  = \floor{\frac{e - \ceil{\frac{s' + d_0' + n_c}{2}} - h_\eta}{2}}.
\end{align*}
Then $ W_{m_{11},n_{11}} $ is given as follows:
\begin{enumerate}[$($a$)$]
  \item If $n_c > 2e$ (black zone), then
  \[
  W^{\odot}_{m_{11},n_{11}} = 2q^{2e - m_{11}^\odot - n_c - \frac{s'}{2} + \frac{h_1}{4} - \floor{\frac{h_1}{2}}} L\left( 2e + 1 \right).
  \]
  \item If $2e - d_0 + 2 \leq n_c \leq 2e$ and $n_c \geq 2e - 2s'$ (plum zone; note $s' \geq 0$), then
  \[
  W^{\odot}_{m_{11},n_{11}} = q^{e - m_{11}^\odot - \ceil{\frac{n_c}{2}} - \ceil{\frac{s'}{2}}} L\left(e + \floor{\frac{n_c}{2}}\right).
  \]
  which, upon simplification, gives
  \[
    W_{m_{11},n_{11}} = q^{e - m_{11} + \frac{d_0 + o_{a_1} - 2s'}{4} - \ceil{\frac{n_c}{2}}} L\left(e + \floor{\frac{n_c}{2}}\right).
  \]
  \item If $2e - d_0/2 - s' + 1 \leq n_c \leq 2e - d_0 + 1$ (purple zone; note that we must be in letter type \ref{type:D} or \ref{type:E}), then
  \[
  W^{\odot}_{m_{11},n_{11}} = 2q^{e - m_{11}^\odot - \ceil{\frac{n_c}{2}} - \ceil{\frac{s'}{2}}} F\left(\floor{\frac{n_c}{2}}, e - d_0' - \floor{\frac{n_c}{2}}, \emptyset\right).
  \]
  which, upon simplification, gives
  \[
    W_{m_{11},n_{11}} = 2 q^{e - m_{11} + \frac{d_0 + o_{a_1} - 2s'}{4} - \ceil{\frac{n_c}{2}}}
    F\left(\floor{\frac{n_c}{2}}, e - d_0' - \floor{\frac{n_c}{2}}, \emptyset\right).
  \]
  \item In letter types \ref{type:C} and \ref{type:D}, if $ 2e + d_0' - s' - 2\square_C + 1 - h_1 \leq n_c \leq 2e - d_0' - s'$ and $n_c \geq \dfrac{2e - s' - d_0'}{3}$ (blue zone), then
  \[
  W^{\odot}_{m_{11},n_{11}} = q^{-m_{11}^\odot + \floor{\frac{d_0' + h_\eta}{2}} + \tilde n^-} F\left(\floor{\frac{n_c - h_\eta}{2}}, e - d_0' - \floor{\frac{n_c - h_\eta}{2}} - \tilde n , \tilde n\right)
  \]
  which, upon simplification, gives
  \[
  W_{m_{11},n_{11}} = q^{-m_{11} + d_0' + \frac{2h_\eta + o_{a_1}}{4} + \tilde n^-} F\left(\floor{\frac{n_c - h_\eta}{2}}, e - d_0' - \floor{\frac{n_c - h_\eta}{2}} - \tilde n , \tilde n\right).
  \]
  \item In letter type \ref{type:A} (green zone),
  \[
    W^{\odot} = (1 + \epsilon_C) q^{-m_{11}^\odot} L\left(n_c - \frac{1}{2}\right)
  \]
  so (as $h_1 = 1$)
  \[
    W = (1 + \epsilon_C) q^{-m_{11} + \frac{d_0 - 1 + o_{a_1}}{4}} L\left(n_c - \frac{1}{2}\right).
  \]
  \item In letter type \ref{type:B}, if $n_c \leq 2e - 2s' - 1$ (green zone),
  \[
    W^{\odot} = (1 + \epsilon_C) q^{-m_{11}^\odot + \floor{s'/2}}L\left(n_c + s'\right),
  \]
  so
  \[
    W = (1 + \epsilon_C) q^{-m_{11} + \frac{d_0 + o_{a_1} + 2 s'}{4}}L\left(n_c + s'\right).
  \]
  \item In letter types \ref{type:C}--\ref{type:E}, if $\square_C - d_0' + h_1/2 < n_c \leq 2e + d_0' - s' - 2\square_C - h_1$ (green zone),
  \begin{align*}
    W^{\odot} &= (1 + \epsilon_C) q^{-m_{11}^\odot + \ell_C} F\bigg(\Big\lfloor\frac{n_c - h_\eta}{2}\Big\rfloor,
    \Big\lceil\frac{s'}{2}\Big\rceil + \ell_C - d_0' + h_\eta + \One_{2\nmid n_c + h_\eta}, \\
    &\quad {} e - \Big\lceil \frac{n_c + h_\eta}{2} \Big\rceil - \Big\lceil \frac{s'}{2} \Big\rceil - \ell_C\bigg).
  \end{align*}
  \item If $n_c < (2e - s' - d_0')/3$ and $n_c \leq \square_C - d_0' + h_1/2$ (red zone; note that we must be in type \ref{type:C} or \ref{type:D}), then
  \begin{align*}
    W_{m_{11},n_{11}}^{\odot} &= \sum_{\floor{\frac{n_c - h_\eta}{2}} \leq k < \tilde{n}} q^{-m_{11}^\odot + \floor{\frac{d_0' + h_\eta}{2}} + k}(1 + \epsilon_C)G_{h_\eta}^\cross\Bigg( k, \ceil{\frac{s' - d_0' + n_c}{2}} + h_\eta, \\
    &\qquad{} e - \ceil{\frac{s' + d_0' + n_c}{2}} - k - h_\eta \Bigg) + \\
    &\quad{} + q^{-m_{11}^\odot + \floor{\frac{d_0' + h_\eta}{2}} + \tilde n^{-} } G_{h_\eta}(\tilde{n}, e - d_0' - 2\tilde{n}, \tilde{n})
  \end{align*}
\end{enumerate}
\end{lem}

\begin{proof}
Specializing the summation lemma \ref{lem:sum_strong} to this splitting type and using the value of $v\left(N\left(\gamma\gamma^{\odot}\right)\right)$ given in Lemma \ref{lem:tfm_conic} yields the sum
\[
  \sum_{\delta \in \cF_0} W_{m_{11},n_{11}}^\odot(\delta) = 2q^{2e - m_{11}^\odot - n_c - \frac{s'}{2} + \frac{h_1}{4} - \floor{\frac{h_1}{2}}},
\]
which yields the desired total in the black zone.

\subsubsection{Plum zone.} Here the support is
\begin{equation*}
  1 + B_{\theta_1}(m_c, n_c) = \left\{[1 + c_1 \pi^{n_c} \theta_1] : c_1 \in \OO_K\right\}.
\end{equation*}
We scale so that $\theta_1 \equiv (1;0)$ mod $\pi^{s' + 1/2}$ and get
\begin{align*}
  1 + \pi^{n_c} \theta_1 &\in \left(1 + \pi^{n_c} (1;0)\right) \cdot \cF_{\min\{n_c + s', 2e\}} \\
  &\subseteq \iota\left(\cF_{\floor{n_c/2}}\right) \cdot \cF_{\min\{n_c + s', 2e\}} \\
  &= \cF_{e + \floor{n_c/2}} \cdot \cF_{\min\{n_c + s', 2e\}} \\
  &= \cF_{e + \floor{n_c/2}}.
\end{align*}
Also, the box $B = 1 + B_{\theta_1}(m_c, n_c)$ is a group, at least for $m_c = 2e$ which we can assume. We claim $[B]$ has full signature $\ZERO.\ZERO^{e+\floor{n_c/2}}.\STAR^{e - \floor{n_c/2}}.\STAR$ by downward induction on $n_c$. In the base case $n_c = 2e$, we only need to get the intimate unit, which we can get by putting an appropriate unit for $c_1$. The induction step only has content if $n_c$ is odd, and then as $c_1$ varies over $\OO_K/\pi\OO_K$, the resulting element $1 + c_1 \pi^{n_c} \theta_1$ ranges through the $q$-many cosets of $\cF_{e + (n_c-1)/2}/\cF_{e + (n_c+1)/2}$, as desired.

Again because $B$ is a group, the thickness is uniform and can be retrieved from the summation lemma.

\subsubsection{Purple zone.} The purple zone is done by the exact same method; now the supplementary boxgroup
\[
  \iota\left(\cF_{\floor{n_c/2}}\right) = T\left(\floor{\frac{n_c}{2}}, e - d_0' - \floor{\frac{n_c}{2}}, \emptyset\right)
\]
appears. 

\subsubsection{Blue zone.} As in the other splitting types, our strategy is to note that
\[
\beta \in 1 + B_{m_c,n_c} \subseteq 1 + B_{m', n_c}
\]
for some $m' \leq m_c$ for which $[1 + B_{m',n_c}]$ is a boxgroup. Here, we find that the gray-blue condition is the most stringent one, and we take
\[
  m' = \floor{e + \frac{n_c + s' - d_0' + 1}{2}}
\]
to get
\[
  [\beta] \in [1 + B_{m',n_c}] = T\left(\floor{\frac{n_c - h_\eta}{2}}, e - d_0' - \floor{\frac{n_c - h_\eta}{2}} - \tilde n , \tilde n\right).
\]
The thickness $W_{m_{11},n_{11}}(\delta)$ is found by noting that the conics are all blue. The summing step then proceeds analogously to the other splitting types.

\subsubsection{Green zone.} We have labeled three cases ``green zone,'' because although the details vary by letter type, they all have in common that $\square_C$ is so low that we get green conics of fixed squareness:
\begin{itemize}
  \item Case \ref{type:A} is the simplest. We have $n_c \in \ZZ + 1/2$, and clearly
  \[
    [\beta] \in [\alpha : \alpha \equiv 1 \mod \pi^{n_c}] = \cF_{n_c - 1/2}.
  \]
  The conics are all tiny and green. Since $\cF_0$ is uncharmed, $\epsilon_C$ is equidistributed on this level space, rendering the summing step easy.
  \item In case \ref{type:B}, we may assume that $\theta_1 \equiv (1;0) \mod \pi^{s' + 1/2}$, and then
  \begin{align*}
    \beta &\in \left(1 + \pi^{n_c}c_1(1;0)\right) \cF_{n_c + s'} \\
    &\subseteq \iota\left(1 + \pi^{n_c} \OO_K\right) \cdot \cF_{n_c + s'}.
  \end{align*}
  The bound $s' < d_0/2 - 1$ defining letter type \ref{type:B} ensures that the second factor dominates, establishing the bounding step. Since
  \[
    n_c + s' \leq 2e - s' - 1
  \]
  and $\hat\omega_C\diamondsuit\heartsuit$ has level $s'$, we get equidistribution of $\epsilon_C$ on the level space, rendering the summing step easy.
  \item Finally, cases \ref{type:C}--\ref{type:E} parallel the green zone of the other splitting types. We take
  \[
    m' = n_c + 2 \ceil{\frac{s'}{2}} - d_0' + \square_C,
  \]
  making the gray-green inequality an equality, which is the most stringent of the boxgroup-defining inequalities, and get the desired bound. The equidistribution of $\epsilon_C$ holds by level considerations, and the summing step is routine.
\end{itemize}

\subsubsection{Red zone.} Considerations of space prevent us from writing out the proof, which is like that in the unramified splitting types with the following changes:
\begin{itemize}
  \item We reduce to the case that $\pi^{h_\eta}\theta_1 = \eta^2$ is a square using Lemma \ref{lem:recentering_1pow21}, and there we will prove the result with the $G$'s replaced by $F$'s. We begin with an arbitrary
  \[
  \beta = 1 + \pi^{n_c} b \theta_1 = 1 + \pi^{n_c - h_\eta} b \eta^2, \quad b \in \OO_K.
  \]
  \item Assume first that $\pi^{n_c - h_\eta} b$ is \emph{not} a square modulo $\pi^{2 \tilde n}$, and let $k$ be the largest integer such that $\pi^{n_c - h_\eta} b$ is a square modulo $\pi^{2k+1}$. Write
  \[
    \pi^{n_c - h_\eta} b = \left(\pi^{\ceil{\frac{n_c - h_\eta}{2}}} a \right)^2 + \pi^{2k + 1} c, \quad \pi \nmid c.
  \]
  \item Let $\zeta = 1 + \pi^{\ceil{\frac{n_c - h_\eta}{2}}} a \eta$. We will show that
  \[
    \frac{\beta}{\zeta^2} \in 1 + B, \quad B = B_{\theta_1}\left(2\ceil{\frac{n_c + h_\eta}{2}} + 2k + 1 + s', 2k + 1 + h_\eta\right),
  \]
  establishing that
  \begin{equation} \label{beta_in_boxgp_1pow21}
    [\beta] \in T\left(k, \ceil{\frac{s' - d_0' + n_c}{2}}, e - \ceil{\frac{s' + d_0' + n_c}{2}} - k \right).
  \end{equation}
  \item Transform
  \begin{equation}
    \frac{\beta}{\zeta^2} = \left(1 - \frac{2\pi^{\ceil{\frac{n_c - h_\eta}{2}}}a \eta}{\zeta^2}\right)\left(1 + \pi^{2k+1+h_\eta} c \theta_1\right) - \frac{\pi^{2\ceil{\frac{n_c + h_\eta}{2}} + 2k + 1} a^2 c \theta_1^2}{\zeta^2}.
  \end{equation}
  \item Check that $\eta \in \OO_K[\theta_1] = \OO_K + B_{\theta_1}(s',0)$, using the last part of Lemma \ref{lem:eta_1pow21} and the relation $s' < 2e - d_0'$ needed for the red zone to be nonempty.
  \item Deduce (using the appropriate bounds) that the summands
  \begin{equation}\label{eq:3_summands_1pow21}
    \frac{2\pi^{\ceil{\frac{n_c - h_\eta}{2}}}a \eta}{\zeta^2}, \quad \pi^{2k+1+h_\eta} c \theta_1, \textand \frac{\pi^{2\ceil{\frac{n_c + h_\eta}{2}} + 2k + 1} a^2 c \theta_1^2}{\zeta^2}
  \end{equation}
  all lie in the requisite box $B$, establishing \eqref{beta_in_boxgp_1pow21}.
  \item To replace $T$ by $T^\cross$ in \eqref{beta_in_boxgp_1pow21}, reexamine the moving case of Lemma \ref{lem:boxgps_1pow21}. When the second box exponent changes from $2k + 2 + h_\eta$ to $2k + 1 + h_\eta$, the new elements all have exact level $2k + h_\eta + d_0'$. Thus it suffices to show that
  \[
    \frac{\beta}{\zeta^2} \notin 1 + B', \quad B' = B_{\theta_1}\left(2\ceil{\frac{n_c + h_\eta}{2}} + 2k + 1 + s', 2k + 2 + h_\eta\right).
  \]
  Since first and third quantities in \eqref{eq:3_summands_1pow21} lie in $B'$ but the second does not, this is straightforward to see.
  \item If $\pi^{n_c - h_\eta}b$ in fact \emph{is} a square modulo $\pi^{2\tilde n}$, then we find by the same method that $[\beta]$ is in the support of the last term of the claimed ring total.
  \item Using bounds on levels, we find that the conics are green for the main sum and blue for the last term, and also that $\epsilon_C$ is equidistributed on each $T^\cross$ appearing in the main sum. The rest of the proof, including the summing step, is completely like the unramified splitting types. \qedhere
\end{itemize}
\end{proof}

\section[The zones when \texorpdfstring{$\cN_{11}$}{N11} is weakly active]{The zones when \texorpdfstring{$\cN_{11}$}{N11} is weakly active (brown, yellow, and lemon)}\label{sec:weak}

We now turn our attention to first vector problems such that $N_{11}$ is weakly active and $m_{11} > 2e$. Note that $s > 0$ induces a distinguished splitting $R = K \cross Q$; in particular, this section is only relevant in splitting types $(111)$, $(12)$, and $(1^2 1)$. The resolvent conditions $\cM_{11}$ and $\cN_{11}$ simplify to
\begin{alignat*}{2}
  \cM_{11}&\colon &\tr(\xi_1^2) &\equiv 0 \mod \pi^{m_{11}} \\
  \cN_{11}&\colon &\xi_1^{(K)} &\equiv 0 \mod \pi^{n_{11}/2}
\end{alignat*}
For $\xi_1$ to satisfy this, its two $Q$-components must be units.
Under the transformation of Lemma \ref{lem:tfm_conic}, the conditions can also be written as
\begin{alignat*}{2}
  \cM_{11}\colon &&\lambda^\diamondsuit(\delta^\odot{\xi_1^\odot}^2) &\equiv 0 \mod \pi^{m_{11}^\odot} \\
  \cN_{11}\colon &&\xi_1^{\odot(K)} &\equiv 0 \mod \pi^{\frac{n_{11}}{2} - v^{(K)}(\gamma\gamma^\odot)} \\
  \iff& &\xi_1^{\odot(K)} &\equiv 0 \mod \pi^{\ceil{\frac{n_{11}}{2} - v^{(K)}(\gamma\gamma^\odot)}}.
\end{alignat*}

For simplicity we let
\[
  n_{11}^\odot = \ceil{\frac{n_{11}}{2} - v^{(K)}(\gamma\gamma^\odot)},
\]
which sometimes differs slightly from the $n^\odot$ of Lemma \ref{lem:N11}.

In each of the three applicable splitting types, we will find a \emph{brown zone} where $n_{11}$ is so high that $\xi_1^{\odot(K)}$ can be taken to be $0$, so the only $\delta^\odot$ for which there is a solution are those where $\delta^{(Q)}$ is the class of an element on the traceless line of $Q$. Let $\clubsuit$ be such an element of the form
\[
  \clubsuit = (\alpha + \bar \alpha; \alpha - \bar\alpha; \bar\alpha - \alpha)
\]
where $\alpha$ is a generator for $\OO_Q$ as an $\OO_K$-module. Remarks are in order:
\begin{itemize}
  \item In splitting type $(111)$, we can take $\clubsuit = (1;1;-1)$.
  \item In splitting type $(12)$, taking $\alpha = (1 + \sqrt{D_0})/2$ gives $\clubsuit = (1; \bar\zeta_2\sqrt{D_0})$.
    \item In splitting type $(1^21)$, taking $\alpha = \pi_Q$ gives $\clubsuit = (\tr \pi_Q; \bar\zeta_2\sqrt{D_0})$.
\end{itemize}
Denote by $Y_{\fP}(\delta^\clubsuit)$ the volume of $\xi'_1$ satisfying the $\cM_{11}$ and $\cN_{11}$ conditions when
\[
  [\delta^\clubsuit] = [\hat\omega_C \clubsuit \delta] = [\diamondsuit\clubsuit \delta^\odot],
\]
that is, $Y_{\fP}$ is a translation of $W_{\fP}$:
\[
  Y_\fP(\delta^\clubsuit) = W_{\fP}\left(\hat\omega_C \clubsuit \delta^\clubsuit\right).
\]
Observe that the dependence on $\theta_1$ has been nullified and, if $m_{11}^\odot > 2e$, that $Y_{\fP}$ is supported on the vanishing locus of the quadratic form
\[
  \epsilon^\clubsuit(\delta^\clubsuit) \coloneqq \epsilon(\clubsuit\delta^\clubsuit).
\]
Let $Y^\odot_{\fP}$ be the corresponding volume of $\xi^\odot_1$.

As in the strong zones, we need a summation lemma.
\begin{lem} \label{lem:sum_weak}
We have
\[
  \sum_{\delta^\clubsuit \in \delta_0^\clubsuit \cF_0} Y_{\fP}(\delta^\clubsuit) = \size{H^0} q^{2e - m_{11}^\odot - n_{11}^\odot + 2v_K\left(\delta^{\odot(Q)}\right)},
\]
where $m_{11}^\odot$, $\gamma$, and $\delta^\odot$ are given by Lemma \ref{lem:tfm_conic} and depend only on the class $\delta^\clubsuit\cF_0 \in H^1/\cF_0$.
\end{lem}
\begin{proof}
Since $\xi_1^{(Q)}$ must be a unit, $\xi_1^{\odot(Q)}$ has fixed valuation
\[
  v^{(Q)}\left(\xi_1^\odot\right) = -v^{(Q)}\left(\gamma\gamma^\odot\right) \eqqcolon r.
\]
Meanwhile, $v\left(\xi_1^{\odot(K)}\right) \geq n_{11}^\odot$ can vary. For $i \geq n^\odot_{11}$, let $Y_{\fP}^{(i)}(\delta^\clubsuit)$ denote the volume of $\xi_1^\odot$ satisfying the first vector problem $\fP$ and having $v_K(\xi^{\odot(K)}) = i$. We have
\[
  Y_{\fP} = \sum_{i = n_{11}^\odot}^\infty Y_{\fP}^{(i)},
\]
since only the measure-zero set where $\xi^{\odot(K)} = 0$ has been dropped. Let
\begin{equation} \label{eq:beta_weak}
  \beta = \frac{\delta^\odot {\xi^\odot}^2}{\pi^{2r}}
\end{equation}
Observe that $2r$ is an integer (by reference to Table \ref{tab:tfm_conic}) so $\beta \in R$. Moreover, as we vary $[\delta^\clubsuit] \in \delta^\clubsuit \cF_0$ and $\xi_1^\odot$ in the set whose volume is $Y_{\fP}^{(i)}$, we get that $\beta$ varies in the region of primitive members of $\OO_R$ such that
\begin{align*}
  v_K(\beta^{(K)}) &= v_K\left(\delta^{\odot(K)}\right) + 2i - 2r \\
  v_K(\beta^Q) &= 0 \\
  \lambda^\diamondsuit(\beta) &\equiv 0 \mod \pi^{m_{11}^\odot - 2r}.
\end{align*}
Conditions of the form $\beta^{(K)} \equiv 0 \mod \pi^{n_{11}}$ and $\lambda^\diamondsuit(\beta) \equiv 0 \mod \pi^{m_{11}}$ cut out a box of volume $q^{-m_{11}-n_{11}}$; so the volume of $\beta$ is
\begin{align*}
  &\left(1 - \frac{1}{q}\right)q^{-\left(v_K\left(\delta^{\odot(K)}\right) + 2i - 2r\right) - \left(m_{11}^\odot - 2r\right)} \\
  &= \left(1 - \frac{1}{q}\right) q^{-2i - m_{11}^\odot - v_K\left(\delta^{\odot(K)}\right) + 4r}.
\end{align*}
Consider the sequence
\[
  \xi' \longmapsto \frac{\xi'}{\left(\pi^i; \pi_Q^{2r}\right)}
  \longmapsto \frac{\xi'^2}{\left(\pi^{2i}; \pi_Q^{4r}\right)}
  \longmapsto \frac{\delta^\odot \xi'^2}{\pi^{2r}} = \beta.
\]
Every term is a primitive vector in $\OO_R$, so we can consider projective volumes. The linear map of dividing by $\left(\pi^i; \pi_Q^{2r}\right)$ scales volumes by $q^{i + 2r}$. Squaring by units scales volumes by $1/\left(\size{H^0} \cdot q^{2e}\right)$ on regions symmetric under multiplication by $R^\cross[2]$, as we noted above in the proof of Lemma \ref{lem:sum_strong}. Finally, multiplication by
\[
  \frac{\delta^\odot \cdot \left(\pi^{2i}; \pi_Q^{4r}\right)}{\pi^{2r}}
\]
scales volumes by
\[
  q^{-v_K\left(\delta^{\odot(K)}\right) - 2v_K\left(\delta^{\odot(Q)}\right) - 2i + 2r}.
\]
So overall, a volume $Y_{\fP}^{(i)}(\delta)$ of $\xi'$ transforms to a volume of
\[
  \frac{1}{\size{H^0}} q^{-2e-i-v_K\left(\delta^{\odot(K)}\right) - 2v_K\left(\delta^{\odot(Q)}\right) + 4r} Y_{\fP}^{(i)}(\delta).
\]
Summing over $[\delta] \in \delta_0\cF_0$,
\[
  \begin{split}
  \lefteqn{\frac{1}{\size{H^0}} q^{-2e-i-v_K\left(\delta^{\odot(K)}\right) - 2v_K\left(\delta^{\odot(Q)}\right) + 4r} \sum_{\delta \in \delta_0\cF_0} Y_{\fP}^{(i)}(\delta)} \\
  &= \left(1 - \frac{1}{q}\right) q^{-2i - m_{11}^\odot - v_K\left(\delta^{\odot(K)}\right) + 4r},
  \end{split}
\]
that is,
\[
  \sum_{\delta \in \delta_0\cF_0} Y_{\fP}^{(i)}(\delta)
  = \size{H^0}\left(1 - \frac{1}{q}\right) q^{2e - i - m_{11}^\odot + 2v_K\left(\delta^{\odot(Q)}\right)}.
\]
Summing over $i \geq n^\odot_{11}$, the right-hand side becomes a geometric series and we get
\[
  \sum_{\delta \in \delta_0\cF_0} Y_{\fP}(\delta)
  = \size{H^0} q^{2e - m^\odot_{11} - n_{11}^\odot + 2v_K\left(\delta^{\odot(Q)}\right)},
\]
as desired.
\end{proof}

\subsection{Unramified}

Here $\clubsuit = (1 ; \bar\zeta_2 \sqrt{D_0})$, where $D_0$ is scaled so that $D_0 \equiv 1 \mod 4$. Observe that $\clubsuit^{(Q)}$ is traceless.

\begin{lem}\label{lem:12_alpha}
There is an $\alpha_0 \in \OO_Q$ such that
\[
  \ker\left(\tr_{\OO_R/\OO_K}\right) = \left\langle\clubsuit (0;1), \clubsuit (1; \alpha_0^2)\right\rangle
\]
as $\OO_K$-modules.
\end{lem}
\begin{proof} The conic
\[
  \cA(\xi) = \tr(\clubsuit \xi^2)
\]
has determinant $1$. Since $\clubsuit \equiv 1 \mod 2$, $\cA$ has maximal squareness $\floor{e/2}$. It has Brauer class $\epsilon(\cA) = 1$, since $\xi = (0;1)$ is a solution. Hence, by Lemma \ref{lem:conic_1}, not all the $\OO_K$-points of $\cA$ lie in a single $1$-pixel. The reduction of $\cA$ modulo $\pi$ consists of $(q+1)$-many $1$-pixels, only one of which has vanishing $K$-coordinate. Hence there is a solution $\xi = (a; \alpha_0)$ with $\pi \nmid a$. Rescaling, we can take $a = 1$.
\end{proof}

\begin{lem}\label{lem:111_weak_zones}
Let $\fP$ be a first vector problem in splitting type $(111)$ or $(12)$ with $\cN_{11}$ weakly active. If $h_1 = 0$, let
\[
\tilde{n} = \floor{\frac{2e - n_{11} + 2}{4}} = \ceil{\frac{e - n_{11}^\odot}{2}}.
\]
Then $ Y_{\fP} $ is given in terms of $n_{11}$ and $n_{11}^\odot = \floor{n_{11}}$ as follows:
\begin{enumerate}[$($a$)$]
  \item If $n_{11} > 2e$ (brown zone), then
  \begin{equation*}
  Y_{\fP} = \begin{cases}
    2q^{e - m_{11} - n_{11}^\odot}  F\left( 0, e, \emptyset \right) & [\delta_0 \hat\omega_C] \in \cF_0 \\
    2q^{e - m_{11} - n_{11}^\odot} \xo F\left( 0, e, \emptyset \right) & [\delta_0 \hat\omega_C] \in (1;\pi;\pi)\cF_0
  \end{cases}
  \end{equation*}
  \item If $ 0 < n_{11} \leq 2e $ (yellow zone), then
  \begin{equation*}
  Y_{\fP} = \begin{cases}
     \begin{aligned}
       \ds \sum_{\smash{0 \leq \ell < \tilde n}} q^{-m_{11} + \ell}(1 + \epsilon^\clubsuit)F^\cross\left( \ell, n^\odot_{11}, e - n^\odot_{11} - \ell \right) \\
       {}+ q^{-m_{11} + \floor{\frac{2e - n_{11}}{4}}} F(\tilde{n}, e - 2\tilde{n}, \tilde{n}),
     \end{aligned} & h_1 = 0 \\
     \ds q^{-m_{11}} (1 + \epsilon^\clubsuit) \xo F\left(0, n^\odot_{11}, e - n^\odot_{11}\right) & h_1 = 1.
  \end{cases}
  \end{equation*}
\end{enumerate}
\end{lem}

\begin{proof}
In the brown zone, the conditions imply that
\[
  (\delta^\odot{\xi^{\odot}}^2)^Q \equiv a \clubsuit \mod \pi^{2e + 1}
\]
for some $a \in \OO_K$, necessarily in $\OO_K^\cross$. So $[\delta^\odot] \in \clubsuit\iota(K) = \clubsuit T(\emptyset, e, \emptyset)$. It's easy to see that all values occur, and the solution volume is constant within the appropriate coarse coset, because the $K$-coordinate of $\beta = \delta^\odot{\xi^{\odot}}^2/\pi^{h_1}$ can range over all of $\pi^{n} \OO_K^\cross$ (for $n \geq 2n^\odot_{11} - h_1$ of the correct parity) while $\beta^{(Q)} \equiv \clubsuit^{(Q)} \mod \pi^{2e + 1}$ remains of constant class. So
\[
  Y_\fP(\delta^\clubsuit) = 2q^{e - m_{11} - n_{11}^\odot}
\]
for each $\delta^\clubsuit$ in the support, as desired.

In particular, $T(\emptyset, e, \emptyset)$ is charmed for $\epsilon^\clubsuit$.

In the yellow zone, the conditions imply that
\[
  \beta^Q \equiv a \clubsuit \mod \pi^{2n_{11}^\odot + h_1}
\]
for some $a$ in $\OO_K$, necessarily in $\OO_K^\cross$. Hence
\[
  [\delta^\clubsuit] \in \begin{cases}
    \epsilon^{\clubsuit-1}(1) \intsec T^\cross\left(\emptyset, n^\odot_{11}, e - n^\odot_{11} \right) & h_1 = 1 \\
    \epsilon^{\clubsuit-1}(1) \intsec T\left(0, n^\odot_{11}, e - n^\odot_{11} \right) & h_1 = 0.
  \end{cases}
\]
Note that the boxgroups have $\ell_1 = n^\odot_{11} \leq \ceil{s/2}$, so they are well defined. Also, $F = G$ in the answer, because everything contains $ \iota(K^\cross) \cdot \cF_e$, which is charmed for $\epsilon^\clubsuit$.

In the case $h_1 = 1$, the conics are all tiny and green, and we get the bound
\begin{align*}
  Y_\fP \leq q^{-m_{11}}(1 + \epsilon^\clubsuit) \xo F\left(0, n^\odot_{11}, e - n^\odot_{11}\right),
\end{align*}
which is exactly as desired. The summing step proceeds routinely, noting that $\epsilon^\clubsuit$ is equidistributed because everything is contained in a non-charmed coarse coset.

In the case $h_1 = 0$, that is, $[\delta^\odot] \in \cF_0$, some further analysis must be done to narrow the support. By Lemma \ref{lem:conic_lift}, we can find a $\xi_1^\odot$ such that $\cM_{11}$ is an equality on the nose and also such that $\xi_1^{\odot(K)} \neq 0$ (to allow recovery of $[\delta^\clubsuit] = [\delta^\odot\diamondsuit\clubsuit\xi'^2]$). Normalize $\delta^\clubsuit$, which hitherto was only a class in $H^1$, by letting
\[
  \delta^\clubsuit = \diamondsuit \clubsuit^{-1} \delta^\odot \in \OO_R^\cross.
\]
Then $\tr(\clubsuit \delta^\clubsuit{\xi_1^\odot}^2) = 0$, so by Lemma \ref{lem:12_alpha} there are $u, b \in \OO_K$ such that
\[
  \delta^\clubsuit{\xi_1^\odot}^2 = u(0;1) + b \pi^{2n_{11}^\odot}(1; \alpha_0^2).
\]
Since $\xi_1^\odot$ is primitive, we have $\pi \nmid u$ and, rescaling, can assume that $u = 1$ since we are curious about the $H^1$-class of the right-hand side. We have $b \neq 0$.

Assume first that $b$ is \emph{not} a square modulo $\pi^{v_K(b) + 2\tilde{n}}$. Write
\[
  b = a^2(1 + \pi^{2k+1}c), \quad \pi \nmid c, \quad 0 \leq k < \tilde n.
\]
We will show that $[\delta^\clubsuit] $ lies in the $k$th term
\[
  T^\cross\left(k, n^\odot_{11}, e - n^\odot_{11} - k \right).
\]
We compute:
\begin{align*}
  [\delta^\clubsuit] &= [(0;1) + b \pi^{2n^\odot_{11}}(1; \alpha_0^2)] \\
  &= [(b ; 1 + b\pi^{2n^\odot_{11}} \alpha_0^2)] \\
  &= [1 + c\pi^{2k+1} ; 1 + a^2\pi^{2n^\odot_{11}} \alpha_0^2 + a^2 c\pi^{2n^\odot_{11} + 2k + 1} \alpha_0^2] \\
  &\equiv [1 + c\pi^{2k+1} ; 1 + a^2\pi^{2n^\odot_{11}} \alpha_0^2 ] \mod \cF_{k + n^\odot_{11}} \\
  &\equiv [1 + c\pi^{2k+1} ; (1 + a \pi^{n^\odot_{11}} \alpha_0)^2 ] \mod \cF_{\ceil{(e + n^\odot_{11})/2}} \subseteq \cF_{k + n^\odot_{11}} \\
  &= [1 + c\pi^{2k+1} ; 1] \\
  &\in T^\cross\left(k, n^\odot_{11}, e - n^\odot_{11} - k \right).
\end{align*}
Also, the conics here are green of squareness $k$. Likewise, if $b$ \emph{is} a square modulo $\pi^{v_K(b) + 2\tilde{n}}$, the same computation shows that $[\delta^\clubsuit] \in T(\tilde{n}, e - 2\tilde{n}, \tilde{n})$. Here the conics are blue, and we have the bounding step.

For the summing step, we note that $\epsilon^\clubsuit$ is equidistributed on each support $T^\cross(\ell, n^\odot_{11}, e - n^\odot_{11} - \ell)$, since it is a union of cosets of $\cF_{e - \ell - 1}$ inside a non-charmed coset of $\cF_{\ell+1}$. The sum is then easy to compute and compare against the total of Lemma \ref{lem:sum_weak}.
\end{proof}

To translate from $Y_\fP$ back to $W_\fP$ is not hard:
\begin{cor}\label{cor:111_weak_zones}
Let $\fP$ be a first vector problem in splitting type $(111)$ or $(12)$ with $\cN_{11}$ weakly active. If $h_1 = 0$, let
\[
\tilde{n} = \floor{\frac{2e - n_{11} + 2}{4}} = \ceil{\frac{e - n_{11}^\odot}{2}}.
\]
Then $ W_{\fP} $ is given in terms of $n_{11}$ and $n_{11}^\odot = \floor{n_{11}}$ as follows:
\begin{enumerate}[$($a$)$]
  \item If $n_{11} > 2e$ (brown zone), then
  \begin{equation*}
    W_{\fP} = \left\{\begin{alignedat}{3}
      &2q^{e - m_{11} - \ceil{n_{11}/2}} F\left( 0, e, \emptyset \right) && [\delta] \in \cF_0, \quad &&\text{$s$ even} \\
      &2q^{e - m_{11} - \ceil{n_{11}/2}} \xo F\left( 0, e, \emptyset \right) \quad && [\delta] \in (1;\pi;\pi)\cF_0, \quad &&\text{$s$ odd} \\
      &2q^{e - m_{11} - \floor{n_{11}/2}} F\left( 0, e, \emptyset \right) && [\delta] \in \cF_0, \quad &&\text{$s$ odd} \\
      &2q^{e - m_{11} - \floor{n_{11}/2}} \xo F\left( 0, e, \emptyset \right) && [\delta] \in (1;\pi;\pi)\cF_0, \quad &&\text{$s$ even}
    \end{alignedat}\right.
  \end{equation*}
  \item If $ 0 < n_{11} \leq 2e $ (yellow zone), then
  \begin{equation*}
    W_{\fP} = \left\{\begin{alignedat}{3}
      &\ds \begin{aligned}
        \sum_{\smash{0 \leq \ell < \tilde n}} q^{-m_{11} + \ell}(1 + \epsilon_C)G^\cross\big( \ell, n^\odot_{11}, e - n^\odot_{11} - \ell \big) \\
         {} + q^{-m_{11} + \big\lfloor{\frac{2e - n_{11}}{4}}\big\rfloor} G(\tilde{n}, e - 2\tilde{n}, \tilde{n}),
      \end{aligned} \quad && h_1 = 0 \\
      &\ds q^{-m_{11}}(1 + \epsilon_C) F\left(0, \floor{\frac{n_{11}}{2}}, e - \floor{\frac{n_{11}}{2}}\right), && h_1 = 1, \quad \text{$s$ odd} \\
      &\ds q^{-m_{11}}(1 + \epsilon_C) \xo F\left(0, \floor{\frac{n_{11}}{2}}, e - \floor{\frac{n_{11}}{2}}\right), && h_1 = 1, \quad \text{$s$ even}
    \end{alignedat}\right.
  \end{equation*}
\end{enumerate}
\end{cor}

\begin{proof}
In the brown zone, the support is contained in $\hat\omega_C \clubsuit \iota(K)$, but using $s > 2e$, this coset is the identity. Intersecting this with the two possible coarse cosets yields the claimed answers.

In the ``long'' yellow zone answer ($h_1 = 0$), the $F$'s in Lemma \ref{lem:111_weak_zones} can be changed to $G$'s, because the answer implies that $T(\tilde n, e - 2\tilde n, \tilde n)$ is charmed for $\epsilon^\clubsuit$. When translating from $Y$ back to $W$, we keep these $G$'s, switching quadratic forms from $\epsilon^\clubsuit$ to $\epsilon_C$.

If $h_1 = 1$ and $s$ is odd, the fact that $\delta = 1$ yields a nonzero ring volume implies that the $\xo G$ simplifies to $F$.

If $h_1 = 1$ and $s$ is even, we claim that the $\xo G$ simplifies to $\xo F$. Note that $\hat\omega_C$ is constrained by $s$: we have $\theta_1 \equiv (1;0;0) \mod \pi^s$, so $\hat\omega_C \equiv (0; \clubsuit^Q) \mod \pi^s$ and
\[
[\hat\omega_C\clubsuit] \in \iota(K^\cross) \cdot (1 + \pi^s \OO_R)
\subseteq 1 + \pi^{n_{11}} \OO_R = T\left(\emptyset, \floor{\frac{n_{11}}{2}}, e - \floor{\frac{n_{11}}{2}}\right).
\]
So the translation $[\delta^\clubsuit] \mapsto [\delta]$ does not affect the support.
\end{proof}
If $h_1 = 0$ and $s$ is even, then $\delta = 1$ must also lie in some term of the sum. Which term it is can be determined using the levels of the conics. If
\[
\ell_C \geq \tilde{n},
\]
then $\delta = 1$ lies in the last term and all $G$'s can be made $F$'s. Otherwise, the $G$ terms for
\[
\ell \leq \ell_C
\]
simplify to $F$'s.

Also, in the very special case that $h_1 = 1$, $s$ is odd, and $n_{11} = 2e$ (on the yellow-brown border), the term
\[
\epsilon_C F(0,e,0)
\]
appearing in the ring volume admits an anomalous simplification. Since the charmed coset of $T(0,e,0)$ is $T^\cross(\emptyset,e,0)$, we have
\begin{equation}
  G(0, e, 0) = \xo F(0,e,0).
\end{equation}
Taking Fourier transforms of both sides,
\begin{equation}
  \epsilon_C F(0,e,0) = F\xo(0,e,0).
\end{equation}
This accounts for the \texttt{Fx-yellow-special} zone in the code.

  \subsection{Splitting type (\texorpdfstring{$1^21$}{1²1})}
As usual, the partially ramified splitting type causes further complications.

Fix a first vector problem with $\cN_{11}$ weakly active. Note that
$m_{11} \geq n_{11} + d_0/2 > d_0/2$, which is good enough for $\xi_1^\odot$ to be meaningful. Also, note that
\begin{equation*}
  n^\odot_{11} = \ceil{\frac{n_{11}}{2} - v^{(K)}(\gamma\gamma^\odot)} =
  \ceil{\frac{2n_{11} - d_0 + h_1}{4}}.
\end{equation*}
If
\[
  n_{11} \leq \frac{d_0 - h_1}{2},
\]
then $\cN_{11}$ is automatic, so the answer $W_{m_{11},n_{11}} = W_{m_{11},0}$ is the same as in the beige zone, to be computed below. Assume we are not in this case. Then the cases $h_1 = 1, 3$ offer no solutions, as there
\[
  \left(\delta^\odot {\xi_1^\odot}^2\right)^{(Q)} \sim \pi_Q^{h_1}
\]
has nonzero $\fI$-value to first order, which cannot be canceled by the $K$-term since $\pi \mid \xi_1^{\odot(K)}$. So we assume that $h_1 \in \{0, 2\}$.

The $\beta$ of Lemma \ref{lem:sum_weak} simplifies to
\[
  \beta = \frac{\delta^\odot {\xi_1^\odot}^2}{\pi^{h_1/2}},
\]
and
\[
  \cN_{11} \iff \beta^{(K)} \equiv 0 \mod \pi^{n_\beta},
\]
where
\[
  n_\beta = 2n^\odot_{11} - \frac{h_1}{2}
\]
is an integer, either even (if $h_1 = 0$) or odd (if $h_1 = 2$). The support of $Y_\fP$ consists of the classes $[\delta^\clubsuit] = [\beta \diamondsuit\clubsuit]$ where $\beta$ ranges over the non-zero-divisors in $\OO_R$ satisfying
\[
  \lambda^\diamondsuit(\beta) = 0, \quad \beta^{(Q)} \sim 1, \quad v(\beta^{(K)}) \geq n_\beta, \quad v(\beta^{(K)}) \equiv n_\beta \mod 2.
\]

\begin{lem} \label{lem:1pow2_1_weak_zones}
Assume that $n_{11} > (d_0 - h_1)/2$, $h_1 \in \{0, 2\}$, and $m_{11} > 2e + (d_0 - h_1)/2$. The volume $W^\odot_{\fP}$ of solutions to the first vector problem depends on the integer
\[
  n_\beta = 2\ceil{\frac{2n_{11} - d_0 + h_1}{4}} - \frac{h_1}{2}
\]
 as follows.
\begin{itemize}
  \item If $n_{\beta} \geq 2e - \ceil{d_0/2} + 1$ (brown zone), then
  \[
    W^\odot_{\fP} = \begin{cases}
      2q^{e - m_{11}^\odot - n_{11}^\odot} F(0, e - d_0', \emptyset), & d_0 - h_1 \equiv 0 \mod 4 \\
      2q^{e - m_{11}^\odot - n_{11}^\odot} \xo F(0, e - d_0', \emptyset), & d_0 - h_1 \equiv 2 \mod 4 \\
      q^{e - m_{11}^\odot - n_{11}^\odot} L(e), & d_0 = 2e + 1
    \end{cases}
  \]
  \item If $d_0/2 \leq n_\beta \leq 2e - d_0/2$ (yellow zone), then, noting that $d_0 = 2d_0'$ is even and that we are in letter type \ref{type:C}, \ref{type:D}, or \ref{type:E}, the answer is as follows:
  \begin{itemize}
    \item If $d_0 - h_1 \equiv 2 \mod 4$, that is, $a_1' \in \ZZ + 1/2$, so $n_\beta \equiv d_0' + 1$ mod $2$, then
    \[
      W^\odot_{\fP} = \begin{cases}
        \ds (1 + \epsilon_C) q^{-m_{11}^\odot + \floor{\frac{d_0' - 1}{2}}} F\bigg(0, \frac{n_\beta - d_0' - 1}{2}, e - \frac{n_\beta + d_0' - 1}{2}\bigg), & \text{\ref{type:C}, \ref{type:E}} \\
        \ds (1 + \epsilon_C) q^{-m_{11}^\odot + \floor{\frac{d_0' - 1}{2}}} \xo F\bigg(0, \frac{n_\beta - d_0' - 1}{2}, e - \frac{n_\beta + d_0' - 1}{2}\bigg), & \text{\ref{type:D}}
      \end{cases}
    \]
    \item If $d_0 - h_1 \equiv 0 \mod 4$, that is, $a_1' \in \ZZ$, so $n_\beta \equiv d_0'$ mod $2$, then
    \begin{align*}
      W^\odot_{\fP} &= \sum_{k=0}^{\tilde n - 1}(1 + \epsilon_C) q^{-m_{11}^\odot + k + \floor{\frac{d_0'}{2}}} G^\cross\left(k, \frac{n_\beta - d_0'}{2}, e - \frac{n_\beta + d_0'}{2} - k\right) + \\
      &\quad + q^{-m_{11}^\odot + \floor{d_0'/2} + \tilde{n}^-} G(\tilde n, e - d_0' - 2\tilde n, \tilde n)
    \end{align*}
    where
    \[
      \tilde n = \floor{\frac{2e - d_0' - n_\beta}{4}} \textand
      \tilde n^- = \floor{\frac{2e - d_0' - n_\beta - 2}{4}}.
    \]
  \end{itemize}
  \item If $0 < n_\beta \leq d_0/2$ (lemon zone), then
  \begin{align*}
    W^\odot_{\fP} &= \sum_{k = 0}^{\floor{\frac{e - n_\beta - 1}{2}}} (1 + \epsilon_C) q^{-m_{11}^\odot + \floor{n_\beta/2} + k} G^\cross(n_\beta + 2k) + q^{-m_{11}^\odot + \floor{e/2}} G(e)
  \end{align*}
\end{itemize}
\end{lem}
\begin{proof} We consider each zone in turn.
\subsubsection{Brown zone.}
Let $n = n_\beta$. Here it suffices to prove that, as $c$ varies over $\OO_K^\cross$, the element
\[
  \beta = (\pi^{n} c ; 1 + \pi_Q \pi^{n} c)
\]
attains all classes in $\iota(\pi^{n} \OO_K^\cross)$ equally often. Now $[\beta] \in \iota(c) \cF_{2e - \ceil{d_0/2} + 1}$, but the claim that all classes occur equally often is less immediate. Nevertheless, as promised, we will do it directly, without reference to an $\cN_{11}$-lemma.

Since
\[
  \iota : \cF_{K, e - \ceil{d_0/2} + 1} \to \cF_{2e - \ceil{d_0/2} + 1}
\]
is an isomorphism of groups, we may consider its inverse $\iota^{-1}$. Also, let
\[
  \sigma : \OO_K^\cross/\left(\OO_K^\cross\right)^2 \to \OO_K^\cross
\]
be a section of the natural projection that has the minimal distance property
\[
  \size{\sigma([x_1]) - \sigma([x_2])} \leq \size{x_1 - x_2}.
\]
Such a $\sigma$ can be made, for instance, by letting $\sigma(x)$ be the product of the elements in the Shafarevich basis decomposition of $x$ in accordance with Proposition \ref{prop:Sh_basis_quartic}. Now consider the map
\begin{align*}
  j : \OO_K^\cross &\to \OO_K^\cross \\
  c &\mapsto c \cdot \sigma\big(\iota^{-1}(1 + \pi_Q\pi^{n} c)\big)
\end{align*}
We claim that $j$ is an \emph{isometry} in the $p$-adic metric.

Let $c_1, c_2 \in \OO_K^\cross$ be given with $\size{c_1 - c_2} = \size{\pi}^{k}$. Let $n = n_\beta$. It suffices to prove that
\[
  \sigma\big(\iota^{-1}(1 + \pi_Q\pi^{n} c_1)\big) \equiv \sigma\big(\iota^{-1}(1 + \pi_Q\pi^{n} c_2)\big) \mod \pi^{k+1}.
\]
However, we have
\begin{equation*}
  \ell\left(\frac{1 + \pi_Q\pi^{n} c_1}{1 + \pi_Q\pi^{n} c_2}\right) \geq \min\{n + k, 2e+1\}
\end{equation*}
and since $\iota$ increases levels by $e$ in the aforementioned region of invertibility,
\[
  \ell\left(\iota^{-1}\left(\frac{1 + \pi_Q\pi^{n} c_1}{1 + \pi_Q\pi^{n} c_2}\right)\right) \geq \min\{n - e + k, e + 1\}
\]
which ensures that
\[
  v\left(\sigma\left(\iota^{-1}\left(\frac{1 + \pi_Q\pi^{n} c_1}{1 + \pi_Q\pi^{n} c_2}\right)\right) - 1\right) \geq 2(n - e + k) > k,
\]
as desired. So $j$ is an isometry. In particular, $j$ is bijective and preserves volumes. So the $\beta_c$, which satisfy
\[
[\beta_c] = \big[\iota\big(\pi^{n} j(c)\big)\big],
\]
are equidistributed in $\iota(\pi^{n} \OO_K^\cross)$. Summing over $k = n + 2i$ and getting the scaling from the summation lemma \ref{lem:sum_weak} yields the claimed brown-zone answer.

\subsubsection{Yellow zone (short answer).} When $n_\beta \equiv d_0' + 1 \mod 2$, the element $\beta^\odot = (c_1 \pi^{n_\beta}; 1 + c_1 \pi^{n_\beta} \pi_Q)$, $2 \mid v(c_1)$, is seen to lie in
\[
  \iota\left(\pi^{h_1/2} \OO_K^\cross\right) \cdot L(n_\beta) = \iota(\pi^{h_1/2}) \cdot T\left(0, \frac{n_\beta - d_0' - 1}{2}, e - \frac{n_\beta + d_0' - 1}{2}\right),
\]
which is exactly the support of the claimed ring total. The conics are all green of level $\floor{(d_0' - 1)/2}$, and $\epsilon^\clubsuit$ is equidistributed because we are in an uncharmed coset of $\cF(d_0')$, so the summing step is routine.

\subsubsection{Yellow zone (long answer).}
 When $n_\beta \equiv d_0' \mod 2$, more work is in order. By analogy with the other splitting types, our intuition is to try to write $\beta^\odot/\heartsuit$ as an approximate sum of squares. Begin with the multiplication law
\[
  \pi_Q^2 = t\pi_Q + u\pi, \quad t \sim \pi^{d_0'}, \quad u \sim 1.
\]
Note that the elements
\[
  \alpha_0 = (0;1) \textand \alpha_1 = (t; u\pi + t\pi_Q)
\]
span $\ker \lambda^\diamondsuit$ over $K$, but over $\OO_K$, they span the index-$q^{d_0'}$ sublattice
\[
  \Lambda = \left\langle\alpha \in \OO_R : \lambda^\diamondsuit(\alpha) = 0 \textand \alpha^{(K)} \equiv 0 \mod \pi^{d_0'}\right\rangle.
\]
Due to the inequality $n_\beta \geq d_0'$, $\beta^\odot$ lies in $\Lambda$ and can be written in the form
\[
  \beta^\odot = c_0\alpha_0 + c_1\alpha_1 = (c_1 t; c_0 + c_1 u\pi + c_1 t\pi_Q).
\]
Here $c_0$ is a unit, and we can scale so that $c_0 = 1$; while $c_1 = b\pi^{n_\beta - d_0'},$ where $b$ has even valuation.

If $b$ is \emph{not} a square modulo $\pi^{v(b) + 2\tilde n}$, write
\[
  b = a^2(1 + \pi^{2k+1} c), \quad \pi\nmid c, \quad 0 \leq k < \tilde n.
\]
We claim $[\beta^\odot]$ belongs to the support of the $k$th term in the claimed ring total:
\begin{align*}
  [\beta^\odot] &= \left[\left(a^2\pi^{n_\beta - d_0'}(1 + \pi^{2k+1}c) t; 1 + a^2\pi^{n_\beta - d_0'} \pi_Q^2 + a^2 \pi^{n_\beta - d_0' + 2k + 1}\pi_Q^2\right)\right] \\
  &= \heartsuit \left[\left(1 + \pi^{2k+1}c; 1 + a^2\pi^{n_\beta - d_0'} \pi_Q^2 + a^2 \pi^{n_\beta - d_0' + 2k + 1}\pi_Q^2\right)\right] \\
  &\equiv \heartsuit \left[\left(1 + \pi^{2k+1}c; (1 + a^2\pi^{n_\beta - d_0'} \pi_Q^2)(1 + a^2 \pi^{n_\beta - d_0' + 2k + 1}\pi_Q^2)\right)\right] \\
  &\equiv \heartsuit \left[\left(1 + \pi^{2k+1}c; (1 + 2a\pi^{(n_\beta - d_0')/2} + a^2\pi^{n_\beta - d_0'} \pi_Q^2)(1 + a^2 \pi^{n_\beta - d_0' + 2k + 1}\pi_Q^2)\!\right)\!\right]\\
  &= \heartsuit \left[\left(1 + \pi^{2k+1}c; (1 + a\pi^{(n_\beta - d_0')/2})^2(1 + a^2 \pi^{n_\beta - d_0' + 2k + 1}\pi_Q^2)\right)\right] \\
  &= \heartsuit \left[\left(1 + \pi^{2k+1}c; (1 + a^2 \pi^{n_\beta - d_0' + 2k + 1}\pi_Q^2)\right)\right] \\
  &= \heartsuit \left[\left(1 + \pi^{2k+1}c; 1 + a^2 \pi^{n_\beta - d_0' + 2k + 1}t\pi_Q + a^2 \pi^{n_\beta - d_0' + 2k + 2} u\right)\right] \\
  &\equiv \heartsuit \left[\left(1 + \pi^{2k+1}c; 1 + a^2 \pi^{n_\beta - d_0' + 2k + 2} u\right)\right] \\
  &= \heartsuit \cdot \iota\left((1 + \pi^{2k+1}c)(1 + a^2 \pi^{n_\beta - d_0' + 2k + 2} u)\right) \mod \cF_{n_\beta + 2k + 1},
\end{align*}
the dropped terms, being of valuation at least $n_\beta + 2k + 1 + 1/2$ by various evident combinations of the known inequalities, may be omitted. We get that
\[
  [\beta^\odot] \in \heartsuit \iota(1 + \pi^{2k+1}\OO_K)\cF_{n_\beta + 2k + 1} = \heartsuit T\left(k, \frac{n_\beta - d_0'}{2}, e - \frac{n_\beta + d_0'}{2} - k\right),
\]
and evidently \emph{not} in $\heartsuit \iota(1 + \pi^{2k+3}\OO_K)\cF_{n_\beta + 2k + 1}$, allowing us to refine $T$ to $T^\cross$.

If $b$ is a square modulo $\pi^{v(b) + 2\tilde n}$, a similar calculation shows that
\[
  [\beta^\odot] \in \heartsuit T\left(\tilde n, e - d_0' - 2\tilde n, \tilde n\right).
\]
Keeping track of the levels and colors of the conics, we get the bounding step, with the $G$'s specified to $F$'s. Then, using the equidistribution of $\epsilon^\clubsuit$, we get the summing step together with the extra information that $T\left(\tilde n, e - d_0' - 2\tilde n, \tilde n\right)$ is charmed for $\epsilon^\clubsuit$.

\subsubsection{Lemon zone.} Here $n^\odot = 0$. The proof consists of
\begin{itemize}
  \item transforming $\cN_{11}$ into a condition on the level of $\delta^\clubsuit$, using the easy inequality that for $k < d_0/2$,
  \[
    v(\fI(\beta^\odot)) \geq k \iff \ell(\beta^\odot) \geq k
  \]
  \item further restricting $\ell(\delta^\clubsuit)$ using the level parity lemma \ref{lem:level_parity}, and
  \item splitting the beige conics that appear when $\ell(\heartsuit \beta^\odot) \geq 2\floor{e/2}$ into generic and special $1$-pixels, which contribute for $h_1 = 0$ and $h_1 = 2$ respectively. The computation is identical to that for the beige zone in splitting type $(1^3)$ (see Lemma \ref{lem:beige_1pow3}). \qedhere
\end{itemize}
\end{proof}

\begin{rem}
The yellow and lemon zones have been defined as overlapping (at $n_\beta = d_0')$ because both colors give the same answer in this case. This flexibility allows us to define the yellow-lemon boundary to be given respectively by the inequalities
\[
  n_{11} > d_0 - \frac{o_1}{2} \textand n_{11} \leq d_0 - \frac{o_1}{2}.
\]
\end{rem}
\section{The beige zone}\label{sec:beige}
We term \emph{beige} the zone where $\cN_{11}$ is inactive but $\cM_{11}$ is sufficiently active that the conic transformation of Lemma \ref{lem:tfm_conic} applies.

\subsection{Unramified}
The following is immediate from the computation of the Igusa zeta function in Lemmas \ref{lem:conic_1} and \ref{lem:conic_pi}.
\begin{lem} \label{lem:beige_ur}
In unramified splitting type, if $m_{11} > 2e$ but $\cN_{11}$ is inactive (beige zone), then the ring volume for $\xi'_1$ is given by
  \begin{align*}
W_{m_{11},0} &= \begin{cases}
\ds q^{-m_{11}} \left[(1 + \epsilon_C) \sum_{0 \leq \ell \leq \floor{\frac{e}{2}} - 1}
q^{\ell}G^{\cross}\left(\ell, 0, e - \ell\right)
+  (\text{core})\right] & [\delta^\odot] \in \cF_0 \\
\ds q^{-m_{11}} (1 + \epsilon_C) \big(\xo G(0,0,e) + \xo\xo G(0,0,e) \big) & [\delta^\odot] \notin \cF_0
\end{cases}
\end{align*}
where
\[
(\text{core}) = \begin{cases}
 q^{\frac{e}{2}} \displaystyle \left( 1 + \frac{1}{q} \right) G\left(\frac{e}{2}, 0, \frac{e}{2} \right) &
e \text{ even} \\
 q^{\frac{e-1}{2}} \displaystyle (1 + \epsilon_C) G\left(\frac{e-1}{2}, 0, \frac{e+1}{2}\right) &
e \text{ odd}
\end{cases}
\]
\end{lem}

In the $[\delta^\odot] \notin \cF_0$ case, the answer would more strictly be written as a restriction to the particular coarse coset specified by the discrete data, but we do not do so, as all $(\size{H^0} - 1)$-many non-charmed coarse cosets admit the same extender indices and will be immediately summed. We know which coset of $\cF_0$ is charmed, and hence:
\[
  \xo G(0,0,e) + \xo\xo G(0,0,e) = \begin{cases}
    \xo F(0,0,e) + \xo\xo F(0,0,e) & \text{$s$ even} \\
    F(0,0,e) + \xo\xo F(0,0,e) & \text{$s$ odd}.
  \end{cases}
\]

\subsection{Splitting type (\texorpdfstring{$1^3$}{1³})}
Here, a little more care is required to deal with the restrictions on the $1$-pixel of $\xi^\odot$ that remain active in the beige zone.

\begin{lem}\label{lem:beige_1pow3}
In splitting type $(1^3)$, a first vector problem with $m_{11} > 2e$ and $n_{11} \leq 0$ (beige zone) has the answer, for $h_1 = 1$,
\begin{align*}
  W_{m_{11}, 0} &= q^{-\tri m_{11}}(1 + \epsilon_C) \sum_{k = 0}^{\ceil{e/2} - 1} q^k \big(G_1(k, 0, e - k) - G_1(k+1, -1, e - k)\big) + {} \\
  &\quad + q^{-\tri m_{11} + \ceil{e/2}} G_1\left(\ceil{\frac{e}{2}}, -\One_{2 \nmid e}, \ceil{\frac{e}{2}}\right)
\end{align*}
and, for $h_1 = -1$,
\begin{align*}
  W_{m_{11}, 0} = q^{-\tri m_{11} + 2}(1 + \epsilon_C) \!\!\!\sum_{k = 0}^{\lfloor{e/2}\rfloor - 1}{\!\!\!} q^k \big(G_{-1}(k, 1, e - k - 1) - G_{-1}(k+1, 0, e - k - 1)\big) \\
  \quad{} + q^{-\tri m_{11} + 2 + \floor{e/2}} G_{-1}\left(\floor{\frac{e}{2}}, \One_{2 \nmid e}, \floor{\frac{e}{2}}\right).
\end{align*}
\end{lem}
\begin{proof}
Since $\cF_e$ is charmed for $\epsilon$, all of the $G$'s can be viewed as selectors for the coset of the indicated boxgroup containing $[\hat\omega_C]$. These $G$'s, happily enough, are the level spaces where the level
\[
  \ell = \ell(\delta^\odot \heartsuit)
\]
in Lemma \ref{lem:level_parity} attains a fixed value. When $\ell < 2\floor{e/2}$ is even, the conic lies in a generic pixel and contributes to the $k = \ell/2$ term of the $h_1 = 1$ answer. When $\ell < 2\floor{e/2}$ is odd, the conic lies in the special pixel and contributes to the $k = (\ell - 1)/2$ term of the $h_1 = -1$ answer. When $\ell \geq 2\floor{e/2}$, the conic lies in multiple pixels and contributes to the trailing terms of one or both answers. As the volume within each $1$-pixel is governed by the $n^\odot = 1$ case of Lemma \ref{lem:conic_1}, the answer can be computed explicitly and agrees with the stated formula.
\end{proof}

  \subsection{Splitting type (\texorpdfstring{$1^21$}{1²1})}
The same method works in splitting type $(1^21)$: we simply state the result.
\begin{lem}
In the beige zone of splitting type $(1^21)$, if $m_{11}^\odot > 2e$, the total $W_{m_{11},0}^{\odot}$ is given as follows:
\begin{itemize}
  \item If $h_1 = 3$, there are no solutions.
  \item If $h_1 = 1$, then
  \begin{align*}
    W_{m_{11},0}^{\odot} &=
    (1 + \epsilon_C)q^{-m_{11}^\odot} G^\cross(-1) \\
    &= \begin{cases}
      (1 + \epsilon_C)q^{-m_{11}^\odot} F(0) & \text{type \ref{type:A}}\\
      (1 + \epsilon_C)q^{-m_{11}^\odot} \left(F(-1) - F(0)\right) & \text{types \ref{type:B}--\ref{type:E}}.
    \end{cases}
  \end{align*}
  \item If $h_1 = 0$ or $h_1 = 2$, then
  \[
    W^\odot_{m_{11},0} = \sum_{k = 0}^{\floor{\frac{e - h_1/2 - 1}{2}}} (1 + \epsilon_C) q^{-m_{11}^\odot + k} G^\cross(2k + h_1/2) + q^{-m_{11}^\odot + \floor{e/2}} G(e).
  \]
\end{itemize}
\end{lem}
The finite case analysis in this lemma is also checked by the accompanying file \verb|quartic_ON_1e21.sage|.

\section{Orthogonality}

We are now ready to prove the following significant property of boxgroups.
\begin{lem}\label{lem:orth}
For any boxgroup $ T_{\theta_1}(\ell_0, \ell_1, \ell_2) $, its orthogonal complement is given by $ T_{\theta_1}(\ell_2, \ell_1, \ell_0) $.
\end{lem}
\begin{rem}
This is an example of an \emph{explicit reciprocity law}, that is, a formula for the Hilbert symbol in a certain region. There is a wide literature on explicit reciprocity laws, but we suspect that this one is new. In our proof, the only fact we use about the Hilbert pairing is that it is the associated bilinear form to $\epsilon$ (and $\epsilon_C$). This enables us to use various facts about $\epsilon$ gleaned in the preceding sections. We will pose first vector problems involving various values of the resolvent datum $\theta_1$ and of the discrete datum $n_{11}$ to plug into the lemmas regarding the ring volumes.
\end{rem}
\begin{proof}[Proof of Lemma \ref{lem:orth}]
Our proof has the following structure:
\begin{itemize}
  \item Replace $\theta_1$ by $\theta_1' = \eta^2$ from the $\eta$-lemma (Lemma \ref{lem:eta_ur}, \ref{lem:eta_1pow3}, or \ref{lem:eta_1pow21}), and use the gray-green inequality to show that the boxgroups are unchanged;
  \item Increase $\ell_1$ by $1$ repeatedly, showing that this only strengthens the claim, until this cannot be done without violating the gray-red or gray-blue inequality.
  \item Exhibit the relevant boxgroups as terms in zone totals in the red and blue zones, and deduce that $\epsilon_C(\alpha) = 1$ in relevant regions.
  \item Deduce values of $\left\langle\alpha, \beta\right\rangle = \epsilon_C(1) \epsilon_C(\beta) \epsilon_C(\alpha) \epsilon_C(\alpha\beta)$.
\end{itemize}

We carry out the proof in detail in the unramified splitting type only, the proof in the other types being very similar.

We first note that if any of $\ell_0, \ell_1, \ell_2$ is the symbol $\emptyset$, the result follows easily from the self-orthogonality of $\iota(K^\cross)$ (if applicable), as mentioned above. So we can assume that the $\ell_i$ are integers. By definition, they must satisfy
 \begin{align}
\ell_0 + \ell_1 + \ell_2 &= e \\
\ell_1 &\leq \ell_0 + \frac{s}{2} + 1 && (\text{the gray-red inequality})\label{eq:orth_red} \\
\ell_1 &\leq \frac{s + \square_C + 1}{2} && (\text{the gray-green inequality}) \label{eq:orth_green} \\
\ell_1 &\leq \ell_2 + \frac{s}{2} + 1, && (\text{the gray-blue inequality})\label{eq:orth_blue}
\end{align}

We first reduce to the case $\square_C = e$. If $s$ is even, this is accomplished, as in the proof of the red zone, by replacing $\theta_1$ by $\theta_1' = \eta^2$ and noting that, by the gray-green inequality, the boxgroups $T_{\theta_1}(\cdots) = T_{\theta_1'}(\cdots)$ are unchanged. If $s$ is odd, we simply replace $\theta_1$ by a $\theta_1' \equiv (1;0;0) \mod \pi^{\max\{2e + 1, s\}}$ whose corresponding $\hat\omega_{C'}$ is in $\cF_{\floor{e/2}}$. (For example, $\theta_1' = (1; 8 \pi^{2s} \bar\zeta_2 \sqrt{D_0}))$ is found to work.) Then since $\theta_1' \equiv \theta_1 \mod \pi^s$ and all boxgroups satisfy the gray-green inequality $\ell_1 \leq (s + 1)/2$, the boxgroups are unchanged. Incidentally, we can also assume as a result of this reduction that $s$ is even.

Now the gray-green inequality is subsumed by the gray-blue and gray-red ones.  Using $\cF_i^\perp = \cF_{e-i}$, the truth of the lemma for a triple $(\ell_0, \ell_1, \ell_2)$, $\ell_1 \geq 1$, implies its truth for the triples $(\ell_0 + 1, \ell_1 - 1, \ell_2)$ and $(\ell_0, \ell_1 - 1, \ell_2 + 1)$. Hence we can run these reductions backward, increasing $\ell_1$ until we reach an obstruction. This usually happens if either the gray-blue or the gray-red inequality becomes an equality, but it can also happen in two special cases, which we dispatch now:
\begin{itemize}
  \item $\ell_0 = \ell_2 = 0$. Here $s \geq 2e - 2$, and the self-orthogonality of $T(0, e, 0)$ follows from that of $T(\emptyset, e, \emptyset) = \iota(K^\cross)$, unless $e = 1$ and $s = 0$, in which case the gray-blue and gray-red inequalities are also equalities.
  \item Both the gray-blue and gray-red inequalities are $1$ away from equality, that is,
  \[
    (\ell_0, \ell_1, \ell_2) = \left(\frac{2e - s}{6}, \frac{e + s}{3}, \frac{2e - s}{6}\right).
  \]
  As we see, this is only possible if $2e \equiv s \mod 3$. This space shows up as the support in the blue zone for the first vector problem
  \[
    m_{11} = 2e + 1, \quad n_{11} - s = \frac{2e - s}{3},
  \]
  right on the blue-red boundary. It is therefore isotropic and, by virtue of its size, maximal isotropic.
\end{itemize}
So we are left with the case that, without loss of generality, the gray-blue inequality is an equality
\[
  \ell_1 = \ell_2 + \frac{s}{2} + 1.
\]
We thus have
\[
  (\ell_0, \ell_1, \ell_2) = \left( e - 2\ell_2 - \frac{s}{2} - 1, \ell_2 + \frac{s}{2} + 1, \ell_2\right).
\]
Let $V = T(\ell_0, \ell_1, \ell_2)$ and $W = T(\ell_2, \ell_1, \ell_0)$ be the claimed orthogonals. By the gray-red inequality, $\ell_0 \geq \ell_2$ and $V \subseteq W$.

If $\ell_0 = \ell_2$, we again have a unique group
\[
  (\ell_0, \ell_1, \ell_2) = \left(\frac{2e - s - 2}{6}, \frac{e + s + 2}{3}, \frac{2e - s - 2}{6}\right)
\]
which shows up as the support in the blue zone for the first vector problem
\[
  m_{11} = 2e + 1, \quad n_c = \frac{2e - s - 2}{3},
\]
right on the blue-red boundary. It is therefore isotropic and, by virtue of its size, maximal isotropic.

So we may assume that $\ell_0 > \ell_2$ and $V \subsetneq W$. The space $V$ shows up as the blue-zone support for the first vector problem
\[
  m_{11} = 2e + 1, \quad n_c = 2e - 4\ell_2 - 1.
\]
Hence $V$ is isotropic. Also, the first vector problem
\[
  m_{11} = 2e + 1, \quad n_c = 2\ell_2 + 1
\]
lies in the red zone. The first term of its answer is a positive multiple of
\[
  (1 + \epsilon_C) G^\cross(W),
\]
but since $\square_C = e$, the $G$ can be replaced by $F$. Our strategy is as follows. Since $W^\cross = T^\cross(\ell_2, \ell_1, \ell_0)$ generates $W$ as a group, it's enough to show that any $\alpha \in W^\cross$ and $\beta \in V$ are orthogonal. We may write
\begin{align*}
\left\langle\alpha, \beta\right\rangle &= \epsilon_C(1) \epsilon_C(\beta) \epsilon_C(\alpha) \epsilon_C(\alpha\beta) \\
&= \epsilon_C(\alpha) \epsilon_C(\alpha\beta).
\end{align*}
Here $\epsilon_C(1)=1$, and $\epsilon_C(\beta)=1$ because $V$ is the blue-zone support on which $\epsilon_C$ is identically $1$. So if $\alpha$ and $\beta$ are \emph{not} orthogonal, then, applying the transformation $\alpha \mapsto \alpha\beta$ if need be, we may assume
\[
  \epsilon_C(\alpha) = 1 \textand \epsilon_C(\alpha\beta) = -1.
\]
Since $\epsilon_C(\alpha) = 1$, $\alpha$ is in the support of the red-zone answer $W_{m_{11} = 2e + 1, n_c = 2\ell_2 + 1}$. So there is a $\psi$ in the box $1 + B_{\theta_1}(\infty, 2\ell_2 + 1)$ representing the class $[\alpha] \in H^1$. Recenter, using Lemma \ref{lem:recenter}, and consider the $\fP'$ with $\theta_1' = \psi^{-1}\theta_1$, $m_{11}' = 2e+1$, and
\[
  n_c' = 2e - 4\ell_2 - 1.
\]
Since
\[
  \ell_{\fP'} = \ell(\alpha) = \ell_2,
\]
this problem is still in the blue zone (right on the blue-green boundary), and we get that $\epsilon_{C'} = 1$ on the boxgroup
\[
  1 + B_{\theta_1'}(\ell_0, \ell_1, \ell_2).
\]
Since $\theta_1' \equiv \theta_1 \mod \pi^{2\ell_2 + 1} \OO_K[\theta_1]$, the subscript can be changed from $\theta_1'$ to $\theta_1$ without changing the boxgroup. So
\[
  1 = \epsilon_{\fP'}(\beta) = \epsilon_\fP(\alpha \beta),
\]
a contradiction. This completes the proof in the unramified splitting types.

\subsubsection{Splitting type (\texorpdfstring{$1^3$}{1³}).}
Here, for boxgroups of type $h = 1$, we get
\begin{itemize}
\item $\ds F\left(e - 2\ell - 1 + h, \ell + 1 - h, \ell\right)$ from $n_{11} = 2e - 4\ell - 1 + 4h/3$ (blue zone)
\item $\ds G^\cross\left(\ell, \ell + 1 - h, e - 2\ell - 1 + h\right)$ from $n_{11} = 2\ell + 1 - 2h/3$ (red zone, or blue zone if $\ell_0 = \ell_2$)
\item exceptionally,
\[
  F\left(\frac{e + h}{3}, \frac{e - 2h}{3}, \frac{e + h}{3}\right),
\]
if the arguments are integers, from $n_{11} = 2e/3$ (blue zone).
\end{itemize}
The blue and red supports in the first two items are exchanged by reversing the outer indices, and the exceptional support is self-reversing. The same recentering argument as in the unramified case shows that the first two are orthogonal; their orders multiply to $\size{H^1}$, while the exceptional support is maximal isotropic.

\subsubsection{Splitting type (\texorpdfstring{$1^21$}{1²1}).}
We must be in letter type \ref{type:C}, \ref{type:D}, or \ref{type:E}, and $d_0 = 2d_0'$ must be even. We again replace $\theta_1$ by $\theta_1' = \eta^2$, where $\eta$ comes from Lemma \ref{lem:eta_1pow21}, and note that this makes $\ell_C$ achieve its maximal possible value $\floor{e/2}$ without changing boxgroups. In particular, we need only consider letter types \ref{type:C} and \ref{type:D}. Now we reduce to the case where $\ell_1$ cannot be increased further. If $\ell_2 = 0$, we reduce to the case where $\ell_2 = \emptyset$: the characteristic function of $T(a,b,0)$ is the sum of that of $T(a,b,\emptyset)$ and the corresponding $\cF_0$-piece, and the same decomposition holds for the claimed orthogonal complement. The case is therefore covered by the remarks following the construction of supplementary boxgroups, and we derive
\begin{itemize}
  \item $F(e - d_0' - 2\ell - 2, \ell + 2, \ell)$ from $n_c = 2e - 2d_0' - 4\ell - 2 - h_\eta$ (blue zone)
  \item $G^\cross(\ell, \ell + 2, e - d_0' - 2\ell - 2)$ from $n_c = 2\ell + h_\eta$ (red zone, or blue zone if $\ell_0 = \ell_2$)
  \item exceptionally,
  \[
    F\left(\frac{e - d_0' - 1}{3}, \frac{e - d_0' + 2}{3}, \frac{e - d_0' - 1}{3}\right)
  \]
  if the arguments are integers, from
  \[
    n_c = \frac{2e - 2d_0' + 1}{3} - h_\eta.
  \]
\end{itemize}
The blue and red supports in the first two items are again exchanged by reversing the outer indices, and the exceptional support is self-reversing. The same recentering argument proves orthogonality; the orders of the first two supports multiply to $\size{H^1}$, while the exceptional support is maximal isotropic.
\end{proof}

\subsection{\texorpdfstring{$E$}{E}-forms of the long answers}
In the red, yellow, and beige zones when $h = 0$, the answer, as announced in \ref{lem:111_strong_zones}, \ref{lem:111_weak_zones}, and \ref{lem:beige_ur}, is a sum of $G$ and $\epsilon_C G$ terms. It is capable of a simplification.

\begin{defn} \label{defn:E}
  If $ F : H^1 \to \CC$ is a weighting, define
  \[
  E_F = q^{-2e}\widehat{\epsilon_C \widehat {F}}.
  \]
  The scale factor is motivated by the fact that, for $f : H^1 \to \CC$,
  \[
    \hat{\hat{f}} = q^{2e} f.
  \]
  If $F({\cdots})$ is defined for some sequence of symbols ``${\cdots}$'', we put
  \[
    E({\cdots}) = E_{F({\cdots})}.
  \]
\end{defn}
For instance, in unramified splitting type, if the $\ell_i$ are integers, we have
\[
E(\ell_0,\ell_1,\ell_2)
= \begin{cases}
G(\ell_0,\ell_1,\ell_2), & \ell_0 \leq \ell_2 \\
q^{\ell_2 - \ell_0} \epsilon_C G(\ell_2,\ell_1,\ell_0), & \ell_2 < \ell_0.
\end{cases}
\]
At $\ell_0=\ell_2$, the two expressions agree.

\begin{lem} \label{lem:to_E}
  \begin{enumerate}[$($a$)$]
    \item \label{to_E:red_yellow} A sum of the form
    \[
    \sum_{\ell = a}^{\ceil{b/2} - 1} q^\ell (1 + \epsilon_C) G^\cross(\ell, e-b, b-\ell) + q^{\floor{b/2}} G\left(\ceil{\frac{b}{2}}, e - 2\ceil{\frac{b}{2}}, \ceil{\frac{b}{2}}\right),
    \]
    where $b \geq 2a$, can be rewritten as
    \[
    \sum_{\ell = a}^{b - a} q^\ell E(\ell, e-b, b-\ell) - \sum_{\ell = a}^{b - a - 1} q^\ell E(\ell + 1, e - b - 1, b-\ell).
    \]
    \item \label{to_E:beige} In unramified splitting types, a sum of the form
    \[
    (1 + \epsilon_C) \sum_{a \leq \ell \leq \floor{\frac{e}{2}} - 1}
    q^{\ell}G^{\cross}\left(\ell, 0, e - \ell\right)
    +  (\text{core}),
    \]
    where $(core)$ denotes the beige-zone core (see Lemma \ref{lem:beige_ur}) and $a < \floor{e/2}$, can be rewritten as
    \[
    \sum_{\ell = a}^{e - a} q^{\ell} E(\ell, 0, e-\ell) - \sum_{\ell = a+1}^{e-a-1} q^{\ell - 1} E(\ell, 0, e-\ell)
    \]
    where the second sum may be empty or may have to be interpreted according to the natural convention
    \[
    \sum_{\ell = 1}^{-1} x_\ell = -x_0.
    \]
    \item In ramified splitting types, a sum of the form
    \[
      \sum_{k = a}^{\floor{\frac{e - h - 1}{2}}} (1 + \epsilon_C) q^{k} G^\cross(2k + h) + q^{\floor{e/2}} G(e),
    \]
    where $h \in \{0,1\}$, where $G(\ell)$ denotes the indicator function of the charmed coset of $\FF_\ell$ for $0 \leq \ell \leq e$, and $G^\cross(\ell) = G(\ell) - G(\ell + 1)$, can be rewritten as
    \[
      \sum_{k = a}^{e - h - a} q^{k} E(2k + h) - \sum_{k = a}^{e - h - a - 1} q^{k} E(2k + h + 1).
    \]
  \end{enumerate}
\end{lem}

\begin{proof}
  The proof is straightforward, converting each $G$ and $\epsilon_C G$ into an $E$ and merging the ranges of summation.
\end{proof}

\begin{cor}\label{cor:E-forms}
For unramified splitting type, we get in the red zone:
\begin{align*}
  W_{m_{11},n_{11}} &= \sum_{\ell = \floor{\frac{n_c}{2}}}^{e - n_{11} + \frac{s}{2}} q^{-m_{11} + \ell} E\left( \ell, \ceil{\frac{n_{11}}{2}}, e - \ceil{\frac{n_{11}}{2}} - \ell \right) \\
  & \quad {} - \sum_{\ell = \floor{\frac{n_c}{2}}}^{e - n_{11} + \frac{s}{2} - 1} q^{-m_{11} + \ell} E\left( \ell + 1, \ceil{\frac{n_{11}}{2}} - 1, e - \ceil{\frac{n_{11}}{2}} - \ell\right),
\end{align*}
and in the yellow zone:
\begin{align*}
W_{m_{11},n_{11}} &=
\sum_{\ell = 0}^{e - \ceil{n_{11}/2}} q^{-m_{11} + \ell}E\left( \ell, \ceil{\frac{n_{11}}{2}}, e - \ceil{\frac{n_{11}}{2}} - \ell \right) \\
& \quad {}- \sum_{\ell = 0}^{e - \ceil{n_{11}/2} - 1} q^{-m_{11} + \ell}E\left( \ell+1, \ceil{\frac{n_{11}}{2}} - 1, e - \ceil{\frac{n_{11}}{2}} - \ell \right),
\end{align*}
and in the beige zone:
\[
  W_{m_{11}} = \sum_{\ell = 0}^e q^{-m_{11} + \ell} E\left(\ell, 0, e-\ell\right) - \sum_{\ell = 1}^{e-1} q^{-m_{11} + \ell - 1} E(\ell, 0, e-\ell).
\]

For splitting type $(1^3)$, we get in the red zone:
\begin{align*}
   W_{m_{11},n_{11}} &= \sum_{\ell = \floor{\frac{\tri n}{2}}}^{e - \tri n + h} q^{-\tri m + 1 - h + \ell} E_h\left( \ell, \ceil{\frac{\tri n}{2}} - h, e - \ceil{\frac{\tri n}{2}} - \ell + h \right) \\
&\quad {} - \sum_{\ell = \floor{\frac{\tri n}{2}}}^{e - \tri n + h - 1} q^{-\tri m + 1 - h + \ell} E_h\left( \ell + 1, \ceil{\frac{\tri n}{2}} - h - 1, e - \ceil{\frac{\tri n}{2}} - \ell + h \right),
\end{align*}
and in the beige zone $(h_1 = 1)$:
\[
  W_{m_{11}} = \sum_{\ell = 0}^{e} q^{-\tri m + \ell} E_{1}(\ell,0,e-\ell) - \sum_{\ell = 0}^{e-1} q^{-\tri m + \ell} E_{1}(\ell+1,-1,e-\ell)
\]
and in the beige zone $(h_1 = -1)$:
\[
  W_{m_{11}} = \sum_{\ell = 0}^{e-1} q^{-\tri m + 2 + \ell} E_{-1}(\ell,1,e-\ell-1) - \sum_{\ell = 0}^{e-2} q^{-\tri m + 2 + \ell} E_{-1}(\ell + 1,0,e-\ell-1).
\]

For splitting type $(1^21)$, we get in the red zone (letter types \ref{type:C} and \ref{type:D}, $d_0$ even):
\begin{align*}
  W_{m_{11},n_{11}}^{\odot} &= \sum_{k = \floor{\frac{n_c - h_\eta}{2}}}^{ e - n_c - \frac{s' + d_0' + h_\eta}{2}} q^{-m_{11}^\odot + \floor{\frac{d_0' + h_\eta}{2}} + k}E_{h_\eta} \bigg( k, \\
  &\quad{} \ceil{\frac{s' - d_0' + n_c}{2}} + h_\eta, e - \ceil{\frac{s' + d_0' + n_c}{2}} - k - h_\eta \bigg) \\
  &\quad{} - \sum_{k = \floor{\frac{n_c - h_\eta}{2}}}^{e - n_c - \frac{s' + d_0' + h_\eta}{2} - 1} q^{-m_{11}^\odot + \floor{\frac{d_0' + h_\eta}{2}} + k}E_{h_\eta} \bigg( k + 1, \\
  &\quad{} \ceil{\frac{s' - d_0' + n_c}{2}} + h_\eta - 1, e - \ceil{\frac{s' + d_0' + n_c}{2}} - k - h_\eta \bigg)
\end{align*}
and in the beige zone:
\begin{align*}
  W_{m_{11}, 0}^{\odot} =
  \sum_{k = 0}^{e - h_1/2} q^{-m_{11}^\odot + k} E(2k + h_1/2) - \sum_{k = 0}^{e - h_1/2 - 1} q^{-m_{11}^\odot + k} E(2k + h_1/2 + 1).
\end{align*}
\end{cor}

We deal with all these $E$ terms as follows. Let
\[
  \ell_E = \begin{cases}
    \floor{\dfrac{\min\{\ell_d, e\} - d_0'}{2}} & \text{splitting type $(1^2 1)$} \\
    \ell_C, & \text{other splitting types.}
  \end{cases}
\]
(Note that $\ell_d$ is defined in Lemma \ref{lem:level_parity}.)
\begin{itemize}
  \item When $\ell_0 \leq \ell_E$ and $\ell_0 \leq \ell_2$, the $E$ came from a $G(\ell_0, \ell_1, \ell_2)$ with $\ell_0 \leq \ell_E$. As we observed in the ``Further remarks'' sections following the red, yellow, and beige zones, such a $G$ is interconvertible with an $F$. Hence such an $E$ will be changed to $F$ if it appears in the final answer (after smearing and applying the $\xi_2$ restrictions: see below).
  \item When $\ell_2 \leq \ell_E$ and $\ell_2 < \ell_0$, the $E$ came from an $\epsilon_C G(\ell_2, \ell_1, \ell_0)$ with $\ell_2 \leq \ell_C$. Such an $\epsilon_C G$ is interconvertible with $\epsilon_C F$ and hence is its own Fourier transform, up to the inevitable factor of $q^e$. We annotate it as $E_{bal}$ (``bal'' for ``balanced''). The same can happen to $\epsilon_C F\xo$ and $\epsilon_C F\xo\xo$ terms, which we accordingly notate as $E_{bal} x$ and $E_{bal} xx$.
  \item Exceptionally, in the $s$ odd case of splitting types $(111)$ and $(12)$, we have
  \[
    E(0, \ell_1, e - \ell_1) = \xo F(0, \ell_1, e - \ell_1)
  \]
  for $0 \leq \ell_1 < s/2$. Similarly, there are two $E = \xo F$ cases in splitting type $(1^2 1)$ notated in the code.
  \item When $\ell_E < \ell_0$ and $\ell_C < \ell_2$, the above transformations do not apply. We annotate the $E$ as $E_{side}$ and note that, for reflection to hold, either
  \begin{itemize}
    \item The $E_{side}$ pairs with its Fourier transform, an $\epsilon_C F$ from the green zone, or
    \item The $E_{side}$ cancels with a like term for a different value of the discrete data. Indeed, decreasing $a_1$ by $1$ increases $n_{11}$ by $2$ in the red, yellow, or beige zone, causing most of the positive terms to reappear with a negative sign.
  \end{itemize}
\end{itemize}

\section{Smeared answers}
When $m_{11}$ is small, we use the smearing lemma (Lemma \ref{lem:smear}). In this subsection, we compute the appropriate smear $\Sm_{r}(W)$ for each of the answers $W$ found in the preceding subsections.

It is convenient to express as much as possible in terms of \emph{sparks} that vanish suddenly as the smear index $r$, or equivalently $m_{11}$, is decreased.

\begin{defn}
  A function $W : H^1 \to \CC$ is a \emph{spark of level $r_0$} if for all $r$, $0 \leq r \leq e' + 1$,
  \[
  \Sm_r(W) = \begin{cases}
    W & r \geq r_0 \\
    0 & r < r_0.
  \end{cases}
  \]
\end{defn}
\begin{lem}
  Let $0 \leq r_0 \leq e' + 1$. A function $W$ is a spark of level $r_0$ if and only if its Fourier transform $\widehat W$ is supported on the set $\cF_{e' - r_0} \backslash \cF_{e' - r_0 + 1}$ of elements of exact level $e' - r_0$ (or $-1/2$, in the case $r_0 = e'+1$).
\end{lem}
\begin{proof}
  Using the familiar Fourier duality between multiplication and convolution, we have the relation
  \begin{equation} \label{eq:Sm}
    \Sm_r(W) = \left( \One_{\cF_{e'-r}} \cdot \widehat{W} \right) \afterhat\, .
  \end{equation}
  So $W$ is a spark of level $r_0$ if and only if
  \[
  \One_{\cF_{\ell}} \cdot \widehat{W} = \begin{cases}
    \widehat{W} & \ell \leq e' - r_0 \\
    0 & \ell > e' - r_0.
  \end{cases}
  \]
  This evidently happens exactly when $\widehat{W}$ is supported on $\cF_{e' - r_0} \backslash \cF_{e' - r_0 + 1}$, as desired.
\end{proof}
Drawing on the repertory of Fourier transforms we computed in Lemma \ref{lem:FT_charm}, as well as the definition of $E$-functions, we get the following.
\begin{lem}\label{lem:sparks}
  Let $\fcr(\ell_0)$ be the minimal level of an element in a boxgroup denoted $T(\ell_0, \ell_1, \ell_2)$, namely
  \begin{itemize}
    \item $\fcr(\ell) = \ell$ in unramified splitting type
    \item $\fcr(\ell) = 2\ell$ in splitting type $(1^3)$ if $h = 1$
    \item $\fcr(\ell) = 2\ell+1$ in splitting type $(1^3)$ if $h = -1$
      \item $\fcr(\ell) = 2\ell + d_0'$ in splitting type $(1^21)$, letter types \ref{type:D} and \ref{type:E}
      \item $\fcr(\ell) = 2\ell + d_0'+1$ in splitting type $(1^21)$, letter type \ref{type:C}.
  \end{itemize}
  The following functions are sparks:
  \begin{enumerate}[$($a$)$]
    \item \label{sparks:green} If $ T(\ell_0,\ell_1,\ell_2) $ is defined and $ \ell_0, \ell_2 \geq \ell_E + 1 $, then
    \[
    W = \epsilon_C F(\ell_0,\ell_1,\ell_2)
    \]
    is a spark of level $e' - \fcr(\ell_E)$. This will be used in the green zone.
    \item \label{sparks:greenAB} In splitting type $(1^21)$, if $F(\ell_0)$ is defined and $\ell_0, e'-\ell_0 > \ell_d$, where $\ell_d = \ell(\hat\omega_C\diamondsuit\heartsuit)$ governs $\ell_C$, then
    \[
      W = \epsilon_C F(\ell_0)
    \]
    is a spark of level $e' - \floor{\ell_d}$. This will be used in the green zone in letter types \ref{type:A} and \ref{type:B}.
    \item \label{sparks:red} If both terms are defined and $\ell_1 \geq 1$, then
    \[
    W = q \cdot E(\ell_0,\ell_1,\ell_2) - E(\ell_0, \ell_1 - 1, \ell_2 + 1)
    \]
    is a spark of level $\fcr(e - \ell_2)$. This will be used in the red and yellow zones.
    \item \label{sparks:beige_ur} In unramified splitting type, if both terms are defined, then
    \[
    W = q^2 \cdot E(\ell_0,0,\ell_2) - E(\ell_0 - 1, 0, \ell_2 + 1)
    \]
    is a spark of level $\fcr(\ell_0)$. This will be used in the beige zone.
    \item \label{sparks:beige_ram} In ramified splitting type, if both terms are defined, then
    \[
    W = q \cdot E(\ell_0) - E(\ell_0 - 1)
    \]
    is a spark of level $\ell_0$. This will be used in the beige zone.
    \item \label{sparks:Gxx} In tame splitting type, expressions of the form
    \[
      G\xo(\ell_0,\ell_1,0), \quad \xo G\xo(0,e,0), \textand G\xo\xo(e,0,0)
    \]
    are sparks of level $e' + 1$.
    \item \label{sparks:Gxx_wild} In splitting type $(1^2 1)$, expressions of the form
    \[
      G\xo(\ell_0,\ell_1,0), \quad \xo G\xo(0,e - d_0' ,0)
    \]
    are sparks of level $e' - d_0' + 1$, while
    \[
      G\xo(e)
    \]
    is a spark of level $e' + 1$.
  \end{enumerate}
\end{lem}

Although a single $F$ is not generally a spark, we do have the relation $\cF_{e' - \fcr(\ell_2)} \subseteq T(\ell_0, \ell_1, \ell_2)$, from which $F(\ell_0, \ell_1, \ell_2)$ and $G(\ell_0, \ell_1, \ell_2)$ are stable under smears of levels $r \geq e' - \fcr(\ell_2)$.

We are now ready to compute explicit answers for the smear. Note that we do \emph{not} try to write $W_{m_{11},n_{11}}$ for each value of $m_{11}$ and $n_{11}$. Instead, we express $W_{m_{11},n_{11}}$ for $m_{11}$ large as a sum of sparks and stable terms whose appearance and disappearance can be coded simply.

One region that we do \emph{not} have to work out is the \emph{gray zone} where $m_{11}$ is so low as to satisfy all the conditions of Lemma \ref{lem:boxgps_ur}, resp.\ \ref{lem:boxgps_1pow3},
  resp.\ \ref{lem:boxgps_1pow21}
for the defining of boxgroups. There the zone total is simply
\begin{equation}
  W_{m_{11},n_{11}} = \frac{\size{H^0} q^{2e - m_{11} - n_c + \frac{d_0}{2} - (s/2 - v(N(\gamma)))}}{\size{T(\ell_0, \ell_1, \ell_2)}} \cdot F(\ell_0, \ell_1, \ell_2),
\end{equation}
where $T(\ell_0,\ell_1,\ell_2)$ is the corresponding boxgroup.

We now consider each zone in turn:
\begin{itemize}
  \item The black and brown zones need no smear, as $m_{11} \geq n_{11} > 2e$ automatically in them.
  \item The purple and blue zones have as answer a single $F$. It is stable under smears as long as $m_{11}$ is above the gray zone, since we computed the support by relaxing $\cM_{11}$ until we hit the gray zone.
  \item The green-zone answer is a sum of type $F + \epsilon_CF$. The $F$ is stable, as we got it by relaxing $\cM_{11}$ until we hit the gray zone. The $\epsilon_CF$ is a spark by Lemma \ref{lem:sparks}\ref{sparks:green}.
  \item In the red, yellow, and beige zones in the charmed coarse coset, the answer is a difference of two series of $E$'s. The $E$'s pair up to form sparks of the types in \ref{lem:sparks}, leaving one singleton (two in the beige zone), a mostly stable $F$.
  In the yellow and beige zones, if the positive sum gets cut down to $1$ or $0$ terms respectively, the negative sum has length $-1$ and must be coded up in a \emph{special zone}. This happens in a few cases, as shown in the code.
\end{itemize}

\section{The average value of a quadratic character on a box}
\label{sec:char_box}
The results in this subsection, coupled with the strong-zone answers in Section \ref{sec:strong}, yield a quick proof of Theorem \ref{thm:char_box}, which at first glance is unrelated to the topic of this memoir.

\begin{proof}[Proof of Theorem \ref{thm:char_box}]
  We may assume that $B$ contains $\pi \OO_{K}$, as enlarging $1 + B$ to $(1 + \pi \OO_K)(1 + B)$ does not change the average of a character that vanishes on $K^\cross$. Now we can take a reduced basis
  \[
  B = \pi \OO_K + \pi^{n} \OO_K\theta_1 + \pi^m \OO_K\theta_2,
  \]
  and observe that $B$ is one of the boxes that came up in Lemma \ref{lem:beta}. The strong-zone total $W_{m,n}(\delta)$ that we have computed in Section \ref{sec:strong} can also be interpreted (up to scaling) as the volume of $\beta \in B$ of class $\delta$. Hence the average in question is
  \[
  I = \frac{\sum_\delta \chi(\delta) W_{m,n}(\delta)}{\sum_\delta W_{m,n}(\delta)}
  = \frac{\widehat W_{m,n}(\chi)}{\widehat W_{m,n}(1)}.
  \]
  We wish to understand the possible values of this as $\chi$ ranges over $H^1$.
  In view of the smearing lemma (Lemma \ref{lem:smear}), as $m$ decreases, the value of the average either stays constant or becomes $0$. Since increasing $m$ only makes the theorem stronger, we can assume that we are in the case of Lemma \ref{lem:111_strong_zones}, \ref{lem:1pow3_strong_zones}, or \ref{lem:1pow21_strong_zones} according to splitting type.

  In the black, plum, purple, and blue zones, $W_{m,n}(\delta) = c \cdot F(T)$ for some subgroup $T \subseteq H^1$, so $I$ takes the value $1$ or $0$ according as $\chi(T) = 1$ or not.

  In the green zone,
  \[
  W_{m,n}(\delta) = c (F(T) + \epsilon_C F(T))
  \]
  for some boxgroup $T = T(\ell_0, \ell_1, \ell_2)$ on which $\epsilon_C$ is equidistributed, so
  \[
  I = \begin{cases}
    \big(F(T^\perp) + G(T^\perp)\big)(\chi), & T^\perp \supseteq T \\
    \ds \left(F(T^\perp) + \epsilon_C \frac{\sqrt{\size{H^1}}}{\size{T}} G(T)\right)(\chi), & T^\perp \subsetneq T.
  \end{cases}
  \]
  In either case, the identity coset $T$ is uncharmed, so the two terms have disjoint supports. In the first case, $I$ is either $1$ or $0$. In the second case, we can additionally get a value of
  \[
  \epsilon_C \frac{\sqrt{\size{H^1}}}{\size{T}} = \pm q^{\ell_0 - \ell_2}.
  \]
  The value of $i = \ell_2 - \ell_0$ evidently satisfies $1 \leq i \leq e$ by our setup of boxgroups.

  Finally, in the red zone, the Fourier transform is easier to compute using the $E$-form (Corollary \ref{cor:E-forms}), which is of the form
  \[
  W_{m,n} = \sum_{\ell = a}^{b - a} q^{c + \ell} E(\ell, e-b, b-\ell) - \sum_{\ell = a}^{b - a - 1} q^{c + \ell} E(\ell + 1, e - b - 1, b-\ell).
  \]
  The Fourier transform, by definition of $E$, is
  \begin{align*}
    \widehat{W}_{m,n} = q^c \epsilon_C \cdot \Bigg(\sum_{\ell=a}^{b-a} q^{b - \ell} F(b - \ell, e - b, \ell) - \sum_{\ell = a}^{b - a - 1} q^{b - \ell - 1} F(b - \ell, e - b - 1, \ell + 1) \Bigg).
  \end{align*}
  Scaling by $\widehat{W}_{m,n}(1) = q^{c+b-a}$, we get
  \[
  I = \epsilon_C \Bigg(\sum_{\ell=a}^{b-a} q^{a - \ell} F(b - \ell, e - b, \ell) - \sum_{\ell = a}^{b - a - 1} q^{a - \ell - 1} F(b - \ell, e - b - 1, \ell + 1) \Bigg).
  \]
  The boxgroups on which the terms are supported form the nested chain
  \[
  \begin{aligned}
  T(b-a,e-b,a)
    &\subset T(b-a,e-b-1,a+1)
     \subset T(b-a-1,e-b,a+1) \\
    &\subset \cdots \subset T(a,e-b,b-a)
  \end{aligned}
  \]
  and appear alternately with positive and negative coefficients. So there are two types of behavior upon plugging in any individual $\chi$:
  \begin{itemize}
    \item $\chi$ lies in an even number of boxgroups in the chain, and they cancel in pairs to yield $I = 0$.
    \item $\chi$ lies in an odd number of boxgroups in the chain, and only the smallest one yields a contribution $I = \pm q^{a - \ell}$. The negative of the exponent satisfies
    \[
    0 \leq \ell - a \leq b - a \leq e.
    \]
  \end{itemize}
  We must exclude the possibility that $I = -1$. This can be done by noting that $I$ is the average value of a character that takes only the values $1$ and $-1$, and that in a small neighborhood of $1 \in 1 + B$, $\chi$ is identically $1$ (continuity of $\chi$ is automatic, because all values near $1$ are squares).
\end{proof}
\chapter{Ring volumes for \texorpdfstring{$\xi'_2$}{xi'2}}
\label{cha:xi2}
\section{Tame splitting types}
\label{sec:xi2_tame}
Fix $\xi_1$ and all the data leading up to it. By Lemma \ref{lem:inactives}, there are only three possibilities for $\bar \xi_2$: either it is unrestricted, in which case the volume is given by the white-zone answer in Section \ref{sec:white}, or it is restricted by $\cM_{12}$ or $\cM_{22}$.

\subsection{\texorpdfstring{$\cM_{12}$}{M12}}
\begin{lem} \label{lem:M12_tame_tfm}
  Assume that $R$ is in a tame splitting type and $\cM_{12}$ is active. Fix $\xi_1$ satisfying the $\cM_{11}$, $\cN_{11}$ conditions, and normalize $\gamma_2$ as in Lemma \ref{lem:gamma_white}. Then $\cM_{12}$ is equivalent to a relation of the form
  \[
  \lambda^\diamondsuit(\alpha_1 \xi'_2) \equiv 0 \mod \pi^{\floor{m_{12}}}
  \]
  where $\alpha_1 \in \OO_R$ is primitive.
\end{lem}
\begin{proof}
  The $\cM_{12}$ condition says that
  \[
  \tr(\xi_1\xi_2) \equiv 0 \mod \pi^{m_{12}}.
  \]
  We have
  \[
  \tr(\xi_1\xi_2) = \tr(\xi_1\gamma_2 \xi_2') = \lambda^\diamondsuit(\diamondsuit \xi_1 \gamma_2 \xi_2').
  \]
  Observe that
  \[
  \diamondsuit \xi_1 \gamma_2 \in \pi^{-a_1' - a_2'} R.
  \]
  Let $r_{12} \in a_1' + a_2' + \ZZ$ be the unique value such that
  \[
  \alpha_1 = \pi^{r_{12}} \diamondsuit \xi_1 \gamma_2
  \]
  is a primitive vector in $\OO_R$. Then
  \[
  \cM_{12} : \lambda^\diamondsuit(\alpha_1 \xi'_2) \equiv 0 \mod \pi^{m_{12} + r_{12}},
  \]
  and the exponent is seen to be an integer. To show that it is $\floor{m_{12}}$, it's enough to prove that
  \[
  -1 < r_{12} \leq 0.
  \]
  We examine the cases.
  \begin{itemize}
    \item If $R$ is unramified and $[\delta \hat\omega_C] \in \cF_0$, then $\diamondsuit$, $\gamma_2$ are units and $\xi_1$ is primitive, so $r_{12} = 0$.
    \item If $R$ is unramified and $[\delta \hat\omega_C] \in (1; \pi; \pi) \cF_0$, then $\diamondsuit$ is a unit. Since $\cM_{11}$ is satisfied, we have $j_0 \neq 1$. Using the na\"ive choice of $\gamma_1$ from Lemma \ref{lem:gamma_white}, which differs from the $\xi_1^\odot$ by which we actually found the ring volume for $\xi_1$, we have $\gamma_1 \sim (\sqrt{\pi}; 1 ; 1)$; moreover, $(\xi_1')^Q$ is primitive. There are then two subcases, $j_0 = 2$ and $j_0 = 3$. In both cases we find that the scaling of $\alpha_1$ is controlled by the $Q$-components and $r_{12} = -1/2$ or $0$ respectively.
    \item If $R$ is totally ramified, then $\diamondsuit \sim \pi_R^2$, $\xi_1 \sim 1$, $\gamma_2 \sim 1$. According as the product $\diamondsuit \xi_1 \gamma_2$ lies in $R$, $\bar\zeta_3 R$, or $\bar\zeta_3^2 R$, we must take $r_{12} = 0$, $-2/3$, or $-1/3$ respectively.
  \end{itemize}
  Thus in all cases $-1 < r_{12} \leq 0$, as desired.
\end{proof}
This allows us to compute the ring volume for $\xi_2$:
\begin{lem}\label{lem:M12_tame}
  If $\cM_{12}$ is active, the solution volume for $\xi'_2$ is
  \[
  q^{-\floor{m_{12}}},
  \]
  except in either of the following cases, when $\cM_{12}$ has no solutions: $m_{12} \in \ZZ$ in splitting type $(1^3)$; or $R$ is unramified, $[\delta\hat\omega_C] \in (1;\pi;\pi)\cF_0$, and $(a_1,a_2,a_3) \equiv (1/2,1/2,0) \mod 1$.
\end{lem}
\begin{proof}
  By Lemma \ref{lem:M12_tame_tfm}, the $\cM_{12}$-condition is given by one of the form
  \[
  \lambda^\diamondsuit(\alpha_1 \xi'_2) \equiv 0 \mod \pi^{\floor{m_{12}}}.
  \]
  Here $\alpha_1$ is primitive, so we have a linear relation modulo $\pi^{\floor{m_{12}}}$ which yields a solution volume of $(1 + 1/q) q^{-\floor{m_{12}}}$. The solutions must be further whittled down using the restrictions on $\xi'_{2}$ in Lemma \ref{lem:gamma_white} as well as the condition that $\xi_1$ and $\xi_2$ be linearly independent mod $\mm_{\bar K}$. Here it's important to note that the solutions to $\cM_{12}$ are distributed equally among $(q + 1)$-many $1$-pixels.

  \begin{itemize}
    \item If $R$ is unramified and $[\delta \hat\omega_C] \in \cF_0$, then $a_1' \in \ZZ$. One of the $1$-pixels is that of $\xi_1$, which violates the linear independence, so we eliminate it.
    \item If $R$ is unramified, $[\delta \hat\omega_C] \in (1;\pi;\pi)\cF_0$, and
    \[
      (a_1, a_2, a_3) \equiv \left(\frac{1}{2}, 0, \frac{1}{2}\right) \mod 1,
    \]
    then $m_{12} \in \ZZ + 1/2$, and the condition that $\pi \nmid (\xi'_2)^{(K)}$ eliminates one $1$-pixel (no more, because $(\xi'_1)^Q$ is primitive).
    \item If $R$ is unramified, $[\delta \hat\omega_C] \in (1;\pi;\pi)\cF_0$, and
    \[
    (a_1, a_2, a_3) \equiv \left(\frac{1}{2}, \frac{1}{2}, 0\right) \mod 1,
    \]
    then $m_{12} \in \ZZ$ and $\gamma_1 = \gamma_2 \sim (\sqrt{\pi} ; 1 ; 1)$. When the $\cM_{12}$-condition
    \[
      \tr(\gamma_1^2 \xi_1' \xi_2') \equiv 0 \mod \pi^{m_{12}}
    \]
    is looked at mod $\pi$, it uniquely determines $\xi_2'^Q \equiv \xi_1'^Q \mod \pi$ (by $\cM_{11}$, the value $\xi_2' = \xi_1'$ is a solution). But then $\xi_2 \equiv \xi_1 \mod \sqrt{\pi}$, violating linear independence. So $\cM_{12}$ is unsatisfiable in this case.
    \item If $R$ is totally ramified, then the condition that $\xi'_2$ be a unit eliminates one $1$-pixel, unless all its solutions are non-units. This happens exactly when $\alpha_1 \sim \pi_R^2$, which is seen to be equivalent to $m_{12} \in \ZZ$.
  \end{itemize}
  Thus, outside the stated exceptional cases, we eliminate one $1$-pixel, leaving a ring volume of $q^{-\floor{m_{12}}}$.
\end{proof}

\subsection{\texorpdfstring{$\cM_{22}$}{M22}}
When $\cM_{22}$ is active, of course $\cM_{11}$ is also, and we normalize both $\gamma_1$ and $\gamma_2$ (that is, $\xi'_1$ and $\xi'_2$) according to Lemma \ref{lem:tfm_conic}. In tame splitting types, we find that the conic
\[
  \cM : \lambda^\diamondsuit (\delta^\odot {\xi_i^\odot}^2) \equiv 0 \mod m^\odot_{ii}
\]
is actually the same conic, but that the $\xi'_i$ cannot even lie in the same $1$-pixel. The following two lemmas detail when this can happen.

\begin{lem}\label{lem:A_M22_1}
Let $\cA$ be a conic of determinant $1$ on a lattice $\Lambda$ over $\OO_K$, and let $m$ be an integer, $1 \leq m \leq e$. Then there are coprimitive $\vec{x}_1, \vec{x}_2 \in \PP(\Lambda)$ satisfying
\begin{equation}\label{eq:A_M22_1}
  \cA(\vec{x}_1) \equiv \cA(\vec{x}_2) \equiv 0 \mod \pi^m
\end{equation}
if and only if the squareness $\square(\cA)$ satisfies
\[
  \square(\cA) \geq m.
\]
Moreover, if $\square(\cA) \geq m$, then for fixed $\vec{x}_1$, the volume of $\vec{x}_2$ satisfying \eqref{eq:A_M22_1} and coprimitive to $\vec{x}_1$ is $q^{-\ceil{m/2}}$, split evenly among $q$-many $1$-pixels.
\end{lem}
\begin{proof}
If coprimitive $\vec{x}_1, \vec{x}_2$ satisfy \eqref{eq:A_M22_1}, we can complete them to a basis $(\vec{x}_1, \vec{x}_2, \vec{x}_3)$ of $\Lambda$. Note that $\pi \nmid \cA(\vec{x}_3)$, or else the determinant could not be $1$. Since $m \leq e$, the cross terms vanish modulo $\pi^m$, and the determinant-one hypothesis forces $\cA(\vec{x}_3)$ to be a square unit modulo $\pi^m$. We may therefore scale $\cA$ so that $\cA(\vec{x}_3) \equiv 1 \mod \pi^m$, and now we see that $\cA$ is a square modulo $\pi^m$.

Conversely, suppose that $\cA$ is a square $\lambda^2$ of a linear form modulo $\pi^m$. Then for $\vec{x} \in \PP(\Lambda)$,
\begin{align*}
  \cA(\vec{x}) &\equiv 0 \mod \pi^m \\
  \iff \lambda(\vec{x})^2 &\equiv 0 \mod \pi^m \\
  \iff \lambda(\vec{x}) & \equiv 0 \mod \pi^{\ceil{m/2}}.
\end{align*}
Since $\lambda \nequiv 0 \mod \pi$, this has solution volume $(1 + 1/q)q^{-\ceil{m/2}}$, split evenly among $(q+1)$-many $1$-pixels. If $\vec{x}_1$ is given, then $\vec{x}_2$ can occupy any $1$-pixel except the one containing $\vec{x}_1$.
\end{proof}

The following lemma limits $m_{22}$:
\begin{lem}\label{lem:m22_le_e} The conditions on $\xi'_2$ can be satisfied only if $m_{22} \leq e$ and
  $m_{22}^{\odot} \leq e$.
\end{lem}
\begin{proof}
For the first part, note that if $m_{22} > e$, then
\[
  m_{12} = \frac{m_{11} + m_{22}}{2} - e + t > \frac{e + e}{2} - e = 0,
\]
so $\cM_{12}$ is active, contradicting Lemma \ref{lem:inactives}. In unramified splitting type, this is the entire content of the lemma, since $m_{22}^{\odot} = m_{22}$. In splitting type $(1^3)$, note that
\[
  m_{11}^\odot = m_{11} - \frac{2h_1}{3} \textand
  m_{22}^{\odot} = m_{22} - \frac{2h_2}{3}
\]
satisfy $m_{11}^\odot + m_{22}^{\odot} = m_{11} + m_{22}$, since $\{h_1, h_2\} = \{1, -1\}$. Also, $m_{11}^\odot > m_{22}^{\odot}$ since $m_{11} - m_{22} = 2(a_2 - a_1) \geq 2/3$. So if $m_{22}^{\odot} > e$, then $m_{11}^\odot > e$ and $\cM_{12}$ is active as above.
\end{proof}

This enables us to compute the ring volume for $\xi'_2$.
\begin{lem}\label{lem:M22}
Let $R$ be in a tame splitting type.
Fix the discrete data of a quartic ring and a $\xi'_1$ satisfying its conditions. The $\cM_{22}$-condition is solvable for $\xi'_2$ if and only if the following conditions are satisfied:
\begin{itemize}
  \item $[\delta \hat\omega_C] \in \cF_0$;
  \item the value of $a_2$ mod $\ZZ$ allows for a $\gamma_2$ and $m_{22}^{\odot}$ according to Lemma \ref{lem:tfm_conic}, and
  \item $m_{22}^{\odot} \leq \min\{2 \ell(\cM) + 1, e\}$.
\end{itemize}
In such cases, the volume of $\xi'_2$ is
\begin{itemize}
  \item $q^{1 - \ceil{m_{22}^{\odot}/2}}$ if $R$ is totally ramified, $h_1 = 1$ and $h_2 = -1$,
  \item $q^{-\ceil{m_{22}^{\odot}/2}}$ otherwise.
\end{itemize}
\end{lem}
\begin{proof}
The necessity of the restrictions on $m_{22}^{\odot}$ is shown by the foregoing lemmas. If they are satisfied, then $\cM_{22}$ transforms to a linear condition with solution volume $(1 + 1/q)q^{-\ceil{m_{22}^{\odot}/2}}$, distributed equally among $(q+1)$-many $1$-pixels. We must check its solutions against the other restrictions on $\xi'_2$:
\begin{itemize}
  \item If $R$ is unramified and $[\delta \hat\omega_C] \in \cF_0$, then $\xi'_1$, $\xi'_2$ are required to be coprimitive, eliminating one of the $1$-pixels.
  \item If $R$ is unramified and $[\delta \hat\omega_C] \in (1;\pi;\pi)\cF_0$ then $j_0$ must be $3$ in order for both $\gamma_1$ and $\gamma_2$ to exist. Then since $m_{22}^{\odot} = m_{22} \geq 1$, we must have two solutions $\xi'_1$, $\xi'_2$ to the transformed conic $\cM$ modulo $\pi$ whose $Q$-components are coprimitive. But $\cM(\xi') \mod \pi$ only depends on the $Q$-component and has only a unique solution in $\PP(\OO_Q/\pi \OO_Q)$, so $\cM_{22}$ is unsatisfiable if active.
  \item If $R$ has splitting type $(1^3)$, then for $\gamma_1, \gamma_2$ to exist, we must have $h_i \neq 0$ so $\{h_1, h_2\} = \{1, -1\}$. Here, coprimitivity between the $\xi_i$ is subsumed by the condition that each $\xi^\odot_i$ lie in its correct domain
  \[
    \xi^\odot_i \sim \pi_R^{1-h_i}.
  \]
  As we noted in the proof of Lemma \ref{lem:eta_1pow3}, the solutions to $\cM$ mod $\pi$ comprise $q$-many $1$-pixels of $\xi^\odot \sim 1$ and one $1$-pixel of $\xi^\odot \sim \pi_R^2$. Hence if $h_1 = -1$ and $h_2 = 1$, we retain $q$ of the $(q+1)$-many $1$-pixels, getting a volume $q^{-\ceil{m_{22}^{\odot}/2}}$. But if $h_1 = 1$ and $h_2 = -1$, then $\xi^\odot_2$ is restricted to one $1$-pixel. This gives a volume of $q^{-1-\ceil{m_{22}^{\odot}/2}}$, but we multiply back by $q^2$ since $\xi_2' = \pi_R^{-2} \xi_2^\odot$. \qedhere
\end{itemize}
\end{proof}
Of the conditions, only
\begin{equation}\label{eq:M22}
m_{22}^{\odot} \leq \min\{2 \ell(\cM) + 1, e\}
\end{equation}
is not trivial to verify. The following solves it:
\begin{lem}
Suppose that the discrete data is fixed in such a way that
\begin{itemize}
  \item $\cM_{22}$ is active,
  \item $[\delta \hat\omega_C] \in \cF_0$ (so the conic has determinant $1$),
  \item the value of $a_2$ mod $\ZZ$ allows for a $\gamma_2$ and $m_{22}^{\odot}$ according to Lemma \ref{lem:tfm_conic}.
\end{itemize}
Also suppose that $\xi_1$ is fixed, satisfying the conditions $\cM_{11}$, $\cN_{11}$ governing it. Then the remaining condition \eqref{eq:M22} can be checked as follows:
\begin{itemize}
  \item In the black, purple, and blue zones, it is automatic.
  \item In the green zone, it is equivalent to
  \[
    m_{22}^{\odot} \leq 2 \ell_C + 1.
  \]
  \item In the red, yellow, and beige zones, it restricts the sum to only use terms $G^\cross (\ell_0,\ell_1,\ell_2)$ with
  \[
    m_{22}^{\odot} \leq 2 \ell_0 + 1.
  \]
\end{itemize}
\end{lem}
\begin{proof}
When the conic $\cM$ is green, its level was computed as part of the finding of the zone total for $\xi'_1$. So it remains to prove that if the conic is black, purple, or blue, \eqref{eq:M22} is satisfied. That the conic is black, purple, or blue implies that
\[
  n_{11} \geq 2e - 4\ell(\cM) - 1.
\]
We already know $m_{22} \leq e$. Suppose that
\[
  m_{22} \geq 2 \ell(\cM) + 2.
\]
Since $m_{12}$ is inactive,
\[
  m_{11} = 2 m_{12} - m_{22} + 2e - 2t \leq 2e - 2\ell(\cM) - 2.
\]
But then
\[
  n_{22} = m_{22} - m_{11} + n_{11} \geq 3,
\]
so $\cN_{22}$ is active and we have a contradiction.
\end{proof}

\section{Splitting type (\texorpdfstring{$1^21$}{1²1})}
\label{sec:xi2_wild}
Extra complications occur when $R = K \cross Q$ is partially wildly ramified. A priori, it is the extender indices $\bar a_1, \bar a_2, \bar a_3$ that govern whether the resolvent conditions $\cM_{ij}$ and $\cN_{ij}$ are active. However, we wish to count rings based on their reduced basis vectors $\rho_i = \pi^{a_i} \xi_i$. By Lemma \ref{lem:a2bar}, $\bar a_2$ is related to the $a_i$:
\begin{equation}
  \label{eq:a2bar_copy}
  \bar a_2 = \min \left\{a_2 + v^{(K)}(\xi_1),\; a_2 + v^{(K)}(\xi_2),\; a_2 + \frac{d_0 - 1}{2},\; a_3\right\}.
\end{equation}
We will say that a situation is of \emph{$\bar{a}_2$-type $\xi_1$, $\xi_2$, $d_0$, or $a_3$} according as the four respective arguments achieve the minimum in \eqref{eq:a2bar_copy}. If more than one argument achieves the minimum, we prefer one $\bar{a}_2$-type over another in the order
\[
  d_0 \succ a_3 \succ \xi_1 \succ \xi_2.
\]
Note that $\bar{a}_2 = a_2$ in flavors that force one of the $v^{(K)}(\xi_i)$-terms to equal $0$ (the ``easy'' flavors $002$, $123$, $213$, $231$, and $321$).

For $1 \leq i \leq 2$, let
\[
\kappa_i = \min\left\{v^{(K)}(\xi_i), \frac{d_0 - 1}{2}\right\},
\]
observing that
\[
\bar a_2 = a_2 + \min\left\{\kappa_1, \kappa_2, \frac{d_0 - 1}{2}, a_3 - a_2\right\}.
\]

We first examine $\kappa_1$.
\begin{lem}\label{lem:kappa_1} For any $\xi_1$ satisfying the $\cM_{11}$ and $\cN_{11}$ conditions, if $m_{11}$ is large enough that $m_{11} > 2e$ and $m_{11}^\odot > 2e$, the value of
  \[
  \kappa_1 = \min\left\{v^{(K)}(\xi_1), \frac{d_0 - 1}{2}\right\}
  \]
  is constrained as follows:
  \begin{enumerate}[$($a$)$]
    \item In the yellow and darker zones, it depends only on the values of $d_0$ and $s'$:
    \begin{itemize}
      \item In letter type \ref{type:A},
      \[
      \kappa_1 = \frac{d_0 - 1}{4}.
      \]
      \item In letter type \ref{type:B},
      \[
      \kappa_1 = \frac{2s' + d_0}{4}.
      \]
      \item In letter types \ref{type:C}--\ref{type:E},
      \[
      \kappa_1 = \frac{d_0 - 1}{2}.
      \]
    \end{itemize}
    \item In the lemon and beige zones, it depends only on what term $G(\ell)$ or $G^\cross(\ell)$ we are in:
    \begin{itemize}
      \item On $\xo G(0)$,
      \[
      \kappa_1 = \frac{d_0 - 1}{4}.
      \]
      \item On $G^\cross(\ell)$, $0 \leq \ell < d_0/2 - 1$,
      \[
      \kappa_1 = \frac{2\ell + d_0}{4}.
      \]
      \item Within $G(\ceil{d_0/2} - 1)$,
      \[
      \kappa_1 = \frac{d_0 - 1}{2}.
      \]
    \end{itemize}
  \end{enumerate}
\end{lem}
\begin{proof}
  The cases have been divided up by the value of $\ell = \ell\left(\delta^\odot\right) = \ell\left(\delta_0 \hat\omega_C \diamondsuit\right)$.

  If $\ell = -1/2$, we must have $h_1 = 1$, $\delta^{\odot(Q)} \sim \pi_Q, \xi_1^{\odot(Q)} \sim 1$. Then
  \[
  \fI\left(\delta^{\odot(Q)}{\xi_1^{\odot(Q)}}^2\right) \sim 1,
  \]
  so to satisfy $\cM_{11}$, we must have
  \[
  v^{(K)}\left(\xi_1^\odot\right) = 0,
  \]
  so by Table \ref{tab:tfm_conic},
  \[
  v^{(K)}\left(\xi_1\right) = \frac{d_0 - h_1}{4} = \frac{d_0 - 1}{4},
  \]
  as desired.

  If $0 \leq \ell < d_0/2 - 1$, then $h_1$ is either $0$ or $2$. It is easy to show that, for any $\alpha \in \OO_Q$,
  \[
  v_K\left(\fI(\alpha^2)\right) \geq \floor{\frac{d_0}{2}}.
  \]
  We have
  \[
  v\left(\fI\left(\delta^{\odot(Q)}{\xi_1^{\odot(Q)}}^2\right)\right) = \ell + \frac{h_1}{2},
  \]
  because the square ${\xi_1^{\odot(Q)}}^2$ is too close to being in $K$ to cancel the main term of $\fI\left(\delta^{\odot(Q)}\right)$. So to satisfy $\cM_{11}$,
  \[
  v^{(K)}\left(\delta^{\odot}{\xi_1^\odot}^2\right) = \ell + \frac{h_1}{2}.
  \]
  Since $\delta^{\odot(K)}$ is a unit, this forces $h_1/2 \equiv \ell \mod 2$ and
  \[
  v^{(K)}\left(\xi_1^{\odot}\right) = \frac{\ell + \frac{h_1}{2}}{2} = \frac{2\ell + h_1}{4}.
  \]
  Hence
  \[
  v^{(K)}\left(\xi_1\right) = \frac{2\ell + h_1}{4} + \frac{d_0 - h_1}{4} = \frac{2\ell + d_0}{4},
  \]
  as desired.

  Finally, if $\ell \geq d_0/2 - 1$, then $h_1$ is either $0$ or $2$. We have
  \[
  v\left(\fI\left(\delta^{\odot(Q)}{\xi_1^{\odot(Q)}}^2\right)\right) \geq \frac{d_0}{2} - 1 + \frac{h_1}{2},
  \]
  so to satisfy $\cM_{11}$,
  \[
  v^{(K)}\left(\delta^{\odot}{\xi_1^\odot}^2\right) \geq \frac{d_0}{2} - 1 + \frac{h_1}{2}.
  \]
  Since $\delta^{\odot(K)}$ is a unit, the $\delta^{\odot}$ factor can be dropped, and
  \[
  v^{(K)}\left(\xi_1\right) \geq \frac{\frac{d_0}{2} - 1 + \frac{h_1}{2}}{2} + \frac{d_0 - h_1}{4} = \frac{d_0 - 1}{2}.
  \]
  So $\kappa_1 = (d_0 - 1)/2$, as desired.
\end{proof}

Now fix the discrete data $a_1, a_2, \bar a_2$. Let
\[
m_{12} = b_2 + 3t - 3e - a_1 - a_2 + s \textand
m_{22} = b_2 + 2t - 2e - 2a_2 + s,
\]
that is, $m_{ij}$ is computed like $\bar m_{ij}$ using $a_i$ in place of $\bar a_i$. Let
\[
  u = \ceil{\bar{a}_2 - a_2 - \frac{o_2}{4}}.
\]

The finite flavor analyses in the next four lemmas are also checked by the accompanying file \verb|quartic_ON_1e21.sage|.
\begin{lem}\label{lem:xi2free} If neither $\cM_{12}$ nor $\cM_{22}$ is active, the volume $Z$ of $\xi_2'$ satisfying the needed restriction on its $(K)$-valuation is as follows:
  \begin{equation}
    \begin{tabular}{r | cc}
      & \multicolumn{2}{c}{---$\bar a_2$-type---} \\
      $Z$ & $\xi_2$ & Others \\ \hline
      Flavor $200$, $u = 0$ & $1$ & $1 + 1/q$ \\
      Other cases & $(1 - 1/q)q^{-u}$ & $q^{-u}$
    \end{tabular}
  \end{equation}
\end{lem}
\begin{proof}
The value of $\bar a_2$ requires us to find $\xi'_2$ satisfying
\[
  v^{(K)}(\xi'_2) = u,
\]
in $\bar{a}_2$-type $\xi_2$, or
\[
  v^{(K)}(\xi'_2) \geq u,
\]
in the other cases. Flavor 200 is exceptional because $\xi_2$ can be either $K$- or $Q$-led. (In flavor 002, the linear independence between $\xi_1$ and $\xi_2$ eliminates that extra possibility.) The computation is routine in all cases.
\end{proof}

\begin{lem}\label{lem:M12_Q}
  Suppose we are in flavor $(o_i)_i = 020, 200, 132, \text{ or } 312$, that $\bar a_2 > a_2$, and $\cM_{12}$ is active but $\cM_{22}$ is not. In $\bar{a}_2$-type $a_3$, we must have
  \begin{equation}\label{eq:M12_Q_a3}
    v^{(K)}(\xi_1) \geq \bar m_{12},
  \end{equation}
  and then $\cM_{12}$ is automatic, that is, $Z$ is the same as in the preceding lemma. Otherwise, the answer depends on which of the three quantities
  \[
    m_{12}, \quad d_0', \quad u' = u + \kappa_1 + \frac{o_2}{4}
  \]
  is smallest:
  \begin{enumerate}[$($a$)$]
    \item If
    \[
      m_{12} \leq u' \textand m_{12} \leq d_0',
    \]
    then $\cM_{12}$ is automatic, that is, $Z$ is the same as in the preceding lemma.
    \item If
    \[
      d_0' < m_{12} \textand d_0' \leq u',
    \]
    then there are solutions only in flavors 020 and 200 (for $d_0$ even) or 132 and 312 (for $d_0$ odd). The solution volume is
    \[
      Z = \begin{cases*}
        \left(1 - \dfrac{1}{q}\right) q^{-u + d_0' - m_{12}}, & $\bar a_2$-type $\xi_2$ \\
        q^{-u + d_0' - m_{12}}, & otherwise.
      \end{cases*}
    \]
    \item If
    \[
      u' < d_0' \textand u' < m_{12},
    \]
    then the solution volume is
    \[
      Z = q^{-u + u' - m_{12}} = q^{-m_{12} + v^{(K)}(\xi_1) + o_2/4},
    \]
    except when $u' < d_0' - 1/2$ in $\bar a_2$-type $\xi_2$, where $Z = 0$.
  \end{enumerate}
\end{lem}

\begin{proof}
The case of $\bar{a}_2$-type $a_3$ is exceptional because there $\bar\xi_2 = \xi_3$, so $\cM_{12}$ is independent of $\xi_2$ and equivalent to
\[
  \tr(\xi_1\xi_3) \equiv 0 \mod \bar m_{12}.
\]
As $\xi_3$ is $K$-led, the $K$-term dominates, yielding \eqref{eq:M12_Q_a3}.

In the other cases, $\cM_{12}$ is equivalent to
\[
  \tr(\xi_1\xi_2) \equiv 0 \mod m_{12}.
\]
As all conditions are linear, the volume of $\xi_2$ is straightforward to compute.
\end{proof}

\begin{lem}\label{lem:M12_K}
Suppose we are in flavor $(o_i)_i = 002, 200, 123, \text{ or } 321$, that $\bar a_2 = a_2$ and that $\cM_{12}$ is active but $\cM_{22}$ is not. Then the volume $Z$ of permissible $\xi_2'$ depends on which of the three quantities
\[
  m_{12}, \quad \kappa_1, \quad d' = \floor{d_0'} + \frac{(o_3 - o_1 - 2)\One_{2\mid d_0} + o_1}{4}
\]
is smallest:
\begin{enumerate}[$($a$)$]
  \item If
  \[
    m_{12} \leq \kappa_1 \textand m_{12} \leq d',
  \]
  then $\cM_{12}$ is automatic and the ring volume is the same as in Lemma \ref{lem:xi2free}.
  \item If
  \[
    d' \leq \kappa_1 \textand d' < m_{12},
  \]
  then the ring volume is
  \[
    Z = q^{d' - m_{12}}.
  \]
  \item If $\kappa_1 < m_{12}$ and $\kappa_1 < d'$, then $Z = 0$.
\end{enumerate}
\end{lem}
\begin{proof}
Since $\bar{a}_2 = a_2$, we must have a $K$-led $\xi_2$. The proof is a straightforward analysis of linear conditions.
\end{proof}

\begin{lem}
Suppose that $\cM_{22}$ is active but $\cM_{12}$ is not. There are only solutions in $\bar{a}_2$-types $\xi_1$, $\xi_2$, and $d_0$ and flavors 020, 200, 132, and 312. The solution volume depends on which of the three quantities
\[
  \frac{m_{22}}{2}, \quad \bar a_2 - a_2, \quad \frac{d_0}{4} + \ell_d
\]
is smallest, where
\[
  \ell_d = \ell(\delta^\odot \heartsuit)
\]
is either constant or sufficiently large on each term of each zone total for $\xi_1$ that we can regard it as a constant:
\begin{enumerate}[$($a$)$]
  \item If
  \[
    \frac{m_{22}}{2} \leq \frac{d_0}{4} + \ell_d \textand \frac{m_{22}}{2} \leq \bar a_2 - a_2,
  \]
  then $\cM_{22}$ is automatic and the solution volume is the same as in Lemma \ref{lem:xi2free}.
  \item If
  \[
    \bar a_2 - a_2 < \frac{m_{22}}{2} \textand \bar a_2 - a_2 \leq \frac{d_0}{4} + \ell_d,
  \]
  then $\cM_{22}$ is satisfiable only in $\bar a_2$-types $\xi_1$ and $d_0$. Its solution volume is
  \[
    Z = q^{-\ceil{\frac{m_{22}}{2} - \frac{o_2}{4}}}.
  \]
  \item If ${d_0}/{4} + \ell_d < \frac{m_{22}}{2}$ and ${d_0}/{4} + \ell_d < \bar a_2 - a_2$, then $Z = 0$.
\end{enumerate}
\end{lem}
\begin{proof}
As in the other splitting types, the condition that $\tr(\xi^2) \equiv 0$ has two linearly independent solutions $\xi = \xi_1$, $\xi_2$ constrains $m_{22}$ to be small enough that $\cM_{22}$ is equivalent to a linear condition
\[
  \lambda(\xi_2') \equiv 0 \mod \pi^{m'}.
\]
Analyzing the leading term of $\lambda$ and comparing with the other linear conditions on $\xi_2$ yields the claimed results.
\end{proof}

\chapter{Further remarks on the code}
\label{cha:code}
In the attached code, we use the computer programs Sage and LattE to compute the generating function of rings. First we count ``zone tuples'' consisting of integer values of the following variables:
\begin{itemize}
  \item $e, t, b_1, b_2, s$.
  \item $\verb|a1f| = \floor{a_1}, \verb|a2f| = \floor{a_2}, o_1, o_2, o_3$. Here we've decomposed
  \[
    a_i = \floor{a_i} + \frac{o_i}{o},
  \]
  where
  \[
    o = \begin{cases}
    2 & \text{in splitting types $(111)$, $(12)$, and $(3)$} \\
    3 & \text{in splitting type $(1^3)$}
     \\ 4 & \text{in splitting type $(1^21)$.}
\end{cases}.
  \]
  The $o_i$ belong to one of a finite number of ``flavors'' coding the classes of the $a_i$ and $b_i$ mod $\ZZ$.

  Note that $\floor{a_3}$ is missing from the variable list, as its value is uniquely determined by the discriminant identity (Lemma \ref{lem:rsv}).
  \item The squareness parameter. In unramified splitting type it is coded as
  \[
    \verb|lCf| = \begin{cases}
      -1, & \square_C = 0\\
      \ell_C, & 0 < \square_C < e \\
      \floor{\dfrac{e}{2}}, & \square_C = e.
    \end{cases}
  \]
  This unambiguously determines $\square_C$. We note that in our answers, $\square_C$ only appears in the answer and bounds of the green zone. When $\square_C = e$, there is no green zone and its bounds are far from being achieved, so we can equate $\floor{\ell_C}$ with $\verb|lCf|$ without affecting the answers. In splitting type $(1^3)$, we take $\ell_C$ itself as our variable. In splitting type $(1^2 1)$, it is most convenient to use $\verb|ld0|$, the $\ell_d$ that appears in Lemma \ref{lem:level_parity} when $\delta = 1$. Note that \verb|ld0| is determined by the other data, except in letter types \ref{type:C} and \ref{type:D}, $d_0$ even. There its parity is known, so we introduce a variable $\ell_m$ and set
  \[
    \verb|ld0| = d_0' + 2\ell_m - h_\eta.
  \]
  \item $k$, an index of summation used to input the answers in the red, yellow, and beige zones.
\end{itemize}

If $\floor{m_{11}/2}$ appears in a zone answer, we add a variable $\verb|m11fl|$ and impose one of the equations
  \[
    m_{11} = 2 \cdot \verb|m11fl| \textor m_{11} = 2 \cdot \verb|m11fl| + 1.
  \]
  Similarly we treat $\floor{n_{11}/2}$, $\floor{n_c/2}$, $\tilde{n}$, $\floor{m_{22}/2}$, $\floor{e/2}$, $\floor{s/2}$. We use the convention that appending \verb|f| to a variable name takes its floor, while appending \verb|fl| or \verb|ceil| divides by $2$ (sometimes translating by an appropriate constant) and \emph{then} takes the floor or ceiling, respectively.

In splitting type $1^2 1$, we use \verb|vK_xi1| and \verb|vK_xi2| to denote variables that either equal the $v^{(K)}(\xi_i)$ or are large enough when these valuations no longer matter (as in the $\kappa_i$ of Lemma \ref{lem:kappa_1}).

The code sets
\[
  \verb|RRR| = \ZZ[E, T, B_1, B_2, S, \verb|A1F|,\ldots]
\]
and represents the lattice-point generating function in $\Frac(\texttt{RRR})$.
Here we use the convention that the value of a lowercase variable appears as an exponent of the corresponding uppercase variable. For instance, the generating function of the three lattice points on the line segment $y = 2x$, $0 \leq x \leq 2$ is
\[
  1 + XY^2 + X^2Y^4.
\]

Because each variable is bounded below in terms of the preceding variables, the power series is formally convergent. Because each zone is delimited by finitely many linear inequalities with $\ZZ$-coefficients, the generating function is a rational function, computed by Barvinok's algorithm as implemented in LattE. The code also sets
\[
  \verb|RINGS_RING| = \ZZ[\verb|E_|,\verb|T_|,\verb|LCF_|,\verb|SFL_|,\verb|B1_|,\verb|B2_|,q,\verb|L0|,\verb|L2|]
\]
and encodes ring totals as substitutions in $\Frac(\texttt{RINGS\_RING})$.
The trailing underscores are to prevent the computer from conflating certain elements of $\verb|RRR|$ with their counterparts in $\verb|RINGS_RING|$, although the reader can think of them as identified. Two new variables $\verb|L0|,\verb|L2|$ (\verb|G0| and \verb|G1| in splitting type $1^2 1$), whose exponents are the $\ell_0$ and $\ell_2$ of a boxgroup, complete the description of a ring total, along with a string \verb|Ftype| that tells the kind of weighting ($F$, $F\xo\xo$, $E_{bal}$, etc.).

\section{The impish case}\label{sec:impish}
Some comments are in order about the impish case of Lemma \ref{lem:impish}, where $q = 2$, $b_1 = b_2$, $C_\ttt = \OO_K + \pi^{b_1}\OO_K^3$. Here, if $m_{ij}$ and $n_{ij}$ had been computed from the formulas of Lemma \ref{lem:rsv}, we would get strange relations such as $n_{ij} > m_{ij}$, violating the activity implications of Figure \ref{fig:inactives} on page \pageref{fig:inactives} and leading to a profusion of exceptional cases. Fortunately, we can compute ring totals much more simply.

Assume that $\bar b_1 = \bar b_2$ (we do not assume that $q = 2$).  Let
\[
  (\pi^{a_1}\xi_{1,\OO}, \pi^{a_2}\xi_{2,\OO}, \pi^{a_3}\xi_{3,\OO})
\]
be a reduced basis for $\sqrt{\delta} I_\OO$. Multiplication by a unit in $\OO_{\bar K}^3$ relates this to the extender basis $(\pi^{a_i}\bar \xi_i)_i$ for $\bar I = \omega \bar I_\OO$; we do not speak of a reduced basis for $\bar I$, since $\omega_C \in \left(\OO_{\bar K}^3\right)^\cross$ makes sense over $\OO_{\bar K}$ but not over $\OO_K$.

Then, with moduli defined using the situation over $\OO_{\bar K}$, we have $m_{ij} = n_{ij}$ and $s = 0$, leading to the following simplifications:
\begin{itemize}
  \item $\cM_{22}$ cannot be active, as then $\cN_{22}$ would be active.
  \item $\cM_{12}$ cannot be active, as then $\cN_{12}$ would be active.
  \item The $\cM_{11}$ and $\cN_{11}$ conditions can be written more simply as
  \[
  \text{All coordinates of } \xi_1^2 \text{ are congruent} \mod \pi^{m_{11}}.
  \]
  This constrains $\delta$ to lie in the level space $\cF_v$, where $v$ depends on $a_1$ and $b_1$ but not $q$, and the volume of $\xi_1$ is independent of $\delta \in \cF_v$. On the other hand, when $q > 2$ and the ring volume is computed by the foregoing method, we get answers in the black, gray, and white zones of the form $c q^{-u} F(v, 0, e-v)$, where the constants $c, u, v$ depend on $a_1$ and $b_1$ but not on $q$. Consequently, the answers thereby derived are valid in the impish case as well. Thus we enter our answers into the code as they have been computed, with no special treatment of the impish case.

\end{itemize}

\appendix

\chapter{Examples of zone totals}
\label{cha:zone_examples}
The following tables serve to illustrate some of the totals computed.

\section{Tame local quartic O-N}
We can now add up the number of rings with given reduced indices $(a_1, a_2, a_3)$ depending on what zone they lie in. We first devise a shorthand for notating subsets of the space $H^1(M_{K^3})$ of quartic algebras with totally split resolvent. We place them on a grid in the following manner:
\begin{center}
  \begin{tabular}{|cc|cc|}
    \hline
    $(1;1;1)$ & $(u;u;1)$ & $(\pi;\pi;1)$ & $(u\pi;u\pi;1)$ \\
    $(1;u;u)$ & $(u;1;u)$ & $(\pi;u\pi;u)$ & $(u\pi;\pi;u)$ \\ \hline
    $(1;\pi;\pi)$ & $(u;u\pi;\pi)$ & $(\pi;1;\pi)$ & $(u\pi;u;\pi)$ \\
    $(1;u\pi;u\pi)$ & $(u;\pi;u\pi)$ & $(\pi;u;u\pi)$ & $(u\pi;1;u\pi)$ \\ \hline
  \end{tabular}
\end{center}
noting that:
\begin{itemize}
  \item The unramified $L$ occupy the upper left quadrant;
  \item Split $L$'s occupy the left column, main diagonal, and top row;
  \item The $\epsilon_C$-charmed coset of $H^1_\ur$ is in the top left or right quadrant, according as $s$ is even or odd.
\end{itemize}
We can then indicate subsets of $H^1(M_{K^3})$ by shading the appropriate cells of a $4 \times 4$ grid. With this notation, the support and thickness of the families of rings with fixed discrete data is shown in Table \ref{tab:zones_odd}.
\begin{table}[htbp]
  \footnotesize
  \setlength{\tabcolsep}{2pt}
  \begin{tabular}{ccccl}
    & & \multicolumn{2}{c}{Support} \\
    Zone & Flavor & $s$ even & $s$ odd & Number of rings with given discrete data \\ \hline
    white & 0 &
    \support{\sh03\sh13\sh02\sh12} &
    \support{\sh33\sh23\sh32\sh22} &
    $ \begin{cases}
      1 & a_1 = a_2 = a_3 \\
      q^{2a_3 - 2a_1}\big(1 + \frac{1}{q} + \frac{1}{q^2}\big) & a_1 = a_2 < a_3 \\
      & \text{ or } a_1 < a_2 = a_3 \\
      q^{2a_3 - 2a_1}\big(1 + \frac{1}{q}\big)\big(1 + \frac{1}{q} + \frac{1}{q^2}\big) & a_1 < a_2 < a_3 \\
    \end{cases} $
    \\
    white & 1,2,3 &
    \support{\sh00\sh01\sh10\sh11\sh20\sh21\sh22\sh23\sh30\sh31\sh32\sh33} &
    \support{\sh00\sh01\sh10\sh11\sh03\sh13\sh02\sh12\sh20\sh21\sh30\sh31} &
    $\begin{cases}
      q^{2a_3 - 2a_1 - 1} & a_1 = a_2 < a_3 \\
      & \text{ or } a_1 < a_2 = a_3 \\
      q^{2a_3 - 2a_1 - 1}\big(1 + \frac{1}{q}\big) & a_1 < a_2 < a_3 \\
    \end{cases}$
    \\
    \begin{tabular}{c}
      beige \\
      ($\cM_{12}$ inactive)
    \end{tabular}
      & 0 &
    \support{\sh03\sh13\sh02\sh12} &
    \support{\sh33\sh23\sh32\sh22} &
    $\begin{cases}
      q^{2a_3 - 2a_1 - m_{11}}\big(1 + \frac{1}{q}\big) & a_1 < a_2 = a_3 \\
      q^{2a_3 - 2a_1 - m_{11}}\big(1 + \frac{1}{q}\big)^2 & a_1 < a_2 < a_3 \\
    \end{cases}$
    \\
    \quickcol{beige \\ ($\cM_{12}$ inactive)} & 2,3 &
    \quickcol{\support{\sh00\sh01\sh30\sh21\sh23\sh33}\\or sim.} &
    \quickcol{\support{\sh30\sh31\sh00\sh11\sh13\sh03}\\or sim.} &
    $2q^{2a_3 - 2a_1 - m_{11} - 1}$
    \\
    \begin{tabular}{c}
beige \\
($\cM_{12}$ active)
\end{tabular} & 0 &
    \support{\sh03\sh13\sh02\sh12} &
    \support{\sh33\sh23\sh32\sh22} &
    $q^{2a_3 - 2a_1 - m_{11} - m_{12}}\big(1 + \frac{1}{q}\big)$
    \\
    \quickcol{beige \\ ($\cM_{12}$ active)} & 2 &
    \quickcol{\support{\sh00\sh01\sh30\sh21\sh23\sh33}\\ or sim.} &
    \quickcol{\support{\sh30\sh31\sh00\sh11\sh13\sh03}\\ or sim.} &
    $2q^{2a_3 - 2a_1 - m_{11} - m_{12} - 1/2}$
    \\
    \begin{tabular}{c}
brown \\
($\cM_{12}$ inactive)
\end{tabular} & 0 &
    \support{\sh03\sh13} &
    \support{\sh23\sh33} &
    $\begin{cases}
      2q^{2a_3 - 2a_1 - m_{11} - \ceil{n_{11}/2}} & a_1 < a_2 = a_3 \\
      2q^{2a_3 - 2a_1 - m_{11} - \ceil{n_{11}/2}}\big(1 + \frac{1}{q}\big) & a_1 < a_2 < a_3 \\
    \end{cases}$
    \\
    \quickcol{beige \\ ($\cM_{12}$ inactive)} & 2,3 &
    \support{\sh23\sh33} &
    \support{\sh03\sh13} &
    $2q^{2a_3 - 2a_1 - m_{11} - \ceil{(n_{11} + 1)/2}}$
    \\
    \begin{tabular}{c}
brown \\
($\cM_{12}$ active)
\end{tabular} & 0 &
    \support{\sh03\sh13} &
    \support{\sh23\sh33} &
    $2q^{2a_3 - 2a_1 - m_{11} - m_{12} - \ceil{n_{11}/2}}$ \\
    \quickcol{beige \\ ($\cM_{12}$ active)} & 2 &
    \support{\sh23\sh33} &
    \support{\sh03\sh13} &
    $2q^{2a_3 - 2a_1 - m_{11} - m_{12} - \ceil{(n_{11} + 1)/2} + 1/2}$ \\
    \begin{tabular}{c}
black \\
($\cM_{12}$ inactive)
\end{tabular} & 0 &
    \support{\sh03} &
    \ignore{Nothing} &
    $\begin{cases}
      4q^{2a_3 - 2a_1 - m_{11} - n_{11} + s/2} & a_1 < a_2 = a_3 \\
      4q^{2a_3 - 2a_1 - m_{11} - n_{11} + s/2}\big(1 + \frac{1}{q}\big) & a_1 < a_2 < a_3 \\
    \end{cases}$
    \\
    \quickcol{beige \\ ($\cM_{12}$ inactive)} & 2,3 &
    \ignore{Nothing} &
    \support{\sh03} &
    $4q^{2a_3 - 2a_1 - m_{11} - n_{11} + (s-1)/2}$
    \\
    \begin{tabular}{c}
black \\
($\cM_{12}$ active)
\end{tabular} & 0 &
    \support{\sh03} &
    \ignore{Nothing} &
    $4q^{2a_3 - 2a_1 - m_{11} - m_{12} - n_{11} + s/2}$
    \\
    \quickcol{beige \\ ($\cM_{12}$ active)} & 2 &
    \ignore{Nothing} &
    \support{\sh03} &
    $4q^{2a_3 - 2a_1 - m_{11} - m_{12} - n_{11} + s/2}$
  \end{tabular}
  \caption{The possible zones for the enumeration of subrings of a quartic ring (tame, split resolvent)}
  \label{tab:zones_odd}
\end{table}

Note that in many cases, the combination of zone and flavor force $a_1 < a_2 < a_3$, so only one ring total is given.

We abbreviate the flavors 000, 011, 101, 110 to 0, 1, 2, and 3 respectively. The disclaimer ``or similar'' is used in the beige zone because the value of $[\hat \omega_C]$ affects which six cells occur in the support.
\input{tables.tex}

\backmatter

\bibliographystyle{amsplain}
\bibliography{Master,me}
\end{document}

%% file: zonedefs_color.tex
\newcommand{\blackzone}{\textcolor{black}{$\blacksquare$~}}
\newcommand{\plumzone}{\textcolor{violet!90!black}{$\blacklozenge$~}}
\newcommand{\purplezone}{\textcolor{violet!70!magenta}{$\CIRCLE$~}}
\newcommand{\bluezone}{\textcolor{blue}{$\blacktriangledown$~}}
\newcommand{\greenzone}{\textcolor{green!70!black}{$\blacktriangleright$~}}
\newcommand{\redzone}{\textcolor{red}{$\blacktriangle$~}}
\newcommand{\brownzone}{\textcolor{brown}{$\blacksquare$~}}
\newcommand{\yellowzone}{\rlap{\textcolor{yellow!90!black}{$\blacktriangledown$~}}$\triangledown$~}
\newcommand{\lemonzone}{\rlap{\textcolor{yellow!30}{$\blacktriangleright$~}}$\triangleright$~}
\newcommand{\whitezone}{{$\lozenge$~}}

%% file: tables.tex
\section{Wild local quartic O-N}

For brevity, the following conventions have been observed throughout the wild cases:
\begin{itemize}
  \item In the red, yellow, and lemon zones, the answer (which is independent of $\ell_C$) is the sum of the entries of the corresponding row.
  \item In the other color zones, the answer is the single entry corresponding to the appropriate values of $\ell_C$ and $n_{11}$.
  \item The invariable factor of $q^{-m_{11}^\odot}$ has been omitted. It is the only dependency on $m_{11}$, provided that $m_{11} > 2e$.
\end{itemize}

\newpage
\subsection{Unramified}
The following are valid in the unramified splitting types ($(111)$, $(12)$, and $(3)$): more specifically, the first table (with $s=0$) applies to all three, the remaining ones to splitting types $(111)$ and $(12)$.

\noindent\begin{minipage}{\textwidth}
  \centering
  \begin{sideways}
      \begin{tabular}{l}
For $e = 8, s = 0$: \blackzone black, \bluezone blue, \greenzone green, \redzone red, and \whitezone white zones.$\vphantom{|}$
\\ \begin{tabular}{r|lllll}
      $n_{11}$ & $\ell_C = 0$ & $1$ & $2$ & $3$ & $\geq4$ \\ \hline
      $17$ &
      \blackzone{$\size{H^0} q^{-1} F(8,\emptyset,\emptyset)$} &
      \blackzone{$\size{H^0} q^{-1} F(8,\emptyset,\emptyset)$} &
      \blackzone{$\size{H^0} q^{-1} F(8,\emptyset,\emptyset)$} &
      \blackzone{$\size{H^0} q^{-1} F(8,\emptyset,\emptyset)$} &
      \blackzone{\tiny $\size{H^0} q^{-1} F(8,\emptyset,\emptyset)$} $\vphantom{0^{0^0}}$ 
      \\
      $16$ & \bluezone{$q^{0} F(8,0,0)$} & \bluezone{$q^{0} F(8,0,0)$} & \bluezone{$q^{0} F(8,0,0)$} & \bluezone{$q^{0} F(8,0,0)$} & \bluezone{$q^{0} F(8,0,0)$} \\
      $15$ & \bluezone{$q^{0} F(7,1,0)$} & \bluezone{$q^{0} F(7,1,0)$} & \bluezone{$q^{0} F(7,1,0)$} & \bluezone{$q^{0} F(7,1,0)$} & \bluezone{$q^{0} F(7,1,0)$} \\
      $14$ & \greenzone{$(1 + \epsilon_C)q^{0} F(7,0,1)$} & \bluezone{$q^{0} F(7,0,1)$} & \bluezone{$q^{0} F(7,0,1)$} & \bluezone{$q^{0} F(7,0,1)$} & \bluezone{$q^{0} F(7,0,1)$} \\
      $13$ & \greenzone{$(1 + \epsilon_C)q^{0} F(6,1,1)$} & \bluezone{$q^{0} F(6,1,1)$} & \bluezone{$q^{0} F(6,1,1)$} & \bluezone{$q^{0} F(6,1,1)$} & \bluezone{$q^{0} F(6,1,1)$} \\
      $12$ & \greenzone{$(1 + \epsilon_C)q^{0} F(6,0,2)$} & \bluezone{$q^{1} F(6,1,1)$} & \bluezone{$q^{1} F(6,1,1)$} & \bluezone{$q^{1} F(6,1,1)$} & \bluezone{$q^{1} F(6,1,1)$} \\
      $11$ & \greenzone{$(1 + \epsilon_C)q^{0} F(5,1,2)$} & \bluezone{$q^{1} F(5,2,1)$} & \bluezone{$q^{1} F(5,2,1)$} & \bluezone{$q^{1} F(5,2,1)$} & \bluezone{$q^{1} F(5,2,1)$} \\
      $10$ & \greenzone{$(1 + \epsilon_C)q^{0} F(5,0,3)$} & \greenzone{$(1 + \epsilon_C)q^{1} F(5,1,2)$} & \bluezone{$q^{1} F(5,1,2)$} & \bluezone{$q^{1} F(5,1,2)$} & \bluezone{$q^{1} F(5,1,2)$} \\
      $9$ & \greenzone{$(1 + \epsilon_C)q^{0} F(4,1,3)$} & \greenzone{$(1 + \epsilon_C)q^{1} F(4,2,2)$} & \bluezone{$q^{1} F(4,2,2)$} & \bluezone{$q^{1} F(4,2,2)$} & \bluezone{$q^{1} F(4,2,2)$} \\
      $8$ & \greenzone{$(1 + \epsilon_C)q^{0} F(4,0,4)$} & \greenzone{$(1 + \epsilon_C)q^{1} F(4,1,3)$} & \bluezone{$q^{2} F(4,2,2)$} & \bluezone{$q^{2} F(4,2,2)$} & \bluezone{$q^{2} F(4,2,2)$} \\
      $7$ & \greenzone{$(1 + \epsilon_C)q^{0} F(3,1,4)$} & \greenzone{$(1 + \epsilon_C)q^{1} F(3,2,3)$} & \bluezone{$q^{2} F(3,3,2)$} & \bluezone{$q^{2} F(3,3,2)$} & \bluezone{$q^{2} F(3,3,2)$} \\
      $6$ & \greenzone{$(1 + \epsilon_C)q^{0} F(3,0,5)$} & \greenzone{$(1 + \epsilon_C)q^{1} F(3,1,4)$} & \greenzone{$(1 + \epsilon_C)q^{2} F(3,2,3)$} & \bluezone{$q^{2} F(3,2,3)$} & \bluezone{$q^{2} F(3,2,3)$} \\
      $5$ & \greenzone{$(1 + \epsilon_C)q^{0} F(2,1,5)$} & \greenzone{$(1 + \epsilon_C)q^{1} F(2,2,4)$} & \redzone{$(1 + \epsilon_C)q^{2} G^\cross(2,3,3)$} & \redzone{$q^{2} G(3,2,3)$} &  \\
      $4$ & \greenzone{$(1 + \epsilon_C)q^{0} F(2,0,6)$} & \greenzone{$(1 + \epsilon_C)q^{1} F(2,1,5)$} & \redzone{$(1 + \epsilon_C)q^{2} G^\cross(2,2,4)$} & \redzone{$q^{3} G(3,2,3)$} &  \\
      $3$ & \greenzone{$(1 + \epsilon_C)q^{0} F(1,1,6)$} & \redzone{$(1 + \epsilon_C)q^{1} G^\cross(1,2,5)$} & \redzone{$(1 + \epsilon_C)q^{2} G^\cross(2,2,4)$} & \redzone{$q^{3} G(3,2,3)$} &  \\
      $2$ & \greenzone{$(1 + \epsilon_C)q^{0} F(1,0,7)$} & \redzone{$(1 + \epsilon_C)q^{1} G^\cross(1,1,6)$} & \redzone{$(1 + \epsilon_C)q^{2} G^\cross(2,1,5)$} & \redzone{$(1 + \epsilon_C)q^{3} G^\cross(3,1,4)$} & \redzone{$q^{3} G(4,0,4)$} \\
      $1$ & \redzone{$(1 + \epsilon_C)q^{0} G^\cross(0,1,7)$} & \redzone{$(1 + \epsilon_C)q^{1} G^\cross(1,1,6)$} & \redzone{$(1 + \epsilon_C)q^{2} G^\cross(2,1,5)$} & \redzone{$(1 + \epsilon_C)q^{3} G^\cross(3,1,4)$} & \redzone{$q^{3} G(4,0,4)$} \\
      $\leq 0$ & \whitezone{$(1 + \epsilon_C)q^{0} G^\cross(0,0,8)$} & \whitezone{$(1 + \epsilon_C)q^{1} G^\cross(1,0,7)$} & \whitezone{$(1 + \epsilon_C)q^{2} G^\cross(2,0,6)$} & \whitezone{$(1 + \epsilon_C)q^{3} G^\cross(3,0,5)$} & \whitezone{$q^{4}\left(1 + \frac{1}{q}\right) G(4,0,4)$}
    \end{tabular}
\end{tabular}
    \vspace{0.1in}
  \end{sideways}
\end{minipage}

\noindent\begin{minipage}{\textwidth}
  \centering
  \begin{tabular}{c}
For $e = 6, s = 14$: \purplezone purple zone.$\vphantom{|}$
  
  \\ \begin{tabular}{r|l}
    $n_c$ & $\ell_C = \text{any}$ \\ \hline$\vphantom{0^{0^0}}$
    $12$ & \purplezone{$2q^{0} F(6,0, \emptyset)$} \\
    $11$ & \purplezone{$2q^{0} F(5,1, \emptyset)$} \\
    $10$ & \purplezone{$2q^{1} F(5,1, \emptyset)$} \\
    $9$ & \purplezone{$2q^{1} F(4,2, \emptyset)$} \\
    $8$ & \purplezone{$2q^{2} F(4,2, \emptyset)$} \\
    $7$ & \purplezone{$2q^{2} F(3,3, \emptyset)$} \\
    $6$ & \purplezone{$2q^{3} F(3,3, \emptyset)$} \\
    $5$ & \purplezone{$2q^{3} F(2,4, \emptyset)$} \\
    $4$ & \purplezone{$2q^{4} F(2,4, \emptyset)$} \\
    $3$ & \purplezone{$2q^{4} F(1,5, \emptyset)$} \\
    $2$ & \purplezone{$2q^{5} F(1,5, \emptyset)$} \\
    $1$ & \purplezone{$2q^{5} F(0,6, \emptyset)$} 
  \end{tabular}
\end{tabular}
  \vspace{0.1in}
\end{minipage}
\vfill
\noindent\begin{minipage}{\textwidth}
  \centering
  \begin{tabular}{c}
For $e = 11, s = 7$: \purplezone purple and \greenzone green zones.$\vphantom{|}$
  
  \\ \begin{tabular}{r|l}
    $n_c$ & $\ell_C = 3$ \\ \hline$\vphantom{0^{0^0}}$$22$ & \purplezone{$2q^{0} F(11,0, \emptyset)$} \\
    $21$ & \purplezone{$2q^{0} F(10,1, \emptyset)$} \\
    $20$ & \purplezone{$2q^{1} F(10,1, \emptyset)$} \\
    $19$ & \purplezone{$2q^{1} F(9,2, \emptyset)$} \\
    $18$ & \purplezone{$2q^{2} F(9,2, \emptyset)$} \\
    $17$ & \purplezone{$2q^{2} F(8,3, \emptyset)$} \\
    $16$ & \purplezone{$2q^{3} F(8,3, \emptyset)$} \\
    $15$ & \greenzone{$(1 + \epsilon_C)q^{3} F(7,4,0)$} \\
    $14$ & \greenzone{$(1 + \epsilon_C)q^{3} F(7,3,1)$} \\
    $13$ & \greenzone{$(1 + \epsilon_C)q^{3} F(6,4,1)$} \\
    $12$ & \greenzone{$(1 + \epsilon_C)q^{3} F(6,3,2)$} \\
    $11$ & \greenzone{$(1 + \epsilon_C)q^{3} F(5,4,2)$} \\
    $10$ & \greenzone{$(1 + \epsilon_C)q^{3} F(5,3,3)$} \\
    $9$ & \greenzone{$(1 + \epsilon_C)q^{3} F(4,4,3)$} \\
    $8$ & \greenzone{$(1 + \epsilon_C)q^{3} F(4,3,4)$} \\
    $7$ & \greenzone{$(1 + \epsilon_C)q^{3} F(3,4,4)$} \\
    $6$ & \greenzone{$(1 + \epsilon_C)q^{3} F(3,3,5)$} \\
    $5$ & \greenzone{$(1 + \epsilon_C)q^{3} F(2,4,5)$} \\
    $4$ & \greenzone{$(1 + \epsilon_C)q^{3} F(2,3,6)$} \\
    $3$ & \greenzone{$(1 + \epsilon_C)q^{3} F(1,4,6)$} \\
    $2$ & \greenzone{$(1 + \epsilon_C)q^{3} F(1,3,7)$} \\
    $1$ & \greenzone{$(1 + \epsilon_C)q^{3} F(0,4,7)$} 
  \end{tabular}
\end{tabular}
  \vspace{0.1in}
\end{minipage}

\noindent\begin{minipage}{\textwidth}
  \begin{sideways}
    \begin{tabular}{l}
For $e = 10, s = 4$: \purplezone purple, \bluezone blue, \greenzone green, and \redzone red zones.$\vphantom{|}$
  
    \\ \begin{tabular}{r|lllll}
      $n_{11}$ & $\ell_C = 0$ & $1$ & $2$ & $3$ & $\geq4$ \\ \hline$\vphantom{0^{0^0}}$$20$ & \purplezone{$2q^{0} F(10,0, \emptyset)$} & \purplezone{$2q^{0} F(10,0, \emptyset)$} & \purplezone{$2q^{0} F(10,0, \emptyset)$} & \purplezone{$2q^{0} F(10,0, \emptyset)$} & \purplezone{$2q^{0} F(10,0, \emptyset)$} \\
      $19$ & \purplezone{$2q^{0} F(9,1, \emptyset)$} & \purplezone{$2q^{0} F(9,1, \emptyset)$} & \purplezone{$2q^{0} F(9,1, \emptyset)$} & \purplezone{$2q^{0} F(9,1, \emptyset)$} & \purplezone{$2q^{0} F(9,1, \emptyset)$} \\
      $18$ & \purplezone{$2q^{1} F(9,1, \emptyset)$} & \purplezone{$2q^{1} F(9,1, \emptyset)$} & \purplezone{$2q^{1} F(9,1, \emptyset)$} & \purplezone{$2q^{1} F(9,1, \emptyset)$} & \purplezone{$2q^{1} F(9,1, \emptyset)$} \\
      $17$ & \purplezone{$2q^{1} F(8,2, \emptyset)$} & \purplezone{$2q^{1} F(8,2, \emptyset)$} & \purplezone{$2q^{1} F(8,2, \emptyset)$} & \purplezone{$2q^{1} F(8,2, \emptyset)$} & \purplezone{$2q^{1} F(8,2, \emptyset)$} \\
      $16$ & \bluezone{$q^{2} F(8,2,0)$} & \bluezone{$q^{2} F(8,2,0)$} & \bluezone{$q^{2} F(8,2,0)$} & \bluezone{$q^{2} F(8,2,0)$} & \bluezone{$q^{2} F(8,2,0)$} \\
      $15$ & \bluezone{$q^{2} F(7,3,0)$} & \bluezone{$q^{2} F(7,3,0)$} & \bluezone{$q^{2} F(7,3,0)$} & \bluezone{$q^{2} F(7,3,0)$} & \bluezone{$q^{2} F(7,3,0)$} \\
      $14$ & \greenzone{$(1 + \epsilon_C)q^{2} F(7,2,1)$} & \bluezone{$q^{2} F(7,2,1)$} & \bluezone{$q^{2} F(7,2,1)$} & \bluezone{$q^{2} F(7,2,1)$} & \bluezone{$q^{2} F(7,2,1)$} \\
      $13$ & \greenzone{$(1 + \epsilon_C)q^{2} F(6,3,1)$} & \bluezone{$q^{2} F(6,3,1)$} & \bluezone{$q^{2} F(6,3,1)$} & \bluezone{$q^{2} F(6,3,1)$} & \bluezone{$q^{2} F(6,3,1)$} \\
      $12$ & \greenzone{$(1 + \epsilon_C)q^{2} F(6,2,2)$} & \bluezone{$q^{3} F(6,3,1)$} & \bluezone{$q^{3} F(6,3,1)$} & \bluezone{$q^{3} F(6,3,1)$} & \bluezone{$q^{3} F(6,3,1)$} \\
      $11$ & \greenzone{$(1 + \epsilon_C)q^{2} F(5,3,2)$} & \bluezone{$q^{3} F(5,4,1)$} & \bluezone{$q^{3} F(5,4,1)$} & \bluezone{$q^{3} F(5,4,1)$} & \bluezone{$q^{3} F(5,4,1)$} \\
      $10$ & \greenzone{$(1 + \epsilon_C)q^{2} F(5,2,3)$} & \greenzone{$(1 + \epsilon_C)q^{3} F(5,3,2)$} & \bluezone{$q^{3} F(5,3,2)$} & \bluezone{$q^{3} F(5,3,2)$} & \bluezone{$q^{3} F(5,3,2)$} \\
      $9$ & \greenzone{$(1 + \epsilon_C)q^{2} F(4,3,3)$} & \greenzone{$(1 + \epsilon_C)q^{3} F(4,4,2)$} & \bluezone{$q^{3} F(4,4,2)$} & \bluezone{$q^{3} F(4,4,2)$} & \bluezone{$q^{3} F(4,4,2)$} \\
      $8$ & \greenzone{$(1 + \epsilon_C)q^{2} F(4,2,4)$} & \greenzone{$(1 + \epsilon_C)q^{3} F(4,3,3)$} & \bluezone{$q^{4} F(4,4,2)$} & \bluezone{$q^{4} F(4,4,2)$} & \bluezone{$q^{4} F(4,4,2)$} \\
      $7$ & \greenzone{$(1 + \epsilon_C)q^{2} F(3,3,4)$} & \greenzone{$(1 + \epsilon_C)q^{3} F(3,4,3)$} & \bluezone{$q^{4} F(3,5,2)$} & \bluezone{$q^{4} F(3,5,2)$} & \bluezone{$q^{4} F(3,5,2)$} \\
      $6$ & \greenzone{$(1 + \epsilon_C)q^{2} F(3,2,5)$} & \greenzone{$(1 + \epsilon_C)q^{3} F(3,3,4)$} & \greenzone{$(1 + \epsilon_C)q^{4} F(3,4,3)$} & \bluezone{$q^{4} F(3,4,3)$} & \bluezone{$q^{4} F(3,4,3)$} \\
      $5$ & \greenzone{$(1 + \epsilon_C)q^{2} F(2,3,5)$} & \greenzone{$(1 + \epsilon_C)q^{3} F(2,4,4)$} & \redzone{$(1 + \epsilon_C)q^{4} G^\cross(2,5,3)$} & \redzone{$q^{4} G(3,4,3)$} &  \\
      $4$ & \greenzone{$(1 + \epsilon_C)q^{2} F(2,2,6)$} & \greenzone{$(1 + \epsilon_C)q^{3} F(2,3,5)$} & \redzone{$(1 + \epsilon_C)q^{4} G^\cross(2,4,4)$} & \redzone{$q^{5} G(3,4,3)$} &  \\
      $3$ & \greenzone{$(1 + \epsilon_C)q^{2} F(1,3,6)$} & \redzone{$(1 + \epsilon_C)q^{3} G^\cross(1,4,5)$} & \redzone{$(1 + \epsilon_C)q^{4} G^\cross(2,4,4)$} & \redzone{$q^{5} G(3,4,3)$} &  \\
      $2$ & \greenzone{$(1 + \epsilon_C)q^{2} F(1,2,7)$} & \redzone{$(1 + \epsilon_C)q^{3} G^\cross(1,3,6)$} & \redzone{$(1 + \epsilon_C)q^{4} G^\cross(2,3,5)$} & \redzone{$(1 + \epsilon_C)q^{5} G^\cross(3,3,4)$} & \redzone{$q^{5} G(4,2,4)$} \\
      $1$ & \redzone{$(1 + \epsilon_C)q^{2} G^\cross(0,3,7)$} & \redzone{$(1 + \epsilon_C)q^{3} G^\cross(1,3,6)$} & \redzone{$(1 + \epsilon_C)q^{4} G^\cross(2,3,5)$} & \redzone{$(1 + \epsilon_C)q^{5} G^\cross(3,3,4)$} & \redzone{$q^{5} G(4,2,4)$} 
    \end{tabular}
\end{tabular}
  \end{sideways}
  \vspace{0.1in}
\end{minipage}

\subsection{Splitting type (\texorpdfstring{$1^3$}{1³})}\ \\[\baselineskip]
\noindent\begin{minipage}{\textwidth}
  \centering
     \begin{tabular}{c}
For $e = 4$, $h = 1$: \blackzone black, \bluezone blue, \greenzone green, and \redzone red zones.$\vphantom{|}$
  
  \\ \begin{tabular}{r|lll}
    $\tri n_{11}$ & $\ell_C = 0$ & $1$ & $2$ \\ \hline$\vphantom{0^{0^0}}9$ & \blackzone{$q^{0} F(4, 0, 0)$} & \blackzone{$q^{0} F(4, 0, 0)$} & \blackzone{$q^{0} F(4, 0, 0)$} \\
    $8$ & \greenzone{$(1 + \epsilon_C)q^{0} F(4,-1,1)$} & \bluezone{$q^{0} F(4,-1,1)$} & \bluezone{$q^{0} F(4,-1,1)$} \\
    $7$ & \greenzone{$(1 + \epsilon_C)q^{0} F(3,0,1)$} & \bluezone{$q^{0} F(3,0,1)$} & \bluezone{$q^{0} F(3,0,1)$} \\
    $6$ & \greenzone{$(1 + \epsilon_C)q^{0} F(3,-1,2)$} & \bluezone{$q^{1} F(3,0,1)$} & \bluezone{$q^{1} F(3,0,1)$} \\
    $5$ & \greenzone{$(1 + \epsilon_C)q^{0} F(2,0,2)$} & \bluezone{$q^{1} F(2,1,1)$} & \bluezone{$q^{1} F(2,1,1)$} \\
    $4$ & \greenzone{$(1 + \epsilon_C)q^{0} F(2,-1,3)$} & \greenzone{$(1 + \epsilon_C)q^{1} F(2,0,2)$} & \bluezone{$q^{1} F(2,0,2)$} \\
    $3$ & \greenzone{$(1 + \epsilon_C)q^{0} F(1,0,3)$} & \redzone{$(1 + \epsilon_C)q^{1} G^\cross(1,1,2)$} & \redzone{$q^{1} G(2,0,2)$} \\
    $2$ & \greenzone{$(1 + \epsilon_C)q^{0} F(1,-1,4)$} & \redzone{$(1 + \epsilon_C)q^{1} G^\cross(1,0,3)$} & \redzone{$q^{2} G(2,0,2)$} \\
    $1$ & \redzone{$(1 + \epsilon_C)q^{0} G^\cross(0,0,4)$} & \redzone{$(1 + \epsilon_C)q^{1} G^\cross(1,0,3)$} & \redzone{$q^{2} G(2,0,2)$} 
  \end{tabular}
\end{tabular}
\end{minipage}
\vfill

\noindent\begin{minipage}{\textwidth}
  \centering
     \begin{tabular}{c}
For $e = 4$, $h = -1$: \blackzone black, \bluezone blue, \greenzone green, and \redzone red zones.$\vphantom{|}$
  
  \\ \begin{tabular}{r|lll}
    $\tri n_{11}$ & $\ell_C = 0$ & $1$ & $2$ \\ \hline$\vphantom{0^{0^0}}9$ & \blackzone{$q^{0} F(4, 0, 0)$} & \blackzone{$q^{0} F(4, 0, 0)$} & \blackzone{$q^{0} F(4, 0, 0)$} \\
    $8$ & \blackzone{$q^{1} F(4, 0, 0)$} & \blackzone{$q^{1} F(4, 0, 0)$} & \blackzone{$q^{1} F(4, 0, 0)$} \\
    $7$ & \bluezone{$q^{1} F(3,1,0)$} & \bluezone{$q^{1} F(3,1,0)$} & \bluezone{$q^{1} F(3,1,0)$} \\
    $6$ & \bluezone{$q^{2} F(3,1,0)$} & \bluezone{$q^{2} F(3,1,0)$} & \bluezone{$q^{2} F(3,1,0)$} \\
    $5$ & \bluezone{$q^{2} F(2,2,0)$} & \bluezone{$q^{2} F(2,2,0)$} & \bluezone{$q^{2} F(2,2,0)$} \\
    $4$ & \greenzone{$(1 + \epsilon_C)q^{2} F(2,1,1)$} & \bluezone{$q^{2} F(2,1,1)$} & \bluezone{$q^{2} F(2,1,1)$} \\
    $3$ & \greenzone{$(1 + \epsilon_C)q^{2} F(1,2,1)$} & \bluezone{$q^{2} F(1,2,1)$} & \bluezone{$q^{2} F(1,2,1)$} \\
    $2$ & \greenzone{$(1 + \epsilon_C)q^{2} F(1,1,2)$} & \bluezone{$q^{3} F(1,2,1)$} & \bluezone{$q^{3} F(1,2,1)$} \\
    $1$ & \redzone{$(1 + \epsilon_C)q^{2} G^\cross(0,2,2)$} & \redzone{$q^{3} G(1,2,1)$} &  \\
    
  \end{tabular}
\end{tabular}
  \vspace{0.1in}
\end{minipage}
\vfill\vfill\newpage

\noindent\begin{minipage}{\textwidth}
\begin{sideways}
  \centering
     \begin{tabular}{c}
For $e = 9$, $h = 1$: \blackzone black, \bluezone blue, \greenzone green, and \redzone red zones.$\vphantom{|}$
  
    \\ \begin{tabular}{r|lllll}
      $\tri n_{11}$ & $\ell_C = 0$ & $1$ & $2$ & $3$ & $4$ \\ \hline$\vphantom{0^{0^0}}19$ & \blackzone{$q^{0} F(9, 0, 0)$} & \blackzone{$q^{0} F(9, 0, 0)$} & \blackzone{$q^{0} F(9, 0, 0)$} & \blackzone{$q^{0} F(9, 0, 0)$} & \blackzone{$q^{0} F(9, 0, 0)$} \\
      $18$ & \greenzone{$(1 + \epsilon_C)q^{0} F(9,-1,1)$} & \bluezone{$q^{0} F(9,-1,1)$} & \bluezone{$q^{0} F(9,-1,1)$} & \bluezone{$q^{0} F(9,-1,1)$} & \bluezone{$q^{0} F(9,-1,1)$} \\
      $17$ & \greenzone{$(1 + \epsilon_C)q^{0} F(8,0,1)$} & \bluezone{$q^{0} F(8,0,1)$} & \bluezone{$q^{0} F(8,0,1)$} & \bluezone{$q^{0} F(8,0,1)$} & \bluezone{$q^{0} F(8,0,1)$} \\
      $16$ & \greenzone{$(1 + \epsilon_C)q^{0} F(8,-1,2)$} & \bluezone{$q^{1} F(8,0,1)$} & \bluezone{$q^{1} F(8,0,1)$} & \bluezone{$q^{1} F(8,0,1)$} & \bluezone{$q^{1} F(8,0,1)$} \\
      $15$ & \greenzone{$(1 + \epsilon_C)q^{0} F(7,0,2)$} & \bluezone{$q^{1} F(7,1,1)$} & \bluezone{$q^{1} F(7,1,1)$} & \bluezone{$q^{1} F(7,1,1)$} & \bluezone{$q^{1} F(7,1,1)$} \\
      $14$ & \greenzone{$(1 + \epsilon_C)q^{0} F(7,-1,3)$} & \greenzone{$(1 + \epsilon_C)q^{1} F(7,0,2)$} & \bluezone{$q^{1} F(7,0,2)$} & \bluezone{$q^{1} F(7,0,2)$} & \bluezone{$q^{1} F(7,0,2)$} \\
      $13$ & \greenzone{$(1 + \epsilon_C)q^{0} F(6,0,3)$} & \greenzone{$(1 + \epsilon_C)q^{1} F(6,1,2)$} & \bluezone{$q^{1} F(6,1,2)$} & \bluezone{$q^{1} F(6,1,2)$} & \bluezone{$q^{1} F(6,1,2)$} \\
      $12$ & \greenzone{$(1 + \epsilon_C)q^{0} F(6,-1,4)$} & \greenzone{$(1 + \epsilon_C)q^{1} F(6,0,3)$} & \bluezone{$q^{2} F(6,1,2)$} & \bluezone{$q^{2} F(6,1,2)$} & \bluezone{$q^{2} F(6,1,2)$} \\
      $11$ & \greenzone{$(1 + \epsilon_C)q^{0} F(5,0,4)$} & \greenzone{$(1 + \epsilon_C)q^{1} F(5,1,3)$} & \bluezone{$q^{2} F(5,2,2)$} & \bluezone{$q^{2} F(5,2,2)$} & \bluezone{$q^{2} F(5,2,2)$} \\
      $10$ & \greenzone{$(1 + \epsilon_C)q^{0} F(5,-1,5)$} & \greenzone{$(1 + \epsilon_C)q^{1} F(5,0,4)$} & \greenzone{$(1 + \epsilon_C)q^{2} F(5,1,3)$} & \bluezone{$q^{2} F(5,1,3)$} & \bluezone{$q^{2} F(5,1,3)$} \\
      $9$ & \greenzone{$(1 + \epsilon_C)q^{0} F(4,0,5)$} & \greenzone{$(1 + \epsilon_C)q^{1} F(4,1,4)$} & \greenzone{$(1 + \epsilon_C)q^{2} F(4,2,3)$} & \bluezone{$q^{2} F(4,2,3)$} & \bluezone{$q^{2} F(4,2,3)$} \\
      $8$ & \greenzone{$(1 + \epsilon_C)q^{0} F(4,-1,6)$} & \greenzone{$(1 + \epsilon_C)q^{1} F(4,0,5)$} & \greenzone{$(1 + \epsilon_C)q^{2} F(4,1,4)$} & \bluezone{$q^{3} F(4,2,3)$} & \bluezone{$q^{3} F(4,2,3)$} \\
      $7$ & \greenzone{$(1 + \epsilon_C)q^{0} F(3,0,6)$} & \greenzone{$(1 + \epsilon_C)q^{1} F(3,1,5)$} & \greenzone{$(1 + \epsilon_C)q^{2} F(3,2,4)$} & \bluezone{$q^{3} F(3,3,3)$} & \bluezone{$q^{3} F(3,3,3)$} \\
      $6$ & \greenzone{$(1 + \epsilon_C)q^{0} F(3,-1,7)$} & \greenzone{$(1 + \epsilon_C)q^{1} F(3,0,6)$} & \greenzone{$(1 + \epsilon_C)q^{2} F(3,1,5)$} & \redzone{$(1 + \epsilon_C)q^{3} G^\cross(3,2,4)$} & \redzone{$q^{3} G(4,1,4)$} \\
      $5$ & \greenzone{$(1 + \epsilon_C)q^{0} F(2,0,7)$} & \greenzone{$(1 + \epsilon_C)q^{1} F(2,1,6)$} & \redzone{$(1 + \epsilon_C)q^{2} G^\cross(2,2,5)$} & \redzone{$(1 + \epsilon_C)q^{3} G^\cross(3,2,4)$} & \redzone{$q^{3} G(4,1,4)$} \\
      $4$ & \greenzone{$(1 + \epsilon_C)q^{0} F(2,-1,8)$} & \greenzone{$(1 + \epsilon_C)q^{1} F(2,0,7)$} & \redzone{$(1 + \epsilon_C)q^{2} G^\cross(2,1,6)$} & \redzone{$(1 + \epsilon_C)q^{3} G^\cross(3,1,5)$} & \redzone{$q^{4} G(4,1,4)$} \\
      $3$ & \greenzone{$(1 + \epsilon_C)q^{0} F(1,0,8)$} & \redzone{$(1 + \epsilon_C)q^{1} G^\cross(1,1,7)$} & \redzone{$(1 + \epsilon_C)q^{2} G^\cross(2,1,6)$} & \redzone{$(1 + \epsilon_C)q^{3} G^\cross(3,1,5)$} & \redzone{$q^{4} G(4,1,4)$} \\
      $2$ & \greenzone{$(1 + \epsilon_C)q^{0} F(1,-1,9)$} & \redzone{$(1 + \epsilon_C)q^{1} G^\cross(1,0,8)$} & \redzone{$(1 + \epsilon_C)q^{2} G^\cross(2,0,7)$} & \redzone{$(1 + \epsilon_C)q^{3} G^\cross(3,0,6)$} & \redzone{$(1 + \epsilon_C)q^{4} G^\cross(4,0,5)$} \\
      $1$ & \redzone{$(1 + \epsilon_C)q^{0} G^\cross(0,0,9)$} & \redzone{$(1 + \epsilon_C)q^{1} G^\cross(1,0,8)$} & \redzone{$(1 + \epsilon_C)q^{2} G^\cross(2,0,7)$} & \redzone{$(1 + \epsilon_C)q^{3} G^\cross(3,0,6)$} & \redzone{$(1 + \epsilon_C)q^{4} G^\cross(4,0,5)$} 
    \end{tabular}
\end{tabular}
    \vspace{0.1in}
  \end{sideways}
\end{minipage}

\noindent\begin{minipage}{\textwidth}
\begin{sideways}
  \centering
     \begin{tabular}{c}
For $e = 9$, $h = -1$: \blackzone black, \bluezone blue, \greenzone green, and \redzone red zones.$\vphantom{|}$
  
    \\ \begin{tabular}{r|lllll}
      $\tri n_{11}$ & $\ell_C = 0$ & $1$ & $2$ & $3$ & $4$ \\ \hline$\vphantom{0^{0^0}}19$ & \blackzone{$q^{0} F(9, 0, 0)$} & \blackzone{$q^{0} F(9, 0, 0)$} & \blackzone{$q^{0} F(9, 0, 0)$} & \blackzone{$q^{0} F(9, 0, 0)$} & \blackzone{$q^{0} F(9, 0, 0)$} \\
      $18$ & \blackzone{$q^{1} F(9, 0, 0)$} & \blackzone{$q^{1} F(9, 0, 0)$} & \blackzone{$q^{1} F(9, 0, 0)$} & \blackzone{$q^{1} F(9, 0, 0)$} & \blackzone{$q^{1} F(9, 0, 0)$} \\
      $17$ & \bluezone{$q^{1} F(8,1,0)$} & \bluezone{$q^{1} F(8,1,0)$} & \bluezone{$q^{1} F(8,1,0)$} & \bluezone{$q^{1} F(8,1,0)$} & \bluezone{$q^{1} F(8,1,0)$} \\
      $16$ & \bluezone{$q^{2} F(8,1,0)$} & \bluezone{$q^{2} F(8,1,0)$} & \bluezone{$q^{2} F(8,1,0)$} & \bluezone{$q^{2} F(8,1,0)$} & \bluezone{$q^{2} F(8,1,0)$} \\
      $15$ & \bluezone{$q^{2} F(7,2,0)$} & \bluezone{$q^{2} F(7,2,0)$} & \bluezone{$q^{2} F(7,2,0)$} & \bluezone{$q^{2} F(7,2,0)$} & \bluezone{$q^{2} F(7,2,0)$} \\
      $14$ & \greenzone{$(1 + \epsilon_C)q^{2} F(7,1,1)$} & \bluezone{$q^{2} F(7,1,1)$} & \bluezone{$q^{2} F(7,1,1)$} & \bluezone{$q^{2} F(7,1,1)$} & \bluezone{$q^{2} F(7,1,1)$} \\
      $13$ & \greenzone{$(1 + \epsilon_C)q^{2} F(6,2,1)$} & \bluezone{$q^{2} F(6,2,1)$} & \bluezone{$q^{2} F(6,2,1)$} & \bluezone{$q^{2} F(6,2,1)$} & \bluezone{$q^{2} F(6,2,1)$} \\
      $12$ & \greenzone{$(1 + \epsilon_C)q^{2} F(6,1,2)$} & \bluezone{$q^{3} F(6,2,1)$} & \bluezone{$q^{3} F(6,2,1)$} & \bluezone{$q^{3} F(6,2,1)$} & \bluezone{$q^{3} F(6,2,1)$} \\
      $11$ & \greenzone{$(1 + \epsilon_C)q^{2} F(5,2,2)$} & \bluezone{$q^{3} F(5,3,1)$} & \bluezone{$q^{3} F(5,3,1)$} & \bluezone{$q^{3} F(5,3,1)$} & \bluezone{$q^{3} F(5,3,1)$} \\
      $10$ & \greenzone{$(1 + \epsilon_C)q^{2} F(5,1,3)$} & \greenzone{$(1 + \epsilon_C)q^{3} F(5,2,2)$} & \bluezone{$q^{3} F(5,2,2)$} & \bluezone{$q^{3} F(5,2,2)$} & \bluezone{$q^{3} F(5,2,2)$} \\
      $9$ & \greenzone{$(1 + \epsilon_C)q^{2} F(4,2,3)$} & \greenzone{$(1 + \epsilon_C)q^{3} F(4,3,2)$} & \bluezone{$q^{3} F(4,3,2)$} & \bluezone{$q^{3} F(4,3,2)$} & \bluezone{$q^{3} F(4,3,2)$} \\
      $8$ & \greenzone{$(1 + \epsilon_C)q^{2} F(4,1,4)$} & \greenzone{$(1 + \epsilon_C)q^{3} F(4,2,3)$} & \bluezone{$q^{4} F(4,3,2)$} & \bluezone{$q^{4} F(4,3,2)$} & \bluezone{$q^{4} F(4,3,2)$} \\
      $7$ & \greenzone{$(1 + \epsilon_C)q^{2} F(3,2,4)$} & \greenzone{$(1 + \epsilon_C)q^{3} F(3,3,3)$} & \bluezone{$q^{4} F(3,4,2)$} & \bluezone{$q^{4} F(3,4,2)$} & \bluezone{$q^{4} F(3,4,2)$} \\
      $6$ & \greenzone{$(1 + \epsilon_C)q^{2} F(3,1,5)$} & \greenzone{$(1 + \epsilon_C)q^{3} F(3,2,4)$} & \greenzone{$(1 + \epsilon_C)q^{4} F(3,3,3)$} & \bluezone{$q^{4} F(3,3,3)$} & \bluezone{$q^{4} F(3,3,3)$} \\
      $5$ & \greenzone{$(1 + \epsilon_C)q^{2} F(2,2,5)$} & \greenzone{$(1 + \epsilon_C)q^{3} F(2,3,4)$} & \redzone{$(1 + \epsilon_C)q^{4} G^\cross(2,4,3)$} & \redzone{$q^{4} G(3,3,3)$} &  \\
      $4$ & \greenzone{$(1 + \epsilon_C)q^{2} F(2,1,6)$} & \greenzone{$(1 + \epsilon_C)q^{3} F(2,2,5)$} & \redzone{$(1 + \epsilon_C)q^{4} G^\cross(2,3,4)$} & \redzone{$q^{5} G(3,3,3)$} &  \\
      $3$ & \greenzone{$(1 + \epsilon_C)q^{2} F(1,2,6)$} & \redzone{$(1 + \epsilon_C)q^{3} G^\cross(1,3,5)$} & \redzone{$(1 + \epsilon_C)q^{4} G^\cross(2,3,4)$} & \redzone{$q^{5} G(3,3,3)$} &  \\
      $2$ & \greenzone{$(1 + \epsilon_C)q^{2} F(1,1,7)$} & \redzone{$(1 + \epsilon_C)q^{3} G^\cross(1,2,6)$} & \redzone{$(1 + \epsilon_C)q^{4} G^\cross(2,2,5)$} & \redzone{$(1 + \epsilon_C)q^{5} G^\cross(3,2,4)$} & \redzone{$q^{5} G(4,1,4)$\!} \\
      $1$ & \redzone{$(1 + \epsilon_C)q^{2} G^\cross(0,2,7)$} & \redzone{$(1 + \epsilon_C)q^{3} G^\cross(1,2,6)$} & \redzone{$(1 + \epsilon_C)q^{4} G^\cross(2,2,5)$} & \redzone{$(1 + \epsilon_C)q^{5} G^\cross(3,2,4)$} & \redzone{$q^{5} G(4,1,4)$\!}
    \end{tabular}
\end{tabular}
  \end{sideways}
  \vspace{0.1in}
\end{minipage}

\newpage
\subsection{Splitting type (\texorpdfstring{$1^21$}{1²1}) (dark zones)}\ \\[\baselineskip]
\noindent\begin{minipage}[t]{.5\textwidth}
For $e = 6$, any $d_0'$, $s' = -1/2$, type (A): \blackzone black and \greenzone green zones.$\vphantom{|}$
  
  \begin{tabular}{r|l}
    $n_c$ & $\ell_C = -1/2$ \\ \hline$\vphantom{0^{0^0}}$
    $25/2$ & \blackzone{$2q^{0} (F + F_\times)(12)$} \\
    $23/2$ & \greenzone{$(1 + \epsilon_C) F(11)$} \\
    $21/2$ & \greenzone{$(1 + \epsilon_C) F(10)$} \\
    $19/2$ & \greenzone{$(1 + \epsilon_C) F(9)$} \\
    $17/2$ & \greenzone{$(1 + \epsilon_C) F(8)$} \\
    $15/2$ & \greenzone{$(1 + \epsilon_C) F(7)$} \\
    $13/2$ & \greenzone{$(1 + \epsilon_C) F(6)$} \\
    $11/2$ & \greenzone{$(1 + \epsilon_C) F(5)$} \\
    $9/2$ & \greenzone{$(1 + \epsilon_C) F(4)$} \\
    $7/2$ & \greenzone{$(1 + \epsilon_C) F(3)$} \\
    $5/2$ & \greenzone{$(1 + \epsilon_C) F(2)$} \\
    $3/2$ & \greenzone{$(1 + \epsilon_C) F(1)$} \\
    $1/2$ & \greenzone{$(1 + \epsilon_C) F(0)$} 
  \end{tabular}
  \vspace{0.1in}
\end{minipage}
\hfill
\noindent\begin{minipage}[t]{0.45\textwidth}
For $e = 6, d_0' \geq 2, s' = 0$, type (B): \blackzone black, \plumzone plum, and \greenzone green zones.$\vphantom{|}$
  
  \begin{tabular}{r|l}
    $n_c$ & $\ell_C = 0$ \\ \hline$\vphantom{0^{0^0}}$
    $13$ & \blackzone{$2q^{-1} (F + F_\times)(12)$} \\
    $12$ & \plumzone{$q^{0} F(12)$} \\
    $11$ & \greenzone{$(1 + \epsilon_C)q^{0} F(11)$} \\
    $10$ & \greenzone{$(1 + \epsilon_C)q^{0} F(10)$} \\
    $9$ & \greenzone{$(1 + \epsilon_C)q^{0} F(9)$} \\
    $8$ & \greenzone{$(1 + \epsilon_C)q^{0} F(8)$} \\
    $7$ & \greenzone{$(1 + \epsilon_C)q^{0} F(7)$} \\
    $6$ & \greenzone{$(1 + \epsilon_C)q^{0} F(6)$} \\
    $5$ & \greenzone{$(1 + \epsilon_C)q^{0} F(5)$} \\
    $4$ & \greenzone{$(1 + \epsilon_C)q^{0} F(4)$} \\
    $3$ & \greenzone{$(1 + \epsilon_C)q^{0} F(3)$} \\
    $2$ & \greenzone{$(1 + \epsilon_C)q^{0} F(2)$} \\
    $1$ & \greenzone{$(1 + \epsilon_C)q^{0} F(1)$} 
  \end{tabular}
  \vspace{0.5in}
\end{minipage}
\hfill
\noindent\begin{minipage}{.4\textwidth}
For $e = 6, d_0' = 13/2, s' = 4$, type (B): \blackzone black, \plumzone plum, and \greenzone green zones.$\vphantom{|}$
  
  \begin{tabular}{r|l}
    $n_c$ & $\ell_C = 2$ \\ \hline$\vphantom{0^{0^0}}$
    $13$ & \blackzone{$2q^{-3} (F + F_\times)(12)$} \\
    $12$ & \plumzone{$q^{-2} F(12)$} \\
    $11$ & \plumzone{$q^{-2} F(11)$} \\
    $10$ & \plumzone{$q^{-1} F(11)$} \\
    $9$ & \plumzone{$q^{-1} F(10)$} \\
    $8$ & \plumzone{$q^{0} F(10)$} \\
    $7$ & \plumzone{$q^{0} F(9)$} \\
    $6$ & \plumzone{$q^{1} F(9)$} \\
    $5$ & \plumzone{$q^{1} F(8)$} \\
    $4$ & \plumzone{$q^{2} F(8)$} \\
    $3$ & \greenzone{$(1 + \epsilon_C)q^{2} F(7)$} \\
    $2$ & \greenzone{$(1 + \epsilon_C)q^{2} F(6)$} \\
    $1$ & \greenzone{$(1 + \epsilon_C)q^{2} F(5)$} 
  \end{tabular}
  \vspace{0.1in}
\end{minipage}
\vfill
\noindent\begin{minipage}{\textwidth}
For $e = 6, d_0' = 3, s' = 2$, type (C): \blackzone black, \plumzone plum, \bluezone blue, \greenzone green, and \redzone red zones.$\vphantom{|}$
  
  \begin{tabular}{r|lll}
    $n_c$ & $\ell_C = 1$ & $2$ & $3$ \\ \hline$\vphantom{0^{0^0}}$
    $13$ & \blackzone{$2q^{-2} (F + F_\times)(12)$} & \blackzone{$2q^{-2} (F + F_\times)(12)$} & \blackzone{$2q^{-2} (F + F_\times)(12)$} \\
    $12$ & \plumzone{$q^{-1} F(12)$} & \plumzone{$q^{-1} F(12)$} & \plumzone{$q^{-1} F(12)$} \\
    $11$ & \plumzone{$q^{-1} F(11)$} & \plumzone{$q^{-1} F(11)$} & \plumzone{$q^{-1} F(11)$} \\
    $10$ & \plumzone{$q^{0} F(11)$} & \plumzone{$q^{0} F(11)$} & \plumzone{$q^{0} F(11)$} \\
    $9$ & \plumzone{$q^{0} F(10)$} & \plumzone{$q^{0} F(10)$} & \plumzone{$q^{0} F(10)$} \\
    $8$ & \plumzone{$q^{1} F(10)$} & \plumzone{$q^{1} F(10)$} & \plumzone{$q^{1} F(10)$} \\
    $7$ & \greenzone{$(1 + \epsilon_C)q^{1} F(3,0,0)$} & \bluezone{$q^{1} F(3,0,0)$} & \bluezone{$q^{1} F(3,0,0)$} \\
    $6$ & \greenzone{$(1 + \epsilon_C)q^{1} F(2,1,0)$} & \bluezone{$q^{1} F(2,1,0)$} & \bluezone{$q^{1} F(2,1,0)$} \\
    $5$ & \greenzone{$(1 + \epsilon_C)q^{1} F(2,0,1)$} & \bluezone{$q^{2} F(2,1,0)$} & \bluezone{$q^{2} F(2,1,0)$} \\
    $4$ & \greenzone{$(1 + \epsilon_C)q^{1} F(1,1,1)$} & \bluezone{$q^{2} F(1,2,0)$} & \bluezone{$q^{2} F(1,2,0)$} \\
    $3$ & \greenzone{$(1 + \epsilon_C)q^{1} F(1,0,2)$} & \greenzone{$(1 + \epsilon_C)q^{2} F(1,1,1)$} & \bluezone{$q^{2} F(1,1,1)$} \\
    $2$ & \greenzone{$(1 + \epsilon_C)q^{1} F(0,1,2)$} & \redzone{$(1 + \epsilon_C)q^{2} G^\cross(0,2,1)$} & \redzone{$q^{2} G(1,1,1)$} \\
    $1$ & \greenzone{$(1 + \epsilon_C)q^{1} F(0,0,3)$} & \redzone{$(1 + \epsilon_C)q^{2} G^\cross(0,1,2)$} & \redzone{$q^{3} G(1,1,1)$} 
  \end{tabular}
  \vspace{0.1in}
\end{minipage}

\noindent\begin{minipage}{\textwidth}
\begin{sideways}
  \centering
     \begin{tabular}{c}
For $e = 9, d_0' = 1, s' = 0$, type (C): \blackzone black, \plumzone plum, \bluezone blue, \greenzone green, and \redzone red zones.$\vphantom{|}$
    \\
    \begin{tabular}{r|lllll}
      $n_c$ & $\ell_C = 0$ & $1$ & $2$ & $3$ & $4$ \\ \hline$\vphantom{0^{0^0}}$
      $19$ & \blackzone{$2q^{-1} (F + F_\times)(18)$} & \blackzone{$2q^{-1} (F + F_\times)(18)$} & \blackzone{$2q^{-1} (F + F_\times)(18)$} & \blackzone{$2q^{-1} (F + F_\times)(18)$} & \blackzone{$2q^{-1} (F + F_\times)(18)$} \\
      $18$ & \plumzone{$q^{0} F(18)$} & \plumzone{$q^{0} F(18)$} & \plumzone{$q^{0} F(18)$} & \plumzone{$q^{0} F(18)$} & \plumzone{$q^{0} F(18)$} \\
      $17$ & \greenzone{$(1 + \epsilon_C)q^{0} F(8,0,0)$} & \bluezone{$q^{0} F(8,0,0)$} & \bluezone{$q^{0} F(8,0,0)$} & \bluezone{$q^{0} F(8,0,0)$} & \bluezone{$q^{0} F(8,0,0)$} \\
      $16$ & \greenzone{$(1 + \epsilon_C)q^{0} F(7,1,0)$} & \bluezone{$q^{0} F(7,1,0)$} & \bluezone{$q^{0} F(7,1,0)$} & \bluezone{$q^{0} F(7,1,0)$} & \bluezone{$q^{0} F(7,1,0)$} \\
      $15$ & \greenzone{$(1 + \epsilon_C)q^{0} F(7,0,1)$} & \bluezone{$q^{1} F(7,1,0)$} & \bluezone{$q^{1} F(7,1,0)$} & \bluezone{$q^{1} F(7,1,0)$} & \bluezone{$q^{1} F(7,1,0)$} \\
      $14$ & \greenzone{$(1 + \epsilon_C)q^{0} F(6,1,1)$} & \bluezone{$q^{1} F(6,2,0)$} & \bluezone{$q^{1} F(6,2,0)$} & \bluezone{$q^{1} F(6,2,0)$} & \bluezone{$q^{1} F(6,2,0)$} \\
      $13$ & \greenzone{$(1 + \epsilon_C)q^{0} F(6,0,2)$} & \greenzone{$(1 + \epsilon_C)q^{1} F(6,1,1)$} & \bluezone{$q^{1} F(6,1,1)$} & \bluezone{$q^{1} F(6,1,1)$} & \bluezone{$q^{1} F(6,1,1)$} \\
      $12$ & \greenzone{$(1 + \epsilon_C)q^{0} F(5,1,2)$} & \greenzone{$(1 + \epsilon_C)q^{1} F(5,2,1)$} & \bluezone{$q^{1} F(5,2,1)$} & \bluezone{$q^{1} F(5,2,1)$} & \bluezone{$q^{1} F(5,2,1)$} \\
      $11$ & \greenzone{$(1 + \epsilon_C)q^{0} F(5,0,3)$} & \greenzone{$(1 + \epsilon_C)q^{1} F(5,1,2)$} & \bluezone{$q^{2} F(5,2,1)$} & \bluezone{$q^{2} F(5,2,1)$} & \bluezone{$q^{2} F(5,2,1)$} \\
      $10$ & \greenzone{$(1 + \epsilon_C)q^{0} F(4,1,3)$} & \greenzone{$(1 + \epsilon_C)q^{1} F(4,2,2)$} & \bluezone{$q^{2} F(4,3,1)$} & \bluezone{$q^{2} F(4,3,1)$} & \bluezone{$q^{2} F(4,3,1)$} \\
      $9$ & \greenzone{$(1 + \epsilon_C)q^{0} F(4,0,4)$} & \greenzone{$(1 + \epsilon_C)q^{1} F(4,1,3)$} & \greenzone{$(1 + \epsilon_C)q^{2} F(4,2,2)$} & \bluezone{$q^{2} F(4,2,2)$} & \bluezone{$q^{2} F(4,2,2)$} \\
      $8$ & \greenzone{$(1 + \epsilon_C)q^{0} F(3,1,4)$} & \greenzone{$(1 + \epsilon_C)q^{1} F(3,2,3)$} & \greenzone{$(1 + \epsilon_C)q^{2} F(3,3,2)$} & \bluezone{$q^{2} F(3,3,2)$} & \bluezone{$q^{2} F(3,3,2)$} \\
      $7$ & \greenzone{$(1 + \epsilon_C)q^{0} F(3,0,5)$} & \greenzone{$(1 + \epsilon_C)q^{1} F(3,1,4)$} & \greenzone{$(1 + \epsilon_C)q^{2} F(3,2,3)$} & \bluezone{$q^{3} F(3,3,2)$} & \bluezone{$q^{3} F(3,3,2)$} \\
      $6$ & \greenzone{$(1 + \epsilon_C)q^{0} F(2,1,5)$} & \greenzone{$(1 + \epsilon_C)q^{1} F(2,2,4)$} & \greenzone{$(1 + \epsilon_C)q^{2} F(2,3,3)$} & \bluezone{$q^{3} F(2,4,2)$} & \bluezone{$q^{3} F(2,4,2)$} \\
      $5$ & \greenzone{$(1 + \epsilon_C)q^{0} F(2,0,6)$} & \greenzone{$(1 + \epsilon_C)q^{1} F(2,1,5)$} & \greenzone{$(1 + \epsilon_C)q^{2} F(2,2,4)$} & \redzone{$(1 + \epsilon_C)q^{3} G^\cross(2,3,3)$} & \redzone{$q^{3} G(3,2,3)$} \\
      $4$ & \greenzone{$(1 + \epsilon_C)q^{0} F(1,1,6)$} & \greenzone{$(1 + \epsilon_C)q^{1} F(1,2,5)$} & \redzone{$(1 + \epsilon_C)q^{2} G^\cross(1,3,4)$} & \redzone{$(1 + \epsilon_C)q^{3} G^\cross(2,3,3)$} & \redzone{$q^{3} G(3,2,3)$} \\
      $3$ & \greenzone{$(1 + \epsilon_C)q^{0} F(1,0,7)$} & \greenzone{$(1 + \epsilon_C)q^{1} F(1,1,6)$} & \redzone{$(1 + \epsilon_C)q^{2} G^\cross(1,2,5)$} & \redzone{$(1 + \epsilon_C)q^{3} G^\cross(2,2,4)$} & \redzone{$q^{4} G(3,2,3)$} \\
      $2$ & \greenzone{$(1 + \epsilon_C)q^{0} F(0,1,7)$} & \redzone{$(1 + \epsilon_C)q^{1} G^\cross(0,2,6)$} & \redzone{$(1 + \epsilon_C)q^{2} G^\cross(1,2,5)$} & \redzone{$(1 + \epsilon_C)q^{3} G^\cross(2,2,4)$} & \redzone{$q^{4} G(3,2,3)$} \\
      $1$ & \greenzone{$(1 + \epsilon_C)q^{0} F(0,0,8)$} & \redzone{$(1 + \epsilon_C)q^{1} G^\cross(0,1,7)$} & \redzone{$(1 + \epsilon_C)q^{2} G^\cross(1,1,6)$} & \redzone{$(1 + \epsilon_C)q^{3} G^\cross(2,1,5)$} & \redzone{$(1 + \epsilon_C)q^{4} G^\cross(3,1,4)$} \\
      &&&&& $\quad {} + q^{4} G(4,0,4)$
    \end{tabular}
  \end{tabular}
  \end{sideways}
  \vspace{0.1in}
\end{minipage}

\noindent\begin{minipage}{\textwidth}
  
For $e = 9, d_0' = 5, s' = 4$, type (C): \blackzone black, \plumzone plum, \bluezone blue, \greenzone green, and \redzone red zones.$\vphantom{|}$
  
  \begin{tabular}{r|lll}
    $n_c$ & $\ell_C = 2$ & $3$ & $4$ \\ \hline$\vphantom{0^{0^0}}$$\vphantom{0^{0^0}}$
    $19$ & \blackzone{$2q^{-3} (F + F_\times)(18)$} & \blackzone{$2q^{-3} (F + F_\times)(18)$} & \blackzone{$2q^{-3} (F + F_\times)(18)$} \\
    $18$ & \plumzone{$q^{-2} F(18)$} & \plumzone{$q^{-2} F(18)$} & \plumzone{$q^{-2} F(18)$} \\
    $17$ & \plumzone{$q^{-2} F(17)$} & \plumzone{$q^{-2} F(17)$} & \plumzone{$q^{-2} F(17)$} \\
    $16$ & \plumzone{$q^{-1} F(17)$} & \plumzone{$q^{-1} F(17)$} & \plumzone{$q^{-1} F(17)$} \\
    $15$ & \plumzone{$q^{-1} F(16)$} & \plumzone{$q^{-1} F(16)$} & \plumzone{$q^{-1} F(16)$} \\
    $14$ & \plumzone{$q^{0} F(16)$} & \plumzone{$q^{0} F(16)$} & \plumzone{$q^{0} F(16)$} \\
    $13$ & \plumzone{$q^{0} F(15)$} & \plumzone{$q^{0} F(15)$} & \plumzone{$q^{0} F(15)$} \\
    $12$ & \plumzone{$q^{1} F(15)$} & \plumzone{$q^{1} F(15)$} & \plumzone{$q^{1} F(15)$} \\
    $11$ & \plumzone{$q^{1} F(14)$} & \plumzone{$q^{1} F(14)$} & \plumzone{$q^{1} F(14)$} \\
    $10$ & \plumzone{$q^{2} F(14)$} & \plumzone{$q^{2} F(14)$} & \plumzone{$q^{2} F(14)$} \\
    $9$ & \greenzone{$(1 + \epsilon_C)q^{2} F(4,0,0)$} & \bluezone{$q^{2} F(4,0,0)$} & \bluezone{$q^{2} F(4,0,0)$} \\
    $8$ & \greenzone{$(1 + \epsilon_C)q^{2} F(3,1,0)$} & \bluezone{$q^{2} F(3,1,0)$} & \bluezone{$q^{2} F(3,1,0)$} \\
    $7$ & \greenzone{$(1 + \epsilon_C)q^{2} F(3,0,1)$} & \bluezone{$q^{3} F(3,1,0)$} & \bluezone{$q^{3} F(3,1,0)$} \\
    $6$ & \greenzone{$(1 + \epsilon_C)q^{2} F(2,1,1)$} & \bluezone{$q^{3} F(2,2,0)$} & \bluezone{$q^{3} F(2,2,0)$} \\
    $5$ & \greenzone{$(1 + \epsilon_C)q^{2} F(2,0,2)$} & \greenzone{$(1 + \epsilon_C)q^{3} F(2,1,1)$} & \bluezone{$q^{3} F(2,1,1)$} \\
    $4$ & \greenzone{$(1 + \epsilon_C)q^{2} F(1,1,2)$} & \greenzone{$(1 + \epsilon_C)q^{3} F(1,2,1)$} & \bluezone{$q^{3} F(1,2,1)$} \\
    $3$ & \greenzone{$(1 + \epsilon_C)q^{2} F(1,0,3)$} & \greenzone{$(1 + \epsilon_C)q^{3} F(1,1,2)$} & \bluezone{$q^{4} F(1,2,1)$} \\
    $2$ & \greenzone{$(1 + \epsilon_C)q^{2} F(0,1,3)$} & \redzone{$(1 + \epsilon_C)q^{3} G^\cross(0,2,2)$} & \redzone{$q^{4} G(1,2,1)$} \\
    $1$ & \greenzone{$(1 + \epsilon_C)q^{2} F(0,0,4)$} & \redzone{$(1 + \epsilon_C)q^{3} G^\cross(0,1,3)$} & \redzone \quickcol{
      $(1 + \epsilon_C)q^{4} G^\cross(1,1,2)\!\!$ \\
     $\quad {} + q^{4} G(2,0,2)\!\!\!$}
  \end{tabular}
  \vspace{0.1in}
\end{minipage}

\noindent\begin{minipage}{\textwidth}
For $e = 6, d_0' = 1, s' = 5$, type (D): \blackzone black, \plumzone plum, \purplezone purple, \bluezone blue, \greenzone green, and \redzone red zones.$\vphantom{|}$
  
  \begin{tabular}{r|lll}
    $n_c$ & $\ell_C = 0$ & $1$ & $2, 3$ \\ \hline$\vphantom{0^{0^0}}$$\vphantom{0^{0^0}}$
    $13$ & \blackzone{$2q^{-4} (F + F_\times)(12)$} & \blackzone{$2q^{-4} (F + F_\times)(12)$} & \blackzone{$2q^{-4} (F + F_\times)(12)$} \\
    $12$ & \plumzone{$q^{-3} F(12)$} & \plumzone{$q^{-3} F(12)$} & \plumzone{$q^{-3} F(12)$} \\
    $11$ & \purplezone{$2q^{-3} F(5,0, \emptyset)$} & \purplezone{$2q^{-3} F(5,0, \emptyset)$} & \purplezone{$2q^{-3} F(5,0, \emptyset)$} \\
    $10$ & \purplezone{$2q^{-2} F(5,0, \emptyset)$} & \purplezone{$2q^{-2} F(5,0, \emptyset)$} & \purplezone{$2q^{-2} F(5,0, \emptyset)$} \\
    $9$ & \purplezone{$2q^{-2} F(4,1, \emptyset)$} & \purplezone{$2q^{-2} F(4,1, \emptyset)$} & \purplezone{$2q^{-2} F(4,1, \emptyset)$} \\
    $8$ & \purplezone{$2q^{-1} F(4,1, \emptyset)$} & \purplezone{$2q^{-1} F(4,1, \emptyset)$} & \purplezone{$2q^{-1} F(4,1, \emptyset)$} \\
    $7$ & \purplezone{$2q^{-1} F(3,2, \emptyset)$} & \purplezone{$2q^{-1} F(3,2, \emptyset)$} & \purplezone{$2q^{-1} F(3,2, \emptyset)$} \\
    $6$ & \bluezone{$q^{0} F(3,2,0)$} & \bluezone{$q^{0} F(3,2,0)$} & \bluezone{$q^{0} F(3,2,0)$} \\
    $5$ & \bluezone{$q^{0} F(2,3,0)$} & \bluezone{$q^{0} F(2,3,0)$} & \bluezone{$q^{0} F(2,3,0)$} \\
    $4$ & \greenzone{$(1 + \epsilon_C)q^{0} F(2,2,1)$} & \bluezone{$q^{0} F(2,2,1)$} & \bluezone{$q^{0} F(2,2,1)$} \\
    $3$ & \greenzone{$(1 + \epsilon_C)q^{0} F(1,3,1)$} & \bluezone{$q^{0} F(1,3,1)$} & \bluezone{$q^{0} F(1,3,1)$} \\
    $2$ & \greenzone{$(1 + \epsilon_C)q^{0} F(1,2,2)$} & \bluezone{$q^{1} F(1,3,1)$} & \bluezone{$q^{1} F(1,3,1)$} \\
    $1$ & \redzone{$(1 + \epsilon_C)q^{0} G^\cross(0,3,2)$} & \redzone{$q^{1} G(1,3,1)$} &    \end{tabular}
  \vspace{0.1in}
\end{minipage}

\noindent\begin{minipage}{\textwidth}
For $e = 6, d_0' = 2, s' = 4$, type (D): \blackzone black, \plumzone plum, \purplezone purple, \bluezone blue, \greenzone green, and \redzone red zones.$\vphantom{|}$
  
  \begin{tabular}{r|lll}
    $n_c$ & $\ell_C = 1$ & $2$ & $3$ \\ \hline$\vphantom{0^{0^0}}$$\vphantom{0^{0^0}}$
    $13$ & \blackzone{$2q^{-3} (F + F_\times)(12)$} & \blackzone{$2q^{-3} (F + F_\times)(12)$} & \blackzone{$2q^{-3} (F + F_\times)(12)$} \\
    $12$ & \plumzone{$q^{-2} F(12)$} & \plumzone{$q^{-2} F(12)$} & \plumzone{$q^{-2} F(12)$} \\
    $11$ & \plumzone{$q^{-2} F(11)$} & \plumzone{$q^{-2} F(11)$} & \plumzone{$q^{-2} F(11)$} \\
    $10$ & \plumzone{$q^{-1} F(11)$} & \plumzone{$q^{-1} F(11)$} & \plumzone{$q^{-1} F(11)$} \\
    $9$ & \purplezone{$2q^{-1} F(4,0, \emptyset)$} & \purplezone{$2q^{-1} F(4,0, \emptyset)$} & \purplezone{$2q^{-1} F(4,0, \emptyset)$} \\
    $8$ & \purplezone{$2q^{0} F(4,0, \emptyset)$} & \purplezone{$2q^{0} F(4,0, \emptyset)$} & \purplezone{$2q^{0} F(4,0, \emptyset)$} \\
    $7$ & \purplezone{$2q^{0} F(3,1, \emptyset)$} & \purplezone{$2q^{0} F(3,1, \emptyset)$} & \purplezone{$2q^{0} F(3,1, \emptyset)$} \\
    $6$ & \bluezone{$q^{1} F(3,1,0)$} & \bluezone{$q^{1} F(3,1,0)$} & \bluezone{$q^{1} F(3,1,0)$} \\
    $5$ & \bluezone{$q^{1} F(2,2,0)$} & \bluezone{$q^{1} F(2,2,0)$} & \bluezone{$q^{1} F(2,2,0)$} \\
    $4$ & \greenzone{$(1 + \epsilon_C)q^{1} F(2,1,1)$} & \bluezone{$q^{1} F(2,1,1)$} & \bluezone{$q^{1} F(2,1,1)$} \\
    $3$ & \greenzone{$(1 + \epsilon_C)q^{1} F(1,2,1)$} & \bluezone{$q^{1} F(1,2,1)$} & \bluezone{$q^{1} F(1,2,1)$} \\
    $2$ & \greenzone{$(1 + \epsilon_C)q^{1} F(1,1,2)$} & \bluezone{$q^{2} F(1,2,1)$} & \bluezone{$q^{2} F(1,2,1)$} \\
    $1$ & \redzone{$(1 + \epsilon_C)q^{1} G^\cross(0,2,2)$} & \redzone{$q^{2} G(1,2,1)$} & 
  \end{tabular}
  \vspace{0.1in}
\end{minipage}

\noindent\begin{minipage}{0.5\textwidth}
For $e = 6, d_0' = 4, s' = 20$, type (D): \blackzone black, \plumzone plum, and \purplezone purple zones.$\vphantom{|}$
  
  \begin{tabular}{r|l}
    $n_c$ & $\ell_C = 2,3$ \\ \hline$\vphantom{0^{0^0}}$$\vphantom{0^{0^0}}$
    $13$ & \blackzone{$2q^{-11} (F + F_\times)(12)$} \\
    $12$ & \plumzone{$q^{-10} F(12)$} \\
    $11$ & \plumzone{$q^{-10} F(11)$} \\
    $10$ & \plumzone{$q^{-9} F(11)$} \\
    $9$ & \plumzone{$q^{-9} F(10)$} \\
    $8$ & \plumzone{$q^{-8} F(10)$} \\
    $7$ & \plumzone{$q^{-8} F(9)$} \\
    $6$ & \plumzone{$q^{-7} F(9)$} \\
    $5$ & \purplezone{$2q^{-7} F(2,0, \emptyset)$} \\
    $4$ & \purplezone{$2q^{-6} F(2,0, \emptyset)$} \\
    $3$ & \purplezone{$2q^{-6} F(1,1, \emptyset)$} \\
    $2$ & \purplezone{$2q^{-5} F(1,1, \emptyset)$} \\
    $1$ & \purplezone{$2q^{-5} F(0,2, \emptyset)$}   \end{tabular}
  \vspace{0.1in}
\end{minipage}
\hfill
\noindent\begin{minipage}{0.4\textwidth}
For $e = 6, d_0' = 13/2, s' = 7$, type (D): \blackzone black and \plumzone plum zones.$\vphantom{|}$
  
  \begin{tabular}{r|l}
    $n_c$ & $\ell_C = 3$ \\ \hline$\vphantom{0^{0^0}}$$\vphantom{0^{0^0}}$
    $13$ & \blackzone{$2q^{-5} (F + F_\times)(12)$} \\
    $12$ & \plumzone{$q^{-4} F(12)$} \\
    $11$ & \plumzone{$q^{-4} F(11)$} \\
    $10$ & \plumzone{$q^{-3} F(11)$} \\
    $9$ & \plumzone{$q^{-3} F(10)$} \\
    $8$ & \plumzone{$q^{-2} F(10)$} \\
    $7$ & \plumzone{$q^{-2} F(9)$} \\
    $6$ & \plumzone{$q^{-1} F(9)$} \\
    $5$ & \plumzone{$q^{-1} F(8)$} \\
    $4$ & \plumzone{$q^{0} F(8)$} \\
    $3$ & \plumzone{$q^{0} F(7)$} \\
    $2$ & \plumzone{$q^{1} F(7)$} \\
    $1$ & \plumzone{$q^{1} F(6)$} 
  \end{tabular}
  \vspace{0.1in}
\end{minipage}

\noindent\begin{minipage}[t]{0.45\textwidth}
For $e = 6, d_0' = 1, s' = 2$, type~(E): \blackzone black, \plumzone plum, \purplezone purple, and \greenzone green zones.$\vphantom{|}$
  
  \begin{tabular}{r|l}
    $n_c$ & $\ell_C = 0$ \\ \hline$\vphantom{0^{0^0}}$$\vphantom{0^{0^0}}$
    $13$ & \blackzone{$2q^{-2} (F + F_\times)(12)$} \\
    $12$ & \plumzone{$q^{-1} F(12)$} \\
    $11$ & \purplezone{$2q^{-1} F(5,0, \emptyset)$} \\
    $10$ & \purplezone{$2q^{0} F(5,0, \emptyset)$} \\
    $9$ & \greenzone{$(1 + \epsilon_C)q^{0} F(4,1,0)$} \\
    $8$ & \greenzone{$(1 + \epsilon_C)q^{0} F(4,0,1)$} \\
    $7$ & \greenzone{$(1 + \epsilon_C)q^{0} F(3,1,1)$} \\
    $6$ & \greenzone{$(1 + \epsilon_C)q^{0} F(3,0,2)$} \\
    $5$ & \greenzone{$(1 + \epsilon_C)q^{0} F(2,1,2)$} \\
    $4$ & \greenzone{$(1 + \epsilon_C)q^{0} F(2,0,3)$} \\
    $3$ & \greenzone{$(1 + \epsilon_C)q^{0} F(1,1,3)$} \\
    $2$ & \greenzone{$(1 + \epsilon_C)q^{0} F(1,0,4)$} \\
    $1$ & \greenzone{$(1 + \epsilon_C)q^{0} F(0,1,4)$} 
\end{tabular}
  \vspace{0.1in}
\end{minipage}
\hfill
\noindent\begin{minipage}[t]{0.45\textwidth}
For $e = 6, d_0' = 1, s' = 20$, type (E): \blackzone black, \plumzone plum, and \purplezone purple zones.$\vphantom{|}$
  
  \begin{tabular}{r|l}
    $n_c$ & $\ell_C = 0$ \\ \hline$\vphantom{0^{0^0}}$$\vphantom{0^{0^0}}$
    $13$ & \blackzone{$2q^{-11} (F + F_\times)(12)$} \\
    $12$ & \plumzone{$q^{-10} F(12)$} \\
    $11$ & \purplezone{$2q^{-10} F(5,0, \emptyset)$} \\
    $10$ & \purplezone{$2q^{-9} F(5,0, \emptyset)$} \\
    $9$ & \purplezone{$2q^{-9} F(4,1, \emptyset)$} \\
    $8$ & \purplezone{$2q^{-8} F(4,1, \emptyset)$} \\
    $7$ & \purplezone{$2q^{-8} F(3,2, \emptyset)$} \\
    $6$ & \purplezone{$2q^{-7} F(3,2, \emptyset)$} \\
    $5$ & \purplezone{$2q^{-7} F(2,3, \emptyset)$} \\
    $4$ & \purplezone{$2q^{-6} F(2,3, \emptyset)$} \\
    $3$ & \purplezone{$2q^{-6} F(1,4, \emptyset)$} \\
    $2$ & \purplezone{$2q^{-5} F(1,4, \emptyset)$} \\
    $1$ & \purplezone{$2q^{-5} F(0,5, \emptyset)$} 
  \end{tabular}
  \vspace{0.1in}
\end{minipage}

\subsection{Splitting type (\texorpdfstring{$1^21$}{1²1}) (light zones)}\ \\[\baselineskip]
\noindent\begin{minipage}{\textwidth}
  \centering
     \begin{tabular}{c}
For $e = 6, d_0' = 1, o_1 = 2$: \brownzone brown, \yellowzone yellow, and \lemonzone lemon zones.$\vphantom{|}$
  
  \\ \begin{tabular}{rc|l}
    $n_\beta$ & $n_{11}$ & \\ \hline$\vphantom{0^{0^0}}$
    $12$ & $12,13$ & \brownzone{$2q^{0} \xo F(0, 5, \emptyset)$} \\
    $10$ & $10,11$ & \yellowzone{$(1 + \epsilon_C)q^{0} \xo G(0,4,1)$} \\
    $8$ & $8,9$ & \yellowzone{$(1 + \epsilon_C)q^{0} \xo G(0,3,2)$} \\
    $6$ & $6,7$ & \yellowzone{$(1 + \epsilon_C)q^{0} \xo G(0,2,3)$} \\
    $4$ & $4,5$ & \yellowzone{$(1 + \epsilon_C)q^{0} \xo G(0,1,4)$} \\
    $2$ & $2,3$ & \yellowzone{$(1 + \epsilon_C)q^{0} \xo G(0,0,5)$} \\
    $0$ & $0,1$ & \lemonzone \quickcol{$(1 + \epsilon_C)q^{0} G^\cross(0) + (1 + \epsilon_C)q^{1} G^\cross(2) + (1 + \epsilon_C)q^{2} G^\cross(4)$ \\ $\quad{} + q^{3} G(6)$}
  \end{tabular}
\end{tabular}
  \vspace{0.1in}
\end{minipage}

\noindent\begin{minipage}{\textwidth}
  \centering
     \begin{tabular}{c}
For $e = 6, d_0' = 1, o_1 = 0$: \brownzone brown and \yellowzone yellow zones.$\vphantom{|}$
  
  \\ \begin{tabular}{rc|l}
    $n_\beta$ & $n_{11}$ & \\ \hline$\vphantom{0^{0^0}}$
    $13$ & $13,14$ & \brownzone{$2q^{-1} F(0, 5, \emptyset)$} \\
    $11$ & $11,12$ & \yellowzone{$q^{-1} G(0,5,0)$} \\
    $9$ & $9,10$ & \yellowzone{$q^{0} G(0,5,0)$} \\
    $7$ & $7,8$ & \yellowzone{$(1 + \epsilon_C)q^{0} G^\cross(0,3,2) + q^{0} G(1,3,1)$} \\
    $5$ & $5,6$ & \yellowzone{$(1 + \epsilon_C)q^{0} G^\cross(0,2,3) + q^{1} G(1,3,1)$} \\
    $3$ & $3,4$ & \yellowzone{$(1 + \epsilon_C)q^{0} G^\cross(0,1,4) + (1 + \epsilon_C)q^{1} G^\cross(1,1,3) + q^{1} G(2,1,2)$} \\
    $1$ & $1,2$ & \yellowzone{$(1 + \epsilon_C)q^{0} G^\cross(0,0,5) + (1 + \epsilon_C)q^{1} G^\cross(1,0,4) + q^{2} G(2,1,2)$} 
  \end{tabular}
\end{tabular}
  \vspace{0.1in}
\end{minipage}

\noindent\begin{minipage}{\textwidth}
  \centering
     \begin{tabular}{c}
For $e = 6, d_0' = 4, o_1 = 0$: \brownzone brown, \yellowzone yellow, and \lemonzone lemon zones.$\vphantom{|}$
  
  \\ \begin{tabular}{rc|l}
    $n_\beta$ & $n_{11}$ & \\ \hline$\vphantom{0^{0^0}}$
    $12$ & $15,16$ & \brownzone{$2q^{0} F(0, 2, \emptyset)$} \\
    $10$ & $13,14$ & \brownzone{$2q^{1} F(0, 2, \emptyset)$} \\
    $8$ & $11,12$ & \yellowzone{$q^{1} G(0,2,0)$} \\
    $6$ & $9,10$ & \yellowzone{$q^{2} G(0,2,0)$} \\
    $4$ & $7,8$ & \yellowzone{$(1 + \epsilon_C)q^{2} G^\cross(0,0,2) + q^{2} G(1,0,1)$} \\
    $2$ & $5,6$ & \lemonzone{$(1 + \epsilon_C)q^{1} G^\cross(2) + (1 + \epsilon_C)q^{2} G^\cross(4) + q^{3} G(6)$} \\
    $0$ & $3,4$ & \lemonzone \quickcol{$(1 + \epsilon_C)q^{0} G^\cross(0) + (1 + \epsilon_C)q^{1} G^\cross(2) + (1 + \epsilon_C)q^{2} G^\cross(4)$ \\ $\quad{} + q^{3} G(6)$}
  \end{tabular}
\end{tabular}
  \vspace{0.1in}
\end{minipage}

\noindent\begin{minipage}{\textwidth}
  \centering
     \begin{tabular}{c}
For $e = 6, d_0' = 4, o_1 = 2$: \brownzone brown, \yellowzone yellow, and \lemonzone lemon zones.$\vphantom{|}$
  
  \\ \begin{tabular}{rc|l}
    $n_\beta$ & $n_{11}$ & \\ \hline$\vphantom{0^{0^0}}$
    $13$ & $16,17$ & \brownzone{$2q^{-1} \xo F(0, 2, \emptyset)$} \\
    $11$ & $14,15$ & \brownzone{$2q^{0} \xo F(0, 2, \emptyset)$} \\
    $9$ & $12,13$ & \brownzone{$2q^{1} \xo F(0, 2, \emptyset)$} \\
    $7$ & $10,11$ & \yellowzone{$(1 + \epsilon_C)q^{1} \xo G(0,1,1)$} \\
    $5$ & $8,9$ & \yellowzone{$(1 + \epsilon_C)q^{1} \xo G(0,0,2)$} \\
    $3$ & $6,7$ & \lemonzone{$(1 + \epsilon_C)q^{1} G^\cross(3) + (1 + \epsilon_C)q^{2} G^\cross(5) + q^{3} G(6)$}\\
    $1$ & $4,5$ & \lemonzone \quickcol{$(1 + \epsilon_C)q^{0} G^\cross(1) + (1 + \epsilon_C)q^{1} G^\cross(3) + (1 + \epsilon_C)q^{2} G^\cross(5)$ \\ $\quad{} + q^{3} G(6)$}
  \end{tabular}
\end{tabular}
  \vspace{0.1in}
\end{minipage}

\noindent\begin{minipage}{\textwidth}
  \centering
     \begin{tabular}{c}
For $e = 6, d_0' = 6 + 1/2, o_1 = 3$: \brownzone brown and \lemonzone lemon zones.$\vphantom{|}$
  
  \\ \begin{tabular}{rc|l}
    $n_\beta$ & $n_{11}$ & \\ \hline$\vphantom{0^{0^0}}$
    $12$ & $17 + 1/2,18 + 1/2$ & \brownzone{$q^{0} F(6)$} \\
    $10$ & $15 + 1/2,16 + 1/2$ & \brownzone{$q^{1} F(6)$} \\
    $8$ & $13 + 1/2,14 + 1/2$ & \brownzone{$q^{2} F(6)$} \\
    $6$ & $11 + 1/2,12 + 1/2$ & \lemonzone{$q^{3} G(6)$}\\
    $4$ & $9 + 1/2,10 + 1/2$ & \lemonzone{$(1 + \epsilon_C)q^{2} G^\cross(4) + q^{3} G(6)$}\\
    $2$ & $7 + 1/2,8 + 1/2$ & \lemonzone{$(1 + \epsilon_C)q^{1} G^\cross(2) + (1 + \epsilon_C)q^{2} G^\cross(4) + q^{3} G(6)$}\\
    $0$ & $5 + 1/2,6 + 1/2$ & \lemonzone \quickcol{$(1 + \epsilon_C)q^{0} G^\cross(0) + (1 + \epsilon_C)q^{1} G^\cross(2)$ \\ $\quad{} + (1 + \epsilon_C)q^{2} G^\cross(4) + q^{3} G(6)$}
  \end{tabular}
\end{tabular}
  \vspace{0.1in}
\end{minipage}

\noindent\begin{minipage}{\textwidth}
  \centering
     \begin{tabular}{c}
For $e = 6, d_0' = 6 + 1/2, o_1 = 1$: \brownzone brown and \lemonzone lemon zones.$\vphantom{|}$
  
  \\ \begin{tabular}{rc|l}
    $n_\beta$ & $n_{11}$ & \\ \hline$\vphantom{0^{0^0}}$
    $13$ & $18 + 1/2,19 + 1/2$ & \brownzone{$q^{-1} F(6)$} \\
    $11$ & $16 + 1/2,17 + 1/2$ & \brownzone{$q^{0} F(6)$} \\
    $9$ & $14 + 1/2,15 + 1/2$ & \brownzone{$q^{1} F(6)$} \\
    $7$ & $12 + 1/2,13 + 1/2$ & \brownzone{$q^{2} F(6)$} \\
    $5$ & $10 + 1/2,11 + 1/2$ & \lemonzone{$(1 + \epsilon_C)q^{2} G^\cross(5) + q^{3} G(6)$} \\
    $3$ & $8 + 1/2,9 + 1/2$ & \lemonzone{$(1 + \epsilon_C)q^{1} G^\cross(3) + (1 + \epsilon_C)q^{2} G^\cross(5) + q^{3} G(6)$}\\
    $1$ & $6 + 1/2,7 + 1/2$ & \lemonzone \quickcol{$(1 + \epsilon_C)q^{0} G^\cross(1) + (1 + \epsilon_C)q^{1} G^\cross(3)$ \\ $\quad{} + (1 + \epsilon_C)q^{2} G^\cross(5) + q^{3} G(6)$}
  \end{tabular}
\end{tabular}
  \vspace{0.1in}
\end{minipage}